\documentclass[francais]{cedram-aif}

\usepackage[all]{xy}
% ...
\usepackage[T1]{fontenc}
\usepackage[applemac]{inputenc} 
\usepackage{color}
\usepackage[francais]{babel} 
\usepackage{textcomp}

\definecolor{darkgreen}{rgb}{0.2,0.5,0.2}

\catcode`\@=11
\def\mysmash{\relax 
  \ifmmode\def\next{\mathpalette\mathsm@mysh}\else\let\next\makesm@mysh
  \fi\next}
\def\makesm@mysh#1{\setbox\z@\hbox{#1}\finsm@mysh}
\def\mathsm@mysh#1#2{\setbox\z@\hbox{$\m@th#1{#2}$}\finsm@mysh}
\def\finsm@mysh{\ht\z@\z@ \dp\z@\z@ \box\z@ }
\def\finsm@myshbot{\dp\z@\z@ \box\z@}
\def\finsm@myshtop{\ht\z@\z@ \box\z@}
\def\finhalfsm@myshbot{\dp\z@0.5\dp\z@ \box\z@}
\def\finhalfsm@myshtop{\ht\z@0.5\ht\z@ \box\z@}
\def\finhalfsm@mysh{\ht\z@0.5\ht\z@ \dp\z@0.5\dp\z@ \box\z@}
\def\fingoodsm@mysh#1#2{\ht\z@#1\ht\z@ \dp\z@#2\dp\z@ \box\z@}
\def\smashtop#1{\begingroup\let\finsm@mysh\finsm@myshtop\mysmash{#1}\endgroup}
\def\smashbot#1{\begingroup\let\finsm@mysh\finsm@myshbot\mysmash{#1}\endgroup}
\def\halfsmashtop#1{\begingroup\let\finsm@mysh\finhalfsm@myshtop\mysmash{#1}\endgroup}
\def\halfsmashbot#1{\begingroup\let\finsm@mysh\finhalfsm@myshbot\mysmash{#1}\endgroup}
\def\halfsmash#1{\begingroup\let\finsm@mysh\finhalfsm@mysh\mysmash{#1}\endgroup}
\def\goodsmash#1#2#3{\begingroup\def\finsm@mysh{\fingoodsm@mysh{#1}{#2}}\mysmash{#3}\endgroup}

\let\goth\mathfrak

\catcode`\:=13
\def:{\ifmmode\string:\else\unskip~\string:\fi}

\def\Rhom{\hbox{{\bfit R}{\it hom}\kern2pt}}
\def\RhomDX{\ifmmode\hbox{{\bfit R}{\it hom}\kern2pt}_\DX\else$\RhomDX$\fi}
\def\RhomDXan{\ifmmode\hbox{{\bfit R}{\it hom}\kern2pt}_\DXan\else$\RhomDX$\fi}
\def\limproj{{\displaystyle\lim_{\leftarrow}}}

\def\DX{{{\cal D}_X}}
\def\DXan{{{\cal D}_{X^{an}}}}

\def\DR{\mathop{\it DR}\nolimits}

\def\codim{\mathord{\rm codim}}

\def\Ker{\mathop{\rm Ker}\nolimits}

\def\dim{\mathop{\rm dim}\nolimits}

\def\log{\mathop{\rm Log}\nolimits}

\def\Im{\mathop{\rm Im}\nolimits}
\def\ext{\mathop{\rm Ext}\nolimits}
\def\hom{\mathop{\rm Hom}\nolimits}

\def\alg{\mathop{\rm alg}\nolimits}

\let\CHI\chi
\def\chi{\raise1pt\hbox{$\CHI$}}

\def\bfF{\mathbf F}
\def\bfR{\mathbf R}
\def\bfL{\mathbf L}

\def\mathbffragile#1{\mathchoice
{\hbox{\boldmath$\textstyle#1$}}
{\hbox{\boldmath$\displaystyle#1$}}
{\hbox{\boldmath$\scriptstyle#1$}}
{\hbox{\boldmath$\scriptstyle#1$}}
}
\def\mathbf{\protect\mathbffragile}

\def\Im{\mathop{\rm Im}\nolimits}

\def\theostyle{\bf}

\def\demo{\ifhmode\par\fi\vskip-\lastskip\vskip1ex\penalty-50
\noindent{\theostyle Démonstration. }\leftmargin0pt\def\theostyle{\sl}\ignorespaces}

\def\enddemo{\gdef\theostyle{\bf}%
\ifmmode\hbox to1pt{\hss$\scriptstyle\blacksquare$}\else\unskip\null\nobreak\hfill\nobreak
${\scriptstyle\blacksquare}$\vskip1em\fi\leftmargin0pt}

\def\endsubdemo{\ifmmode\hbox to0pt{\hss$\scriptstyle\square$}\else
\unskip\null\nobreak\hfill\nobreak$\scriptstyle\square$\vskip1,5ex\fi}

\def\demor#1{\ifhmode\par\fi\vskip-\lastskip\vskip1em\penalty-50
\noindent{\bf Démonstration #1. }\ignorespaces}

\def\ps{$p$-adiques }

\def\MLS{\mathop{\rm MLS}\nolimits}

\def\actimes#1#2{{\setbox1=\hbox{$\scriptscriptstyle#2$}%
\setbox3=\hbox{$#1$}%
\setbox4=\hbox{$\scriptstyle#1$}%
\dimen0=0,5\wd3
\setbox5=\hbox{$\scriptscriptstyle\nwarrow$}%
\setbox6=\hbox{$\otimes$}%
\advance\dimen0 by -\wd5
\ifdim\dimen0<\wd1\dimen0=\wd1\fi
\setbox0=\hbox{$\scriptscriptstyle\nwarrow$\kern-1pt
\vbox{\offinterlineskip\parindent0pt\hsize\dimen0
\hbox to\hsize{\hss\copy1\hss}%
\vrule height-0,95pt depth1,2pt
width\dimen0}$\,\vcenter{\kern0,5\ht3\offinterlineskip
\parindent0pt\hbox{$\scriptstyle#1$}}$}
{#1}\vtop{\kern0pt\hbox{\kern-0,5\wd3\box0}}\kern-0,5\wd4\kern-0,5\wd6}{\otimes}\,}
\def\Lotimes{\mathop{\mathop{\otimes}\limits^\bfL }}

\def\flechedroite#1{\mathrel{\hbox to#1{\rightarrowfill}}}
\def\flechegauche#1{\mathrel{\hbox to#1{\leftarrowfill}}}

\def\double#1{\,\vcenter{\offinterlineskip\hbox{$\scriptstyle#1$}\kern1pt\hbox{$\scriptstyle#1$}}\,}
\def\iso{\stackrel{\vtop to0pt{\kern0pt\hbox{$\scriptscriptstyle\Sim$}\vss}}{\to}}

\newdimen\scriptwd
\newdimen\scriptdp
\def\fbigotimes{\mathop{\hbox{\footnotesize$\bigotimes$}}}

\def\otimesdisplay_#1{\vtop{\mathsurround0pt
\offinterlineskip
\hbox to\scriptwd{\hss\footnotesize$\bigotimes$\hss}
\kern\scriptdp
\hbox to\scriptwd{\hss$\scriptstyle#1$\hss}}}

\def\Lotimesdisplay_#1{\vtop{\mathsurround0pt
\offinterlineskip
\hbox to\scriptwd{\hss\footnotesize$\stackrel\bfL \bigotimes$\hss}
\kern\scriptdp
\hbox to\scriptwd{\hss$\scriptstyle#1$\hss}}}

\def\Otimes_#1{\mathchoice
{\otimesdisplay_{#1}}
{\otimesdisplay_{#1}}
{\otimes_{#1}}
{\otimes_{#1}}}

\def\LOtimes_#1{\mathchoice {\Lotimesdisplay_{#1}}
{\Lotimesdisplay_{#1}}
{\stackrel\bfL \otimes_{#1}}
{\stackrel\bfL \otimes_{#1}}}

\def\Lotimes{\stackrel\bfL \otimes}

\def\cTor{\mathop{\mathcal T\mkern-4mu\it or}\nolimits}
\def\cHom{\mathop{\cal H\!\it om}\nolimits}
\def\cDiff{\mathop{\cal D\it iff}\nolimits}
\def\alg{\mathop{\rm alg}\nolimits}
\def\diff{\mathord{\rm diff\,}}
\def\Prefais{\mathop{\mathbf{Prefais}}\nolimits}
\def\Fais{\mathop{\mathbf{Fais}}\nolimits}
\def\Inv{\mathop{\rm Inv}\nolimits}
\def\fonct{\mathrel{\rightsquigarrow}}
\def\Sp{\mathop{\rm Sp}\nolimits}

\def\dotuteur{\vrule height\dimen0 width0pt\endgroup}
\def\tuteur{\begingroup\afterassignment\dotuteur\dimen0=}

\let\mydownarrow\downarrow
\def\vardownarrow{\begingroup\afterassignment\dovardownarrow\dimen0=}
\def\dovardownarrow{\hbox{$\left\mydownarrow\vcenter
to\dimen0{}\right.$}\endgroup}

\def\downarrow{\vardownarrow15pt\vcenter to20pt{} }

\let\myuparrow\uparrow
\def\varuparrow{\begingroup\afterassignment\dovaruparrow\dimen0=}
\def\dovaruparrow{\left\myuparrow\vcenter
to\dimen0{}\right.\endgroup}

\def\uparrow{\varuparrow15pt\vcenter to20pt{} }

\def\varrightarrow#1#2{\mathop{\hbox to#1{\hss\hbox to
#2{\rightarrowfill}\hss}}\limits}

\let\varto\varrightarrow

\def\varrightarrows#1#2#3{\mathrel{\hbox to#2{\hss$\vcenter{\offinterlineskip
\count11=#1\loop
\smash{\raise-2pt\hbox to#3{\rightarrowfill}}
\advance\count11 by-1 \ifnum\count11 > 0  \kern3.5pt \repeat}$\hss}}}

\def\vegal#1#2#3{\llap{$\scriptstyle#1$\kern2pt}\left|\!\vcenter to#3{}\right|\rlap{\kern2pt$\scriptstyle#2$}}

\def\fixepreskip{\abovedisplayshortskip\abovedisplayskip}
\def\fixepostskip{\belowdisplayshortskip\belowdisplayskip}
\def\preskip{\afterassignment\fixepreskip\abovedisplayskip}
\def\postskip{\afterassignment\fixepostskip\belowdisplayskip}

\let\too\longrightarrow
\def\longhookrightarrow{\lhook \joinrel \too}
\def\inj{\mathord{\rm inj}}
\def\surj{\mathord{\rm surj}}
\def\Coker{\mathop{\rm Coker}\nolimits}

{\catcode`\@=11
\gdef\closeoutglossary{\immediate\closeout\@glossaryfile}}

\def\trans{{\,}^{t}\!}

\def\inffrag{\mathop{\rm inf}\nolimits}
\def\inf{\protect\inffrag}

\let\qed\relax

\def\rmD{\mathop{\bf {}D}\nolimits}
\def\Db{\rmD\sp{b}}
\def\Dplus{\rmD\sp{+}}
\def\Dmoins{\rmD\sp{-}}
\def\dR{\mathord{dR}}
\def\Tr{\mathop{\rm Tr}\nolimits}
\def\Ab{\mathord{\rm Ab}}
\def\sing{\mathord{\rm sing}}
\def\Adj{\mathord{\rm Adj}}
\def\aff{\mathord{\rm aff}}
\def\Ouv{\mathop{\cal Ouv}\nolimits}
\def\Retr{\mathop{\rm Retr}\nolimits}
\def\txtMod{Mod}
\def\RmMod{\mathop{\rm \txtMod}\nolimits}
\def\Modfragile{\mathchoice
{\hbox{\rm-\txtMod}}
{\hbox{\rm-\txtMod}}
{\hbox{\scriptsize\rm-\txtMod}}
{\hbox{\tiny\rm-\txtMod}}}
\def\Moddfragile{\mathchoice
{\hbox{\rm \txtMod-}}
{\hbox{\rm \txtMod-}}
{\hbox{\scriptsize\rm \txtMod-}}
{\hbox{\rm\tiny \txtMod-}}}

\def\Mod{\protect\Modfragile}
\def\Modd{\protect\Moddfragile}
\def\indiced#1{\vcenter{\rlap{$\scriptstyle#1$}}}

\def\SmAff{\mathop{\rm SmAff}\nolimits}
\def\Sms{\mathop{\rm Sms}\nolimits}
\def\MLS{\mathop{\rm MLS}\nolimits}
\def\DDD{I\mkern-6muD}
\def\EEE{I\mkern-6muE}
\def\dag{\sp{\dagger}}
\def\smalldag{\sp{\smalldagger}}

\def\daginf{\dag\sb{\inf}}
\def\smalldaginf{\smalldag\sb{\inf}}

\def\Ydaginf{Y\daginf}
\def\Xdaginf{X\daginf}

\def\Udaginf{{U}\daginf}

\let\smalldagger\dagger

\def\makeletterdaginf#1{
\expandafter\edef\csname #1daginf\endcsname
{\noexpand\protect\csname #1daginfchoose\endcsname}
\expandafter\def\csname #1daginfchoose\endcsname
{\mathchoice
{#1\daginf}
{#1\daginf}
{\goodsmash{0.8}1{#1\daginf}}
{\goodsmash{0.8}1{#1\smalldaginf}}}
\expandafter\edef\csname cal#1dag\endcsname
{\noexpand\protect\csname cal#1dagchoose\endcsname}
\expandafter\def\csname cal#1dagchoose\endcsname{
\mathchoice
{\cal #1\dag}
{\cal #1\dag}
{\cal #1\dag}
{\cal #1\smalldag}
}
\expandafter\edef\csname cal#1daginf\endcsname
{\noexpand\protect\csname cal#1daginfchoose\endcsname}
\expandafter\def\csname cal#1daginfchoose\endcsname{
\mathchoice
{\cal #1\daginf}
{\cal #1\daginf}
{\cal #1\daginf}
{\cal #1\smalldaginf}
}}

\makeletterdaginf X
\makeletterdaginf Y
\makeletterdaginf W
\makeletterdaginf M
\makeletterdaginf N
\makeletterdaginf V
\makeletterdaginf I
\makeletterdaginf U
\makeletterdaginf D
\makeletterdaginf Z
\makeletterdaginf T

\def\calDdag{\cal D\dag}
\def\calD{\cal D}

\let\mycal\mathcal
\def\calfragile#1{\let\exitcal\relax
            \if#1H{\myrsfs H}\else
             \if#1D{\myrsfs D}\else
              \if#1X{\myrsfs X}\else
               \if#1Y{\myrsfs Y}\else
                \if#1T{\myrsfs T}\else
\def\exitcal{\mycal#1}\fi\fi\fi\fi\fi\exitcal}

\def\cal{\protect\calfragile}

\def\Spec{\mathop{\rm Spec}\nolimits}

\def\newdisplayline#1#2#3{\hskip#1\postdisplaypenalty10000
$$$$\hskip#2\abovedisplayshortskip#3}

\let\get\leftarrow
\def\set#1/{\{#1\}}

\def\goodSup#1#2#3#4{\hbox to#1{\vbox{\vss\hbox{\kern#3$\scriptstyle#4$}\kern#2}\hss}}
\def\goodSub#1#2#3#4{\hbox to#1{\vtop{\kern#2\hbox{\kern#3$\scriptstyle#4$}\vss}\hss}}
\def\goodsub#1#2#3#4{\hbox to#1{\vtop to0pt{\kern#2\hbox{\kern#3$\scriptstyle#4$}\vss}\hss}}

\def\sub#1#2{\goodsub{#1}{4pt}{-4mm}{#2}}
\def\Sub#1#2{\goodSub{#1}{4pt}{-4mm}{#2}}

\def\SubX{\goodSub{10pt}{4pt}{-4mm}}
\def\SubV{\goodSub{5pt}{4pt}{-4mm}}
\def\SubO{\goodSub{0pt}{4pt}{-4mm}}

\newbox\myleqnobox
\newdimen\ecartlignes\ecartlignes=2pt
\def\deuxlignes#1#2#3\\#4\endlignes{{\global\setbox\myleqnobox=\hbox{}\def\leqno##1{\global\setbox\myleqnobox=\hbox{$##1$ }}%
\advance\hsize by-\leftmargin%
\setbox100=\hbox{$\displaystyle#3{}$}
\setbox101=\hbox to\hsize{\hss$\displaystyle{}#4$\kern#2}
\setbox100=\hbox to\hsize{\kern#1\box\myleqnobox\box100\hss}
\vtop{
\hbox{\hskip\leftmargin\box100}
\kern\ecartlignes
\hbox{\hskip\leftmargin\box101}
}}}
\def\DeuxLignes{\deuxlignes{0pt}{0pt}}
\def\troislignes#1#2#3\\#4\\#5\endlignes{{\global\setbox\myleqnobox=\hbox{}\def\leqno##1{\global\setbox\myleqnobox=\hbox{$##1$ }}%
\advance\hsize by-\leftmargin%
\setbox100=\hbox{$\displaystyle#3{}$}
\setbox101=\hbox to\hsize{\hss\kern#1$\displaystyle{}#4{}$\hss}
\setbox102=\hbox to\hsize{\hss$\displaystyle{}#5$\kern#2}
\setbox100=\hbox to\hsize{\kern#1\box\myleqnobox\box100\hss}
\vtop{
\hbox{\hskip\leftmargin\box100}
\kern\ecartlignes
\hbox{\hskip\leftmargin\box101}
\kern\ecartlignes
\hbox{\hskip\leftmargin\box102}
}}}
\def\TroisLignes{\troislignes{0pt}{0pt}}
\def\myeqno#1#2:{\hbox to#1{$#2\string:$\hss}}
\def\varmyeqno#1#2#3:{\vtop to0pt{\kern#2\hbox to#1{$#3\string:$\hss}\vss}}

\def\*{\mathord{*}}
\newdimen\pathunit
\newdimen\pathwd
\newcount\xpath
\newcount\ypath
\newcount\tpath
\newcount\cpath
\def\addendpath{}

\def\hline#1{\tpath=\xpath\advance\tpath#1 \drawhline\xpath\tpath\ypath\xpath\tpath}
\def\vline#1{\tpath=\ypath\advance\tpath#1 \drawvline\xpath\ypath\tpath\ypath\tpath}

\def\drawhline#1#2#3{\dimen0=#3\pathunit \dimen1=#1\pathunit \dimen2=#2\pathunit
\ifdim \dimen1>\dimen2 \dimen1=\dimen2 \dimen2=#1\pathunit\fi
\advance\dimen2 by -\dimen1\relax
\hbox to0pt{\kern\dimen1
          \vtop to0pt{\kern-\dimen0\kern-0.5\pathwd
                       \hrule height0.5\pathwd depth0.5\pathwd width\dimen2
                       \kern-0.5\pathwd\kern\dimen0
                      }\kern-\dimen2\kern-\dimen1}}

\def\drawvline#1#2#3{\dimen0=#1\pathunit \dimen1=#2\pathunit \dimen2=#3\pathunit
\ifdim \dimen1>\dimen2 \dimen1=\dimen2 \dimen2=#2\pathunit\fi
\advance\dimen2 by -\dimen1\relax
\hbox to0pt{\kern\dimen0\kern-0.5\pathwd
          \vtop to0pt{\kern-\dimen2\kern-\dimen1
                      \kern-0.5\pathwd
                      \advance\dimen2 by\pathwd
                      \hrule height\dimen2 depth0pt width\pathwd
                      \kern-0.5\pathwd
                      \kern\dimen1
                      }\kern-0.5\pathwd\kern-\dimen0}}

\def\putendpath#1{\raise-2.5pt\hbox to0pt{\hss\vtop to0pt
{\kern-\ypath\pathunit\smash{\kern\xpath\pathunit#1\kern-\xpath\pathunit}\kern\ypath\pathunit}\hss}}

\def\hvpath#1{\raise2.45pt\hbox{\cpath0 \xpath0 \ypath0 \tpath0 \pathspool#1 {} 
\ifx\addendpath\rien\else\putendpath\addendpath\global\let\endpath\rien\fi}}

\def\vhpath#1{\raise2pt\hbox{\cpath1 \xpath0 \ypath0 \tpath0 \pathspool#1 {} 
\ifx\addendpath\rien\else\putendpath\addendpath\global\let\endpath\rien\fi}}

\def\rien{}

\def\pathspool#1 {\def\donne{#1}\ifx\donne\rien\let\nextspoolpath\relax\else
\ifodd\cpath\vline{#1}\else\hline{#1}\fi\advance\cpath1
\let\nextspoolpath\pathspool\fi
\nextspoolpath}

\def\decale#1#2{\kern-#1#2\kern#1}

\def\mathrigid#1 {\thinmuskip=#1
\medmuskip=#1
\thickmuskip=#1 }

\let\emph\em
\def\em#1{{\emph#1\/}}

\let\Rm\relax
\newenvironment{liste}{
\begin{enumerate}\let\Item\item
\def\item##1){\Item[\rm##1)]}}{\end{enumerate}}

\def\paragraphe{\par\bigskip\noindent\refstepcounter{subsection}{\bf \thesubsection. }\ignorespaces}
\def\sousparagraphe{\par\medskip\noindent\refstepcounter{subsubsection}{\bf \thesubsubsection. }\ignorespaces}

{\catcode`\@=11
\gdef\matrix#1{\null\,\vcenter{\normalbaselines\m@th
    \ialign{\hfil$##$\hfil&&\quad\hfil$##$\hfil\crcr
      \mathstrut\crcr\noalign{\kern-\baselineskip}
      #1\crcr\mathstrut\crcr\noalign{\kern-\baselineskip}}}\,}

\gdef\closeoutglossary{\immediate\closeout\@glossaryfile}
}

\def\resetdisplay{\displayindent=\leftmargin
\advance\displaywidth-\leftmargin}

\def\textfontiii{\textfont3=}

\def\displaysize#1{\par\begingroup\abovedisplayskip-0\parskip\csname #1\endcsname\rm\scriptwd=2em
\setbox0=\hbox{$$}\relax
\expandafter\textfontiii\csname #1cmex\endcsname
\setbox0=\hbox{$$}\relax}

\long\def\inhibe#1\endinhibe{\relax}
\let\Bbb\mathbb

\let\mytimes\times
\def\dtimes{\mathop{\smash{\mytimes}}\limits}

\def\goodtimesk#1#2{\mkern#1\mathchoice
{\dotimesk\textstyle\scriptstyle}
{\dotimesk\textstyle\scriptstyle}
{\dotimesk\scriptstyle\scriptscriptstyle}
{\dotimesk\scriptscriptstyle\scriptscriptstyle}\mkern#2}

\def\timesk{\goodtimesk{4mu}{5mu}}
\def\dotimesk#1#2{\vtop{\offinterlineskip
\hbox to0pt{\hss$#1\times$\hss}\kern1pt
\hbox to0pt{\hss$#2k$\hss}}}

\catcode `\@=11
\def\x@donoparitem{\@noparitemfalse
   \global\setbox\@labels\hbox{\hskip -\leftmargin
                               \unhbox\@labels 
                                \hskip \labelsep\hskip1em}\if@minipage\else
  \@tempskipa\lastskip
  \vskip -\lastskip \advance\@tempskipa\@outerparskip
  \advance\@tempskipa -\parskip \vskip\@tempskipa\fi}
  
  \def\@setaltabstracta{%
\ifx\@empty\@alttitle\else
\goodbreak\vbox{\begin{center}%
\large\slshape  \@alttitle  
\end{center}}\nobreak
\fi
  \ifvoid\altabstractbox
  \else
     \skip@5\p@ \advance\skip@-\lastskip
     \advance\skip@-\baselineskip 
     \vskip\skip@\nobreak
    \box\altabstractbox
    \prevdepth\z@ 
  \fi
}

\def\rscript#1{\mkern0mu\rlap{$\mathsurround0pt\scriptstyle#1$}}
\def\lscript#1{\llap{$\mathsurround0pt\scriptstyle#1$}\mkern0mu}

\def\hmsmash#1{\hbox to0pt{\hss$#1$\hss}}
\def\hdecale#1#2{\hskip#1#2\hskip-#1}
\def\sto{\raise-1.5pt\hbox{$\scriptstyle\vec{\kern2pt}$\kern2pt}}

\def\Inf{\mathop{\rm Inf}\nolimits}

\DeclareFontFamily{T1}{smallsy}{}
\DeclareFontShape{T1}{smallsy}{m}{n}{
  <->
  s*[0.7]cmsy10
}{}

\DeclareMathAlphabet {\smallsy} {T1} {smallsy} {m} {n}
\def\smalldagger{\smallsy{\char"79}}

\DeclareFontFamily{T1}{myrsfs}{}
\DeclareFontShape{T1}{myrsfs}{m}{n}{
  <->
  rsfs10
}{}

\DeclareMathAlphabet {\myrsfs} {T1} {myrsfs} {m} {n}

\def\comb#1#2{\big(\,{\vrule height10pt depth5pt width0pt}^{#1}_{#2}\,\big)\;}

\let\mBig\big

\NumberTheoremsIn{subsection}%

\equalenv{notation}{nota}
\equalenv{notations}{notas}
\equalenv{exemple}{exem}
\equalenv{exemples}{exems}
\equalenv{Rema}{rema}
\equalenv{Remas}{remas}

\title
[Topos infinit\'esimal $p$-adique]
{\boldmath Sur le Topos infinit\'esimal $p$-adique\\ d'un sch\'ema lisse I $(\sp{*})$}

\alttitle{Infinit\'esimal $p$-adic topos of a smooth scheme I}

\author{\firstname{Alberto} \lastname{Arabia}}
\address{CNRS et Université Paris 7 - Denis Diderot\\ 
175, rue de Chevaleret\\Paris, 75013}
\email{arabia@math.jussieu.fr}

\catcode`\@=11 \def\@@and{et}
\author{\firstname{Zoghman} \lastname{Mebkhout}}
\address{CNRS et Université Paris 7 - Denis Diderot\\ 
175, rue de Chevaleret\\Paris, 75013}

\email{mebkhout@math.jussieu.fr}

\keywords{
algèbres dag-adiques,
cohomologie de de Rham $p$-adique,
complexe de de Rham $p$-adique,
équations différentielles $p$-adiques,
factorisation $p$-adique de la fonction Zéta,
fonctorialité,
groupe des automorphismes,
module de transfert,
module spécial,
opérateurs différentiels $p$-adiques,
opérations cohomologiques,
relèvements plats,
schémas dag-adiques,
site infinitésimal,
suite de Gysin,
topos infinitésimal.
}
  
\altkeywords{}

\subjclass{
11E95, 12H25, 13Dxx, 13Fxx, 13Jxx, 13N10, 14Axx, 14Fxx, 16Exx, 18Fxx, 18Gxx.}

\dedicatory{dédié à Alexandre Grothendieck\\ à l'occasion de son quatre-vingtième anniversaire}

\makeglossary

\begin{document}


\begin{abstract}\pretolerance800 Afin de disposer des opérations cohomologiques aussi  	\\souples que possible pour la cohomologie de de Rham $p$-adique à coefficients, le but principal de ce mémoire est de résoudre
intrinsèquement \em{du point de vue cohomologique} le problème des \em{relèvements des schémas lisses et de leurs morphismes de la caractéristique $p>0$ à la
caractéristique nulle} ce qui a été l'une des difficultés centrales 
de la théorie de la cohomologie de de Rham des schémas algébriques en caractéristique positive depuis le début. Nous montrons que, bien que
les schémas lisses et leurs morphismes ne se relèvent pas en général du point de vue géométrique, tout se passe comme si c'était bien le cas du
point de vue cohomologique, ce qui est conforme à la Théorie des Motifs de Grothendieck. On en déduit la factorisation $p$-adique de la fonction
Zêta d'une variété algébrique lisse sur un corps fini, éventuellement ouverte, qui est le résultat test de nos méthodes.

Soit $V$ un anneau de valuation discrète complet d'inégales caractéristiques $(p, 0)$ de corps résiduel $k$ et de corps de fractions $K$.
On définit la cohomologie de de Rham $p$-adique d'un schéma lisse sur $k$, à coefficients qui sont des espaces vectoriels sur $K$
et on définit les opérations cohomologiques pour un morphisme de schémas lisses sur $k$. On montre que l'on obtient en particulier un foncteur
contravariant entre la catégorie de tous les schémas lisses et séparés sur $k$ et 
la catégorie dérivée de la catégorie des espaces vectoriels sur $K$. 
On montre la suite exacte de Gysin pour tout couple de schémas 
lisses, ce qui permet en particulier de définir la classe de cohomologie d'un cycle dans le cas d'un corps de base parfait. On montre le lemme de Poincaré-Künneth sur une base lisse.
\end{abstract}

\begin{altabstract}
In order to have cohomological operations for de Rham $p$-adic cohomology with coefficients 
as manageable as possible, the main purpose of this paper is to solve intrinsically
and \em {from a cohomological point of view the lifting problem 
of smooth schemes and their morphisms from characteristic $p>0$ to 
characteristic zero} which has been one of the fundamental difficulties 
in the theory of de Rham cohomology of algebraic schemes in positive characteristic since the beginning. We show that although 
smooth schemes and morphisms fail to lift geometrically, 
it is as if this was the case within the cohomological point of view, which is consistent with the theory of Grothendieck Motives. We deduce the $p$-adic factorization of the 
Zeta function of a smooth algebraic variety, possibly open, over a finite field, which is a key testing result of our methods. 

Let $V$ be a complete discrete valuation ring of unequal characteristics $(p,0)$ with residue field $k$ and fraction field $K$. 
We define the $p$-adic de Rham cohomology of  smooth scheme over $k$ with coefficients which are vector spaces over $K$,
and we define cohomological operations for morphisms of smooth schemes over $k$. We show that we obtain in particular a contravariant functor between the category of all separated smooth schemes over $k$ and 
the derived category of the category of vector spaces over $K$. We give the Gysin exact sequence for every pair of smooth schemes, allowing in particular to define the cohomology class of a cycle in the case of a perfect base field. We prove the Poincaré-Künneth lemma for smooth base.
\end{altabstract}

\begingroup
\catcode`\@=11
\long\def\@footnotetext#1{\insert\footins{%
    \reset@font\footnotesize
    \interlinepenalty\interfootnotelinepenalty
    \splittopskip \footnotesep 
    \splitmaxdepth \dp\strutbox \floatingpenalty \@MM
    \hsize\columnwidth \@parboxrestore
\hskip 2em (*) #1\@finalstrut\strutbox}}

\footnotetext{\`A para\^\i tre dans les Annales de l'Institut Fourier nro. {\bf 60} fascicule 6.
Manuscrit soumis le 10 octobre 2007, révisé le 18 juin 2009, accepté le 21 juillet 2009.}
\endgroup

\null\vskip-1cm\maketitle

\vspace*{-0.5cm}
\def\contentsname{\large\bf Table des matières}

{\footnotesize
\baselineskip=1\baselineskip plus1pt minus 1pt
\setcounter{tocdepth}{2}

\tableofcontents

}
  
\baselineskip=1\baselineskip plus0.5pt minus 0.5pt
\section{Introduction}

 Cet article est le premier d'une série d'articles qui s'inscrivent dans la recherche d'une théorie cohomologique des coefficients $p$-adiques pour la
catégorie des variétés algébriques en dimensions supérieures. Cette théorie est loin  d'avoir atteint son but malgré d'intenses efforts; seul le cas des courbes est bien compris, grâce à la théorie des équations différentielles
$p$-adiques ([C-M$_1$], [C-M$_2$], [C-M$_3$], [C-M$_4$]). 

\paragraphe
Afin de disposer des opérations cohomologiques aussi souples que possible, le but principal de ce mémoire est de résoudre
intrinsèquement \em{du point de vue cohomologique le problème des relèvements des schémas  lisses et  de leurs  morphismes de la caractéristique $p>0$ à la
caractéristique nulle} ce qui a été l'une des difficultés centrales de la théorie depuis qu'elle existe. Nous allons montrer que, bien que
les  schémas  lisses et leurs  morphismes ne se relèvent pas en général du point de vue géométrique, tout se passe comme si c'était bien le cas du
point de vue cohomologique ce qui est conforme à la Théorie des Motifs de Grothendieck. 

\bigskip
Nous rappelons rapidement pour la commodité du lecteur la situation de la théorie  pour situer ce mémoire.

\paragraphe
Le point de départ de la théorie de de Rham $p$-adique est la démonstration de B. Dwork [D] de la rationalité de la fonction Zêta d'une variété algébrique sur un
corps fini, par des méthodes $p$-adiques non cohomologiques. P.~Monsky et G.~Washnitzer ont obtenu ([M-W$_0$], [M-W], [M$_1$], [M$_2$]) une expression cohomologique de
cette fonction pour les variétés algébriques affines lisses par voie $p$-adique, parallèle à celle obtenue par A. Grothendieck par voie
$\ell$-adique  pour toutes les variétés algébriques sur les corps finis ([G$_2$]), en construisant une théorie cohomologique de type de de Rham ([M-W]) en passant   de la caractéristique
$p>0$ à la caractéristique zéro à l'aide de la théorie des schémas formels affines de Grothendieck ([EGA I]).
En retour, leurs travaux  ont inspiré à Grothendieck la construction  du site Infinitésimal et du site Cristallin d'un schéma au-dessus d'un autre ([G$_3$], [G$_4$]).
De plus, Grothendieck  a montré que pour une variété algébrique lisse sur un corps de caractéristique nulle la cohomologie de son  site infinitésimal à valeurs dans le
faisceau structural est canoniquement isomorphe à sa cohomologie de de Rham ([G$_4$]), ce qui constitue une idée très très originale  dans le calcul différentiel algébrique.

\paragraphe
P. Monsky et G. Washnitzer ont laissé ouvert le problème de l'extension du 
foncteur de de Rham $p$-adique à la catégorie des 
variétés algébriques éventuellement non affines. Cette problématique  a déjà été discutée par  Monsky dans son exposé au congrès de Nice de 1970 ([M$_3$]).
P. Berthelot a proposé une théorie cohomologique pour toutes les variétés algébriques via la géométrie analytique rigide ([B$_3$]). Cette  construction est de nature
globale et fait intervenir une compactification, ce qui présente de sérieux inconvénients à ce stade; par exemple, il n'est pas du tout clair que cette définition soit fonctorielle. 
\footnote{Note du rapporteur: ``{\sl In order to get the functoriality, another construction is necessary. There is some recent, unpublished work of Le Stum in this direction.}''}

\paragraphe
Dans l'article fondamental de recherche [M-N$_1$], qui a servi de transition, on  a proposé de revenir au point de vue originel, de nature {\bf locale et algébrique pour la topologie de Zariski en caractéristique $p>0$}, en définissant la cohomologie de de
Rham
$p$-adique pour les variétés algébriques lisses
admettant un relèvement plat comme schéma formel faible au sens de Meredith ([Mr]) via la théorie des $\calDdag $-modules sur un tel schéma ([M-N$_1$],
[M-N$_2$]).  Cette définition   considérée sur le corps des
fractions  d'un anneau de valuation discrète complet d'inégales caractéristiques coïncide canoniquement avec  la définition  de Monsky-Washnitzer pour les variétés 
{\bf affines} lisses. 

\paragraphe
Mais il
n'était pas établi que ces espaces sont indépendants du relèvement pour les variétés {\bf non affines} ni que cette définition  se prolonge à  la catégorie de toutes
les  variétés algébriques lisses sur un corps de caractéristique $p>0$. Cependant, 
ce point de vue, qui s'est avéré être le bon,  nous a conduits à  la théorie des équations différentielles $p$-adiques, qui a
prise sur les questions de la monodromie  arithmétique et  dépasse de loin la seule question  de la cohomologie des variétés algébriques ([C-M$_1$],
[C-M$_2$], [C-M$_3$], [C-M$_4$]).

\paragraphe
Le but de cet article est,  d'une part, de montrer qu'effectivement ces espaces sont indépendants à isomorphisme \em{\bf canonique} près du relèvement et que
leur définition se prolonge  à la catégorie de {\bf tous}  les schémas lisses et, d'autre part, de définir  {\bf les opérations
cohomologiques pour la cohomologie de de Rham
$p$-adique à coefficients} pour un morphisme de schémas  lisses sur un corps de caractéristique positive, mais sans {\bf aucune
hypothèse  restrictive} sur les morphismes pour en déduire en particulier la {\bf fonctorialité} de la cohomologie de de Rham $p$-adique et l'expression cohomologique $p$-adique de la fonction Zêta d'une variété algébrique lisse sur un corps fini. 

\paragraphe
La théorie de de Rham en caractéristique nulle et ses 
nombreuses applications nous enseignent
que   la situation précédente d'un morphisme qui n'est soumis à aucune restriction, de lissité, de platitude  ou de propreté, entre schémas lisses sur le corps de base sur lesquels on  dispose du calcul différentiel et du calcul infinitésimal sur le corps des nombres complexes,  est la situation  essentielle         même  pour le cas
singulier et donne lieu aux problèmes les plus intéressants. Le cas des courbes, sur lequel nous disposons aujourd'hui des résultats les plus complets grâce précisément à la monodromie $p$-adique, montre déjà que la théorie de de Rham arithmétique est hautement non triviale comparée aux autres théories cohomologiques existantes, ce qui explique en partie son retard. Il était évident qu'on ne pouvait pas espérer un progrès significatif dans la théorie des coefficients $p$-adiques sans 
la structure complète de la monodromie locale d'un point singulier d'une équation différentielle $p$-adique d'ordre supérieure ([C-M$_3$], Introduction).

\paragraphe
Nous allons décrire les différentes étapes qui  nous ont conduits à ce résultat. Soit $V$ un anneau de valuation discrète complet d'inégales
caractéristiques
$(p, 0)$, d'idéal maximal $\goth m$, de corps résiduel
$k$ et corps de fractions $K$. Soit $X/k$ un schéma lisse sur $k$, on définit naturellement son site infinitésimal $p$-adique $\Xdaginf $ 
comme la catégorie, munie de sa topologie de Zariski naturelle, dont les objets sont  les schémas  $\dagger$-adiques  $\calUdag :=(U, \cal O_{\calUdag /V})$, au sens de Meredith ([Mr]), {\bf
plats} sur
$V$, où
$U$ est un ouvert 
de l'espace topologique $X$. L'hypothèse de lissité sur $k$ assure que le site $ \Xdaginf $ a beaucoup d'ouverts et de morphismes. Le site $ \Xdaginf $ est naturellement annelé par le
faisceau de
$V$-algèbres
$ \cal O_{\Xdaginf /V}$ et le faisceau de
$K$-algèbres 
$ \cal O_{\Xdaginf /K}:= \cal O_{\Xdaginf /V}\otimes_VK$. La cohomologie
infinitésimale $ H^{\bullet}_{\inf}(X/K, \cal O_{\Xdaginf /K})$ du site $ \Xdaginf $ à valeurs
dans le faisceau structural
$ \cal O_{\Xdaginf /K}$ est un invariant intrinsèque du schéma  $X$ et une question à laquelle nous avons cru pendant longtemps est qu'elle
fournisse les bons nombres de Betti $p$-adiques, ce qui a retardé d'autant notre compréhension du problème. En fait il n'en est rien, contrastant vivement avec le cas de la caractéristique nulle. Pour  voir ceci,
supposons pour simplifier que
$X$ est une variété affine lisse sur le corps $k$ et soit  $A\dag$ une $V$-algèbre faiblement complète de type fini 
au sens de Monsky-Washnitzer ([M-W]) plate sur $V$ qui relève l'algèbre des fonctions de $X$ et $\calXdag := (X,\cal O_{\calXdag /V})$ son schéma
$\dagger$-adique  affine  associé ([Mr]). On peut considérer le groupe 
$G_{A\dag}$ des {\bf automorphismes} de $V$-algèbre de $A\dag$ qui se réduisent à l'identité modulo $\goth m$. Ce groupe se localise en un faisceau de
groupes
$\cal G_{\calXdag }$ sur $X$. Nous montrons dans le théorème \ref{coh-inf} que la cohomologie infinitésimale $ H^{\bullet}_{\inf}(X/K, \cal O_{\Xdaginf /K})$ s'identifie canoniquement à la cohomologie
de groupe
$ H^{\bullet}(\cal G_{\calXdag },\cal O_{\calXdag /K})$ du $\cal G_{\calXdag }$-module $\cal O_{\calXdag /K}$ qui est donc de nature galoisienne, ce qui donne un moyen de 
calcul. Or, nous construisons dans le théorème \ref{lem-poi}  pour l'espace affine $X=\Spec(k[x_1,\dots,x_n]),$  un cocycle qui n'est pas un cobord et donc
une classe de cohomologie non triviale  de
$ H^1(\cal G_{\calXdag },\cal O_{\calXdag /K})$,  montrant que le lemme de Poincaré n'a pas lieu pour cette cohomologie. La cohomologie infinitésimale $p$-adique
d'une variété affine lisse à valeurs dans le faisceau structural \em{ne coïncide pas} avec sa cohomologie de de Rham
$p$-adique et son invariance par rapport au relèvement n'est d'aucun secours pour le problème qui nous concerne ici.

\paragraphe
Nous revenons alors au point de vue introduit dans l'article  de recherche [M-N$_1$]. Soient $R$ un anneau commutatif noethérien, $I$ un idéal  et  $R_1:=R/I$. Dans
les articles [M-N$_1$] et [M-N$_2$], on a défini le faisceau
$\calDdag_{\calXdag /R}$ des opérateurs différentiels sur un schéma 
$\dagger$-adique $\calXdag  =(X,
\cal O_{\calXdag /R})$ par réduction modulo les puissances de l'idéal $I$ en imposant une condition de
croissance du type $\dagger$ sur l'ordre des opérateurs différentiels. Par construction, le faisceau
$\calDdag_{\calXdag /R}$ opère à gauche sur le faisceau $\cal O_{\calXdag /R}$ et on peut donc considérer les modules  de cohomologie de de Rham
$\dagger$-adiques
$\ext^{\bullet}_{\calDdag_{\calXdag /R}}(\cal O_{\calXdag /R}, \cal O_{\calXdag /R})$ définis comme foncteurs dérivés du foncteur covariant
exact à gauche: 
$$\cal M\dag_{\calXdag }\fonct \hom_{\calDdag_{\calXdag /R}}(\cal O_{\calXdag /R}, \cal M\dag_{\calXdag })$$ dans la catégorie
 des $\calDdag_{\calXdag /R}$-modules à gauche $ \calDdag_{\calXdag /R}\Mod$ et la cohomologie globale apparaît naturellement comme
l'hypercohomologie {\bf d'un complexe de Zariski}.

\paragraphe
En particulier, si $R$ est l'anneau $V$ on peut considérer  le faisceau $\calDdag_{\calXdag /V}$, le faisceau $\calDdag_{\calXdag /K}:= \calDdag_{\calXdag /R}\otimes_VK$ et les $K$-espaces vectoriels $\ext^{\bullet}_{\calDdag_{\calXdag /K}}(\cal O_{\calXdag /K}, \cal O_{\calXdag /K})$. Il  est 
montré dans [M-N$_1$]  que ces derniers espaces vectoriels 
coïncident canoniquement avec les espaces de
cohomologie  définis par Monsky-Washnitzer dans le cas d'un relèvement 
plat d'une variété \em{affine} lisse et en vertu de leur théorème ([M-W], Thm. 5.6)
ils sont donc indépendants du relèvement à isomorphisme canonique près. On disposait dès le départ de nombreuses indications notamment grâce aux foncteurs de cohomologie locale pour espérer que ces espaces soient indépendants à isomorphisme
canonique près du relèvement d'une variété algébrique lisse sur $k$ quand il existe.

\paragraphe
Sur un schéma $\dagger$-adique $\calXdag =(X,\cal O_{\calXdag /R})$  on peut considérer le faisceau  $\cal G_{\calXdag }$ des 
endomorphismes de $R$-{\bf algèbres}
du 
faisceau structural 
$\cal O_{\calXdag /R}$ qui se réduisent à l'identité modulo $I$. Il se trouve de façon tout à fait remarquable (Thm. \ref{inc-grp}) 
que le faisceau $\cal G_{\calXdag }$ est un sous-faisceau de semi-groupes  pour la structure {\bf multiplicative} du faisceau d'anneaux $\calDdag_{\calXdag/R}$, autrement dit un élément $g$ du semi-groupe précédent est un opérateur \em{différentiel d'un type très particulier}. Si le faisceau
structural $\cal O_{\calXdag /R}$ est $R$-plat, le faisceau $\cal G_{\calXdag }$ est un faisceau de groupes.  C'est là le lien essentiel
entre la géométrie algébrique en caractéristique $p>0$ et le calcul différentiel en caractéristique zéro. Par exemple, c'est ce plongement qui nous 
a permis de construire une classe de cohomologie non triviale dans la  cohomologie infinitésimale $p$-adique de l'espace affine à valeurs dans le faisceau structural.

\paragraphe
Sur le site infinitésimal $ \Xdaginf $ d'un $R_1$-schéma lisse, défini comme dans 1.5, et qui a aussi beaucoup d'ouverts et de morphismes,    le faisceau
de groupes $ \cal G_{\Xdaginf }$  des automorphismes du faisceau structural qui se réduisent à l'identité modulo $I$ est aussi défini,  ainsi que le  faisceau  d'algèbres $ \calDdag_{\Xdaginf /R}$ des opérateurs
différentiels, et l'on a  une inclusion $ \cal G_{\Xdaginf }\hookrightarrow\calDdag_{\Xdaginf /R}$ pour la structure multiplicative. Un faisceau de
$R_{\Xdaginf }$-modules sur le site $ \Xdaginf $ est muni par construction   d'une action \em{\bf géométrique à gauche} $(\sharp)$ du faisceau de
groupes $\Rm
\cal G_{\Xdaginf }$ (Prop. \ref{act-gro}). En particulier, le faisceau
$ \cal G_{\Xdaginf }$ opère sur le faisceau  $ \calDdag_{\Xdaginf /R}$   par automorphismes intérieurs (Coro. \ref{act-dif}). 
On dispose  de la catégorie $ \calDdag_{\Xdaginf /R}\Mod$ \glossary{$ \calDdag_{\Xdaginf /R}\Mod$} des $ \calDdag_{\Xdaginf /R}$-modules à gauche sur le site infinitésimal. Par construction,  un  
$ \calDdag_{\Xdaginf /R}$-module à gauche sur le site $\calMdaginf $ est muni  de {\bf deux actions} à gauche du faisceau de groupes  $ \cal G_{\Xdaginf }$,  
l'action géométrique $(\sharp)$  et l'action  {\bf différentielle à gauche} $(\diamond)$ qui se fait  
à travers $ \calDdag_{\Xdaginf /R}$.
Nous arrivons à la notion essentielle et centrale de {\bf module spécial}.

\begin{defi*}[\ref{def-spe}]Nous dirons qu'un $ \calDdag_{\Xdaginf /R}$-module à gauche $\calMdaginf $ sur le site  infinitésimal $\dagger$-adique $ \Xdaginf $ est \em{\bf spécial}  si
l'action géométrique $(\sharp)$ du faisceau de groupes 
$ \cal G_{\Xdaginf }$ sur $\calMdaginf $ {\bf se fait à travers}  $ \calDdag_{\Xdaginf /R}$.
\end{defi*} 
Autrement dit, sur un $ \calDdag_{\Xdaginf /R}$-module à gauche {\bf spécial} les deux actions à gauche du faisceau de groupes 
$ \cal G_{\Xdaginf }$,
géométrique $(\sharp)$ et différentielle~$(\diamond)$,  {\bf coïncident}. Ainsi, le $ \calDdag_{\Xdaginf /R}$-module à gauche $ \cal O_{\Xdaginf /R}$ \em{est spécial} par
construction,  mais  le
$\Rm
\calDdag_{\Xdaginf /R}$-module à gauche $ \calDdag_{\Xdaginf /R}$ \em{n'est pas spécial} parce que l'action géométrique $(\sharp)$ est  $(g,P)\mapsto gPg^{-1}$
et que l'action différentielle $(\diamond)$ est  $(g, P)\mapsto gP$.

On note $ (\calDdag_{\Xdaginf /R}, \Sp)\Mod$\glossary{$ (\calDdag_{\Xdaginf /R}, \Sp)\Mod$} la catégorie des $ \calDdag_{\Xdaginf /R}$-modules à gauche spéciaux sur le site $\Xdaginf$.
Nous démontrons le théorème suivant, qui légitime cette
notion.

\begin{theo*}[\ref{can}]Soit $\calXdag $ un relèvement $\dagger$-adique 
plat d'un $R_1$-schéma lisse $X$. Alors, il existe
un foncteur \em{\bf canonique de prolongement}  :
$$  \calDdag_{\calXdag /R}\Mod\to (\calDdag_{\Xdaginf /R}, \Sp)\Mod\,,\leqno P_{\calXdag }:$$ qui est une équivalence de catégories et 
un inverse {\bf canonique} au foncteur image inverse  naturel de restriction
$ R_{\calXdag }$.
\end{theo*} 
Ainsi, un $\calDdag_{\calXdag /R}$-module à gauche se \em{propage canoniquement} en un $ \calDdag_{\Xdaginf /R}$-module à 
gauche \em{spécial}  et, réciproquement, un module $ \calDdag_{\Xdaginf /R}$
à gauche \em{spécial} se restreint canoniquement en un $\calDdag_{\calXdag /R}$-module à gauche. Dans sa lettre [G$_3$] Grothendieck écrivait :
\og \em{Un Cristal possède deux propriétés caractéristiques: la rigidité, et la faculté de croître dans un voisinage approprié}\fg. 
 L'isomorphisme  canonique $(\diamond)$  de restriction du   théorème \ref{iso-can} des modules spéciaux correspond à la rigidité des cristaux et la propagation sur le site infinitésimal du théorème  \ref{can} des modules spéciaux correspond à la croissance des cristaux dans un voisinage approprié. Mais alors que la catégorie des cristaux ne semble fournir une bonne théorie que pour les variétés algébriques  propres et lisses sur $k$, la catégorie des
modules spéciaux semble fournir une bonne théorie pour les variétés {\bf ouvertes}, comme l'illustre le théorème 
\ref{fac-zet} de l'expression cohomologique $p$-adique de la fonction Zêta d'une variété algébrique lisse sur un corps fini. 
\begin{coro*}[\ref{can-cor}]La catégorie de
$\calDdag_{\calXdag /R}$-modules à gauche ne dépend pas à équivalence \em{\bf canonique} près du relèvement plat $\calXdag $   du schéma $X$ lisse sur $R_1$ quand il existe. \end{coro*}Ce théorème 
résout pour les besoins de la théorie des coefficients $p$-adiques le problème des relèvements des \em{schémas lisses}  parce que la catégorie des $\calDdag_{\Xdaginf /R}$-modules à gauche {\bf spéciaux} $ (\calDdag_{\Xdaginf /R}, \Sp)\Mod$ est définie pour tout $R_1$-schéma $X$ lisse  et qu'elle en est un
{\bf invariant intrinsèque}. 

\paragraphe
Nous montrons que la catégorie des modules à gauche spéciaux $ (\calDdag_{\Xdaginf /R},\Sp)\Mod$ est une sous-catégorie pleine, abélienne de la catégorie $ \calDdag_{\Xdaginf /R}\Mod$ des modules à gauche (Prop. \ref{cat-abe}), admet suffisamment d'objets injectifs (Thm. \ref{obj-inj})  d'objets plats (Thm. \ref{exi-pla})  et est un champ (Thm.  \ref{champ}) bien que ce n'est pas une catégorie de modules. Le théorème qui suit 
est corollaire de l'équivalence de catégories précédente.

\begin{theo*}[\ref{coh-deR}]Soit $\calXdag $ un relèvement $\dagger$-adique plat d'un $R_1$-schéma  $X$ lisse. Alors,
il existe des isomorphismes {\bf canoniques} de
$R$-modules:
$$ \ext^{\bullet}_{\calDdag_{\calXdag/R}}(\cal O_{\calXdag /R}, \cal O_{\calXdag /R})\simeq\ext^{\bullet}_{\calDdag_{\Xdaginf /R},\Sp}(\cal
O_{\Xdaginf /R}, \cal
O_{\Xdaginf /R}).$$
\end{theo*}
Les modules de droite sont définis comme les foncteurs dérivés dans la catégorie abélienne $ (\calDdag_{\Xdaginf /R}, \Sp)\Mod$ du foncteur
covariant exact à gauche:  $$ \calMdaginf \fonct \hom_{\calDdag_{ \Xdaginf /R},\Sp}(\cal O_{ \Xdaginf /R}, \calMdaginf )$$ des morphismes
$\calDdag_{ \Xdaginf /R}$-linéaires. Les modules précédents apparaissent donc comme les modules de cohomologie du complexe: $$ \DR(X/R,\calMdaginf ):= \bfR \hom_{\calDdag_{ \Xdaginf /R},\Sp}(\cal O_{ \Xdaginf /R}, \calMdaginf ).$$ L'isomorphisme du théorème précédent montre, d'une part, que la cohomologie de de Rham $\dagger$-adique d'un relèvement ne dépend
pas à isomorphisme \em{canonique} près du relèvement et, d'autre part, que la cohomologie de de Rham $\dagger$-adique  du site infinitésimal prolonge le
foncteur de de Rham sur $X$ à tous les schémas lisses, 
{\bf sans hypothèse d'existence de relèvement}. 

\paragraphe
La catégorie  $\RmMod(R_X)$ des faisceaux de $R$-modules sur $X$ pour la {\bf topologie de Zariski} apparaît naturellement comme une
sous-catégorie pleine de la catégorie $\RmMod(R_{\Xdaginf })$ des $ R_{\Xdaginf }$-modules sur le site infinitésimal $ \Xdaginf $, à savoir la catégorie des
faisceaux de  $ R_{\Xdaginf }$-module sur lesquels l'action géométrique $(\sharp)$ du groupe
$ \cal G_{\Xdaginf }$ est \em{triviale} (Thm. \ref{equ-fai}). Or par construction, si $\calNdaginf $ et $\calMdaginf $ sont deux $ \calDdag_{\Xdaginf /R}$-modules à
gauche {\bf spéciaux}, l'action géométrique du groupe  $ \cal G_{\Xdaginf }$ sur le faisceau (Prop. \ref{fai-hom}) de $ R_{\Xdaginf }$-modules 
$$\Rm
\cHom_{\calDdag_{\Xdaginf /R},\Sp}(\calNdaginf , \calMdaginf )= \cHom_{\calDdag_{\Xdaginf /R}}(\calNdaginf , \calMdaginf )$$   est {\bf triviale} (Coro. \ref{hom-spe}) et  ce faisceau est donc un faisceau pour la \em{topologie de
Zariski} de $X$. En particulier, le faisceau 
$ \cHom_{\calDdag_{\Xdaginf /R}}(\cal O_{\Xdaginf /R}, \cal
M\daginf )$ des sections horizontales de $\calMdaginf $ est un faisceau de Zariski. 
Il en résulte le foncteur exact de de Rham local de catégories triangulées
$$ \Dplus ((\calDdag_{\Xdaginf /R},\Sp)\Mod)\to  \Dplus (\RmMod(R_X))\leqno \dR(X/R,- ):$$ qui à $\calMdaginf $
associe son  complexe de de Rham $\dagger$-adique local: $$ \dR(X/R,\calMdaginf ):=  \bfR \cHom_{\calDdag_{\Xdaginf /R},\Sp}(\cal O_{\Xdaginf /R}, \calMdaginf ).$$  On déduit alors le caractère {\bf local pour la topologie de Zariski } de la définition de la cohomologie de de Rham $p$-adique.

\begin{theo*}[\ref{loc-deR}]La cohomologie de de
Rham
$\dagger$-adique d'un complexe spécial $\calMdaginf $ est  l'hypercohomologie du complexe:   $$ \bfR \cHom_{\calDdag_{\Xdaginf /R},\Sp}(\cal O_{\Xdaginf /R}, \calMdaginf )$$
\em{\bf de  Zariski}. 
\end{theo*}De plus, pour tout fermé
$Z$ de l'espace topologique $X$ la cohomologie de de Rham
$\dagger$-adique  à support dans $Z$ est définie comme d'habitude comme l'hypercohomologie du complexe de cohomologie locale à support dans $Z$. Cela
ramène de nombreuses questions concernant la cohomologie de de Rham
$p$-adique à des questions de nature locale pour la
\em{topologie de Zariski} qui est finalement le point de  vue classique et qui justifie s'il en est encore besoin celui des
schémas $\dagger$-adiques. En tout état de cause, cette localisation pour la topologie de Zariski apparaît comme un point essentiel dans notre théorie.

\paragraphe
Pour aller plus loin il nous faut considérer une situation un peu plus générale. Soit $R\to S$  une extension d'anneaux commutatifs ;
on peut considérer les faisceaux d'anneaux sur le site $\Xdaginf$ obtenus par changement de base:
 $$ \cal O_{\Xdaginf /S}:=\cal O_{\Xdaginf /R}\otimes_RS, \qquad \calDdag_{\Xdaginf /S}:=\calDdag_{\Xdaginf /R}\otimes_RS$$ sur lesquels  le faisceau de groupes $\cal G_{\Xdaginf}$ agit tout aussi
bien.  On a alors la catégorie des $ \calDdag_{\Xdaginf /S}$-modules à gauche $ \calDdag_{\Xdaginf /S}\Mod$ et sa sous-catégorie abélienne pleine $ (\calDdag_{\Xdaginf /S},\Sp)\Mod$ des modules à gauche spéciaux. Tous les résultats précédents restent valables
dans cette situation. Nous notons $ \rmD^*((\calDdag_{\Xdaginf /S},\Sp)\Mod)$
la catégorie dérivée de la catégorie abélienne des $\calDdag_{\Xdaginf /S}$-modules à gauche spéciaux  $ (\calDdag_{\Xdaginf /S},\Sp)\Mod$. Si $\calMdaginf $ un est complexe de la catégorie $ \rmD^{+}((\calDdag_{\Xdaginf /S},\Sp)\Mod)$, on définit \og la cohomologie de de Rham $\dagger$-adique de $X$ à coefficients dans $\calMdaginf $\fg, et nous notons pour simplifier
$\Rm
H^{\bullet}_{\DR}(X/S,\calMdaginf )$, la cohomologie du complexe $$\DR(X/S,\calMdaginf ):=\bfR \hom_{\calDdag_{ \Xdaginf /S},\Sp}(\cal O_{ \Xdaginf /S}, \calMdaginf )$$  qui consiste en 
des $S$-modules. On a donc par définition:
$$\Rm
H^{\bullet}_{\DR}(X/S, \calMdaginf ):=\ext^{\bullet}_{\calDdag_{\Xdaginf /S},\Sp}(\cal
O_{\Xdaginf /S}, \calMdaginf ).$$

En particulier, on peut considérer l'extension $V\to K$. Nous appelons alors \og cohomologie de de Rham $p$-adique de $X$\fg, que nous notons pour simplifier 
$\Rm
H^{\bullet}_{\DR}(X/K)$, la cohomologie $\Rm
H^{\bullet}_{\DR}(X/K, \cal O_{\Xdaginf/K})$. Nous appelons 
\og $i$-ème nombre  de Betti $p$-adique de $X$\fg, et nous le notons $B_{p,i}(X)$, la dimension du $K$-espace vectoriel $\Rm
H^{i}_{\DR}(X/K)$. La suite spectrale locale-globale (Thm. \ref{glo-loc}) ramène la finitude des nombres $B_{p,i}(X)$ d'une variété algébrique, schéma séparé de type fini sur $k$,  lisse sur $k$,  au théorème de finitude des nombres de Betti $p$-adiques des variétés algébriques affines lisses sur $k$ [Me$_2$], ce qui montre qu'on est enfin sur la bonne voie qui est en même temps le point de départ des opérations cohomologiques.

\paragraphe
C'est pour la catégorie $ (\calDdag_{\Xdaginf /V},\Sp)\Mod$ des modules spéciaux qu'on dispose  intrinsèquement des opérations cohomologiques pour un morphisme et nous
construisons dans cet article les trois opérations fondamentales et le foncteur de dualité $\DDD^\vee_{\Xdaginf /V}(-)$ (Def. \ref{def-dua}), à partir desquels l'on peut  construire toutes les autres. 

\begin{theo*}[\ref{com-chl}]Soient  $Z$  un fermé de l'espace topologique $X$ et $\calMdaginf $ un complexe borné à gauche de $\calDdag_{\Xdaginf /S}$-modules à gauche spéciaux. Alors, on a un triangle distingué de
la catégorie $ \Dplus ((\calDdag_{\Xdaginf /S}, \Sp)\Mod)$: 
$$\bfR \Gamma_Z(\calMdaginf )\to\calMdaginf \to\bfR j^{\inf}_*j_{\inf}^{-1}\calMdaginf \to\,,$$ qui se restreint au triangle usuel sur tout ouvert
du site $ \Xdaginf $.
\end{theo*}

\begin{theo*}[\ref{fon-inv}]Soit   $f: Y\to X$ un morphisme de schémas lisses sur $k$. Il existe un foncteur image inverse qui est un foncteur exact de catégorie triangulées:
$$  \Dmoins ((\calDdag_{\Xdaginf /V},\Sp)\Mod)\to \Dmoins ((\calDdag_{\Ydaginf /V},\Sp)\Mod).\leqno f^*_{\diff} :$$ 
\end{theo*}
Le foncteur image inverse se définit aussi sur l'extension $V\to K$.

\medskip

Pour le foncteur image directe nous faisons pour simplifier dans cet article la  restriction de séparation  sur la base, et
nous définissons le foncteur image directe  sur le corps des fractions $ K$. 

\begin{theo*}[\ref{fon-dir}]Soit  $f: Y\to X$ un morphisme de schémas    lisses  sur $k$ tel que $X$ est séparé sur $k$. Il existe un foncteur image directe qui est un foncteur exact de catégorie triangulées:
$$ \Dplus((\calDdag_{\Ydaginf /K},\Sp))\Mod\to \Dplus ((\calDdag_{\Xdaginf /K},\Sp)\Mod).\leqno f^{\diff}_*:$$ 
\end{theo*}

Ainsi, si la base est un
point et   $Y$ est une variété algébrique lisse sur $k$, l'image directe 
$ f_*^{\diff} \cal O_{\Ydaginf /K}$ est un complexe de $ \Dplus (K)$ dont la cohomologie est la cohomologie de de Rham $p$-adique de la source. Le
théorème de finitude [Me$_2$] montre alors que ce complexe est à cohomologie de dimension {\bf finie}.
Dans le cas général, si $f$ est un morphisme de variétés algébriques lisses sur $k$, 
{\bf  le complexe spécial}  ${ f_*^{\diff}} \cal
O_{\Ydaginf /K}$  sur $X$ est l'analogue $p$-adique du complexe  de Gauss-Manin d'un morphisme entre variétés algébriques lisses sur un corps de caractéristique nulle,
dont on s'attend naturellement à ses propriétés de finitude. Si par exemple le morphisme $f$ est
propre et lisse, on s'attend à ce que ses faisceaux de cohomologie soient des fibrés $p$-adiques dont  les rangs sont donnés par les nombres de Betti $p$-adiques
des fibres, et {\bf alors et seulement alors}, si $X$ est une courbe,  la théorie des  équations différentielles $p$-adiques ([C-M$_1$],
[C-M$_2$], [C-M$_3$], [C-M$_4$]) définit leurs monodromies aux points à l'infini.

Le foncteur image directe $i_*^{\diff}$ est défini sur $V$ pour toute immersion fermée $i: Y\to X$ de schémas lisses sur $k$. Pour définir le morphisme de Gysin et le morphisme classe de cohomologie d'un cycle nous étudions les compatibilités entre les foncteurs $i_*^{\diff}, i^*_{\diff},  \bfR \Gamma_Y$.

\begin{theo*}[\ref{adj},\ref{coh-loc},\ref{adj-imm}] Soit $i: Y\to X$ une immersion fermée de schémas lisses sur $k$. Il existe
un morphisme de foncteurs de la catégorie $\Dplus((\calDdag_{\Xdaginf /V},\Sp)\Mod)$: 
$$i_*^{\diff}i^*_{\diff}\to \bfR \Gamma_Y [\codim_XY]\leqno \rm Adj^i_*:$$ et un morphisme de foncteurs de la catégorie $\rmD((\calDdag_{\Ydaginf /V},\Sp)\Mod)$ qui est un isomorphisme:
$$Id \simeq i^*_{\diff}i_*^{\diff} [-\codim_XY]\leqno \rm Adj^*_i:$$ où $Id$ est le foncteur identique. De plus le foncteur $\rm Adj^i_*$ induit des isomorphismes $\rm Adj^i_*(\cal
O_{\Xdaginf /V})$ de la catégorie $\Dplus((\calDdag_{\Xdaginf /V},\Sp)\Mod)$:
$$i_*^{\diff}\cal
O_{\Ydaginf /V}\simeq i_*^{\diff}i^*_{\diff}(\cal
O_{\Xdaginf /V})\simeq \bfR \Gamma_Y(\cal O_{\Xdaginf /V}) [\codim_XY].$$
\end{theo*}

\paragraphe
Ces foncteurs coïncident  avec les foncteurs qu'on peut construire chaque fois qu'on peut relever le triplet $(Y,X, f)$. Cela 
est plutôt rare, même pour un morphisme de courbes propres et lisses. Ces théorèmes  résolvent,  du point de vue cohomologique,  le problème des relèvements des \em{morphismes quelconques entre  schémas lisses},
parce que les foncteurs
$f^*_{\diff}$ et $f_*^{\diff}$ sont des {\bf invariants intrinsèques} du morphisme $f$. Ces foncteurs sont munis de morphismes de compatibilités avec les foncteurs images inverse et directe topologiques agissant sur les complexes de de Rham locaux (Thm. \ref{dr-dir}, Thm. \ref{ima-dr}). Le foncteur dualité est compatible avec le foncteur image directe pour une immersion fermée (Thm. \ref{dua-imm}) et avec le foncteur image inverse pour une projection (Thm. \ref{dua-pro}).
En particulier, toute construction géométrique en caractéristique $p>0$ se relève {\bf intrinsèquement} du point de vue cohomologique en caractéristique zéro.

\paragraphe
Un autre  intérêt de la théorie précédente est de fournir  des invariants sur l'anneau des entiers $V$ et pas seulement sur le corps des
fractions
$K$, bien qu'il faille modifier la définition du foncteur image directe, ce que les méthodes de la géométrie analytique rigide ne peuvent pas faire. Cependant, les méthodes de cet article fournissent déjà une bonne théorie pour la cohomologie de de Rham $p$-adique d'échelon nul dans le cas où l'indice de ramification $e$ de $V$ est strictement plus petite
que $p-1$. Cela permet de redémontrer, comme test, que si $e<p-1$,  la cohomologie de de Rham algébrique d'un relèvement propre et lisse
sur $V$ d'une variété propre et lisse sur $k$ est indépendante du relèvement à isomorphisme canonique près  (Coro. \ref{alg-ada}). De même, les méthodes de cet article montrent 
la fonctorialité sur $V$ de la cohomologie de de Rham $p$-adique d'échelon nul sur $V$ d'une variété algébrique lisse sur $k$, toujours sous la condition $e<p-1$ (Thm. \ref{fon-ent}).

\paragraphe Pour illustrer nos méthodes on démontre les théorèmes suivants, qui sont autant de  tests pour une théorie cohomologique des schémas algébriques:

\begin{theo*}[\ref{fon-deR}] Soit $ \Sms (k)$ la catégorie des schémas lisses et séparés sur $k$ alors la correspondance  $\DR(-/K)$
$$\matrix{
\Sms (k)&\fonct&\Dplus (K)\hfill\cr
X&\fonct& \DR(X/K):=\bfR \hom\SubX{\calDdag_{\Xdaginf /K},\Sp}(
\cal O_{\Xdaginf /K},\cal O_{\Xdaginf /K})
}$$
est un {\bf foncteur contravariant} qui étend le foncteur de Monsky-Wash\-nitzer {\rm [M-W]} et qui induit un
isomorphisme $\DR(f/K)$  dans le cas où le morphisme $f$ de schémas est la projection de l'espace affine $A^n_k\times_k X$ de base $X$ sur $X$.\glossary{$ X\mapsto \DR(X/K)$}
\end{theo*}
Comme application de la formule locale des traces de Monsky [M$_2$],  du théorème de finitude des nombres de Betti $p$-adiques  [Me$_2$] et de la  fonctorialité  précédente on trouve la factorisation $p$-adique de la fonction Zêta
 d'un variété algébrique non singulière sur un corps fini, qui nous a servi de guide, de test et qui justifie nos méthodes:

\begin{theo*}[\ref{fac-zet}]Supposons  que $k$ est un corps fini $\Bbb F_q$ à $q$ éléments et que $f=fr$ est le morphisme de Frobenius relatif à $\Bbb F_q$ d'une variété algébrique 
$X$ lisse purement de dimension $\dim X$. Les endomorphismes de Frobenius $\bfF^{\bullet}$ de la cohomologie de de Rham $p$-adique $ H^{\bullet}_{\DR}(X/K)$ définis par la fonctorialité du théorème précédent 
sont \em{bijectifs}, et l'on a la factorisation $p$-adique de la fonction Zêta de $X$:
$$Z(X,t)= \prod_{i, impair} P_{p,i}(X)\mBig/\prod_{i, pair} P_{p,i}(X),$$ où $P_{p,i}(X):=\det(1-q^{\dim X}({\bfF}^{i})^{-1}t)\in K[t]$, 
$0\leq i\leq 2\dim X$, est le polynôme caractéristique 
de l'endomorphisme $q^{\dim X}({\bfF}^{i})^{-1}$ agissant sur  $ H^{i}_{\DR}(X/K)$.
\end{theo*}

  Cela fournit la première démonstration, à notre connaissance,  de l'expres\-sion cohomologique  $p$-adique de la fonction Zêta d'une variété algébrique
non singulière sur un corps fini par des polynômes, sans restriction  sur les  singularités à l'infini ou de relèvement (par exemple, pour toutes les variétés quasi-projectives lisses à distance finie). Auparavant, on connaissait, d'une part,  le cas affine  [M$_2$],  complété de façon essentielle par le théorème de finitude   [Me$_2$] où intervient déjà de façon cruciale  la théorie des équations différentielles $p$-adiques (parce que les variétés affines et leurs morphismes se relèvent) et, d'autre part,  le cas propre et lisse par la théorie cristalline [B$_2$] (parce qu'il n'y pas de singularités du tout ni à distance finie ni à l'infini).

\begin{theo*}[\ref{sui-gys}] Soit $i: Y\to X$ une immersion fermée de schémas lisses et séparés sur $k$ de complémentaire $j: U\to X$. Alors, on a la 
suite exacte de Gysin:
$$ \cdots \to H_{\DR}^{\bullet -2\codim_XY}(Y/K)\to H_{\DR}^{\bullet }(X/K)\to H_{\DR}^{\bullet }(U/K)\to \cdots,$$ qui étend la suite de Gysin de l'article {\rm [Me$_2$]}
dans le cas affine. En particulier, la classe $c\ell(Y)$ de cohomologie de $Y$ est bien définie comme un élément de $ H_{\DR}^{2\codim_XY}(X/K)$ dont la restriction à $U$ est nulle. On en déduit lorsque le corps de base est parfait
un morphisme \og classe de cohomologie\fg\ gradué:
$$ C^\bullet(X)\to H_{\DR}^{2\bullet}(X/K)\leqno c\ell:$$ entre le groupe gradué par la codimension des cycles de $X$, et le $K$-espace gradué de cohomologie de de Rham $p$-adique de $X$. \end{theo*} 
Cela fournit, à notre connaissance, la première construction du morphisme \og classe de cohomologie d'un cycle\fg\ en théorie de de Rham $p$-adique sans restriction sur les singularités dans le cas d'un corps de base parfait et donne donc lieu aux problèmes classiques des rapports entre classes de cohomologie et cycles d'autant plus que la cohomologie de de Rham $p$-adique est munie de l'action de Frobenius. Rappelons que la théorie cristalline
fournit un morphisme classe de cohomologie pour les cycles non singuliers [B$_2$].

\paragraphe
C'est à l'intérieur de la catégorie $ \rmD^*((\calDdag_{\Xdaginf /K},\Sp)\Mod)$ qu'on pourra définir à l'aide des foncteurs $ f^*_{\diff}$, $f_*^{\diff}$, $\bfR \Gamma_Z$, $\bfR j^{\inf}_*j_{\inf}^{-1}$ des catégories pleines de coefficients ayant les bonnes propriétés
de finitude, mais nous n'aborderons pas cette question dans le présent article. Le cas des courbes nous indique clairement où se trouvent les difficultés et
est en même temps un guide précieux. 
On
dispose déjà de la catégorie de base des fibrés
$p$-adiques
$ \MLS (\cal O_{\Xdaginf /K})$ à savoir la sous-catégorie pleine de la catégorie $ (\calDdag_{\Xdaginf /K},\Sp)\Mod$ des modules à gauche spéciaux
qui sont localement libres de type fini sur le faisceau $ \cal O_{\Xdaginf /K}$ introduite dans  [Me$_2$] dans le cas d'un relèvement. 

\paragraphe
Si
$X$ est une courbe lisse, nous avons défini dans [C-M$_1$], [C-M$_2$], [C-M$_3$] et [C-M$_4$]
la {\bf monodromie} des objets  de la catégorie $ \MLS (\cal O_{\Xdaginf /K})$ en un  point singulier à l'infini, 
en particulier nous avons défini leurs exposants $p$-adiques et leur polygone 
de Newton $p$-adique.
C'est là le fond du problème et le progrès décisif dans la théorie $p$-adique. Nous avons en particulier résolu dans l'article fondamental [C-M$_2$] le problème du prolongement analytique multiforme des solutions locales des équations différentielles sur le plan $p$-adique qui est totalement discontinu  sans lequel rien n'est possible comme nous l'avons mis en évidence dans ([M-N$_1$], introduction). Dans 
l'article [C-M$_4$]  nous avons montré la formule d'Euler-Poincaré globale pour les fibrés $p$-adiques appartenant à une sous-catégorie de la catégorie
$ \MLS (\cal O_{\Xdaginf /K})$ en imposant des restrictions sur leurs exposants $p$-adiques et aux exposants  $p$-adiques de leur module des endomorphismes aux points à l'infini, ce qui montre que c'est le bon point de
vue:  ces restrictions sont satisfaites pour les fibrés $p$-adiques qui proviennent de la géométrie, ce qui donne lieu à la théorie $p$-adique des fonctions $L$ sur les corps finis [C-M$_4$].

\paragraphe
Dans la théorie $p$-adique les propriétés de finitude  sont des questions hautement non triviales, comme l'illustrent le théorème de finitude des
nombres de Betti [Me$_2$] et la formule d'Euler-Poincaré ([C-M$_1$], [C-M$_4$]). Pour  démontrer ces théorèmes on a été amené à développer en profondeur la théorie des équations différentielles $p$-adiques  
([C-M$_1$], [C-M$_2$], [C-M$_3$], [C-M$_4$]), qui est de toute façon au   coeur de  la théorie des coefficients
$p$-adiques et qu'on ne pourra pas éviter par des considérations formelles. La  théorie  des équations différentielles $p$-adiques permet de  définir la
catégorie  des  coefficients
$p$-adiques sur les courbes et surtout de définir leurs monodromies aux points singuliers, qui sont  à la base de  leurs propriétés de finitude, 
en particulier de leur stabilité pour une {\bf immersion ouverte} ([M-N$_1$], Thm. 4.5.4,  [C-M$_1$],  [C-M$_2$], [C-M$_3$], [C-M$_4$]). Les trois premiers article montrent le théorème de l'indice et le dernier article montre l'existence d'un réseau.
La
théorie des coefficients
$p$-adiques sur les courbes a déjà de
nombreuses  applications ([Me$_2$]). 

\paragraphe
Il nous reste à étendre ces résultats en dimensions supérieures. Pour attaquer ce problème on dispose maintenant 
du couple
$ \cal G_{\Xdaginf }\hookrightarrow
\calDdag_{\Xdaginf /V}$  qui est le cadre naturel de la géométrie  différentielle   en inégales caractéristiques. En  surmontant  intrinsèquement les
difficultés liées aux relèvements aussi bien pour les variétés que pour leurs morphismes, ce couple 
permet de raisonner géométriquement pour la topologie de Zariski des schémas en caractéristiques positives. C'est un des outils  de base pour traiter ces questions, les foncteurs $ f^*_{\diff}$, $f_*^{\diff}$, $\bfR \Gamma_Z$, $\bfR j^{\inf}_*j_{\inf}^{-1}$ qu'il  définit 
jouant un rôle essentiel et la théorie trouve un  cadre vraiment naturel permettant de formuler les bons énoncés.

\paragraphe\label{site-infinitesimal-formel}On peut considérer le site $ X^{\wedge}_{\inf}$ \glossary{$ X^{\wedge}_{\inf}$}infinitésimal formel dont les objets $\cal U^\wedge= (U, \cal O_{\cal U^\wedge/R})$
sont les relèvements formels
$R$-plats d'ouverts  de $X$ qui est annelé par le faisceau des fonctions formelles $ \cal O_{X^\wedge_{\inf}/R}$\glossary{$ \cal
O_{X^\wedge_{\inf}/V}$}. On définit de même le faisceau de groupes
$ \cal G_{X^\wedge_{\inf}}$ \glossary{$ \cal
G_{X^\wedge_{\inf}}$} des automorphismes d'algèbres qui se réduisent à l'identité modulo $I$. Localement, un élément de ce groupe est un
opérateur différentiel à coefficients formels défini
de même par réduction modulo les puissances de $ I$ en imposant une condition de croissance de type $\dagger$ sur l'ordre des opérateurs différentiels et  le faisceau $ \calDdag_{X^\wedge_{\inf}/R}$ \glossary{$ \calDdag_{X^\wedge_{\inf}/R}$}des opérateurs différentiels à coefficients formels  garde un sens sur le site infinitésimal formel 
$ X^{\wedge}_{\inf}$ et l'on a un plongement $ \cal G_{X^\wedge_{\inf}}\hookrightarrow \calDdag_{X^\wedge_{\inf}/R}$ pour la structure multiplicative. On définit alors la catégorie $ (\calDdag_{X_{\inf}^\wedge /V},\Sp)\Mod$ \glossary{$ (\calDdag_{X_{\inf}^\wedge /V},\Sp)\Mod$} des $ \calDdag_{X^\wedge_{\inf}/R}$-modules à gauche spéciaux qui est canoniquement équivalente à la
catégorie des $\calDdag_{\cal X^\wedge/R}$-modules \glossary{$\calDdag_{\cal X^\wedge/R}$} à gauche d'un relèvement formel plat lorsqu'il existe. La
cohomologie de de Rham $p$-adique $$ H_{\DR}^{\bullet }(X/K, \cal O_{X^\wedge_{\inf}/K}):=\ext^{\bullet}_{\calDdag_{X^\wedge_{\inf} /K},\Sp}(\cal
O_{X^\wedge_{\inf} /K}, \cal
O_{X^\wedge_{\inf} /K})$$ du site  infinitésimal formel d'une variété algébrique $X$ lisse sur $k$ doit être canoniquement isomorphe à la
cohomologie cristalline de $X$ sur $K$ de façon fonctorielle prolongeant l'isomorphisme de comparaison de Berthelot-Ogus [B-O] 
dans le cas d'un relèvement  lisse. Tous les développements qui ne font pas intervenir des propriétés de finitude de cet article se transposent avec des démonstrations plus simples en remplaçant le couple 
$ \cal G_{\Xdaginf }\hookrightarrow\calDdag_{\Xdaginf /R}$ 
par le  couple  $ \cal G_{X^\wedge_{\inf}}\hookrightarrow\calDdag_{X^\wedge_{\inf}/R}$. 
\glossary{$ \cal G_{X^\wedge_{\inf}}\hookrightarrow\calDdag_{X^\wedge_{\inf}/R}$} Le
lecteur non averti prendra garde à ce que la cohomologie de de Rham $p$-adique
du site infinitésimal formel, tout comme celle du site cristallin, ne fournit pas les bons nombres de Betti $p$-adiques pour les variétés {\bf ouvertes}, ce qui explique la nécessité du passage du site infinitésimal formel
au site infinitésimal 
$\dagger$-adique précisément. D'autre part, on dispose d'un foncteur naturel exact de la catégorie $ (\calDdag_{\Xdaginf /V},\Sp)\Mod$ vers la catégorie $ (\calDdag_{X^\wedge_{\inf} /V},\Sp)\Mod$, ce qui induit un morphisme de comparaison de la cohomologie de de Rham $p$-adique du site $\dagger$-adique
à la cohomologie de de Rham $p$-adique du site formel. On prendra garde aussi à ce que
la suite exacte de Gysin n'a pas lieu dans le cas formel.

\paragraphe
Comme le lecteur peut le constater cette théorie de nature essentiellement algébro-géométrique est un cadre naturel si l'on veut développer une théorie des coefficients $p$-adiques
pour les variétés algébriques.  
On peut penser et espérer qu'à terme elle fournira une théorie cohomologique des coefficients $p$-adiques disposant de tous ses outils essentiels, ce qui est déjà le cas  pour les courbes grâce à la théorie des équations différentielles $p$-adiques ([C-M$_1$],
[C-M$_2$], [C-M$_3$], [C-M$_4$]). Nous avons maintenant  de nombreuses indications allant dans ce sens en dimensions supérieures. 
Elle soulève de nombreuses questions  
et ouvre de nouveaux horizons.
Elle provient  
de trois courants d'idées dont elle est la synthèse. La notion d'algèbre faiblement complète [M-W] et  sa localisation [Mr], la théorie du site
infinitésimal [G$_4$] et bien entendu   la théorie des $\calD_X$-modules des années dix-neuf cent soixante-dix comme théorie des Coefficients de de Rham en Géométrie Algébrique au sens de Grothendieck, qui a fourni la bonne théorie en
caractéristique zéro [Me$_1$].
 Ces trois courants ont historiquement des motivations bien distinctes, mais se complètent harmonieusement dans cette théorie. 


\paragraphe
Nous n'utilisons dans le présent article aucun résultat des articles fondateurs du point de vue conceptuel et technique ([M-W], [M$_1$]), dont les résultats les plus importants apparaissent ici comme des cas particuliers. En particulier, nous retrouvons avec nos méthodes comme test préliminaire (Coro. \ref{fon-aff}) la fonctorialité de la cohomologie de de Rham $p$-adique des variétés affines lisses de l'article [M-W] qui est le point clef de la théorie de ces auteurs.

\paragraphe
La Table des matières donne une idée du contenu de notre article qui comporte deux parties. Les chapitres 1 à 9 développent les fondements de la théorie des 
modules spéciaux sur le site infinitésimal $\dagger$-adique d'un schéma lisse et de la cohomologie de de Rham $\dagger$-adique. Les chapitres 10 à 15 développent les opérations 
cohomologiques pour les catégories des modules spéciaux sur les sites infinitésimaux $p$-adiques 
et les propriétés de  fonctorialité de la cohomologie de de Rham $p$-adique. Les résultats de cet article reposent sur les propriétés de {\bf finitude} des opérateurs différentiels $p$-adiques ([M-N$_1$], [M-N$_2$], [M-N$_3$], [Me$_3$]), sur le théorème des relèvements
des algèbres lisses et de leurs morphismes [A] et sur le critère d'affinité $\dagger$-adique [A-M$_2$].

\paragraphe
À l'occasion de cet article, nous souhaitons remercier Alexandre Grothendieck, 
qui a incité avec force en mai 1983 le deuxième auteur à travailler dans le sujet. Nous espérons que nous avons contribué à la promotion de ses idées parmi les nouvelles générations. 
Tout particulièrement son idée de Cristal en tout genre trouve son aboutissement pour les variétés éventuellement ouvertes  dans la notion de Module Spécial sur le Site Infinitésimal $p$-adique d'un schéma, bien que comme nous l'avons déjà signalé la théorie des Cristaux et la théorie des $\calD_X$-modules ont au départ des motivations complètement indépendantes.
Nous remercions aussi nos collègues Gilles Christol et Luis Narvaez pour leurs contributions respectives qui ont fait progresser le sujet et  l'ont rendu crédible. Les résultats de cet article ont été présentés à la Conférence  Satellite de Géométrie Algébrique du  C.I.M. de Madrid  2006, Segovia (Espagne) 16-19 août  2006.  

\section{Notations, conventions et rappels}
Une difficulté dans ce domaine est la cohérence des notations et de la terminologie. Nous avons essayé d'être aussi soigneux et simples que possible, autant du point de
vue des notations que de la terminologie. En particulier, nous dirons de façon générique \og cohomologie de de Rham
$\dagger$-adique\fg\  pour  cohomologie
de type de de Rham que nous considérons dans cet article et plus spécialement  \og cohomologie de de Rham $p$-adique\fg\ lorsque l'anneau de base est de
valuation d'inégales caractéristiques $p>0$, étant entendu que tous les espaces   de cohomologie de type de Rham naturels et intéressants doivent être 
canoniquement isomorphes par des théorèmes de comparaison convenables.  Afin d'éviter des répétitions, nous utiliserons dans tout  cet article les notations suivantes.

\begin{notation}\begin{liste}\rm
\item 1) {\spaceskip1ex plus 1ex
$R$ anneau commutatif ayant un élément unité, $I$ idéal de $R$, $R_s:= R/I^s, s\geq 1$; en particulier $R_1:= R/I$.}
\item 2) $V,\goth m, k, K, e$  anneau de valuation discrète complet d'inégales caractéristiques $(p,0)$ d'idéal maximal $\goth m$, de corps résiduel $k$ et 
de corps de
fractions $K$ et d'indice de ramification absolu $e:= v_{\goth m}(p)$.

\item 3) $X/R$ schéma  sur $R$.
\item 4) De façon générale, si $f: Y\to X$ est un morphisme d'espaces annelés sur $R$, nous notons $\calD_{Y\to X/R}$ le faisceau des
opérateur différentiels défini par récurrence sur
l'ordre dans ([EGA IV$_4$], $\S 16$) et noté $\cDiff_{R}(f^{-1}\cal O_{X/R}, \cal O_{Y/R})$. 
Les notations $\calD_{Y\to X/R}, \,\,\, \calD_{X/R}:= \cDiff_{R}(\cal O_{X/R}):=\cDiff_{R}(\cal O_{X/R}, \cal O_{X/R})$ et leurs  variantes sont aujourd'hui universelles.
\item 4) Si $X/R$ est un schéma  lisse sur $R$, un système de coordonnées locales $x_1,\dots,x_n$ est constitué de sections locales du faisceau $\cal
O_{X/R}$ telles que leurs différentielles forment une base locale du  fibré des $R$-formes différentielles ([EGA IV$_4$], $\S 16$).
\item 5) $\alpha\in \Bbb N^n, \alpha!:= \alpha_1!\cdots\alpha_n!$, $|\alpha|:= \alpha_1+\cdots +\alpha_n$,
$\comb\alpha\beta:=\alpha!/\beta!(\alpha-\beta)!$, $x^\alpha:= x_1^{\alpha_1}\cdots x_n^{\alpha_n}$.
\item 6) $\Delta^\alpha:= \Delta_x^\alpha := \Delta_{x_1}^{\alpha_1}\cdots\Delta_{x_n}^{\alpha_n}:= \Delta_{1}^{\alpha_1}\cdots\Delta_{n}^{\alpha_n}$ la suite des opérateurs différentiels associés à
$x=(x_1,\dots,x_n)$ et à $\alpha=(\alpha_1,\dots,\alpha_n)\in \Bbb N^n$; ces opérateurs sont caractérisés par 
$\Delta^\alpha(x^\beta)= (^\beta_\alpha)x^{\beta-\alpha}$ ([EGA IV$_4$], Thm. 16.11.2).
\item 7) $i: Y\hookrightarrow X$ immersion fermée de complémentaire $j: U\hookrightarrow X$, $p: Y\times X\to X$ projection du produit sur le second
facteur et $q: Y\times X\to Y$ projection du produit sur le premier facteur, $f: Y\to X$ morphisme de
$Y$ dans
$X$.
\item 8) $r:=r_{W,U}: W\hookrightarrow U$ l'inclusion d'un ouvert dans un autre; nous noterons $r$ pour $r_{W,U}$ quand il n'y a pas de risque de
confusion.
\item 9) Si $f: Y\to X$ est un morphisme, nous notons par $f_*, f^{-1}$ les images directe et inverse topologiques usuelles, par $f_*^{\diff}$ une image
directe différentielle et par $f^*_{\diff}$ une image inverse différentielle.
\item 10) Si $R$ est noethérien et $A$ est une $R$-algèbre, on note $A\dag$ son complété $\dagger$-adique (complété faible) pour la topologie $I$-adique de $R$ [M-W],
c'est donc par construction une algèbre sur le complété séparé $\hat R= R\dag$ de $R$ qui est séparée pour la topologie $I$-adique. On appelera \og$R$-algèbre $\dagger$-adique\fg\ une  algèbre qui coïncide
avec son complété $\dagger$-adique. On pose $A_s:= A/I^s$.
\item 11) Si $A\dag$ est une $R$-algèbre $\dagger$-adique, on note $\Omega_{A\dag/R}:= \Omega^{sep}_{A\dag/R}$ le $A\dag$-module des formes différentielles
{\bf séparées} de $A\dag$.
\item 12)Si $\cal A_X$ est un faisceau d'anneaux sur un espace topologique $X$, on note $ \cal A_X\Mod$ la catégorie des $ \cal A_X$-modules
à gauche et $\Rm
\Modd\,\cal A_X$ la catégorie des $\cal A_X$-modules à droite. Mais si $\cal A_X$ est commutatif, on note $\RmMod(\cal A_X)$ la catégorie des $ \cal A_X$-modules.
\item 13) Si $\Ab$ est une catégorie abélienne, on note $ \rmD^*(\Ab)$ sa catégorie 
dérivée (le signe $*$ pouvant être vide ou prendre les valeurs $+,-,b$ [V]). Mais nous notons $ \rmD^*(\cal A_X)$ au lieu de 
$ \rmD^*(\RmMod(\cal A_X))$ dans le cas d'un faisceau d'anneaux commutatifs.
 
\end{liste}
\end{notation}
Nous compléterons  cette liste de notations au fur et à mesure  et nous dressons leur liste dans l'index des notations.

Comme le lecteur pourra le constater, la théorie de la catégorie homotopique $\rm K^*$ est insuffisante pour le présent article, alors que  la théorie de  la catégorie dérivée $\rmD^*$ [V] est essentielle aussi bien pour les résultats que pour les démonstrations. 

\begin{defi}\label{parfait} Nous rappelons qu'un complexe  de modules sur un faisceau d'anneaux sur un espace topologique ou sur un site est \em {parfait},  selon Grothendieck, (relativement à la catégorie des modules localement libres de rang fini) si localement il admet une résolution de longueur finie dont les termes sont des modules libres de rang fini.\end{defi}

Nous utilisons le critère de platitude locale sous la forme  suivante ([EGA III$_1$], Chap. 0, 10.2.1), ([Bo], Bourbaki, \em{Alg. comm.}, Chap. III, \S 5, n° $2$, Thm. 1).
\begin{prop}\label{cri-pla}Soient deux  $V$-algèbres noethériennes $D$ et $ D'$,  telles que l'idéal $\goth m$ est contenu dans le radical de $D'$. Si $M$ est un $(D, D')$-bimodule qui est un $D'$-module  de type fini et que les réductions modulo $\goth m^s$ de $M$ sont $D_s$-plates pour  tout $s\geq1$,   alors $M$ est un $D$-module à gauche  {\bf plat}.
\end{prop} 

Par définition d'une $V$-algèbre, $V$ est contenu dans le centre de $D$. Les hypothèses de la proposition entraînent que
$M$ est  $D$-idéalement séparé pour la topologie $\goth m$-adique, c'est-à-dire pour tout idéal $J$ à droite de $D$ le $D'$-module à droite $J\otimes_D M$ est séparé pour la topologie $\goth m$-adique.

\subsection{Faisceau donné localement sur un espace topologique}
Si $X$ est un espace topologique, rappelons qu'un faisceau d'ensembles donné localement est la donnée d'un recouvrement de $X$ par des ouverts $U_i$, $i\in J$, d'un
faisceau d'ensembles
$\cal F_i$ sur chaque ouvert, d'un isomorphisme $\alpha_{ij} : \cal F_i|U_{ij}\simeq \cal F_j|U_{ij}$  satisfaisant aux conditions de cocycles $\alpha_{ik}=
\alpha_{ij}\circ\alpha_{jk}$ sur $U_{ijk}:= U_i\cap U_j\cap U_k$. On rappelle le théorème classique qui pour Grothendieck est le cas le plus simple de la
méthode de la  descente.

\begin{theo} Soit $U_i$, $i\in J$ un recouvrement de $X$. Il existe un foncteur \em{canonique} de la catégorie des faisceaux d'ensembles donnés localement sur le
recouvrement
$U_i$, $i\in J$ dans la catégorie des faisceaux d'ensembles sur $X$ qui est un inverse canonique au foncteur de restriction.
\end{theo}

\bigskip

Nous utilisons de façon essentielle  le critère d'affinité d'un schéma faiblement complet [A-M$_2$] et le théorème du symbole total des opérateurs
différentiels $p$-adiques ([M-N$_1$], [M-N$_3$]), que nous allons rappeler rapidement pour la commodité du lecteur.

\subsection{Schémas $\dagger$-adiques, critère d'affinité $\dagger$-adique, produits fibrés}
Soit $R$ un anneau commutatif noethérien muni de la topologie définie par un idéal $I$. La notion de schéma faiblement complet pour la topologie
$I$-adique de $R$ a été définie par D. Meredith [Mr]. Rappelons que c'est un espace  annelé $\calXdag = (X, \cal O_{\calXdag/R})$\glossary{$\calXdag = (X, \cal O_{\calXdag /R})$} localement isomorphe à un schéma $\dagger$-adique affine de type fini. Un schéma
$\dagger$-adique affine se construit ([Mr]) à partir d'une $R$-algèbre
$A\dag$ faiblement complète de type fini, qui est noethérienne en vertu du théorème de Fulton,   exactement comme se construit un schéma formel affine adique noethérien  à partir d'un anneau adique noethérien [EGA I]. Par construction, si
$(X,
\cal O_{\calXdag /R})$ est un schéma faiblement complet, 
c'est un espace localement annelé dont le faisceau structural $\cal O_{\calXdag /R}$ est cohérent et  sa réduction $(X, \cal O_{\calXdag /R}/I):= (X, \cal O_{\calXdag /R}\otimes_{R_X}{R_{1X}})$ est
un
$R_1$-schéma  localement de type fini. En particulier, l'espace topologique $X$ est localement noethérien.  On dira que $\calXdag :=  (X, \cal O_{\calXdag /R})$ est un relèvement $\dagger$-adique de
sa réduction.  Nous utilisons systématiquement la terminologie schéma
{\bf $\dagger$-adique} au lieu  de schéma {\bf faiblement complet} qui peut prêter à confusion puisque complet n'implique pas faiblement complet, le couple  $R\supset I$ étant sous-entendu quand il n'y a
pas de risque de confusion précisément. 

\bigskip
Dans l'article  [A-M$_2$] nous avons démontré le critère d'affinité $\dagger$-adique:

\begin{theo}\label{cri-aff} Un schéma $\dagger$-adique $\calXdag = (X, \cal O_{\calXdag /R})$ est affine {\bf si et seulement si} le $R_1$-schéma
$(X, \cal O_{\calXdag /R}/I)$ est affine.
\end{theo}

Ce critère admet de nombreux corollaires importants et permet de raisonner avec la topologie de Zariski.

\begin{coro}L'espace localement annelé induit sur un ouvert affine par un schéma $\dagger$-adique est un schéma $\dagger$-adique affine.
\end{coro}

Ce corollaire produit beaucoup d'ouverts affines du site infinitésimal $\dagger$-adique en remplaçant l'hypothèse \og{\sl être 
un ouvert affine $\dagger$-adique\/}\fg\  par l'hypothèse plus pratique et plus souple  
\og{\sl être un ouvert affine \/}\fg\ au sens de la théorie des schémas.

\begin{notation}Si $\calUdag $ \glossary{$\calUdag $}est un schéma $\dagger$-adique et $W$ est un ouvert de l'espace topologique $U$, on note $\cal
U\dag|W$\glossary{$\cal
U\dag\mid W$} le schéma
$\dagger$-adique induit par $\calUdag $ sur $W$.
\end{notation}

\begin{coro}\label{acy}Soient  $\calXdag = (X, \cal O_{\calXdag /R})$ un schéma  $\dagger$-adique  et $\cal F$ un $\cal O_{\calXdag /R}$-module cohérent.
Si 
$U$ est un ouvert affine, la cohomologie $ H^i(U,\cal F)$ est nulle pour $i>0$.
\end{coro}

\demo C'est une conséquence directe du critère d'affinité et du théorème d'acyclicité ([A-M$_2$], Thm. 4.2.2), qui étend au-dessus de $R$ le théorème
d'acyclicité
de  Meredith [Mr]. En particulier, si $X$ est séparé, un recouvrement par des ouverts affines est un recouvrement de Leray pour tout $\cal O_{\calXdag /R}$-module cohérent.
\enddemo
On en déduit le critère cohomologique suivant.

\begin{coro}\label{aff-coh}Soit un schéma  $\dagger$-adique $\calXdag = (X, \cal O_{\calXdag /R})$. Les conditions suivantes sont équivalentes: \begin{liste}
\item 1) le schéma $\calXdag $  est $\dagger$-adique affine,  
\item 2) sa réduction  est un schéma affine, 
\item 3) pour tout  $\cal O_{\calXdag /R}$-module cohérent $\cal F$,
le module  $ H^i(X,\cal F)$ est nul pour $i>0$.
\end{liste}
\end{coro}

\demo
Il suffit de montrer que 3) implique 1). Mais  un  $\cal O_{X/R_1}$-module cohérent $\cal F_1$ est aussi   un  $\cal O_{\calXdag /R}$-module cohérent  et est donc cohomologiquement trivial par hypothèse. En vertu du critère d'affinité algébrique de Serre, le schéma $(X, \cal O_{X/R_1})$ est affine, et en vertu du critère d'affinité  $\dagger$-adique,  le schéma
$\calXdag $ est $\dagger$-adique affine.
\enddemo
\begin{Rema} La méthode de l'article [A-M$_2$] montre aussi le critère d'affinité dans le cas des schémas formels
et donc l'analogue du  critère cohomologique précédent. Nous ignorons si ce résultat, qui ne semble pas figurer  dans [EGA], était connu auparavant.\end{Rema} 

Dans le même article,   nous utilisons le critère d'affinité uniquement pour les ouverts affines déjà contenus dans les ouverts $\dagger$-adiques affines pour démontrer l'existence de produits fibrés ([A-M$_2$], Thm. 6.2.3):

\begin{theo}\label{pro-fib}La catégorie des  schémas $\dagger$-adiques sur un anneau noethérien $R$ muni d'une topologie définie par un idéal $I$ admet des produits
fibrés.
\end{theo}
L'existence du produit fibré dans le cas des schémas formels est démontré dans [EGA I].
Nous utilisons aussi le résultat suivant démontré dans ([A-M$_2$], 3.3.7).

\begin{theo}\label{eqi-aff}La catégorie des schémas $\dagger$-adiques affines est équivalente à la catégorie des  $R$-algèbres  $\dagger$-adiques  de
type fini.
\end{theo}

En particulier, un morphisme de schémas $\dagger$-adiques  affines provient d'un morphisme des algèbres associées.

Nous ne considérons dans le présent  article que les algèbres $\dagger$-adiques (topologiquement) {\bf de type fini}. Nous dirons  
algèbre $\dagger$-adique pour algèbre $\dagger$-adique topologiquement de type fini  et schéma $\dagger$-adique pour schéma $\dagger$-adique
localement topologiquement de type fini.  Nous rappelons le théorème ([M-W], Thm. 3.2):

\begin{theo}\label{rel-pla} Soit $u\dag : A\dag\to B\dag$ un morphisme de $R$-algèbres
$\dagger$-adiques. 
\begin{liste}
\item 1) Si la réduction modulo $I$ de $u\dag$ est surjective, alors $u\dag$ est surjective.

\item 2) Si 
 $B\dag$ une $R$-algèbre $\dagger$-adique \em{plate} et si la réduction modulo $I$ de $u\dag$ est un isomorphisme, alors $u\dag$ est un isomorphisme.\end{liste}
\end{theo}

Nous en déduisons le théorème suivant.

\begin{theo}\label{iso-pla}Soient $\cal X\dag_1$ et $\cal X\dag_2$ deux schémas $\dagger$-adiques   et $f\dag$ un morphisme
$\cal X\dag_1\to\cal X\dag_2$ dont la réduction $f$ modulo $I$  est un isomorphisme. Si $\cal X\dag_1$ est $R$-plat,
le morphisme $f\dag$ est un isomorphisme.
\end{theo}

\demo Il s'agit de montrer que le morphisme de faisceaux $R$-algèbres $f^{-1}\cal O_{\cal X\dag_2/R}\to\cal O_{\cal X\dag_1/R}$ est un isomorphisme.
La question est donc locale. On peut donc supposer que $X_1$ et $X_2$ sont  affines. En vertu du critère d'affinité \ref{cri-aff}, les schémas $\dagger$-adiques,
$\cal X\dag_1$ et $\cal X\dag_2$ sont donc affines d'algèbres $A\dag_1$   et $A\dag_2$; de plus, $A\dag_1$ est plate sur $R$. En vertu de l'équivalence \ref{eqi-aff},
le morphisme
$f\dag$ provient d'un morphisme
d'algèbres $A\dag_2\to A\dag_1$ dont la réduction est un isomorphisme. D'après le théorème \ref{rel-pla}, ce morphisme d'algèbres  est un isomorphisme. Enfin, en vertu de
l'équivalence
\ref{eqi-aff}, le morphisme $f\dag$ est un isomorphisme.
\enddemo
\begin{Rema} Le critère d'affinité dans la démonstration du théorème précédent n'est pas nécessaire si on se rappelle
que, sur un schéma $\dagger$-adique affine, l'espace annelé induit  sur un ouvert principal
est lui-même un schéma $\dagger$-affine et que les ouverts principaux forment une base de la topologie.\end{Rema}
\subsection{Algèbres $\dagger$-adiques lisses et schémas $\dagger$-adiques lisses}

\looseness-1
Nous rappelons les propriétés essentielles de $R$-lissité [A]. On dit qu'une $R$-algèbre $A$ est lisse si le morphisme structural
$\Spec A\to \Spec R$  est lisse au sens de la théorie des schémas ([EGA IV$_4$], $\S 17$); c'est une condition locale à la source qui implique que $\Spec A$ est  localement de présentation finie sur $\Spec R$, donc $A$ est de présentation finie sur $R$.
Dans le théorème des relèvements algébriques on ne fait \em{aucune hypothèse} sur le couple
$I\subset R$.

\begin{theo}\label{rel-alg}Si $R$ est un anneau commutatif avec élément unité,  $I\subset R$ un idéal et  $A_1$ une $R_1$-algèbre lisse, alors $A_1$  admet un relèvement {\bf lisse} $A$ sur
$R$.
\end{theo}
Le  premier résultat dans cette direction est dû à Grothendieck dans le cas d'un idéal $I$ nilpotent qui est aussi l'un des tous premiers succès de la théorie des schémas et des éléments nilpotents.   Le cas   d'un couple $I\subset R$ hensélien  est dû à R. Elkik. On rencontre bien entendu 
naturellement des couples $I\subset R$  non henséliens où $R$ est le plus souvent noethérien.

\bigskip
Supposons maintenant que $R$ est noethérien.

\begin{theo}\label{rel-mor}Soient $f_1: A_1\to B_1$ un morphisme de $R_1$-algèbres de type fini, $B\dag$,$ A\dag$ des relèvements $\dagger$-adiques de $B_1,
A_1$. Supposons que
$A_1$ est lisse sur
$R_1$ et que 
$A\dag$
est plate sur $R$.  
\begin{liste}\item 1) Le morphisme $f_1$ admet alors un relèvement $f\dag: A\dag\to B\dag$.
\item 2) Soient $f\dag: A\dag\to C\dag$ un morphisme de $R$-algèbres,  $v\dag: B\dag\to C\dag$  un morphisme de $R$-algèbres  
\em{surjectif} et  $u_1$
un morphisme de
$R_1$-algèbres
$A_1\to B_1$ qui factorise $f_1$:
$f_1= v_1\circ u_1$.
Il existe alors un morphisme de $R$-algèbres $u\dag: A\dag\to B\dag$ qui se réduit modulo $I$ à $u_1$ et qui  factorise $f\dag$: $f\dag = v\dag\circ
u\dag$.
\end{liste}
\end{theo}

La propriété $2)$ dans le théorème précédent exprime que l'algèbre $A\dag$ est très lisse dans la terminologie de Monsky-Washnitzer [M-W].

\begin{defi}Nous dirons qu'une algèbre $\dagger$-adique $A\dag$ est une algèbre $\dagger$-adique lisse sur $R $, sous-entendu pour la topologie $I$-adique,
si l'une des conditions suivantes et équivalentes a lieu:
\begin{liste}
\item 1) l'algèbre  $A\dag$ est plate sur $R$ et sa réduction $A_1$ est lisse sur $R_1$, \item 2) l'algèbre $A_s$ est lisse sur $R_s$ pour tout
$s\geq 1$, \item 3) l'algèbre $A\dag$ est très lisse.
\end{liste}
\end{defi}
Grothendieck dit formellement lisse pour la condition 2) et 
Monsky-Washnitzer disent très lisse pour la condition 3). Comme toutes ces définitions sont équivalentes, il est légitime d'appeler une telle algèbre $\dagger$-adique lisse ou même simplement lisse, quand il n'y a pas de risque de confusion.  Le complété $\dagger$-adique $A\dag$  d'un relèvement algébrique lisse $A$ d'une 
$R_1$-algèbre lisse $A_1$ est donc  une algèbre
$\dagger$-adique lisse.

Nous rappelons qu'on dit qu'un schéma est non singulier, ou est régulier,  si tous ses anneaux locaux sont réguliers. Un schéma lisse sur un corps $k$ est non singulier; de plus,  si le corps de base est parfait, un schéma non singulier localement de type fini sur $k$ est lisse sur $k$ ([EGA IV$_4$], $\S 17.15$).

Nous rappelons  que si $A\dag$ est une $R$-algèbre $\dagger$-adique, le $A\dag$-module $\Omega_{A\dag/R}$ des formes différentielles {\bf séparées} est de type fini ([M-W], Thm. 4.5), et que  si $A\dag$ est $\dagger$-adique lisse, c'est un $A\dag$-module projectif
([M-W], Thm. 4.6). De plus, la réduction modulo $I$ de $\Omega_{A\dag/R}$ est canoniquement isomorphe à $\Omega_{A_1/R_1}$ ([M-W], Thm. 4.4).

\begin{defi}On dit qu'un schéma $\dagger$-adique  $\calXdag  = (X,\cal O_{\calXdag /R})$ est \og $\dagger$-adique lisse\fg\ si son  faisceau structural est
plat sur
$R$ et que sa réduction
$X_1 := (X,\cal O_{\calXdag /R}/I)$\glossary{$X_s := (X,\cal O_{\calXdag /R}/I^s), X_1:= X$} est lisse sur $R_1$.
\end{defi}
En utilisant le corollaire \ref{acy}, on voit qu'un schéma $\dagger$-adique $\calXdag  = (X,\cal O_{\calXdag /R})$ est $\dagger$-adique lisse si et
seulement si son  algèbre au-dessus de tout ouvert affine est une algèbre $\dagger$-adique lisse. Un schéma $\dagger$-adique lisse est plat,
mais si la réduction est lisse  un schéma $\dagger$-adique plat est  $\dagger$-adique lisse.  Quand on sait à l'avance que 
{\bf  la réduction  est lisse}, comme c'est souvent le cas dans cet article,   il vaut mieux  dire \og plat\fg\ au lieu de $\dagger$-adique \og lisse\fg\ pour les schémas $\dagger$-adiques pour être précis et éviter toute confusion.

Si $\calXdag $ est un schéma  $\dagger$-adique, le $\cal O_{\calXdag /R}$-module $\Omega_{\calXdag /R}$ des formes différentielles séparées est cohérent;  et si $\calXdag $ est $\dagger$-adique lisse, c'est un $\cal O_{\calXdag /R}$-module localement libre de type fini.

Dans cet article nous dirons pour abréger relèvement au lieu relèvement $\dagger$-adique, et donc relèvement plat au lieu
de relèvement $\dagger$-adique plat.

\subsection{Le théorème du symbole total d'un opérateur différentiel $p$-adique}
Soit $\calXdag = (X, \cal O_{\calXdag /R})$ un schéma $\dagger$-adique (non nécessairement lisse), on définit le faisceau des opérateurs différentiels 
$\calDdag_{\calXdag /R}$\glossary{$\calDdag_{\calXdag /R}$}
([M-N$_1$], [M-N$_2$]):

\begin{defi}Le faisceau $\calDdag_{\calXdag /R}$ des opérateurs différentiels $\dagger$-adiques est le sous-faisceau du faisceau des endomorphismes du
faisceau structural
$\cal End_R(\cal O_{\calXdag /R})$ des sections $P$ telles que pour tout $s\geq 1$ leur  réduction modulo $I^s$ est un opérateur différentiel sur le
$R_s$-schéma
$X_s:= (X,
\cal O_{\calXdag /R}/I^s)$ dont l'ordre est localement borné par une fonction linéaire en $s$.
\end{defi}
Par définition, le faisceau $\calDdag_{\calXdag /R}$ contient comme sous-faisceau d'algèbres le faisceau $\calD_{\calXdag /R}$ des opérateurs
différentiels d'ordre localement fini, défini par récurrence sur l'ordre comme le faisceau des opérateurs différentiels de l'espace annelé $\calXdag = (X, \cal O_{\calXdag/R})$ [M-N$_2$]. Si $\calXdag $ est $\dagger$-adique, lisse le faisceau $\calD_{\calXdag /R}$ est une $\cal O_{\calXdag/R}$-algèbre filtrée par des $\cal O_{\calXdag/R}$-modules localement libres de type fini ([M-N$_2$], appendice A).  On peut étendre  ce résultat au faisceau 
$\calDdag_{\calXdag /R}$ dans le cas d'un couple $(V, \goth m)$.
 Si
$A\dag$ est une algèbre $\dagger$-adique, on note
$D\dag_{A\dag/R}$ l'anneau des sections globales des opérateurs différentiels du schéma $\dagger$-adique affine associé.
Dans l'article [M-N$_3$]  qui complète l'article de recherche [M-N$_1$] on démontre le théorème du symbole total en deux parties (Thm. 5.0.2, Thm. 6.0.7), 
la première partie définit le morphisme symbole total qui est alors injectif, la seconde partie montre qu'il est surjectif:

\begin{theo}\label{sym-tot} Soient $\calXdag $ un schéma $\dagger$-adique lisse sur $V$ et $U$ un ouvert affine de $X$ d'algèbre $\dagger$-adique $A\dag$ 
munie d'éléments $x_1,\dots,x_n$ dont les  différentielles $dx_1,\dots,dx_n$ forment une base du module des formes différentielles {\bf séparées}
$\Omega_{A\dag/V}$. Soient $P$ un opérateur différentiel de l'anneau $D\dag_{A\dag/V}:= \Gamma(U, \calDdag_{\calXdag /V})$ et  
$a_\alpha$, $\alpha\in
\Bbb N^n$ la suite d'éléments de
$A\dag$ définis par:
$$a_\alpha := \sum_{0\leq\beta\leq\alpha}\comb\alpha\beta(-x)^\beta P(x^{\alpha-\beta}).
$$\glossary{$a_\alpha :=
\sum_{0\leq\beta\leq\alpha}\comb\alpha\beta(-x)^\beta P(x^{\alpha-\beta})$}Alors, l'opérateur
$P$ est égal à la série: 
$$P(x,\Delta_x):= \sum_\alpha a_\alpha\Delta^\alpha_x.$$ De plus,  l'application \og symbole total\fg\ qui à un opérateur différentiel $P$ associe son symbole total: $$P\mapsto \sigma
_P(x,\xi):= \sum_{\alpha}a_\alpha\xi^{\alpha}$$ est un \em{isomorphisme} de \em{\bf $A\dag$-modules à gauche} entre l'anneau des opérateur différentiels
$D\dag_{A\dag/V}$ et l'algèbre
$(A\dag[\xi_1,\dots,\xi_n])\dag$ complétée $\dagger$-adique de l'algèbre $A\dag[\xi_1,\dots,\xi_n]$ pour la topologie $\goth m$-adique de $V$.\glossary{$P(x,\Delta_x):= \sum_\alpha
a_\alpha\Delta^\alpha_x$}\glossary{$\sigma
_P(x,\xi):= \sum_{\alpha}a_\alpha\xi^{\alpha}$}
\end{theo}

En fait, dans les articles ([M-N$_1$], [M-N$_3$]) on ne disposait pas du critère d'affinité et on imposait  en plus que le schéma $\dagger$-adique
induit
$\calXdag |U:=(U, \cal O_{\calXdag /V}|U)$ soit $\dagger$-affine. Le critère d'affinité lève cette restriction gênante.

L'isomorphisme précédent {\bf ne respecte pas}, bien entendu, la structure multiplicative et dépend des coordonnées locales. Le théorème du symbole total permet cependant de produire beaucoup
d'opérateurs différentiels parce que les éléments de l'algèbre
$(A\dag[\xi_1,\dots,\xi_n])\dag$ sont faciles à décrire. Le théorème du symbole total admet de nombreuses applications importantes. Par exemple, on a la trivialité cohomologique, qui jouera un rôle essentiel dans cet article:

\begin{coro}\label{acy-dif}Soit $\calUdag $ un schéma $\dagger$-adique affine lisse sur $V$ muni de coordonnées globales. Alors, les $V$-modules
$ H^{\bullet}(U, \calDdag_{\calUdag /V})$ sont nuls en degrés positifs.
\end{coro}

\demo Soient $A\dag$ l'algèbre de $\calUdag $, $A\dag\to (A\dag[\xi_1,\dots,\xi_n)])\dag$ l'extension naturelle,  $ f\dag:  (\cal T^*\cal U)\dag\to \calUdag $
la projection associée de schémas $\dagger$-adiques et $f: T^*U\to U$ la projection du fibré cotangent de $U$. En vertu du corollaire  \ref{acy} et du théorème du
symbole total, les images directes
$\bfR ^if_*\cal O_{(\cal T^*\cal U)\dag/V}$ sont nulles en degrés positifs et l'on a l'isomorphisme de faisceaux $V$-modules:
$$ \calDdag_{\calUdag /V}\simeq f_*\cal O_{(\cal T^*\cal
U)\dag/V}.$$ Il en résulte  les isomorphismes de $V$-modules:
$$ H^{i}(U, \calDdag_{\calUdag /V})\simeq H^{
i}(T^*U, \cal O_{(\cal T^*\cal U)\dag/V})$$ et l'on est ramené au théorème d'acyclicité \ref{acy} pour les schémas $\dagger$-adiques affines.
\enddemo
Cette démonstration figure  déjà dans l'article [M-N$_3$] dans le cas d'un schéma $\dagger$-adique affine $U$ et le critère d'affinité permet de supposer que $U$ est affine.
\begin{Rema} Tous les résultats de cet article qui reposent sur  le théorème du symbole total  ne sont démontrés que pour un couple $(V, \goth m)$ jusqu'à nouvel ordre.\end{Rema}

\begin{Rema} Dans la définition précédente on peut partir d'un $R$-schéma formel $\cal X^\wedge= (X,\cal O_{\cal X^\wedge/R})$ [EGA I] et définir le faisceau $\calDdag_{\cal
X^\wedge/R}$  comme sous-faisceau du faisceau des endomorphismes du faisceau structural et  en imposant par réduction modulo $I^s$ des conditions de croissances de
type
$\dagger$ aux ordres des opérateurs différentiels. Le théorème du symbole total dans le cas d'un  schéma formellement lisse sur un anneau noethérien $R$ a
lieu et établit de même un isomorphisme de
$\hat A$-modules à gauche entre
$D\dag_{\hat A/R}$ et le complété $\dagger$-adique
$(\hat A[\xi_1,\dots,\xi_n)])\dag$ pour la topologie $I$-adique de $\hat A$ mais la démonstration est élémentaire [M-N$_3$].
De ce point de vue la situation est meilleure.\end{Rema}

\section{Le site infinitésimal $\dagger$-adique $ \Xdaginf $ et son topos}
Soit $R$ un anneau noethérien muni de la topologie $I$-adique définie par un idéal $I\subset R$  et soit $(X, \cal O_{X/R_1})$ un $R_1$-schéma {\bf lisse} ([EGA IV$_4$], $\S 17$).

\subsection{Le site infinitésimal $\dagger$-adique}
Nous allons définir le site infinitésimal $\dagger$-adique de $X$  pour la topologie $I$-adique de $R$ que
nous noterons $ \Xdaginf $\glossary{$ \Xdaginf $}, le couple $I\subset R$ étant sous-entendu. Nous ne considérons pas dans cet article de façon essentielle les changements de base du couple
$(R, I)$.

\begin{defi}
\begin{liste}\item 1) Un objet de la catégorie $ \Xdaginf $ est un ouvert  $U$ de l'espace topologique sous-jacent à $X$ muni d'un relèvement $\calUdag $
\em{\bf plat} sur $R $, qui est donc un schéma 
$\dagger$-adique
$\calUdag  := (U, \cal O_{\calUdag /R})$ $R$-plat,  muni d'un morphisme surjectif 
$$\cal O_{\calUdag /R}\to \cal O_{U/R_1}\to0$$ de noyau $I\cal O_{\calUdag /R}$ où $(U,  \cal O_{U/R_1})$ est le schéma induit sur $U$.
\item 2) Un morphisme de la catégorie $ \Xdaginf $ est une inclusion d'ouverts  $r: W\hookrightarrow U$ munie d'un morphisme
$r^{-1}\cal O_{\calUdag /R}\to\cal O_{\calWdag /R}$ de faisceaux de $R$-algèbres  rendant le diagramme suivant commutatif:
$$\matrix{r^{-1}\cal O_{\calUdag /R}&\too&\cal O_{\calWdag /R}\cr
\downarrow&&\downarrow\cr
r^{-1}\cal O_{ U/R_1}&\too&\cal O_{ W/R_1}\,.}$$
\item 3) Un recouvrement d'un ouvert $ \calUdag $ est une famille de morphismes $\set
r\dag_\gamma:\cal U\dag_\gamma\to \calUdag , \gamma\in J/$ telle que
la famille $\set U_\gamma,
\gamma\in J/$ est un recouvrement de
$U$.
\end{liste}
\end{defi}
En vertu du théorème \ref{pro-fib}, la catégorie des schémas $\dagger$-adiques  sur $R$ admet des produits fibrés, et le produit fibré de deux
ouverts du site $ \Xdaginf $ est plat sur $R$. La catégorie $ \Xdaginf $ est stable par produit fibré:  le changement de base
d'un recouvrement est un recouvrement. Le composé d'un recouvrement par des recouvrements est un recouvrement.  La catégorie $ \Xdaginf $
se trouve munie d'une topologie de Grothendieck qui devient donc un site: le Site Infinitésimal $\dagger$-adique $ \Xdaginf $ du schéma $X$ pour la topologie $I$-adique de $R$. Le théorème des 
relèvements affines \ref{rel-alg} montre que pour $X$ lisse sur $R_1$ le site $ \Xdaginf $ a beaucoup d'ouverts et de morphismes et le critère d'affinité \ref{cri-aff} produit davantage  d'ouverts et de morphismes du site $ \Xdaginf $, ce qui lui donne une grande souplesse.

\begin{notation} Nous notons $\calUdag $ un objet du site infinitésimal $\dagger$-adique $ \Xdaginf $ et 
$r\dag: \calWdag \to \calUdag $ un morphisme du site $ \Xdaginf $.
Si
$r\dag:
\calWdag \to \calUdag $ est un morphisme  du site $ \Xdaginf $ nous noterons $r: W\to U$ l'inclusion topologique et 
$r^*$ le morphisme structural 
$r^{-1}\cal O_{\calUdag /R}\to\cal O_{\calWdag /R}$ de faisceaux de $R$-algèbres  sur $W$ qui est donc un morphisme local. 
\end{notation}

\begin{Rema}\begin{liste}\item 1) Un morphisme $r\dag: \calWdag \to \calUdag $ du site $ \Xdaginf $ se factorise par construction :$$\calWdag \to \calUdag |W\to \calUdag $$
où
$\calUdag |W\to \calUdag $ est le morphisme d'inclusion canonique. Le morphisme $\calWdag \to \calUdag |W$ est un {\bf isomorphisme}, en vertu du
théorème \ref{iso-pla}. En fait, c'est cet {\bf isomorphisme} d'espaces localement annelés qui est le point essentiel dans ce qui va suivre.\item 2) 
On aurait pu  dans la définition du site précédent, partir d'un schéma localement de type fini et plat sur $R_1$ et imposer 
aux objets d'être plats sur $R$, alors les morphismes du site sont des isomorphismes lorsque la réduction est une égalité.
\item 3) On aurait pu  dans la définition du site précédent partir 
d'un schéma localement de type fini  sur $R_1$, 
n'imposer aucune condition restrictive aux  relèvements $\dagger$-adiques mais imposer aux morphismes $\calWdag \to \calUdag |W$ d'être des
isomorphismes. 
\item 4) Mais aujourd'hui, nous  ne savons pas à quelles  conditions, en dehors du cas lisse, on obtient au-dessus d'un point singulier un site dont 
l'ensemble des objets est non vide  dans la situation du 2) 
ou dont l'ensemble des morphismes est non vide dans la situation du 3). 
C'est pour cette raison-là que nous supposerons  dans le présent  article 
que $X$ {\bf est lisse sur $R_1$}, bien que formellement cela n'est pas indispensable, et que certains développements se transposent.
\end{liste}
\end{Rema}

\begin{Rema}À cause de la première partie de  la remarque précédente, il est facile de voir que la catégorie $ X^{\dagger}_{\inf}$ admet des produit fibrés sans
invoquer le théorème général que la catégorie des schémas  $\dagger$-adiques admet des produits fibrés.
\end{Rema}

\begin{Rema}Pour tout $s\geq1$ on définit de même le site $ X^s_{\inf}$ dont les objets sont les relèvements $R_s$-plats des ouverts de $X$.
On a alors des foncteurs canoniques des catégories sous-jacentes  $\Xdaginf\to X^{s+1}_{\inf}\to X^{s}_{\inf}$, pour tout $s\geq1$, puisqu'un ouvert du site $\Xdaginf$ induit par réduction un ouvert du site  $X^{s+1}_{\inf}$ qui induit un ouvert du site $X^{s}_{\inf}$. Les méthodes de cet article  montrent que ces foncteurs  définissent des morphismes de topos abéliens.
\end{Rema}

\subsection{Le topos infinitésimal $\dagger$-adique}
Soient $\Prefais( { \Xdaginf })$, resp. $\Fais( X^{\dagger}_{\inf})$ les catégories des préfaisceaux, resp. des faisceaux, à
valeurs dans la catégorie des ensembles sur le site infinitésimal $\dagger$-adique $ \Xdaginf $.

\bigskip

\begin{defi}Soit $\calUdag $ un ouvert du site  $ \Xdaginf $. On a alors un foncteur de restriction:
$$R_{\calUdag }: \Prefais( { \Xdaginf })\to \Prefais(U)$$ de la catégorie des préfaisceaux d'ensembles sur le site infinitésimal
dans la catégorie des préfaisceaux d'ensembles sur l'espace topologique $U$ qui à un préfaisceau $\cal P_{\inf}$ associe le préfaisceau
$\cal P_{\calUdag }$ défini par:
$$\cal P_{\calUdag }(W):= \cal P_{\inf}(\calUdag |W)\,,$$ où $\calUdag |W$ est le schéma induit sur $W$ par $\calUdag $, qui est bien un objet
du site infinitésimal. Si $W'$ est un ouvert de $W$, le morphisme de restriction
$$\cal P_{\calUdag }(W)\to\cal P_{\calUdag }(W')$$ provient de l'inclusion d'espaces annelés  
$\calUdag |W'\subset\calUdag |W$.
\end{defi}

\begin{notation}Si $\cal F_{\inf}$ est un préfaisceau sur le site $ \Xdaginf $ et $\calUdag $ un ouvert, on note $\cal F_{\calUdag }:= R_{\calUdag }(\cal
F_{\inf})$ la restriction du préfaisceau $\cal F_{\inf}$ à l'ouvert $\calUdag $, c'est donc un préfaisceau sur $U$ pour la topologie de
Zariski.
\end{notation}

\begin{prop}\label{fai-str} Un préfaisceau d'ensembles $\cal F_{\inf}$ sur le site $ \Xdaginf $ est
un faisceau si et seulement si pour tout ouvert 
$\calUdag $  sa restriction $\cal F_{\calUdag }$ est un faisceau sur $U$.
\end{prop}

\demo  Soit $\cal F_{\inf}$ un préfaisceau sur le site infinitésimal et supposons que pour tout ouvert 
$\calUdag $  sa restriction $\cal F_{\calUdag }$ est un faisceau sur $U$. 
Soient  $ 
\{\phi_\alpha: \cal W\dag_\alpha\to U\dag\}$
un recouvrement   de
$\calUdag $  et
$\cup_\alpha  W_\alpha$ le recouvrement de $U$ induit. Notons $\cal U\dag_\alpha := \calUdag |W_\alpha$ le schéma $\dagger$-adique induit par $\cal
U\dag$ sur
$W_\alpha$.  Pour chaque $\alpha$ le morphisme  $\phi_\alpha:  \cal W\dag_\alpha\to \calUdag $ se factorise canoniquement 
par un morphisme : $\cal W\dag_\alpha\to \cal U\dag_\alpha$, qui est automatiquement un isomorphisme.
D'où des isomorphismes  ${\cal F_{\inf}}(\cal U\dag_\alpha)\simeq  {\cal F_{\inf}}(\cal W\dag_\alpha)$
qui fournissent    un  diagramme commutatif:
$$\matrix{ {\cal F}_{\inf}(\calUdag )&\to&\prod_{\alpha} {\cal F_{\inf}}(\cal U\dag_\alpha)&\rightrightarrows&
\prod_{\alpha,\beta} {\cal F}_{\inf}(\cal U\dag_{\alpha,\beta})\cr
\Vert&&\downarrow &&\downarrow\cr
{\cal F}_{\inf}(\calUdag )&\to&\prod_{\alpha} {\cal F_{\inf}}(\cal W\dag_\alpha)&\rightrightarrows& \prod_{\alpha,\beta} {\cal F}_{\inf}(\cal
W\dag_{\alpha,\beta})\,.}$$  Les colonnes sont bijectives. Mais la ligne supérieure est exacte par hypothèse et donc la ligne inférieure est exacte.
\enddemo

\begin{defi} On définit  ${\mathbf{Fam}}(\Ouv(\Xdaginf ))$ comme la catégorie des familles  des couples $(\calUdag , \cal F_{\calUdag })$ indexés   
par les ouverts $\calUdag $ 
munis d'un faisceau 
d'ensemble $\cal F_{\calUdag }$ sur l'espace topologique $U$ et des morphismes 
$(\sharp)_{r\dag}:f^{-1}\cal F_{\calUdag }\to \cal F_{\calWdag }$ indexés par les morphismes $r\dag$ tels que:
\begin{liste} 
\item 1) les morphismes $(\sharp)_{r\dag}$ sont transitifs, 
\item 2) on a l'égalité $\cal F_{\calUdag |W}= \cal F_{\calUdag }|W$,
\item 3) pour l'inclusion canonique
$r\dag: \calUdag |W\hookrightarrow \calUdag $ le morphisme $(\sharp)_{r\dag}$ est l'identité.
\end{liste}
\end{defi} On a alors un foncteur 
canonique de restrictions:$$ \Fais( X^{\dagger}_{\inf})\longrightarrow{\mathbf{Fam}}(\Ouv(\Xdaginf ))\leqno(\rm Rest):$$ qui à un faisceau $\cal F_{\inf}$ associe 
ses restrictions $\cal F_{\calUdag }:= R_{\calUdag }(\cal F_{\inf})$.
\glossary{$ \Fais( X^{\dagger}_{\inf}),{\mathbf{Fam}}(\Ouv(\Xdaginf ))$}

\begin{prop}\label{fai-fam} Le foncteur canonique précédent $(\rm Rest)$ est une équivalence de catégories. Les morphismes $(\sharp)_{r\dag}$ sont de plus des 
isomorphismes.
\end{prop}

\demo 
En effet,  une famille de faisceaux d'ensembles $\cal F_{\calUdag }$ indexée
par les ouverts $\calUdag $ et qui  appartient à la catégorie ${\mathbf{Fam}}(\Ouv(\Xdaginf ))$
définit un préfaisceau sur le site infinitésimal. La restriction de ce préfaisceau à tout ouvert du site est un faisceau de Zariski 
 sur l'ouvert topologique sous-jacent par hypothèse et, en vertu
de la proposition précédente,
ce préfaisceau est un faisceau. On obtient comme cela un inverse canonique au foncteur $(\rm Rest)$. La compatibilité avec la composition des morphismes montre que les morphismes $(\sharp)_{r\dag}$ sont des isomorphismes.
\enddemo

Cette équivalence canonique   produit beaucoup de faisceaux sur le site.\begin{exemples}\begin{liste}\item 1)
La famille $\cal O_{\calUdag /R}$ des faisceaux structuraux des objets du site infinitésimal a  la propriété de restriction par construction. Elle définit
donc 
un faisceau de
$R$-algèbres, que nous noterons \glossary{$ \cal O_{\Xdaginf /R}$}\smash{$ \cal O_{\Xdaginf /R}$}: le faisceau structural du site infinitésimal
$\dagger$-adique $ \Xdaginf $.
\item 2) Pour tout $j\geq0$ la famille
$\Omega^{j}_{\calUdag /R}$ des formes $j$-différentielles séparées a la propriété de restriction par construction. Elle
définit un faisceau de $ \cal O_{\Xdaginf /R}$-modules, que nous noterons \glossary{$\Omega^{j}_{\Xdaginf /R}$}$\Omega^{j}_{\Xdaginf /R}$: le
faisceau des formes
$j$-différentielles séparées du site infinitésimal $ \Xdaginf $.
\item 3) La famille
$\Omega^{\bullet}_{\calUdag /R}$ des complexes de de Rham des formes différentielles séparées a  la propriété de restriction par construction, la
différentielle $d$ commute aux restrictions \glossary{$(\sharp)_{r\dag}$}$(\sharp)_{r\dag}$. Elle définit un complexe de $R$-modules que nous noterons
\glossary{$ \Omega^{\bullet}_{\Xdaginf /R}$}\smash{$ \Omega^{\bullet}_{\Xdaginf /R}$}: le complexe de de Rham du site infinitésimal $ \Xdaginf $.
\item  4) La
famille $\cal T_{\calUdag /R}$  \glossary{$ \cal T_{\Xdaginf /R}$} des fibrés tangents, duaux des fibrés des $1$-formes différentielles séparées 
des objets du site infinitésimal a la propriété de restriction par construction. Elle définit donc 
un faisceau de
$ \cal O_{\Xdaginf /R}$-modules, que nous noterons $ \cal T_{\Xdaginf /R}$: le fibré tangent du site infinitésimal
$\dagger$-adique $ \Xdaginf $.
\end{liste}
\end{exemples}

\vskip1em
Si $r\dag: \calWdag \to \calUdag $ est un morphisme du site infinitésimal et si $P$ est un opérateur différentiel de 
$\calDdag_{\calUdag /R}$,
on définit $(\sharp)_{r\dag}(P)$
comme le composé
$$\cal O_{\calWdag /R}\stackrel{r^{*-1}}{\to}r^{-1}\cal O_{\calUdag /R}\stackrel{r^{-1}P}{\too}r^{-1}\cal O_{\cal
U\dag/R}\stackrel{r^*}{\to}\cal O_{\calWdag /R}$$ dont la réduction modulo $I^s$ est un opérateur différentiel d'ordre égal à celui de la réduction
modulo $I^s$ de $P$. C'est donc un opérateur différentiel de $\calDdag_{\calWdag /R}$ et on obtient 
le morphisme de restriction $$r^{-1}\calDdag_{\calUdag /R}\to\calDdag_{\calWdag /R}\leqno(\sharp)_{r\dag}:$$ 
compatible aux
compositions des morphismes du site. En vertu de la proposition précédente  \ref{fai-fam}, on obtient un faisceau sur le site:

\begin{defi}La famille  qui  à un ouvert $\calUdag $ associe le faisceau $\calDdag_{\calUdag /R}$ munie des morphismes de restriction $(\sharp)_{r\dag}$
définit un faisceau de
$ \cal O_{\Xdaginf /R}$-algèbres  
\glossary{$ \calDdag_{\Xdaginf /R}$} que nous noterons $ \calDdag_{\Xdaginf /R}$: le faisceau des opérateurs différentiels
$\dagger$-adiques sur le site $\Xdaginf $. 
\end{defi}

\begin{prop}La famille  qui  à un ouvert $\calUdag $ associe le faisceau $\calD_{\calUdag /R}$ des opérateurs différentiels d'ordre localement fini définit un sous-faisceau du
$ \cal O_{\Xdaginf /R}$-algèbres \glossary{$ \calD_{\Xdaginf /R}$}$ \calD_{\Xdaginf /R}$ de faisceau $ \calDdag_{\Xdaginf /R}$. 
\end{prop}

\demo
Il s'agit de montrer que les restrictions $(\sharp)_{r\dag}$ laissent stables les opérateurs différentiels d'ordre fini, mais c'est un conséquence du fait
que
les morphismes d'algèbres $r^*$ sont des isomorphismes. En fait, les restrictions $(\sharp)_{r\dag}$ conservent l'ordre des opérateurs différentiels.
\enddemo 

\begin{Rema}\begin{liste}\item 1) Les  raisonnements que nous avons faits pour les faisceaux d'ensembles s'appliquent tout aussi bien pour les
faisceaux à valeurs dans une catégorie abélienne ou pour les modules sur un faisceau d'anneaux sur le site infinitésimal. Nous utiliserons librement les résultats précédents dans ces contextes. \item 2) Par exemple, si $\cal
A_{\Xdaginf }$ est faisceau d'anneaux sur le site infinitésimal, un faisceau  $\cal F_{\inf}$ de $\cal A_{\Xdaginf }$-modules est la donnée pour tout ouvert $\calUdag $
d'un faisceau de 
$\cal A_{\calUdag }$-modules $\cal F_{\calUdag }$ tel que $\cal F_{\calUdag |W}= \cal F_{\calUdag }|W$ et pour tout morphisme $r\dag$ d'un
morphisme de restriction compatible à la composition 
$$\resetdisplay
r^{-1}\cal F_{\calUdag }\to\cal F_{\calWdag }\leqno(\sharp)_{r\dag}:$$ qui est {\bf $(\sharp)_{r\dag}$-linéaire}:
$$(\sharp)_{r\dag}(P(m))=(\sharp)_{r\dag}(P)((\sharp)_{r\dag}(m))$$ pour une section $P$ de $r^{-1}\cal A_{\calUdag }$ et une section $m$ de
$r^{-1}\cal F_{\calUdag }$ et qui se réduit à l'application identique dans le cas de l'inclusion  $\calUdag |W\subset \calUdag $.
\end{liste}
\end{Rema}

\begin{notation}Si $\cal
A_{\Xdaginf }$ est un faisceau d'anneaux sur le site, on note $ \cal
A_{\Xdaginf }\Mod$ la catégorie des $\cal
A_{\Xdaginf }$-modules à gauche sur le site et $ \Modd\cal
A_{\Xdaginf }$ la catégorie des $\cal
A_{\Xdaginf }$-modules à droite sur le site.
\end{notation}

\subsection{Le site infinitésimal $\dagger$-adique affine}
Dans le cas où le schéma $X$ est \em{séparé} il est plus commode d'utiliser le site infinitésimal $\dagger$-adique affine.

\begin{defi}Soit un $R_1$-schéma $X$ lisse et \em{séparé}. On note \glossary{$ X^{\dagger,\aff}_{\inf}$}$ 
X^{\dagger,\aff}_{\inf}$ la sous-catégorie pleine de la catégorie
$ \Xdaginf $ des relèvements des ouverts affines de $X$. Muni de la topologie induite $ X^{\dagger,\aff}_{\inf}$ devient un site, le site infinitésimal
$\dagger$-adique affine de $X$.
\end{defi}
Soit $\Fais( X^{\dagger,\aff}_{\inf})$ la catégorie des faisceaux d'ensembles sur le site affine.

\begin{theo}Soit $X$ un $R_1$-schéma  lisse et \em{séparé}. Alors, 
le foncteur naturel de restriction:
$$\Fais( X^{\dagger}_{\inf})\to \Fais( X^{\dagger,\aff}_{\inf})$$ est une équivalence de catégories.
\end{theo}

\demo Soient $\cal F$ un faisceau d'ensembles sur le site infinitésimal affine et  $\calUdag $ un ouvert du
site infinitésimal. Soit $U_\alpha, \alpha\in I$ un recouvrement par des ouverts affines   de $U$.
En vertu du critère d'affinité \ref{cri-aff},
si $U_\alpha$ est un ouvert affine  de $U$, l'espace annelé $\cal U\dag_\alpha:= (U_\alpha,
\cal O_{\calUdag /R}|U_\alpha)$ est un ouvert du site infinitésimal affine. Soient alors les deux  morphismes de restriction: 
$$\prod_{\alpha\in I}\cal F(\cal U\dag_\alpha)\varrightarrows2{15pt}{15pt} \prod_{\alpha,\beta\in I}\cal F(\cal U\dag_{\alpha,\beta}
)$$  de noyau  (ou lieu des coïncidences) noté
$\Ker(I)$. On définit:
$$   \tilde{\cal F}(\calUdag ):= \Ker(I).$$ Cet ensemble ne dépend pas du recouvrement  $U_\alpha, \alpha\in I$. En effet, soit  $W_\beta,
\beta\in J$ un autre recouvrement; alors pour tout $\alpha\in I$ les ouverts affines $U_{\alpha,\beta}:= U_\alpha\cap W_\beta$ forment un 
recouvrement  de 
$U_\alpha$. Si $s_\beta, \beta\in J$ est une famille de sections de $\cal F$ qui se recollent, leurs restrictions aux ouverts $U_{\alpha, \beta}$
forment une famille de sections qui se recollent et par hypothèse définissent une section $s_\alpha$ au-dessus de l'ouvert $\cal U\dag_\alpha$. La
famille $s_\alpha, \alpha\in I$ se recolle et définit un élément de l'ensemble $\Gamma(\calUdag ,   \tilde{\cal F})$, ce qui fournit une
identification
entre l'ensemble construit à partir du recouvrement $U_\alpha, \alpha\in I$ et celui construit  à partir du recouvrement $W_\beta, \beta\in J$.

Soit
$r\dag :
\calWdag \to\calUdag $ un morphisme du site et soit un recouvrement $U_\alpha, \alpha\in I$ de $U$ tel que pour un sous-ensemble $J\subset I$ les
ouverts
$W_\alpha:= U_\alpha, \alpha\in J$ forment un recouvrement de $W$. On obtient  pour tout $\alpha\in J$ un morphisme
du site affine
$\cal W\dag_\alpha\to\cal U\dag_\alpha$. Si $s_\alpha, \alpha\in I$ est une famille de sections qui se recollent, alors leurs restrictions
aux ouverts affines $\cal W\dag_\alpha$ est une famille de sections qui se recollent.
D'où un morphisme de restriction:
$$  \tilde{\cal F}(\calUdag )\to    \tilde{\cal F}(\calWdag )$$ et on obtient de façon évidente un préfaisceau d'ensembles  sur le
site infinitésimal. 

Il s'agit de montrer que ce préfaisceau est un faisceau. Soit
$\cup_{\alpha\in I}\cal U\dag_\alpha$ un recouvrement de l'ouvert $\calUdag $ et $\cup_{\gamma\in I_\alpha}U_{\alpha,\gamma}$ un recouvrement de l'ouvert
$U_\alpha$ par des ouverts affines. Le morphisme  $\cal U\dag_{\alpha}\to \calUdag |U_\alpha$ induit un morphisme du site affine: 
$$\cal U\dag_{\alpha}|U_{\alpha,\gamma}\to \calUdag |U_{\alpha,\gamma}$$ qui est un isomorphisme en vertu de \ref{iso-pla} et induit une bijection
d'ensembles
$$\cal F(\cal U\dag_{\alpha}|U_{\alpha,\gamma})\to  \cal F(\calUdag |U_{\alpha,\gamma}).$$ On obtient un diagramme commutatif:
$$\matrix{\prod_{\alpha}\tilde {\cal F}(\cal U\dag_\alpha)&\varrightarrows2{0pt}{12ex} & \prod_{\alpha,\beta}
\tilde
{\cal F}(\cal U\dag_{\alpha,\beta})\cr
\downarrow&&\downarrow\cr
\prod_{\alpha,\gamma}\cal F(\cal U\dag_\alpha|U_{\alpha,\gamma})&\rightrightarrows  &\prod_{\alpha,\gamma,\beta,\delta}\cal F(\cal
U\dag_{\alpha,\beta}|U_{\alpha,\gamma,\beta,\delta})\cr
\uparrow&&\uparrow\cr
\prod_{\alpha,\gamma}\cal F(\calUdag |U_{\alpha,\gamma})&\rightrightarrows & \prod_{\alpha,\gamma,\beta,\delta}\cal F(\cal
U\dag|U_{\alpha,\gamma,\beta,\delta})}$$ 
dont les morphismes verticaux des deux dernières lignes sont des isomorphismes. Le noyau de la
dernière ligne
est l'ensemble $\Gamma(\calUdag ,\tilde {\cal F})$. Il suffit de montrer que le morphisme induit sur les noyaux des deux premières lignes est
bijectif. Il est injectif de façon évidente. Un élément du noyau de la deuxième ligne définit une famille de sections de $\tilde{\cal F}$ au-dessus
des ouverts $\cal U\dag_\alpha$ qui se recollent. Le morphisme est surjectif.
Cela montre que le diagramme:
$$\tilde {\cal F}(\calUdag )\too
\prod_{\alpha}\tilde {\cal F}(\cal U\dag_\alpha)\varrightarrows2{3ex}{3ex} \prod_{\alpha,\beta}\tilde {\cal F}(\cal
U\dag_{\alpha,\beta})$$ est exact et le  préfaisceau $\tilde{\cal F}$ est un faisceau. Le foncteur  de restriction du théorème est
essentiellement surjectif.

\bigskip
De même, on voit facilement qu'un morphisme de faisceaux d'ensembles sur le site infinitésimal est uniquement déterminé par sa restriction aux faisceaux
d'ensembles sur le site infinitésimal affine. Le morphisme de restriction du théorème est pleinement fidèle.
\enddemo

\section{Le faisceau  $ \cal G_{\Xdaginf }$ des automorphismes  du faisceau structural qui se réduisent à l'identité modulo $I$} 

\subsection{Le faisceau  d'ensembles de transfert $\cal G_{\calYdag \to\calXdag }$}
\glossary{$\cal G_{\calYdag \to\calXdag }$}Soient $\calYdag $, $\calXdag $ des relèvements  des
$R_1$-schéma  $Y, X$ et  $f: Y\to X$ un morphisme de $R_1$-schémas.

On peut considérer  le faisceau d'ensembles sur $Y$ des
morphisme d'espaces d'annelés de $\calYdag $ dans $\calXdag $ qui se réduisent modulo $I$ au morphisme $f$. 
Mais de façon équivalente il est plus commode pour pour nous de faire la définition suivante qui fait le lien avec le calcul différentiel des chapitres 
qui suivent.

\begin{defi}Le faisceau $\cal G_{\calYdag \to\calXdag }$ d'ensembles sur $Y$ {\bf de transfert}  est le sous-faisceau du
faisceau 
$\cHom_{R}(f^{-1}\cal O_{\calXdag /R},
\cal O_{\calYdag /R})$ des morphismes de {\bf $R$-algèbres}  qui se réduisent modulo l'idéal $I$ au morphisme structural de $R_1$-algèbres
$f_1^*: f^{-1}\cal O_{ X/R_1}\to\cal O_{Y/R_1}$.
\end{defi}
Une section globale du faisceau d'ensembles  $\cal G_{\calYdag \to\calXdag }$ est par construction un morphisme d'espaces localement annelés.
\begin{notation}\begin{liste}\item 1) Dans la notation $\cal G_{\calYdag \to\calXdag }$ le morphisme
$f$ de $R_1$-schémas est sous-entendu quand il n'y a pas de risque de confusion.
\item 2) Notons  $\cal G_{\calXdag }$ \glossary{$\cal G_{\calXdag }$}le faisceau de
transfert correspondant au morphisme identique de $X$.
\item 3) Notons $R[\cal G_{\calYdag \to\calXdag }]$ \glossary{$R[\cal G_{\calYdag \to\calXdag }]$}le faisceau  de $R$-modules  engendré par le
$R$-module libre des sections   du faisceau de transfert 
$\cal G_{\calYdag \to\calXdag }$.
\item 4) Le module $R[\cal G_{\calYdag \to\calXdag }]$ est naturellement un 
$(R[\cal G_{\calYdag }], f^{-1}R[\cal G_{\calXdag }])$-bimodule.
\end{liste}
\end{notation}

Le cas le plus important est le cas où $Y$ et $X$ sont affines d'algèbres $B, A$ et $\calYdag $, $\calXdag $ des relèvements d'algèbres
$B\dag, A\dag$. En vertu de l'équivalence \ref{eqi-aff},
les sections globales du faisceau $\cal G_{\calYdag \to\calXdag }$ correspondent aux morphismes de $R$-algèbres
$A\dag\to B\dag$ qui se réduisent modulo $I$ au morphisme structural $f: A\to B$.

\begin{notation} Si $A\dag$ est une $R$-algèbre, on note $G_{A\dag}$ \glossary{$G_{A\dag}$}l'ensemble des morphismes de $R$-algèbres $A\dag\to A\dag$
qui se réduisent modulo $I$ à l'identité.
\end{notation} 

\begin{prop}Si $A\dag$ est $R$-plate alors l'ensemble $G_{A\dag}$ est un groupe. Plus généralement, si $\calXdag $ est plat sur $R$, le faisceau
$\cal G_{\calXdag }$ d'ensembles est un faisceau de groupes.
\end{prop}

\demo C'est une conséquence du théorème \ref{rel-pla} et du théorème \ref{iso-pla}.
\enddemo

\begin{defi}Si $r\dag$ est un morphisme du site infinitésimal $ \Xdaginf $, le morphisme de faisceaux $r^* : r^{-1}\cal O_{\calUdag /R}\to \cal O_{\calWdag /R}$
est un isomorphisme de faisceaux de $R$-algèbres, en vertu de \ref{iso-pla}. On définit le morphisme de restriction
$$r^{-1}\cal G_{\calUdag }\to\cal G_{\calWdag }\leqno (\sharp)_{r\dag}:$$ par $g\mapsto r^*\circ g\circ r^{*-1}$. La famille $\cal G_{\cal
U\dag}$ et les morphismes $(\sharp)_{r\dag}$ ont la propriété \ref{fai-fam}. Ils définissent un  faisceau de groupes sur le site infinitésimal.
\end{defi}

\begin{notation}On note $ \cal G_{\Xdaginf }$ \glossary{$ \cal G_{\Xdaginf }$} le faisceau de groupes ainsi défini;  c'est le faisceau des
automorphismes de $R$-algèbres du faisceau structural
$ \cal O_{\Xdaginf /R}$ qui se réduisent modulo $I$ à l'identité.
\end{notation}

\begin{prop}\label{act-gro}Un préfaisceau $\cal P_{\inf}$, resp. un faisceau $\cal F_{\inf}$, d'ensembles  sur le site infinitésimal $ \Xdaginf $ est un $ \cal G_{\Xdaginf }$-préfaisceau, resp. faisceau, 
à {\bf gauche} d'ensembles.
\end{prop}

\demo Soient $\calUdag $ un ouvert du site et $ \calUdag \to \calUdag $ un morphisme du site $X\daginf$ défini par $g: \cal O_{\calUdag /R}\to  \cal
O_{\calUdag /R}$ qui est donc un isomorphisme. La restriction
$(\sharp)_g$ induite par $g$ du préfaisceau $\cal P_{\inf}(\calUdag )\to \cal P_{\inf}(\calUdag )$ est donc bijective. Pour deux morphismes $g_1, g_2$
on a $(\sharp)_{g_2g_1}= (\sharp)_{g_2}(\sharp)_{g_1}$. Mais  par construction,
$g$ est une section globale du faisceau  $ \cal G_{\Xdaginf }$ au-dessus de l'ouvert $\calUdag $. Cela définit une action à gauche de $\cal G_{\Xdaginf }(\calUdag )$ sur $\cal P_{\inf}(\calUdag )$. Cette action induit
une action  à gauche de  $\cal G_{\calUdag }$ sur la restriction $\cal P_{\calUdag }$.

Si $r\dag : \calWdag \to \calUdag $ est un morphisme du site, on a un diagramme commutatif:
$$\matrix{\cal P_{\inf}(\calUdag |W)&\varto{0pt}{10ex}^{(\sharp)_g}
&\cal P_{\inf}(\calUdag |W)\cr\downarrow 
(\sharp)_{r\dag}&&\downarrow 
(\sharp)_{r\dag}\cr
\cal P_{\inf}(\calWdag )&
\smash{\varto{0pt}{10ex}^{(\sharp)_{r^*gr^{*-1}}}}
&\cal P_{\inf}(\calWdag )}$$ qui montre que l'action commute aux restrictions.
Les mêmes arguments valent dans le cas du faisceau $\cal F_{\inf}$. Par construction, si $\calUdag $ est un ouvert du site et si 
$\cal F_{\inf}$ est faisceau d'ensemble sur le site sa restriction $R_{\calUdag }(\cal F_{\inf})$ est un $\cal G_{\calUdag }$-module à gauche.
\enddemo

\sousparagraphe{}Un peu plus généralement, si $f: Y\to X$ est un morphisme de $R_1$-schémas lisses et $\cal A_{\calYdag }$ est un faisceau de $R$-algèbres sur $Y$, on peut considérer le 
$\cal A_{\calYdag }$-module à gauche $\cal A_{\calYdag }[\cal G_{\calYdag \to\calXdag }]$ engendré  par le  module libre des sections   du
faisceau de transfert 
$\cal G_{\calYdag \to\calXdag }$. Si de plus $\cal A_{\calXdag }$ est un $\cal G_{\calXdag }$-module à gauche
le module $\cal A_{\calXdag }[\cal G_{\calXdag }]$ est naturellement un faisceau d'algèbres où le produit est le produit croisé défini 
par $P_1g_1P_2g_2:=P_1\rho_{g_1}(P_2)g_1g_2 $, où $\rho$ est la représentation de $\cal G_{\calXdag }$ sur $\cal A_{\calXdag }$. 

 En particulier,
si $ \cal A_{\Xdaginf }$ est un faisceau de $R$-algèbres sur le site $\Xdaginf $, 
pour tout ouvert $\calUdag $  l'algèbre $\cal A_{\calUdag }$ est un $\cal G_{\calUdag }$-module à gauche et le module 
$\cal A_{\calUdag }[\cal G_{\calUdag }]$ est muni du produit croisé. Si $\cal F_{\inf}$ est un $ \cal A_{\Xdaginf }$-module à gauche,
sa restriction $\cal F_{\calUdag }$ est un $\cal A_{\calUdag }[\cal G_{\calUdag }]$-module à gauche.
Le faisceau de transfert\glossary{$\cal A_{\calYdag }[\cal
G_{\calYdag \to\calXdag }]$}
$\cal A_{\calYdag }[\cal G_{\calYdag \to\calXdag }]$ est naturellement un $(\cal A_{\calYdag }[\cal G_{\calYdag }], f^{-1}\cal A_{\calXdag}[\cal G_{\calXdag }])$-bimodule où l'action à droite est définie  par $r\dag P := (\sharp)_{r\dag}(P)r\dag$.

\subsection{La catégorie $ \cal A_{\Xdaginf }\Mod$ des  $\Rm\cal A_{\Xdaginf }$-modules à gauche sur le site $\dagger$-adique}\glossary{$\cal
A_{\Xdaginf }$}Soit $ \cal A_{\Xdaginf }$ un faisceau de $R$-algèbres sur le site infinitésimal. On dispose donc de la catégorie  $ \cal
A_{\Xdaginf }\Mod$ des $ \cal A_{\Xdaginf }$-modules à gauche.
Soit $\cal F_{\inf}$ un faisceau de $ \cal A_{\Xdaginf }$-modules à gauche sur le site $ \Xdaginf $. Si $\calUdag $ est un ouvert du site, alors
sa restriction $\cal F_{\calUdag }$ est un faisceau de Zariski de $ \cal A_{\calUdag }$-modules à gauche sur
$U$, sur le lequel le faisceau de groupes  $\cal G_{\calUdag }$ agit à gauche. C'est donc un module sur l'algèbre du groupe $ \cal A_{\calUdag }[\cal G_{\calUdag }]$. Nous allons donner une définition équivalente des objets de la catégorie 
$ \cal A_{\Xdaginf }\Mod$, qui nous servira de transition
à la catégorie des  modules spéciaux.

\begin{prop}\label{defi'}La donnée d' un $ \cal A_{\Xdaginf }$-module à gauche 
$\cal F_{\inf}$ sur le site $ \Xdaginf $ est la donnée,  pour tout ouvert
$
\calUdag $ du site 
$\Xdaginf $, d'un faisceau $\cal F_{\calUdag }$ de
$\cal A_{\calUdag }[\cal G_{\calUdag }]$-modules à gauche, et la donnée, pour 
 pour tout couple $(\calWdag ,\calUdag )$ d'objets du site $ \Xdaginf $, 
 avec $r: W\hookrightarrow U$, d'un morphisme de
$\cal A_{\calWdag }[\cal G_{\calWdag }]$-modules à gauche 
\glossary{$(\sharp)$}\glossary{\vadjust{\kern-6pt}$\cal A_{\calWdag }[\cal G_{\calWdag \rightarrow \cal
U\dag}]\Otimes_{r^{-1}\cal A_{\calUdag }[\cal G_{\calUdag }]}r^{-1}\cal F_{\calUdag }\simeq \cal F_{\calWdag }$
\hbox to4.1cm{}}
$$\cal A_{\calWdag }[\cal G_{\calWdag \rightarrow \calUdag }]\otimes_{r^{-1}\cal A_{\calUdag }[\cal G_{\calUdag }]}r^{-1}\cal F_{\cal
U\dag}\simeq \cal F_{\calWdag }.\leqno (\sharp):$$ 
En outre, ces données doivent satisfaire les conditions:
\begin{liste}\item 1) pour un couple $(\cal U\dag|W, \cal U\dag)$, avec  $W\hookrightarrow U$,  on a $\cal F_{\calUdag |W}= \cal F_{\calUdag }|W$, et le morphisme $(\sharp)$ coïncide avec le morphisme canonique,
\item 2) pour un triplet $(\cal W'{}\dag ,\calWdag ,\calUdag)$, 
avec $r': W'\hookrightarrow W$ et $r: W\hookrightarrow
U$, le diagramme suivant est commutatif:
$$\let\bigotimes\otimes
\def\quad{\hskip0.1ex}\scriptwd0.8em\mathrigid0mu
\scriptspace0.pt\def\ss{\sb}
\hss\matrix{\cal A\ss{\cal W'{}\dag }[\cal G_{\cal
W'{}\dag \sto
\calWdag }]\Otimes_{r'^{-1}\cal A\ss{\calWdag }[\cal G_{\calWdag }]}r'^{-1}\mBig(\cal A\ss{\cal
W\dag}[\cal G_{\calWdag \sto \calUdag }]\Otimes_{r^{-1}\cal A\ss{\calUdag }[\cal
G_{\calUdag }]} r^{-1}\cal F_{\calUdag }\mBig)&\to &\cal A\ss{\cal W'{}\dag }[\cal
G_{\cal W'{}\dag \sto
\calUdag }]\Otimes_{(r'\circ r)^{-1}\cal A\ss{\calUdag }[\cal G_{\calUdag }]} (r'\circ r)^{-1}\cal F_{\calUdag }\cr
\downarrow&&\downarrow\cr
\cal A\ss{\cal W'{}\dag }[\cal G_{\cal W'{}\dag \rightarrow \calWdag }]\Otimes_{r'^{-1}\cal A\ss{\calWdag }[\cal G_{\calWdag }]} r'^{-1}\cal
F_{\calWdag }&\varto{0pt}{20ex} &\cal F_{\cal W'{}\dag }\,.}\hss$$
\end{liste}
\end{prop}
\demo
Soit $\cal F_{\inf}$ un 
$ \cal A_{ \Xdaginf }$-module à gauche.  En vertu des propositions \ref{fai-fam} et \ref{act-gro}, $\cal F_{\inf}$ est la donnée d'une famille de $\cal A_{\calUdag }[\cal G_{\calUdag }]$-modules
à gauche $\cal F_{\calUdag }$ paramétrée par les ouverts $\calUdag $, et de plus, pour chaque morphisme du site $r\dag: \calWdag \to \calUdag $, d'un morphisme de
restriction
$$ r^{-1}\cal F_{\calUdag }\to \cal F_{\calWdag }.\leqno(\sharp)_{r\dag} :$$

Soient  $\calWdag $, $\calUdag $ deux ouverts tels que $W\subset U$,  et $r^*$ une section globale
au-dessus de $W$ du faisceau de transfert $\cal G_{\calWdag \rightarrow 
\calUdag }$ qui définit une restriction $(\sharp)_{r\dag}: r^{-1}\cal F_{\calUdag }\to \cal F_{\calWdag }$. On a  par linéarité
un morphisme de
$\cal A_{\calWdag }[\cal G_{\calUdag }]$-modules à gauche
$$\cal A_{\calWdag }[\cal G_{\calWdag \rightarrow \calUdag }]\otimes_{r^{-1}\cal A_{\calUdag }[\cal G_{\calUdag }]}r^{-1}\cal F_{\cal
U\dag}\simeq \cal F_{\calWdag }\leqno (\sharp):$$  qui a les propriétés de la proposition. La transitivité des morphismes $(\sharp)_{r\dag}$ entraîne la
transitivité des morphismes $(\sharp)$.

Réciproquement, une famille de $\cal A_{\calUdag }[\cal G_{\calUdag }]$-modules à gauche munie de morphismes 
$(\sharp)$ ayant les propriétés de la proposition définit une famille  munie de morphismes $(\sharp)_{r\dag}$
compatibles à la compositions des morphismes ayant les propriétés de la proposition \ref{fai-fam} et 
définit un objet de la catégorie $\RmMod(\cal A_{\Xdaginf })$.
\enddemo

\begin{theo}\label{rig-spe}Soient  $ \cal A_{\Xdaginf }$ un faisceau d'anneaux sur le site $ \Xdaginf $, $\cal F_{\inf}$ un $ \cal
A_{\Xdaginf }$-module à gauche  et
$\calWdag $ et $\calUdag $
deux ouverts du site
tels que $W\subset U$. Alors, le morphisme de $\cal A_{\calWdag }[\cal G_{\calWdag }]$-modules à gauche 
$(\sharp)$ de la proposition précédente 
$$\cal A_{\calWdag }[\cal G_{\calWdag \rightarrow \calUdag }]\otimes_{r^{-1}\cal A_{\calUdag }[\cal G_{\calUdag }]}r^{-1}\cal F_{\cal
U\dag}\simeq \cal F_{\calWdag }\leqno (\sharp):$$  est un {\bf isomorphisme}.
\end{theo}

\demo
Soit $T\subset W$ un ouvert au-dessus duquel le morphisme d'inclusion  se relève en une section $r^*_T$ de $\cal G_{\calWdag \rightarrow \calUdag }$. Alors, le morphisme $(\sharp)$ est un isomorphisme au-dessus de $T$. En effet, si $w$ est une section de $\cal F_{\calWdag }$
au-dessus de $T$, elle est l'image de $r^*_T\otimes (\sharp)^{-1}_{r^*_T}( w)$ par le morphisme $(\sharp)$, ce qui montre qu'il est surjectif. Si l'image
par le morphisme $(\sharp)$ d'une section le forme $r^*_T\otimes u$ pour une section $u$ de $r^{-1}\cal F_{\calUdag }$ est nulle, nécessairement
$u$ est nulle, mais $r^*_T$ engendre $\cal A_{\calWdag }[\cal G_{\calWdag \rightarrow \calUdag }]$ comme $r^{-1}\cal A_{\calUdag }[\cal G_{\calUdag }]$-module  libre à droite au-dessus de
$T$. Cela montre que le morphisme $(\sharp)$ est injectif.

Soit $T$ un voisinage affine d'un point  de $W$. En vertu du critère d'affinité, les schémas $\dagger$-adiques $\calWdag |T$ et $\calUdag |T$ 
sont affines et, 
en vertu du 
théorème \ref{rel-mor}, le morphisme d'inclusion se relève en une section $r^*_T$ de $\cal G_{\calWdag \rightarrow \calUdag }$ au-dessus de $T$. Le morphisme $(\sharp)$ est un isomorphisme au voisinage des tous
les points de $W$.
\enddemo

\begin{defi}Soient $X$ un $R_1$-schéma lisse, $ \cal A_{\Xdaginf }$ un faisceau d'anneaux sur le site $ \Xdaginf $,  $\calXdag $  un relèvement
 plat sur $R$ de  $X$ et $\cal F_{\calXdag }$ un
$\cal A_{\calXdag }[\cal G_{\calXdag }]$-modules à gauche. Si $\calUdag $ est un ouvert du site de $X$ on définit le $\cal A_{\calUdag }[\cal G_{\calUdag }]$-module à gauche:
$$P_{ \cal A_{\Xdaginf },\calXdag }(\cal F_{\calXdag })_{\calUdag }:= \cal A_{\calUdag }[\cal G_{\calUdag \rightarrow \calXdag}]\otimes_{r^{-1}\cal A_{\calXdag }[\cal G_{\calXdag }]}r^{-1} 
\cal F_{\calXdag }.$$\glossary{$P_{ \cal A_{\Xdaginf },\calXdag }$}Si $\calWdag , \calUdag $ sont deux ouverts tels que $W\subset U$, on définit le morphisme
$$\cal A_{\calWdag }[\cal G_{\calWdag \rightarrow \calUdag }]\otimes_{r^{-1}\cal A_{\calUdag }[\cal G_{\calUdag }]}r^{-1}P_{\cal A_{\Xdaginf }, \calXdag }(\cal F_{\calXdag })_{\calUdag }\to P_{\cal A_{\Xdaginf }, \calXdag }(\cal F_{\calXdag })_{\calWdag }\leqno (\sharp):$$ provenant, par produit tensoriel,
du morphisme naturel provenant de la composition:
$$\cal A_{\calWdag }[\cal G_{\calWdag \rightarrow \calUdag }]\otimes_{r^{-1}\cal A_{\calUdag }[\cal G_{\calUdag }]}r^{-1}\cal A_{\cal
U\dag}[\cal G_{\calUdag \rightarrow \calXdag }]\to  \cal A_{\calWdag }[\cal G_{\calWdag \rightarrow \calXdag }].$$ La famille
$P_{\cal A_{\Xdaginf }, \calXdag }(\cal F_{\calXdag })_{\calUdag }$ et les morphismes $(\sharp)$ satisfont les conditions de la proposition précédente. Ils
définissent donc un
$\cal A_{\Xdaginf }$-module à gauche $P_{\cal A_{\Xdaginf }, \calXdag }(\cal F_{\calXdag })$ sur le site infinitésimal.

On obtient 
le foncteur prolongement
$$ \cal A_{\calXdag }[\cal G_{\calXdag }]\Mod\to
\cal A_{ \Xdaginf }\Mod.\leqno  P_{\cal A_{\Xdaginf }, \calXdag }: $$
\end{defi}
\begin{lemm} Le foncteur $ P_{\cal A_{\Xdaginf }, \calXdag }$ est exact.\end{lemm}
\demo En effet  le $r^{-1}\cal A_{\calXdag }[\cal G_{\calXdag }]$-module à droite $\cal A_{\calUdag }[\cal G_{\calUdag \rightarrow \calXdag }]$ est localement libre de rang $1$.
\begin{exemple}\begin{liste}\samepage\item 1) Le prolongement $ P_{R_{\Xdaginf }, \calXdag }(\cal O_{\calXdag /R})$ du faisceau structural est le faisceau
structural
$\Rm
\cal O_{\Xdaginf /R}$ du site.\item 2) Le complexe de de Rham $\Omega^{\bullet}_{\calXdag /R}$ est un complexe de $R[\cal G_{\calXdag }]$-modules à gauche
et son prolongement $ P_{R_{\Xdaginf }, \calXdag }(\Omega^{\bullet}_{\calXdag /R})$ est, par construction, le complexe de de Rham $ \Omega^{\bullet}_{
\Xdaginf /R}$ du site.
\end{liste}
\end{exemple}

\begin{theo}\label{equ-god}Soit  $\calXdag $  un relèvement \em{plat} sur $R$ d'un  schéma  $X$ lisse  sur $R_1$. Alors, le foncteur
prolongement 
$ P_{\cal A_{\Xdaginf }, \calXdag }$ est une équivalence de catégories
entre les catégories $ \cal A_{\calXdag }[\cal G_{\calXdag }]\Mod$ et  $ \cal A_{ \Xdaginf }\Mod$ 
qui est un inverse \em{\bf canonique} du foncteur de restriction  naturel $R_{\calXdag }$.
\end{theo}

\demo
Il est évident que si $\cal F_{\calXdag }$ est un $\cal A_{\calXdag }[\cal G_{\calXdag }]$-module à gauche, la restriction à $\calXdag $ 
du  $ \cal A_{\Xdaginf }$-module à gauche 
$ P_{\cal A_{\Xdaginf },\calXdag }\cal (\cal F_{\calXdag})$ coïncide avec le $ \cal A_{\calXdag }[\cal G_{\calXdag }]$-module à gauche $\cal F_{\calXdag }$. On a donc l'isomorphisme canonique de
foncteurs $ R_{\calXdag }\circ P_{\cal A_{\Xdaginf }, \calXdag }\simeq  Id$. Soit $\cal F_{\inf}$ un $ \cal A_{\Xdaginf }$-module à gauche et
notons
$\cal F_{\calXdag }$ sa restriction à
$\calXdag $. Si $\calUdag , r: U\hookrightarrow X,$ est un ouvert du site, la valeur du module $P_{\cal A_{\Xdaginf }, \calXdag }(\cal F_{\calXdag}) $ au-dessus de $\calUdag $ est :
$$P_{ \cal A_{\Xdaginf },\calXdag }(\cal F_{\calXdag })_{\calUdag }:= \cal A_{\calUdag }[\cal G_{\calUdag \rightarrow \calXdag}]\otimes_{r^{-1}\cal A_{\calXdag }[\cal G_{\calXdag }]}r^{-1} 
\cal F_{\calXdag }$$ qui est, en vertu du  théorème précédent, canoniquement isomorphe au  $\cal A_{\calUdag }[\cal G_{\calUdag /R}]$-module
$\cal F_{\calUdag }$. Si $\calWdag , \calUdag $ sont deux ouverts tels que $W\subset U$, le diagramme suivant:
$$\hss\def\quad{\hskip4pt}\scriptwd1.3em\mathrigid2mu\scriptspace0.5pt
\matrix{\cal A_{\calWdag }[\cal G_{\calWdag \rightarrow \calUdag }]\Otimes_{r^{-1}\cal A_{\calUdag }[\cal G_{\calUdag }]}
r^{-1}\cal A_{\calUdag }[\cal G_{\calUdag \rightarrow \calXdag }]\Otimes_{r^{-1}\cal A_{\calXdag }[\cal G_{\calXdag }]}r^{-1} 
\cal F_{\calXdag }&\simeq &\cal A_{\calWdag }[\cal G_{\calWdag \rightarrow \calUdag }]\Otimes_{r^{-1}\cal A_{\calXdag }[\cal G_{\calXdag }]}r^{-1}\cal F_{\calUdag }\cr
\downarrow&&\downarrow\cr
\cal A_{\calWdag }[\cal G_{\calWdag \rightarrow \calXdag }]\Otimes_{r^{-1}\cal A_{\cal X}[\cal G_{\calXdag }]}r^{-1} 
\cal F_{\calXdag }&\varto{0pt}{22ex}^\simeq& \cal F_{\cal
W\dag}}\hss$$
est commutatif par la transitivité des morphismes $(\sharp)$ pour le triplet $(\calWdag , \calUdag , \calXdag)$. On a
donc canoniquement l'isomorphisme de foncteurs 
$$ P_{\cal A_{\Xdaginf }, \calXdag }\circ R_{\calXdag }= Id\,.$$
Autrement dit, les foncteurs $P_{\cal A_{\Xdaginf }, \calXdag }$ et $R_{\calXdag }$ sont canoniquement inverses l'un de l'autre.
\enddemo

\begin{coro} La catégorie $ \cal A_{\calXdag }[\cal G_{\calXdag }]\Mod$ des  $\cal A_{ \calXdag }[\cal G_{\calXdag }]$-mo\-dules à gauche ne dépend  pas à équivalence \em{canonique}
près du relèvement $\calXdag $ de $X$ dans le cas d'un relèvement  plat d'un schéma lisse.
\end{coro}

En fait, si $\cal X\dag_1$ et $\cal X\dag_2$ sont deux relèvements lisses de $X$, on a un diagramme commutatif d'équivalences de catégories:
$$\let\quad\relax
\matrix{ \cal A_{\cal X\dag_1}[\cal G_{\calXdag }]\Mod&&
\varto{0pt}{14ex}&& \cal A_{\cal X_2\dag}[\cal
G_{\calXdag }]\Mod\cr
&\searrow&&\swarrow\cr
&& \cal A_{\Xdaginf }\Mod}$$ où les foncteurs obliques sont les prolongements  canoniques et le foncteur horizontal est le foncteur
image inverse ou directe (c'est comme on veut) $$ \cal F_{\cal X\dag_1}\rightarrow \cal A_{\cal X\dag_2}[\cal G_{\cal X\dag_2\rightarrow\cal X\dag_1}]
\otimes_{\cal A_{\cal X\dag_1}[\cal G_{\cal X\dag_1}]}\cal F_{\cal X\dag_1}.$$

\bigskip
La catégorie $ \cal A_{ \Xdaginf }\Mod$ des $ \cal A_{ \Xdaginf }$-modules  à
gauche a toutes les propriétés d'une catégorie de modules sur un
espace topologique  dans le cas d'un schéma  lisse. Par exemple:

\begin{coro}\label{rec-fai}Soit $X$ un schéma  lisse sur $R_1$ et soit $(\cal U_i, i\in I)$ une famille d'ouverts  du site $\Xdaginf $
telle que les ouverts  $U_i,i\in I,$  recouvrent  $X$. Alors, le foncteur qui à un faisceau $\cal F_{\inf}$ de $ \cal A_{ \Xdaginf }$-modules  sur le
site associe ses restrictions
$\cal F_{\cal U_i}, i\in I,$ est une équivalence de catégories entre la catégorie $ \cal A_{ \Xdaginf }\Mod$ et la catégorie des faisceaux donnés localement: les
familles
$\cal F_{\cal U_i}$ de $\cal A_{\calUdag }[\cal G_{\cal U\dag_i}]$-modules, $i\in I,$ munis d'isomorphismes de recollement satisfaisant les conditions de
cocycle.
\end{coro}

\demo Un isomorphisme de recollement entre $\cal F_{\cal U_i}$ et $\cal F_{\cal U_j}$ au-dessus de $U_i\cap U_j$ est un isomorphisme de recollement de faisceaux de
$ \cal A_{(U_i\cap U_j)\daginf}$-modules
sur le site $\dagger$-adique
de $U_i\cap U_j$, ce qui est équivalent, en vertu de ce qui précède, à un isomorphisme de recollement à leur restriction à tout relèvement de $U_i\cap U_j$. Les
conditions de cocycle se traduisent en conditions de cocycles pour les recollements de faisceaux sur des espaces topologiques. On est ramené au recollement
de faisceaux de modules sur les espaces topologiques et de leurs morphismes, sans à avoir à invoquer le recollement plus délicat de faisceaux sur un site.
\enddemo

\subsection{Le faisceau $\cHom_{\cal A_{ \Xdaginf }}(\cal F_{1\inf}, \cal F_{2\inf})$  de  deux  $ \cal A_{ \Xdaginf }$-modules à gauche}

\glossary{$\cHom_{\cal A_{ \Xdaginf }}(\cal F_{1\inf}, \cal F_{2\inf})$}Si $\cal F_{1\inf}$ et $\cal F_{2\inf}$ sont deux  $ \cal A_{ \Xdaginf }$-modules à gauche, où $ \cal A_{ \Xdaginf }$ est un faisceau de $R$-algèbres, nous allons définir le faisceau
$\cHom_{\cal A_{ \Xdaginf }}(\cal F_{1\inf}, \cal F_{2\inf})$ qui est un faisceau sur le site $\Xdaginf$ de $R$-modules.
Si $\calUdag $ est un ouvert, le faisceau  $\cHom_{\cal A_{ \calUdag  }}(\cal F_{1\,\calUdag }, \cal F_{2\,\calUdag })$ est un faisceau de Zariski sur $U$ de $R$-modules. Si $r\dag: \calWdag  \to \calUdag $ est un morphisme du site $\Xdaginf$,
le morphisme de restriction $(\sharp)_{r\dag}: r^{-1}\cal F_{1\,\calUdag }\to \cal F_{1\,\calWdag }$ est inversible.
On définit la restriction
$$r^{-1}\cHom_{\cal A_{ \calUdag  }}(\cal F_{1\,\calUdag }, \cal F_{2\,\calUdag })\to \cHom_{\cal A_{ \calWdag  }}(\cal F_{1\,\calWdag }, \cal F_{2\,\calWdag })\leqno (\sharp)_{r\dag}: $$ par $\phi\mapsto (\sharp)_{r\dag}\circ \phi\circ  ((\sharp)_{r\dag})^{-1}.$
\begin{prop}\label{fai-hom} La famille  $(\calUdag , \cHom_{\cal A_{ \calUdag  }}(\cal F_{1\,\calUdag }, \cal F_{2\,\calUdag }))$ munie
des morphismes de restriction précédents $(\sharp)_{r\dag}$ définit   un faisceau, noté
$\cHom_{\cal A_{ \Xdaginf }}(\cal F_{1\inf}, \cal F_{2\inf})$, qui est un faisceau   de $R$-modules sur le site $\Xdaginf$.
\end{prop}

\demo En effet, cette famille est de façon évidente  un objet  de la catégorie ${\mathbf{Fam}}(\Ouv(\Xdaginf ))$.
\enddemo
\begin{Rema} Le foncteur $\cal F_{\inf}\mapsto \cHom_{\cal A_{ \Xdaginf }}(
\cal A_{ \Xdaginf }, \cal F_{\inf})$ est canoniquement isomorphe par construction au foncteur identique de la catégorie des $\cal A_{ \Xdaginf }$-modules.\end{Rema}
\begin{Rema} Par contre, le lecteur prendra garde que les sections globales du faisceau $\cHom_{\cal A_{ \calUdag  }}(\cal F_{1\,\calUdag }, \cal F_{2\,\calUdag })$ {\bf ne sont pas}, en général, les morphismes entre $\cal F_{1\inf}, \cal F_{2\inf}$ au-dessus de $U$, et que $\rm Hom_{\cal A_{ \Udaginf }}(\cal A_{ \Udaginf },\cal F_{\inf})$ n'est pas égal, en général,  à $\Gamma(\calUdag ,\cal F_{\inf})=\Gamma(U,R_{\calUdag }(\cal F_{\inf}))$. C'est là un point important dans cette théorie.
\end{Rema} 
\subsection{Le faisceau $\cal F_{1\inf}\otimes_{\cal A_{ \Xdaginf }} \cal F_{2\inf}$  d'un   $ \cal A_{ \Xdaginf }$-module à gauche et d'un   $ \cal A_{ \Xdaginf }$-module à droite}

\glossary{$\cal F_{1\inf}\otimes_{\cal A_{ \Xdaginf }} \cal F_{2\inf}$}Si $\cal F_{1\inf}$  est un $ \cal A_{ \Xdaginf }$-modules à droite,  $\cal F_{2\inf}$ un   $ \cal A_{ \Xdaginf }$-module à gauche et $\calUdag $ est un ouvert du site, le faisceau $\cal F_{1\,\calUdag }\otimes_{\cal A_{ \calUdag  }} \cal F_{2\,\calUdag }$ est un faisceau de Zariski sur $U$ de $R$-modules.
 Si $r\dag: \calWdag  \to \calUdag $ est un morphisme du site $\Xdaginf$,
on définit la restriction:
$$r^{-1}\cal F_{1\,\calUdag }\otimes_{\cal A_{\calUdag  }} \cal F_{2\,\calUdag }\to \cal F_{1\,\calWdag }\otimes_{\cal A_{ \calWdag  }} \cal F_{2\,\calWdag }\leqno (\sharp)_{r\dag}: $$ 
comme le produit tensoriel des restrictions.
\begin{prop}\label{fai-ten}\tolerance1000 La famille  $(\calUdag , \cal F_{1\,\calUdag }\otimes_{\cal A_{\calUdag }} \cal F_{2\,\calUdag })$ munie
des morphismes de restriction précédents $(\sharp)_{r\dag}$ définit   un faisceau, noté
$\cal F_{1\inf}\otimes_{\cal A_{ \Xdaginf }} \cal F_{2\inf}$, qui est un faisceau
 de $R$-modules sur le site $\Xdaginf$.
\end{prop}

\demo En effet, cette famille est de façon évidente  un objet de la catégorie ${\mathbf{Fam}}(\Ouv(\Xdaginf ))$.
\enddemo
\subsection{Les foncteurs de prolongement et de restriction}
Soient $\cal A_{\Xdaginf }$ un faisceau d'anneaux sur le site et   $\calUdag $ un objet du site $\Xdaginf $. On a un foncteur de restriction: 
$$ \cal A_{\Xdaginf }\Mod\to \cal A_{\calUdag }[\cal G_{\calUdag }]\Mod\leqno \Rm
R_{\calUdag }:$$ qui à $\cal
F_{\inf}$ associe $\cal F_{\calUdag }$. Nous allons définir un foncteur prolongement:
$$ \cal A_{\calUdag }[\cal G_{\calUdag }]\Mod\to \cal A_{\Xdaginf /R}\Mod.\leqno P_{\cal
A_{\Xdaginf }, \calUdag }:$$ Soient $\cal
W\dag$ un ouvert du site, $r: W\cap U\hookrightarrow W$ l'inclusion canonique et $\cal F_{\calUdag }$
un $\cal A_{\calUdag }[\cal G_{\calUdag }]$-module. On définit:
$$P_{\cal A_{\Xdaginf }, \calUdag }(\cal F_{\calUdag })(\calWdag ):= r_* \mBig(\cal A_{\calWdag |W\cap U}[\cal G_{\calWdag |W\cap U\to \cal
U\dag}]\Otimes_{\cal A_{\calUdag }[\cal G_{\calUdag }]\sb{|W\cap U}}\cal F_{\calUdag }|W\cap U\mBig).$$ 

La famille $P_{\cal A_{\Xdaginf }, \cal
U\dag}(\cal F_{\calUdag })(\calWdag )$ définit de façon évidente un faisceau
$P_{\cal A_{\Xdaginf }, \calUdag }(\cal F_{\calUdag })$ de $\cal A_{\Xdaginf }$-modules à gauche sur le site $\Xdaginf $.

\bigskip\begin{prop}Soit  un  $R_1$-schéma lisse $X$. Alors, le foncteur prolongement $ P_{\cal A_{\Xdaginf }, \calUdag }$ est un adjoint à droite du
foncteur restriction
$ R_{\calUdag }$:
$$ \hom_{\cal A_{\calUdag }[\cal G_{\calUdag }]}(R_{\calUdag }?, ?)\simeq \hom_{\cal A_{\Xdaginf }}(?, P_{\cal A_{\Xdaginf }, \cal
U\dag}?).$$
\end{prop}

\demo
En fait, pour un ouvert $U$ de $X$ on a un foncteur naturel de restriction :
$$\cal A_{\Xdaginf }\Mod\to \cal A_{\Udaginf }\Mod$$ qui admet un adjoint à droite $ P_{\cal
A_{\Xdaginf }, \Udaginf }$. Cet adjoint associe à un 
$\cal A_{\Udaginf }$-module à gauche 
$\cal F_{\Udaginf }$ sur le site
$\Udaginf $   le $\cal A_{\Xdaginf }$-module à gauche $\cal F_{\Xdaginf }$ sur le site
$ \Xdaginf $  dont la valeur  sur un ouvert $\calWdag $ est le faisceau $r_*\cal F_{\Udaginf }(\calWdag |W\cap U)$. L'adjonction de la proposition résulte alors de l'équivalence
du théorème \ref{equ-god}.
\enddemo


De la même façon, on dispose du résultat suivant.  

\begin{prop}Soit  un  $R_1$-schéma lisse $X$, le foncteur de 
restriction
$ R_{\calUdag  }$  admet  un adjoint à gauche:
$$ \hom_{\cal A_{ \Xdaginf }}(P_{\calUdag  !}?, ?)\simeq \hom_{\cal A_{\calUdag }[\cal G_{\calUdag }]}(?, R_{\calUdag }?).\leqno  P_{\cal A_{\Xdaginf }, \calUdag  !} :$$
\end{prop}

\demo Si $r: W\cap U\hookrightarrow W$, on définit:
$$P_{\cal A_{\Xdaginf }, \calUdag  !}(\cal F_{\calUdag })(\calWdag ):= r_! \mBig(\cal A_{\calWdag |W\cap U}[\cal G_{\calWdag |W\cap U\to \cal
U\dag}]\Otimes_{\cal A_{\calUdag }[\cal G_{\calUdag }]_{|W\cap U}}\cal F_{\calUdag |W\cap U}\mBig)\,,$$ où $r_!$ est le foncteur prolongement par zéro. La
famille
$P_{\cal A_{\Xdaginf },
\calUdag  !}(\cal F_{\calUdag })(\calWdag )$ définit de façon évidente un faisceau
$P_{\cal A_{\Xdaginf }, \calUdag  !}(\cal F_{\calUdag })$ de $\cal A_{\Xdaginf }$-modules  à gauche sur le site $\Xdaginf $.
Remarquons que le foncteur $P_{\cal A_{\Xdaginf }, \calUdag  !}$ est exact.
\enddemo

\begin{coro} Soit  un  $R_1$-schéma lisse $X$. Alors, les foncteurs $R_{\calUdag }$ et  $P_{\cal A_{\Xdaginf }, \calUdag }$   transforment injectifs en injectifs.
\end{coro}

\demo En effet, ils sont tous les deux adjoints à droite d'un  foncteur  qui est exact.
\enddemo

\section{L'inclusion $ \cal G_{\Xdaginf }\hookrightarrow \calDdag_{\Xdaginf /R}$}\glossary{$ \cal
G_{\Xdaginf }\hookrightarrow \calDdag_{\Xdaginf /R}$}

\begin{theo}\label{inc-grp}Soit $\calXdag $ un schéma $\dagger$-adique  sur $R$. Le faisceau  $\cal G_{\calXdag }$ est alors un
sous-faisceau de semi-groupes  pour la structure multiplicative du faisceau d'anneaux $\calDdag_{\calXdag /R}$.
\end{theo}

\demo
Par définition, les deux faisceaux $\cal G_{\calXdag }$ et $\calDdag_{\calXdag /R}$ sont des sous-faisceaux du faisceau des endomorphismes
$\cal End_{R}(\cal O_{\calXdag /R})$. Il suffit de montrer que si $g$ est une section du faisceau $\cal G_{\calXdag }$, sa réduction $g_s$ modulo
$I^s$  comme section de $\cal End_{R_s}(\cal O_{ X_s/R_s})$ est un opérateur différentiel d'ordre borné par une fonction linéaire en $s$. Nous allons
montrer qu'en fait
$g_s$ est un opérateur différentiel d'ordre $s-1$ avec  $s\geq1$. Par définition d'un opérateur différentiel d'ordre $s-1$, avec  $s\geq1$, il suffit de montrer que si
$a_1,\dots,a_s$ sont des sections  du faisceau $\cal O_{X_s/R_s}$, alors le $s$-commutateur $[\dots[g_s, a_s],\dots,a_1]$ est un endomorphisme  nul.
Mais on a  l'égalité d'endomorphismes 
$$[\dots[g_s, a_s],\dots,a_1]=(g_s(a_s)-a_s)\cdots (g_s(a_1)-a_1)g_s,$$ 
et comme pour tout élément  $a$ l'élément 
$(g_s(a)-a)$ appartient à l'idéal engendré par 
$I$,  l'endomorphisme $[\dots[g_s, a_s],\dots,a_1]$ appartient à l'idéal engendré par $I^s$ et est donc nul  comme élément de $\cal End_{R_s}(\cal O_{ X_s/R_s})$.
\enddemo

\begin{coro}Si $P$ est un opérateur différentiel sur un schéma  $\dagger$-adique plat  $\calXdag $ et $g$ 
est un élément du groupe $\cal G_{\calXdag }$, l'endomorphisme $gPg^{-1}$ est un opérateur différentiel.
\end{coro}


Un peu plus généralement, on dispose du théorème suivant.

\begin{theo}\label{inc-grp'}Soient $Y\to X$ un morphisme de $R_1$-schéma, $\calYdag $ et $\calXdag $ deux  schémas $\dagger$-adiques  sur $R$ qui
relèvent $Y$ et $X$, et 
$u\dag$ et
$v\dag$ deux  morphismes de $\cal G_{\calYdag \to\calXdag }$. Alors, pour tout $s\geq1$ la
réduction $v_s$ modulo
$I^s$ de $v\dag$ est un opérateur différentiel d'ordre $s-1$ de la réduction $u_s$ modulo
$I^s$ de $u\dag$:   $\calD_{ Y_s\stackrel{u_s}{\to} X_s/R_s}:= \cDiff_{R_s}(u_s^{-1}\cal O_{X_s/R_s}, \cal O_{Y_s/R_s})$
\rm ([EGA IV$_4$], \S 16).
\glossary{$\calD_{ Y_s\stackrel{u_s}{\to} X_s/R_s}$}
\end{theo}

\demo La question est locale. Par définition d'un opérateur différentiel, comme élément du faisceau  $\cHom_{R_s}(u_s^{-1} \cal O_{ X_s/R_s}, \cal O_{
Y_s/R_s} )$ d'ordre $s-1$, $s\geq1$, il suffit de montrer que si
$a_1,\dots,a_s$ sont des sections  du faisceau $\cal O_{ X_s/R_s}$ alors le $s$-commutateur $[\dots[v_s, u_s(a_s)],\dots,u_s(a_1)]$, comme élément du faisceau  $\cHom_{R_s}(u_s^{-1} \cal O_{
X_s/R_s},
\cal O_{ Y_s/R_s} ),$  est
un homomorphisme  nul. Mais on a l'égalité d'homomorphismes: $$[\dots[v_s, u_s(a_s)],\dots,u_s(a_1)]=(v_s(a_s)-u_s(a_s))\cdots (v_s(a_1)-u_s(a_1))v_s\,,$$ 
et comme
pour tout élément
$a$, l'élément
$v_s(a)-u_s(a)$ appartient à l'idéal engendré par 
$I$, l'homomorphisme $[\dots[v_s, u_s(a_s)],\dots,u_s(a_1)]$ appartient à l'idéal engendré par $I^s$ et est donc nul  comme élément du faisceau $\cHom_{R_s}(u_s^{-1}
\cal O_{ X_s/R_s}, \cal O_{ Y_s/R_s} )$.
\enddemo

\begin{coro}\label{act-dif}Soit  un  $R_1$-schéma lisse $X$, le faisceau $ \cal G_{\Xdaginf }$ est un sous-faisceau de groupes pour la structure multiplicative du faisceau d'anneaux 
$ \calDdag_{\Xdaginf /R}$.
\end{coro}

\demo En effet, par construction,  pour tout ouvert $\calUdag $ du site, l'inclusion $\cal
G_{\calUdag }\hookrightarrow \calDdag_{\calUdag /R}$ commute aux restrictions pour tout morphisme du site.
\enddemo

\begin{coro}\label{dev-grp} Soient $V$ un anneau de valuation discrète complet, et $ A\dag$ une $V$-algèbre  $\dagger$-adique lisse  munie
d'éléments 
$x_1,\dots,x_n$ tels que
$\{dx_1,\dots,dx_n\}$ est une base du module des formes différentielles séparées $\Omega_{A\dag/V}$.
L'application qui à un élément  $g$ du groupe
$G_{A\dag} $ associe
$$\delta(g)=(\delta_1(g),\dots,\delta_n(g)):= ((g-1)(x_1),\dots, (g-1)(x_n))$$ est une {\bf bijection} du groupe $G_{A\dag}$ sur $(\goth m A\dag)^n$.
L'application inverse associe à
$a=(a_1,\dots,a_n)$ l'opérateur différentiel
\glossary{$g:= \sum_{\alpha\in
\Bbb N^n}a^\alpha\Delta_x^\alpha, a^\alpha:=a_1^{\alpha_1}\cdots a_n^{\alpha_n}$
\hbox to4cm{}}$$\theta_a:= \sum_{\alpha\in \Bbb N^n}a^\alpha\Delta_x^\alpha,\quad \hbox{avec }a^\alpha:=a_1^{\alpha_1}\cdots a_n^{\alpha_n} .$$
\end{coro}

\demo
En effet, comme $g$ est un opérateur différentiel en vertu du théorème \ref{inc-grp}, il est égal en vertu de la première partie du théorème du symbole total \ref{sym-tot} à la série
$$\sum_{\alpha\in \Bbb N^n}a_\alpha\Delta_x^\alpha\,,\quad\hbox{avec }
a_\alpha = \sum_{0\leq\beta\leq\alpha}\comb\alpha\beta(-x)^\beta g(x^{\alpha-\beta}).$$ Mais comme $g$ est un morphisme d'algèbres, on a l'égalité:
$$a_\alpha = (g(x_1)-x_1)^{\alpha_1}\cdots(g(x_n)-x_n)^{\alpha_n}\,,$$ ce qui montre que $g$ est complètement déterminé par $$\delta(g)= \big(g(x_1)-x_1,\dots,
g(x_n)-x_n\big).$$ Réciproquement, en vertu de la deuxième partie du théorème du symbole total \ref{sym-tot}, si $a_1,\dots,a_n$ sont des éléments de l'idéal $\goth m A\dag$, la série:
$$\theta_a= \sum_{\alpha\in \Bbb N^n}a^\alpha\Delta_x^\alpha$$ est un opérateur différentiel qui est un morphisme d'algèbres, en vertu de la formule de
Leibniz, et qui se réduit à l'identité modulo $\goth m$. 
\enddemo

\begin{Rema}On remarquera que les applications $\delta_i$ sont des dérivations intérieures:
$$\delta_i(g_1g_2)= \delta_i(g_1)+g_1\delta_i(g_2).$$
\end{Rema}

Le corollaire \ref{dev-grp} produit  \em{localement} beaucoup d'éléments du groupe $G_{A\dag}$:

\begin{coro}Soient $V$ un anneau de valuation discrète complet, $ A\dag$ une $V$-algèbre  $\dagger$-adique lisse  munie
d'éléments 
$x_1,\dots,x_n$ tels que
$\set dx_1,\dots,dx_n/$ est une base du module des formes différentielles séparées  $\Omega_{A\dag/V}$, et $\calXdag$  le $V$-schéma $\dagger$-adique associé à l'algèbre
$A\dag$. Alors, le groupe $G_{A\dag}$, qui est aussi l'ensemble des
sections globales du module de transfert  $\cal G_{\calXdag \to\calXdag }$ en vertu de l'équivalence \ref{eqi-aff}, est \em{non trivial}.
\end{coro}

\begin{Rema}\`A partir de là,  le théorème du symbole total donne des conditions nécessaires et suffisantes de recollement  pour qu'une famille
d'éléments définissant des sections locales du faisceau $\cal G_{\calXdag }$ définissent une sections globale, donnant ainsi un moyen éventuel 
de construire 
des éléments non triviaux du groupe $G_{A\dag}$ d'une algèbre $\dagger$-adique lisse $A\dag$ globale.
\end{Rema}

\section{La catégorie $ (\calDdag_{\Xdaginf /R},\Sp)\Mod$ des Modules à gauche spéciaux sur le site infinitésimal et la cohomologie de de Rham $\dagger$-adique}\glossary{$ (\calDdag_{\Xdaginf /R},\Sp)\Mod$} Soit $X$ un $R_1$-schéma lisse. En vertu du paragraphe précédent, 
on a sur le site infinitésimal $ \Xdaginf $    le faisceau  $ \calDdag_{\Xdaginf /R}$ des opérateurs différentiels, qui est un faisceau
de $ \cal O_{\Xdaginf /R}$-algèbres. On dispose donc de la catégorie $ \calDdag_{\Xdaginf /R}\Mod$ des $ \calDdag_{\Xdaginf /R}$-modules à
gauche.

\bigskip
Nous allons définir une sous-catégorie 
pleine  fondamentale
de la catégorie 
$ \calDdag_{\Xdaginf /R}\Mod$, à savoir la 
catégorie
$ (\calDdag_{\Xdaginf /R}, \Sp)\Mod$ 
des modules à gauche \em{\bf spéciaux}; c'est elle qui
va nous permettre de définir la cohomologie de de Rham $\dagger$-adique d'un
schéma  lisse sur $R_1$ ainsi que les opérations cohomologiques pour cette cohomologie. Nous commençons par définir cette catégorie à l'aide des morphismes de transfert ce qui est essentiel pour les opérations cohomologiques et nous montrons que cette définition est équivalente à celle de l'introduction.

\subsection{Le module de transfert $\calDdag_{\calWdag \rightarrow \calUdag /R}$  pour une immersion ouverte}
Soit $r\dag: \calWdag \to \calUdag $ un morphisme du site $ \Xdaginf $. Pour tout $s\geq1$ le morphisme $r\dag$ donne naissance à un morphisme
de schémas $r_s: W_s\rightarrow U_s$ sur $R_s$. On peut considérer le faisceau des opérateurs différentiels $R_s$-linéaires $$\calD_{
W_s\stackrel{r_s}{\rightarrow}  U_s/R_s}:=
\cDiff_{R_s}(r^{-1}\cal O_{U_s/R_s},\cal O_{W_s/R_s})$$ de
$r^{-1}\cal O_{U_s/R_s}$ dans $\cal O_{W_s/R_s}$, qui est par définition  un faisceau d'anneaux filtré ([EGA IV$_4$], $\S 16$).

\begin{defi}Soit $r\dag: \calWdag \to\calUdag $ un morphisme  du site $ \Xdaginf $. On définit le
module  de transfert:
$$\calDdag_{\calWdag \rightarrow \calUdag /R}$$\glossary{$\calDdag_{\calWdag \rightarrow \calUdag /R}$}comme le
sous-module de 
$\cHom_{R}(r^{-1}\cal O_{\calUdag /R},\cal O_{\calWdag /R})$
des morphismes de faisceaux de $R$-\em{modules} qui se réduisent modulo $I^s$ 
à un opérateur différentiel de $\calD_{
W_h\stackrel{r_s}{\rightarrow}  U_s/R_s}$  dont le degré est localement borné par une fonction linéaire en $s$. 
\end{defi}

Le module de transfert $\calDdag_{\calWdag \rightarrow \calUdag /R}$ est 
un $(\calDdag_{\calWdag /R}, r^{-1}\calDdag_{\calUdag /R})$-bimodule de façon naturelle.

\bigskip

\begin{theo}\label{ouv-tra}
Soit    $r\dag: \calWdag \to\calUdag $ un morphisme  du site
$ \Xdaginf $. Le module de transfert
$\calDdag_{\calWdag \rightarrow
\calUdag /R}$ {\bf ne dépend pas} du morphisme d'algèbres $r^*: r^{-1}\cal O_{\calUdag /R}\rightarrow \cal O_{\calWdag /R}$ choisi et est
un
$r^{-1}\calDdag_{\calUdag }$-module à droite localement libre engendré par le faisceau   de transfert $\cal G_{\calWdag \rightarrow \calUdag }$.
\end{theo}

\demo
Soit $P$ une section du module de transfert construit sur $r^*$.
En vertu du théorème  \ref{iso-pla}, le morphisme $r^*$ est inversible : $r^{*-1} :\cal O_{\calWdag /R}\simeq r^{-1}\cal O_{\calUdag /R}.$
En vertu du théorème \ref{inc-grp'}, la réduction
$r^{*-1}_{s}$ est un opérateur différentiel d'ordre $s-1$.
Le morphisme $Q:= r^{*-1}\circ P$ est un morphisme
de 
$\cal End_{R}(\cal O_{\calUdag |W/R})$ dont la réduction modulo $I^s$ est un opérateur différentiel  composé de deux opérateurs
différentiels et dont le degré
est une fonction localement  bornée par une fonction linéaire en $s$,   somme de deux fonctions localement bornées par des fonctions linéaires en $s$. C'est donc par définition un
opérateur différentiel
$Q$
de $\Gamma(W, \calDdag_{\calUdag |W/R})$. On a alors la factorisation:
$$P= r^*\circ Q.$$

Si $r'^*$ est un autre relèvement de $r$, on a la factorisation 
$r^*=r'^*\circ g$, où $g$ est une section globale du faisceau de groupes
$\cal G_{\cal
U\dag|W}$ et est donc un opérateur différentiel de $\calDdag_{\calUdag |W/R}$.
L'égalité: $$P= r^*\circ Q=r'^*\circ g\circ Q$$ montre que si $P$ est un opérateur de transfert construit sur $r^*$, alors c'est aussi un opérateur de transfert 
construit sur $r'^*$. Comme $r^*$ est inversible, cela montre la dernière assertion. D'où le théorème \ref{ouv-tra}.
\enddemo

\begin{defi}Soient un  $R_1$-schéma lisse $X$  et  $\calWdag , \calUdag $ un couple d'objets du site infinitésimal avec  $r: W\hookrightarrow
U$. Le module de transfert
$\calDdag_{\calWdag \rightarrow
\calUdag /R}$ est défini comme recollement 
des modules de transfert \em{locaux} induits par des morphismes de schémas $\dagger$-adiques qui relèvent l'inclusion. En vertu du
théorème précédent, le module de transfert est un sous-$(\calDdag_{\calWdag /R}, r^{-1}\calDdag_{\calUdag /R})$-bimodule bien défini de $\cHom_{R}(r^{-1}\cal O_{\calUdag /R},\cal O_{\calWdag /R})$.
\end{defi}
Cette définition impose que $X$ soit lisse pour {\bf garantir l'existence de relèvements locaux}. 

\subsection{La catégorie des $\Rm\calDdag_{\Xdaginf /R}$-modules à gauche spéciaux $ (\calDdag_{\Xdaginf /R}, \Sp)\Mod$}
Le module de transfert  pour une immersion ouverte permet de définir la catégorie des $\calDdag_{\Xdaginf /R}$-modules à gauche spéciaux.

\begin{defi}\label{mod-gau}Soit  $X$ un  $R_1$-schéma lisse.
Un $\calDdag_{\Xdaginf /R}$-module à gauche {\rm spécial} $\calMdaginf $ \glossary{$\cal
M\daginf $}sur le site
$\Xdaginf $ est la donnée, pour tout ouvert $\calUdag $ du site 
$\Xdaginf $, d'un 
$\calDdag_{\calUdag /R}$-module à gauche $\cal M\dag_{\cal U\dag}$, 
et la donnée, pour tout  couple $(\calWdag , \calUdag )$ d'objets du site 
$\Xdaginf $,   avec $r:W\hookrightarrow U$,
d'un morphisme de $\calDdag_{\calWdag /R}$-modules à gauche:
$$\calDdag_{\calWdag \rightarrow \calUdag /R}\otimes_{r^{-1}\calDdag_{\calUdag /R}}r^{-1}\cal M\dag_{\calUdag /R}\to \cal
M\dag_{\calWdag /R}\,.\leqno(\diamond):$$ 

En outre, ces données doivent satisfaire les conditions:
\begin{liste}\item 1) pour un couple $(\calUdag |W , \calUdag )$,  avec $r:
W\hookrightarrow U$, on a $\cal M\dag_{\calUdag |W}=\cal M\dag_{\calUdag }|W$ et le morphisme $(\diamond)$ coïncide avec le morphisme canonique,
\item 2) pour un triplet  $(\cal W'{}\dag ,\calWdag , \calUdag )$, avec $r': W'\hookrightarrow
W$ et $r: W\hookrightarrow U$,  le diagramme suivant est commutatif:
$$\hss\resetdisplay\scriptdp=6pt
\scriptspace0.5pt\def\quad{\hskip1pt}\mathrigid0mu
\def\smt{\let\bigotimes\otimes\scriptwd1em}
\matrix{\smt\calDdag_{\cal W'{}\dag \sto
\calWdag /R}\Otimes_{r'^{-1}\calDdag_{\calWdag /R}}r'^{-1}(\calDdag_{\calWdag \sto \cal
U\dag/R}\Otimes_{r^{-1}\calDdag_{\calUdag /R}} r^{-1}\cal M\dag_{\calUdag })&\longrightarrow &\smt\calDdag_{\cal
W'{}\dag \sto
\calUdag /R}\Otimes_{(r\circ r'{})^{-1}\calDdag_{\calUdag /R}} (r\circ r'{})^{-1}\cal M\dag_{\calUdag }\cr
\noalign{\kern-3pt}\downarrow&&\downarrow\cr
\calDdag_{\cal W'{}\dag \rightarrow \calWdag /R}\bigotimes_{r'{}^{-1}\calDdag_{\calWdag /R}} r{'}^{-1}\cal
M\dag_{\calWdag }&\decale{-0.5cm}{\varto{10pt}{15ex}} &\cal M\dag_{\cal W'{}\dag }\,.}$$
\end{liste}
\end{defi}

Munis des morphismes naturels, les $\calDdag_{\Xdaginf /R}$-modules à gauche spéciaux  forment une catégorie.

\begin{notation}On note $ (\calDdag_{\Xdaginf /R}, \Sp)\Mod$ la catégorie des 
$\calDdag_{\Xdaginf /R}$-modules à gauche spéciaux.
\end{notation}

\begin{prop}Soit  $X$ un  $R_1$-schéma lisse. 
Un $\calDdag_{\Xdaginf /R}$-module à gauche spécial $\calMdaginf $ est un
$R_{\Xdaginf }$-module sur le site
$\Xdaginf $ et  pour
tout  couple d'objets $(\calWdag , \calUdag )$  avec  $r: W\hookrightarrow U$   
le morphisme 
\glossary{$(\diamond)$}\glossary{$\calDdag_{\calWdag \rightarrow \calUdag /R}\Otimes_{r^{-1}\calDdag_{\calUdag /R}}
r^{-1}\cal
M\dag_{\calUdag }\rightarrow \cal M\dag_{\calWdag }$}$$\calDdag_{\calWdag \rightarrow \calUdag /R}\otimes_{r^{-1}\calDdag_{\calUdag /R}}
r^{-1}\cal
M\dag_{\calUdag }\rightarrow \cal M\dag_{\calWdag }\leqno(\diamond):$$ 
prolonge
le morphisme de restriction géométrique:
$$R[\cal G_{\calWdag \rightarrow \calUdag }]\otimes_{r^{-1}R[\cal G_{\calUdag }]}
r^{-1}\cal
M\dag_{\calUdag }\rightarrow \cal M\dag_{\calWdag }.\leqno(\sharp):$$ En particulier, l'action du groupe $\cal G_{\calUdag }$ sur $\cal
M\dag_{\calUdag }$ se fait à travers celle de $\calDdag_{\calUdag /R}$.
\end{prop}

\demo En vertu du théorème \ref{inc-grp},  un module $\calDdag_{\calUdag /R}$-module à gauche est un $R[\cal G_{\calUdag }]$-module à gauche. 
En vertu du théorème \ref{inc-grp'}, toute section locale du module $\cal G_{\calWdag \rightarrow \calUdag }$ est un opérateur différentiel. 
Il en résulte un morphisme canonique:
$$R[\cal G_{\calWdag \rightarrow \calUdag }]\to \calDdag_{\calWdag \rightarrow \calUdag /R}$$ 
qui induit un morphisme:
$$R[\cal G_{\calWdag \rightarrow \calUdag }]\otimes_{r^{-1}R[\cal G_{\calUdag }]}
r^{-1}\cal
M\dag_{\calUdag }\to\calDdag_{\calWdag \rightarrow \calUdag /R}\otimes_{r^{-1}\calDdag_{\calUdag /R}}
r^{-1}\cal
M\dag_{\calUdag }.$$ 
Les conditions de la proposition \ref{defi'} sont satisfaites et 
un  $\calDdag_{\Xdaginf /R}$-module à gauche
spécial est bien un
$R_{\Xdaginf }$-module   sur le site
$\Xdaginf $ dont l'action du groupe $\cal G_{\calUdag }$ sur $\cal
M\dag_{\calUdag }$ se fait à travers celle de $\calDdag_{\calUdag /R}$. 
\enddemo

\begin{prop}Soit  $X$ un  $R_1$-schéma lisse. La catégorie
des modules à gauche spéciaux  $ (\calDdag_{\Xdaginf /R}, \Sp)\Mod$  est  une
sous-catégorie  pleine de la catégorie $ \calDdag_{\Xdaginf /R}\Mod$ des $\calDdag_{\Xdaginf /R}$-modules à gauche. 
\end{prop}

\demo
Si $\calMdaginf $ est un $ \calDdag_{\Xdaginf /R}$-module à gauche spécial et $r\dag$ est un morphisme $\calWdag \to \calUdag $,
alors on a un morphisme de restriction:
$$\Gamma(U,\cal M\dag_{\calUdag })\to\Gamma(W,\cal M\dag_{\calUdag |W})\to\Gamma( W,\cal M\dag_{\calWdag })\,,$$ qui
provient du morphisme $(\diamond)$ et qui commute par construction avec l'action des opérateurs différentiels.
\enddemo
\noindent Réciproquement, on a ce qui suit.

\begin{prop}\label{def-spe}Soit  $X$ un  $R_1$-schéma lisse. Un $ \calDdag_{\Xdaginf /R}$-module à gauche est {\bf spécial
si et seulement si} l'action géométrique $(\sharp) $ de $\cal G_{\Xdaginf }$ se
fait à travers l'action de $ \calDdag_{\Xdaginf /R}$.
\end{prop}

\demo Par construction, l'action de $\cal G_{\Xdaginf}$ sur un $ \calDdag_{\Xdaginf /R}$-module à gauche  spécial se fait à
travers celle de $ \calDdag_{\Xdaginf /R}$.  Soit $\calMdaginf $ un $ \calDdag_{\Xdaginf /R}$-module à gauche dont
l'action géométrique  de $\cal G_{\Xdaginf}$ se fait à
travers celle de $ \calDdag_{\Xdaginf /R}$, et soit
$(\calWdag , \calUdag )$ un couple d'ouverts du site avec  $W\subset U$, on a alors le
morphisme canonique: 
$$R[\cal G_{\calWdag \rightarrow \calUdag }]\otimes_{r^{-1}R[\cal G_{\calUdag }]}r^{-1}\cal M\dag_{\cal
U\dag}\to \cal M\dag_{\calWdag }.\leqno (\sharp):$$ et aussi un autre morphisme canonique: 
$$R[\cal G_{\calWdag \rightarrow \calUdag }]\otimes_{r^{-1}R[\cal G_{\calUdag }]}r^{-1}\cal M\dag_{\cal
U\dag}\to\calDdag_{\calWdag \rightarrow \calUdag /R}\otimes_{r^{-1}\calDdag_{\calUdag /R}}
r^{-1}\cal
M\dag_{\calUdag }\,,$$ et il s'agit de voir que le morphisme $(\sharp)$ se factorise à travers un morphisme de $\calDdag_{\cal
W\dag/R}$-modules à gauche
$$\calDdag_{\calWdag \rightarrow \calUdag /R}\otimes_{r^{-1}\calDdag_{\calUdag /R}}
r^{-1}\cal
M\dag_{\calUdag }\rightarrow \cal M\dag_{\calUdag }.\leqno(\diamond):$$ Si $r^*$ est une section locale de $\cal G_{\calWdag \rightarrow \calUdag }$, elle engendre $\calDdag_{\calWdag \rightarrow \calUdag /R}$ comme $r^{-1}\calDdag_{\calUdag /R}$-module à
droite, on fait alors correspondre:
$$\mBig(\sum r^*P_i\otimes m_{i\calUdag }\mBig)\longmapsto
r^*\mBig(\sum P_im_{i\calUdag }\mBig),$$ 
d'où une application bien définie qui est un morphisme de 
$\calDdag_{\calWdag /R}$-modules à gauche. Cette application ne dépend pas du générateur $r^*$ parce que précisément l'action de $\cal G_{\calWdag }$ se fait
à travers celle de $\calDdag_{\calWdag /R}$. La famille $\cal M\dag_{\calUdag }$ munie des morphismes $(\diamond)$ ainsi définis a toutes
les propriétés de la définition d'un $ \calDdag_{\Xdaginf /R}$-module à gauche  spécial.
\enddemo

\begin{Rema}La proposition précédente montre que la définition simple des modules spéciaux de l'introduction coïncide avec la définition à l'aide des modules de transfert.
Mais c'est cette dernière qui est essentielle pour les opérations cohomologiques des chapitres qui suivent.
\end{Rema}

\begin{theo}\label{iso-can} Soient $X$ un  $R_1$-schéma lisse et $\calMdaginf $ un $\calDdag_{\Xdaginf /R}$-module à gauche
{\bf spécial}. Alors, les morphismes  de restriction
$$\calDdag_{\calWdag \rightarrow \calUdag /R}\otimes_{r^{-1}\calDdag_{\calUdag /R}}
r^{-1}\cal
M\dag_{\calUdag }\rightarrow \cal M\dag_{\calWdag }\leqno(\diamond):$$ 
sont des  {\bf isomorphismes} de $\calDdag_{\cal W\dag/R}$-modules à gauche.
\end{theo}

\demo  La question est locale et puisque on est dans le cas lisse il existe toujours une section locale $r^*$ du faisceau $\cal G_{\calWdag \rightarrow \calUdag }$ qui est \em{inversible}. Si $m_{\calWdag }$ est une section locale du faisceau $\cal M\dag_{\calWdag }$, elle provient de
la section locale
$r^*\otimes r^{*-1}m_{\calWdag }$ du faisceau 
$\calDdag_{\calWdag \rightarrow \calUdag /R}\otimes_{{r^{-1}\calDdag_{\calWdag /R}}}
{r^{-1}\cal
M\dag_{\calWdag }}$. Le morphisme $(\diamond)$ du théorème  est surjectif. Le faisceau $\calDdag_{\calWdag \rightarrow \calUdag /R}$ est un 
${r^{-1}\calDdag_{\calUdag }}$-module à droite engendré localement par $r^*$. Toute section locale du produit tensoriel est de la forme 
$r^*\otimes m_{\calUdag }$ dont l'image est nulle si et seulement si $m_{\calUdag }$ est nulle. Le morphisme $(\diamond)$ du théorème  est injectif. 
\enddemo

\begin{Rema}Dans le cas lisse,
si $\calMdaginf $ est un $\calDdag_{\Xdaginf /R}$-module 
spécial  de  valeur  $\cal M\dag_{\calUdag }$ sur un
ouvert
$\calUdag $, alors sa valeur  est  canoniquement: 
$$\calDdag_{\calWdag \rightarrow \calUdag /R}\otimes_{\calDdag_{\calUdag /R}}
\cal
M\dag_{\calUdag }$$ sur tout  ouvert $\calWdag $ qui relève $U$, et c'est en cela qu'il est spécial. Autrement dit, un module
spécial est canoniquement  déterminé par sa restriction à un relèvement particulier.
\end{Rema}

\begin{defi}Soit  $\calXdag $  un relèvement  plat  sur $R$ d'un schéma $X$  lisse sur $R_1$ et soit $\cal M\dag_{\calXdag }$ un 
$\calDdag_{\calXdag /R}$-module à gauche. On définit le prolongement $P_{\calXdag }(\cal M\dag_{\calXdag })$ de 
$\cal M\dag_{\calXdag }$
comme  le $\calDdag_{\Xdaginf /R}$-module  à gauche {\rm spécial} 
dont la valeur  sur un ouvert
$\calUdag $, $r: U\hookrightarrow X,$ est   le $\calDdag_{\calUdag  /R}$-module  à gauche: 
$$\calDdag_{\calUdag \rightarrow \calXdag /R}\otimes_{r^{-1}\calDdag_{\calXdag /R}} r^{-1}\cal M\dag_{\calXdag }.$$ 
\end{defi}

Par construction, le prolongement  $P_{\calXdag }(\cal M\dag_{\calXdag })$ d'un  $\calDdag_{\calXdag /R}$-modules à gauche $\cal M\dag_{\calXdag }$
est un 
$ \calDdag_{\Xdaginf }$-module \em{spécial}, les conditions de transitivité étant automatiques. On obtient ainsi un
foncteur\glossary{$P_{\calXdag }$}\glossary{$R_{\calXdag }$}
$$ \calDdag_{\calXdag /R}\Mod\to (\calDdag_{\Xdaginf /R}, \Sp)\Mod.\leqno  P_{\calXdag }:$$ 
\begin{lemm}Le foncteur prolongement  $P_{\calXdag }$ est exact.\end{lemm}
\demo
En effet,   $\calDdag_{\calUdag \rightarrow \calXdag /R}$ est un $r^{-1}\calDdag_{\calXdag /R}$-module  à droite localement libre de rang $1$.

\begin{theo}\label{can}Soit  $\calXdag $  un relèvement plat sur $R$  d'un schéma $X$  lisse sur $R_1$. Le foncteur prolongement  $\Rm
P_{\calXdag }$ est une équivalence de catégories,  qui est un inverse \em{\bf canonique} du
foncteur de restriction $R_{\calXdag }$.
\end{theo}

\demo Il est évident que le composé $R_{\calXdag }\circ P_{\calXdag }$ est canoniquement isomorphe au foncteur identique de la catégorie  $\calDdag_{\calXdag /R}\Mod$. Soit $\calMdaginf $ un 
$\calDdag_{\Xdaginf /R}$-module  à gauche  spécial et soit $\calUdag $ un ouvert du site. Par définition de module spécial
on a un isomorphisme canonique: 
$$\calDdag_{\calUdag \rightarrow \calXdag /R}\otimes_{r^{-1}\calDdag_{\calXdag /R}} r^{-1}\cal M\dag_{\calXdag }\simeq
\cal M\dag_{\calUdag }\leqno(\diamond):$$ de $\calDdag_{\calUdag /R}$-modules à gauche. Si $\calWdag , \calUdag $ sont deux ouverts avec 
$W\subset U$, on a le diagramme:
$$\scriptspace0.5pt\scriptwd1.8em\def\quad{\hskip5pt}\mathrigid3mu
\matrix{\calDdag_{\calWdag \rightarrow \calUdag /R}\Otimes_{r^{-1}\calDdag_{\cal
U\dag/R}}r^{-1} \calDdag_{\calUdag \rightarrow \calXdag /R}\Otimes_{r^{-1}\calDdag_{\calXdag/R}}r^{-1} 
\cal M\dag_{\calXdag }&\simeq &\calDdag_{\calWdag \rightarrow \calUdag /R}\Otimes_{r^{-1}\calDdag_{\calUdag /R}}r^{-1}\cal
M\dag_{\calUdag }\cr
\downarrow&&\downarrow\cr
\calDdag_{\calWdag \rightarrow \calXdag /R}\Otimes_{r^{-1}\calDdag_{\calXdag /R}}r^{-1} 
\calDdag_{\calXdag }&\varto{0pt}{20ex}^\simeq& \cal M\dag_{\cal
W\dag},}$$ qui est commutatif en vertu de la condition de transitivité des morphismes $(\diamond)$ pour un module spécial. Cela  exprime que le composé
$P_{\calXdag }\circ R_{\calXdag }$ est canoniquement isomorphe au foncteur identique de la catégorie $ (\calDdag_{\Xdaginf /R}, \Sp)\Mod$.
\enddemo
\begin{exemple} 
Si $\calUdag $ est un ouvert, la valeur du prolongement  $P_{\calXdag }(\cal O_{\calXdag /R})$ en $\calUdag $ est le faisceau 
$\cal O_{\calUdag /R}$ et donc $P_{\calXdag }(\cal O_{\cal X/R})\simeq \cal O_{\Xdaginf/R}$, alors que la valeur du prolongement  $P_{\calXdag }(\calDdag_{\cal X/R})$ 
en $\calUdag $ 
est le faisceau 
$\calDdag_{\calUdag \to \calXdag /R}$ et donc $P_{\calXdag }(\calDdag_{\cal X/R})$ {\bf n'est pas} isomorphe à $ \calDdag_{\Xdaginf/R}$.\end{exemple}

\begin{coro}\label{can-cor} Soit  un  $R_1$-schéma lisse $X$. Alors, la catégorie $ \calDdag_{\calXdag /R}\Mod$ des $\calDdag_{\calXdag /R}$-modules à gauche ne dépend  pas {\bf à
équivalence canonique} près du relèvement  $\dagger$-adique plat $\calXdag $ de $X$ s'il existe.
\end{coro}

En fait, si $\cal X\dag_1$ et $\cal X\dag_2$ sont deux relèvements plats de $X,$ on a un diagramme commutatif d'équivalences
de catégories:
$$\let\quad\relax\matrix{ \calDdag_{\cal X\dag_1/R}\Mod &&\varto{0pt}{19ex}&& \calDdag_{\cal X\dag_2/R}\Mod \cr
&\searrow&&\swarrow&\cr  
&& (\calDdag_{\Xdaginf /R}, \Sp)\Mod}$$ où les foncteurs obliques sont les prolongements 
canoniques et le foncteur horizontal est le foncteur image inverse ou directe:  $$ \cal M\dag_{\cal
X_1\dag}\mapsto\calDdag_{\cal X\dag_1\rightarrow\cal X\dag_2/R}
\otimes_{\calDdag_{\cal X\dag_1/R}}\cal M\dag_{\cal X\dag_1}$$ où le module de transfert $\calDdag_{\cal X_1\dag\rightarrow\cal X\dag_2/R}$
est défini par le théorème \ref{ouv-tra}.

\begin{prop}\label{cat-abe}Soit $X$ un  $R_1$-schéma lisse. 
La catégorie 
des modules à gauche spéciaux $ (\calDdag_{\Xdaginf /R},\Sp)\Mod$ est
une sous-catégorie abélienne de la catégorie
$ \calDdag_{\Xdaginf /R}\Mod$.
\end{prop}

\demo
En effet, si $u: \calMdaginf \to \calNdaginf $ est un morphisme de la catégorie des modules à gauche spéciaux $ (\calDdag_{\Xdaginf /R},\Rm
\Sp)\Mod$ son noyau
$\Ker(u)$ et son conoyau
$\Coker(u)$ sont des $R_{\Xdaginf }$-modules sur le site. Si $\calWdag , \calUdag $ est un couple d'objets  du site avec  $r: W\hookrightarrow
U$, 
la suite exacte de
$\calDdag_{\calUdag /R}$-modules à gauche
$$0\to \Ker(u)_{\calUdag }\to \cal M\dag_{\calUdag }\to \cal N\dag_{\calUdag }\to \Coker(u)_{\calUdag }\to0$$ donne naissance par platitude (le  $r^{-1}\calDdag_{\calUdag /R}$-module à droite $\calDdag_{\calWdag \rightarrow \calUdag /R}$ étant localement libre de rang $1$)  à la suite exacte de $\calDdag_{\calWdag /R}$-modules à gauche:
$$
\ecartlignes4pt\DeuxLignes  
0
\to
\calDdag_{\calWdag \rightarrow \calUdag /R}\Otimes_{r^{-1}\calDdag_{\calUdag /R}}r^{-1}\Ker(u)_{\calUdag }
\to 
\calDdag_{\calWdag \rightarrow \calUdag /R}\Otimes_{r^{-1}\calDdag_{\calUdag /R}}r^{-1}\cal M\dag_{\calUdag }\to\hskip2cm 
\\
\to
\calDdag_{\calWdag \rightarrow \cal
U\dag/R}\Otimes_{r^{-1}\calDdag_{\calUdag /R}}r^{-1}\cal N\dag_{\calUdag }
\to
\calDdag_{\calWdag \rightarrow \calUdag /R}\Otimes_{r^{-1}\calDdag_{\calUdag /R}}r^{-1}\Coker(u)_{\calUdag }
\to
0\,,
\endlignes$$
qui montre que le noyau et le conoyau du morphisme $u$ sont des $ \calDdag_{\Xdaginf /R}$-modules à gauche spéciaux.
\enddemo

\begin{Rema}Le théorème \ref{can} et son corollaire précédent constituent le succès le plus remarquable du point de vue du site infinitésimal  $ 
\Xdaginf $. De plus, les isomorphismes de restriction  $(\diamond)$  des modules spéciaux du théorème \ref{iso-can} correspondent  à la rigidité des cristaux de Grothendieck [G$_3$], et la propagation des modules spéciaux  du théorème \ref{can} correspond à la croissance des cristaux dans un voisinage approprié. Mais alors que la catégorie des cristaux ne semble fournir une bonne théorie que pour les variétés algébriques  propres et lisses sur $k$, la catégorie des
modules spéciaux semble fournir  une bonne théorie cohomologique pour les variétés éventuellement ouvertes, comme l'illustre la {\bf première démonstration} du  théorème \ref{fac-zet}  
 de l'expression cohomologique $p$-adique de la fonction Zêta d'une variété algébrique lisse sur un corps fini.
\end{Rema}

\sousparagraphe{}Soit $\calXdag $ un relèvement  plat de $X$, on a alors  deux foncteurs 
$ P_{\mkern-3mu\calDdag_{\mkern-6mu \Xdaginf /\!R}, \calXdag }\!$ et $ P_{\calXdag }$ et nous allons étudier leur rapport. On a deux morphismes 
canoniques de faisceaux d'anneaux: l'injection canonique, et la surjection canonique qui à un élément de l'algèbre du groupe $\sum_\alpha P_\alpha.g_\alpha$ associe l'opérateur
différentiel $\sum_\alpha P_\alpha g_\alpha$:
$$\calDdag_{\calXdag /R}\stackrel{\inj}{\longhookrightarrow}\calDdag_{\calXdag /R}[\cal G_{\calXdag }]
\stackrel{\surj}{\too}
\calDdag_{\calXdag /R}.$$ 
Le morphisme composé $surj\circ \inj$ est l'application identique, alors que $\inj\circ 
\surj$
n'est pas l'application identique. Si bien qu'on a deux
foncteurs canoniques de restrictions des scalaires:
$$\calDdag_{\calXdag /R}\Mod\stackrel{\surj^*}{\too}\calDdag_{\calXdag /R}[\cal G_{\calXdag }]\Mod\stackrel{\inj^*}{\too}\calDdag_{\calXdag /R}\Mod.$$ 
Le foncteur composé $\inj^*\circ \surj^*$ est le foncteur identique, mais le  foncteur  $\surj^*\circ \inj^*$ {\bf n'est pas} le foncteur identique.

\begin{prop} Soient  un  $R_1$-schéma lisse $X$ et $\calXdag $ et un relèvement  plat. On a un isomorphisme canonique de foncteurs de la catégorie $\Rm\calDdag_{\calXdag /R}\Mod$
vers la catégorie $\Rm\calDdag_{\Xdaginf /R}\Mod$:
$$ P_{\calDdag_{\Xdaginf /R}, \calXdag }\circ \surj^*\simeq  P_{\calXdag }.$$
\end{prop}

\demo  Soient $\cal M\dag_{\calXdag }$ un $\calDdag_{\calXdag /R}$-module à gauche et $\calUdag $ un ouvert du site $\Xdaginf $. On a alors un morphisme canonique:
$$\calDdag_{\calUdag /R}[\cal G_{\calUdag \to \calXdag }]\otimes_{r^{-1}\calDdag_{\calXdag /R}[\cal G_{\calXdag }]} r^{-1}\cal M\dag_{\calXdag }\to\calDdag_{\calUdag \rightarrow \calXdag /R}\otimes_{r^{-1}\calDdag_{\calXdag /R}} r^{-1}\cal M\dag_{\calXdag }\,,$$ qui commute aux restrictions, et il s'agit de voir que c'est un isomorphisme 
de $\calDdag_{\calUdag /R}$-modules à gauche. Le problème est de nature locale. Soit $r\dag$ un relèvement local de l'inclusion
$U\to X$. Alors, $\calDdag_{\calUdag /R}[\cal G_{\calUdag \to \calXdag }]$ est engendré par $r^*$ comme $r^{-1}\calDdag_{\calXdag /R}[\cal G_{\calXdag }]$-module à droite  et $\calDdag_{\calUdag\to\calXdag /R}$ est engendré par $r^*$ comme $r^{-1}\calDdag_{\calXdag /R}$-module à droite en vertu du théorème 	\ref{ouv-tra}, ce qui montre que le morphisme est à la fois surjectif et injectif. 
\enddemo

\begin{Rema}Le foncteur $ P_{ \calXdag }\circ \inj^*$ n'est pas isomorphe au foncteur
$ P_{\calDdag_{\Xdaginf /R}, \calXdag }$.
\end{Rema}

\begin{Rema} Remarquons aussi que le faisceau structural du site infinitésimal $\dagger$-adique $ \cal O_{\Xdaginf /R}$  est, par construction,
un
$ \calDdag_{\Xdaginf /R}$-module à gauche spécial et est le prolongement   du faisceau structural $\cal O_{\calXdag /R}$ d'un relèvement.
Mais le faisceau
$ \calDdag_{\Xdaginf /R}$ n'est pas spécial,  parce que l'action  de $\cal G_{\Xdaginf }$  par automorphismes intérieurs est distincte de
l'action différentielle à
gauche. Autrement dit, {\bf  il n'existe pas de morphisme canonique}$$ \calDdag_{\calWdag \rightarrow \cal
U\dag/R}\to \calDdag_{\calWdag /R}\leqno (\diamond):$$
 en général.
\end{Rema}

\begin{Rema} Bien que la catégorie $ \calDdag_{\Xdaginf /R}\Mod$ ait  bien un sens dans le cas singulier, elle ne semble pas avoir des propriétés
intéressantes. Dans le cas singulier la théorie en caractéristique nulle nous enseigne de procéder autrement et à notre avis cela ne posera pas de difficulté une fois surmonté tous les problèmes que posent le cas d'un morphisme, sans hypothèse restrictive, de schémas lisses. 
Aussi, dans tout cet article le lecteur  ne  perdra  rien d'intéressant
en  supposant que le schéma $X$ est  lisse sur $R_1$.
\end{Rema}
Si $X$ est lisse sur $R_1$ la catégorie de modules spéciaux est un champ:

\begin{theo}\label{champ}Soit  $X$ un  $R_1$-schéma lisse  alors  la catégorie $ (\calDdag_{\Xdaginf /R}, \Sp)$-modules à gauche spéciaux est un champ sur $X$, c'est-à-dire
ses  objets et ses morphismes sont de nature locale sur $X$.
\end{theo}

\demo C'est une conséquence directe de l'équivalence fondamentale \ref{can}.
\enddemo

\subsection{La catégorie des modules spéciaux a suffisamment d'objets injectifs}
Nous allons montrer que la catégorie des modules spéciaux  a suffisamment d'injectifs, bien que ce ne soit pas  une catégorie de modules, ce qui nous permettra 
de dériver les foncteurs exacts à gauche et, en particulier, de définir la cohomologie de de Rham
$\dagger$-adique en toute généralité.

\begin{theo}\label{obj-inj}Soit  un  $R_1$-schéma lisse $X$.\\
La catégorie $ (\calDdag_{\Xdaginf /R},\Sp)\Mod$ a suffisamment
d'objets injectifs.
\end{theo}

\demo C'est une conséquence directe du théorème \ref{can}.
Soit  $\calUdag $ un ouvert et $\cal M\dag_{\calUdag }$ un $\calDdag_{\calUdag /R}$-module à gauche, on définit le foncteur prolongement $P_{\cal
U\dag}$ par:
$$P_{\calUdag }(\cal M\dag_{\calUdag })(\calWdag ):= r_*\mBig(\calDdag_{\calWdag |W\cap U\to\cal
U\dag/R}\otimes_{r^{-1}\calDdag_{\calUdag /R}}r^{-1}\cal M\dag_{\calUdag }\mBig)$$ où $\calWdag $ est un objet du site $X\daginf$ et $r: W\cap U\hookrightarrow W$ l'inclusion canonique.
On obtient un adjoint à droite du foncteur restriction $R_{\calUdag }$ (exact à gauche)
et donc transforme injectif en injectif.

Si maintenant
$\calMdaginf $ est un $\calDdag_{\Xdaginf /R}$-module à gauche spécial, $\calUdag $ un ouvert et $R_{\calUdag }(\calMdaginf )\hookrightarrow
\cal I\dag_{\calUdag }$ un plongement dans un $\calDdag_{\calUdag /R}$-module à gauche injectif alors 
$$P_{\calUdag }\circ R_{\calUdag }(\calMdaginf )\hookrightarrow P_{\calUdag }(\cal I\dag_{\calUdag })$$ est un plongement dans un $\calDdag_{\Xdaginf /R}$-module à gauche spécial et injectif. Si
$\set\cal U\dag_i,i\in I/$ est une famille d'ouverts du site tels que 
$\set U_i,i\in I/$ est un recouvrement de $X$, le morphisme
produit 
$$\calMdaginf \too \prod_{i\in I}P_{\cal U\dag_i}(\cal I\dag_{\cal U\dag_i})$$ est un plongement dans un $\calDdag_{\Xdaginf /R}$-module à gauche spécial et injectif comme objet de la catégorie $ (\calDdag_{\Xdaginf /R}, \Sp)\Mod$.
\enddemo

\begin{Rema}
{\bf Attention}: Un  $\calDdag_{\Xdaginf /R}$-module à gauche spécial qui est injectif comme objet de la catégorie $ (\calDdag_{\Xdaginf /R}, \Sp)\Mod$ {\bf n'est pas} en général injectif comme objet de la catégorie $ \calDdag_{\Xdaginf /R}\Mod$. C'est là une différence essentielle, voir \ref{con-exe}. C'est un autre point important dans cette théorie.
\end{Rema}

\subsection{La cohomologie de de Rham $\dagger$-adique}
Si $X$ est un schéma  lisse  sur $R_1$, on note $ (\calDdag_{\Xdaginf /R},\Sp)\Mod$ la catégorie des $ \calDdag_{\Xdaginf /R}$-modules à gauche spéciaux sur le site infinitésimal $ \Xdaginf $ de $X$.

La catégorie $ (\calDdag_{\Xdaginf /R},\Sp)\Mod$ est abélienne et admet suffisamment de 
$ \calDdag_{\Xdaginf /R}$-modules {\bf spéciaux injectifs}, en vertu du théorème précédent. On note $ \rmD^*((\calDdag_{\Xdaginf /R},\Sp)\Mod)$ \glossary{$ \rmD^*((\calDdag_{\Xdaginf /R},\Sp)\Mod)$}sa catégorie dérivée. On peut donc
dériver tout foncteur covariant exact à gauche de la catégorie $ (\calDdag_{\Xdaginf /R},\Sp)\Mod$, et en particulier le foncteur: $$\Rm
\calMdaginf \mapsto \hom_{\calDdag_{\Xdaginf /R},\Sp}(
\cal O_{\Xdaginf /R},\calMdaginf )=\hom_{\calDdag_{\Xdaginf /R}}(
\cal O_{\Xdaginf /R},\calMdaginf ),$$ dont on notera $ \ext^{\bullet}_{\calDdag_{X\dag_{\inf/R},\Sp}}( \cal O_{\Xdaginf /R},\cal
M\daginf )$ \glossary{$ \ext^{\bullet}_{\calDdag_{X\dag_{\inf/R},\Sp}}( \cal O_{\Xdaginf /R},\cal
M\daginf )$}les foncteurs dérivés.

\begin{defi}Soient $X$  un schéma  lisse  sur $R_1$ et $\calMdaginf $ un complexe de la catégorie $ \Dplus ((\calDdag_{\Xdaginf /R},\Sp)\Mod)$. On définit la cohomologie de de Rham $\dagger$-adique de $X$ à coefficients dans $\calMdaginf $ comme les $R$-modules
$$ H^\bullet_{\DR}(X/R,\calMdaginf ):= \ext^{\bullet}_{\calDdag_{X\dag_{\inf/R}},\Sp}(\cal O_{\Xdaginf /R},\calMdaginf )\,.$$ 
En
particulier, on définit la cohomologie de de Rham $\dagger$-adique  de $X$ comme: 
\glossary{$ H^\bullet_{\DR}(X/R, \calMdaginf )$}
$$ H^\bullet_{\DR}(X/R):=  H^\bullet_{\DR}(X/R,\cal O_{\Xdaginf /R}):= \ext^{\bullet}_{\calDdag_{\Xdaginf /R},\Sp}(\cal O_{\Xdaginf /R},\cal
O_{\Xdaginf /R}).$$
\end{defi}
Les modules  de cohomologie
$ H^\bullet_{\DR}(X/R,\calMdaginf )$ sont donc les modules de cohomologie du complexe: \glossary{$ \bfR \hom_{\calDdag_{\Xdaginf /R},\Sp}(
\cal O_{\Xdaginf /R},\calMdaginf )$}\glossary{$ \bfR \cHom_{\calDdag_{\Xdaginf /R},\Sp}(
\cal O_{\Xdaginf /R},\calMdaginf )$} \glossary{$\DR(X/R, \calMdaginf )$}
$$\DR(X/R, \calMdaginf ):=  \bfR \hom_{\calDdag_{\Xdaginf /R},\Sp}(
\cal O_{\Xdaginf /R},\calMdaginf ),$$ 
vu comme foncteur dérivé à droite du foncteur :
$$ \calMdaginf \mapsto\hom_{\calDdag_{\Xdaginf /R}}(
\cal O_{\Xdaginf /R},\calMdaginf )$$ de la catégorie des modules spéciaux qui a suffisamment d'injectifs dans la catégorie des
$R$-modules. De plus, l'action de $R$ sur ces modules se fait à travers le morphisme canonique $R\to \hat R$ de $R$ dans son complété séparé pour la topologie $I$-adique.

\begin{theo}\label{coh-deR}Soient $\calXdag $ un relèvement  plat  d'un $R_1$-schéma lisse  $X$ et $\cal M\dag_{\calXdag }$ un complexe de $ \Dplus (\calDdag_{\calXdag/R})$. Il existe des isomorphismes \em{\bf canoniques} de $R$-modules:
$$ \ext^{\bullet}_{\calDdag_{\calXdag /R}}(\cal O_{\calXdag /R},\cal M\dag_{\calXdag })\simeq\ext^{\bullet}_{\calDdag_{\Xdaginf /R},\Sp}(
\cal O_{\Xdaginf /R},P_{\calXdag }(\cal M\dag_{\calXdag })).$$ En particulier, il existe des isomorphismes canoniques:
$$ \ext^{\bullet}_{\calDdag_{\calXdag /R}}(\cal O_{\calXdag /R},\cal O_{\calXdag /R})\simeq\ext^{\bullet}_{\calDdag_{\Xdaginf /R},\Sp}(
\cal O_{\Xdaginf /R},\cal O_{\Xdaginf /R})\,.$$
\end{theo}

\demo Soit $\cal M\dag_{\calXdag }$ un complexe de $ \Dplus (\calDdag_{\calXdag /R})$ et $\cal I\dag_{\calXdag }$ une résolution $\calDdag_{\calXdag /R}$-injective de $\cal M\dag_{\calXdag }$ alors, en vertu du théorème \ref{can}, son prolongement canonique  $P_{\calXdag }(\cal I\dag_{\calXdag })$ est une
résolution
injective de $P_{\calXdag }(\cal M\dag_{\calXdag })$ dans la catégorie $ (\calDdag_{\Xdaginf /R},\Sp)\Mod$. D'autre part, le morphisme de complexes de $R$-modules:
$$ \hom_{\calDdag_{\calXdag /R}}(\cal O_{\calXdag /R},\cal I\dag_{\calXdag })\simeq\hom_{\calDdag_{\Xdaginf /R},\Sp}(
\cal O_{\Xdaginf /R},P_{\calXdag }(\cal I\dag_{\calXdag })),$$ qui est un isomorphisme,  induit en particulier un
isomorphisme de leurs cohomologies,  réalisant ainsi  les isomorphismes du théorème.
\enddemo

En particulier, la cohomologie de de Rham $\dagger$-adique d'un schéma  lisse coïncide  avec la cohomologie de de Rham d'un relèvement quand il
existe, et la cohomologie de de Rham $\dagger$-adique d'un relèvement lisse ne dépend pas à isomorphisme canonique près du relèvement et coïncide donc avec la définition de l'article de recherche [M-N$_1$]. On a alors résolu notre problème.

\bigskip
Un peu plus généralement, si $R\to S$ est un morphisme d'anneaux, on note   $$ \cal O_{\Xdaginf /S}:=  \cal O_{\Xdaginf /R}\otimes_RS,  \calDdag_{\Xdaginf /S}:=  \calDdag_{\Xdaginf /R}\otimes_RS$$ les faisceaux obtenus par changement de base.
Le faisceau de groupes $\cal G_{\Xdaginf }$ opère à gauche sur  les faisceaux $ \cal O_{\Xdaginf /S},  \calDdag_{\Xdaginf /S}$ et l'on dispose  de la catégorie $ (\calDdag_{\Xdaginf /S}, \Sp)\Mod$ des $ \calDdag_{\Xdaginf /S}$-modules à gauche spéciaux sur le site infinitésimal $ \Xdaginf $ de $X$ qui est abélienne et qui a suffisamment d'injectifs. On peut donc considérer la théorie précédente sur S.

\begin{defi}
Soient
$X$  un schéma  lisse  sur $R_1$ et $\calMdaginf $ un complexe de $ \Dplus ((\calDdag_{\Xdaginf /S},\Sp)\Mod)$. On définit la cohomologie de de Rham $\dagger$-adique de $X$ à coefficients   $\calMdaginf $ relativement à $S$ comme les $S$-modules
$$ H^\bullet_{\DR}(X/S,\calMdaginf ):= \ext^{\bullet}_{\calDdag_{X_{\inf/S}},\Sp}(\cal O_{\Xdaginf /S},\calMdaginf ).$$ 
En
particulier, on définit la cohomologie de de Rham $\dagger$-adique  de $X$ comme 
$$ H^\bullet_{\DR}(X/S):=  H^\bullet_{\DR}(X/S,\cal O_{\Xdaginf /S}):= \ext^{\bullet}_{\calDdag_{\Xdaginf /S},\Sp}(\cal O_{\Xdaginf /S},\cal
O_{\Xdaginf /S}).$$
\end{defi}
Le cas le plus important pour l'instant est l'extension $V\to K$, mais les extensions $V\to V_s, s\geq 1$ sont aussi intéressantes et semblent donner lieu à des questions non sans intérêt.
\subsection{Le foncteur de de Rham $\dagger$-adique local}
Soit $X$ un $R_1$-schéma   lisse et soit $\RmMod(R_X)$ la catégorie des faisceaux de $R$-module sur $X$ muni de la topologie de
Zariski. Nous allons définir un foncteur prolongement des faisceaux de $R$-modules:
$$ \RmMod(R_X)\rightarrow  \RmMod(R_{\Xdaginf }).\leqno P_X: $$ Si $\cal F_X$ un faisceau de $R$-modules sur $X$, on lui associe le faisceau $\cal
F_{\Xdaginf }$ dont la valeur    sur un ouvert  $\calUdag $ du site $\Xdaginf $   est le faisceau $\cal F_U$
et dont la valeur sur un morphisme du site $\calWdag \rightarrow \calUdag $
est  le morphisme de restriction $\cal F_U|W\to \cal F_W$. 
\begin{theo}\label{equ-fai}Le foncteur $P_X$\glossary{$P_X$} est une équivalence de catégories de la catégorie $\RmMod(R_X)$ dans la sous-catégorie des
$ R_{\Xdaginf }$-modules
sur le site infinitésimal $X\daginf$ sur lesquels l'action du faisceau de groupes $ \cal G_{\Xdaginf }$ est {\bf triviale}.
\end{theo}

\demo
Un $ R_{\Xdaginf }$-module  sur le site dont l'action de $ \cal G_{\Xdaginf }$ est triviale est un $R_{\Xdaginf }$-module $\cal F_{\inf}$
tel que pour tout ouvert $\calUdag $  l'action du faisceau $\cal G_{\calUdag }$ sur le $R$-module  $\cal F_{\calUdag }$ est triviale. Soit un 
couple d'ouverts 
$(\calUdag ,\calWdag)$ avec $r: W\subset  U$ et $W$ affine. Il existe alors un isomorphisme {\bf canonique} de $R$-modules : $$r^{-1}\cal F_{\calUdag }\simeq \cal F_{\cal
W\dag}.$$ En effet, en vertu du critère d'affinité $\dagger$-adique \ref{cri-aff} l'ouvert $\calUdag |W$ est affine et, en vertu du théorème des relèvements \ref{rel-mor}, l'inclusion $r: W\subset U$ admet des relèvements $r\dag$. Il existe des morphismes de restriction
$(\sharp)_{r\dag}: r^{-1}\cal F_{\calUdag }\simeq \cal F_{\cal
W\dag} $. Si $r\dag_1, r\dag_2$ sont deux relèvements, les restrictions $(\sharp)_{r_1\dag}, (\sharp)_{r_2\dag}$ diffèrent
par une restriction $(\sharp)_{g}$, où $g$ est une section de $\cal G_{\cal
W\dag}$ qui opère trivialement par hypothèse. Deux relèvements induisent la {\bf même} restriction. Pour tout couple  d'ouverts 
$\calUdag ,\calWdag $ avec $r: W\subset  U$ on obtient
un morphisme canonique
$r^{-1}\cal F_{\calUdag }\to \cal F_{\cal
W\dag} $, qui est défini localement et qui est un isomorphisme.
À un tel $R_{\Xdaginf }$-module   $\cal F_{\inf}$ et à un ouvert \em{affine} $U$ on associe le $R_U$-module  défini par 
$$\cal F_U:= \limproj_{\calUdag }\cal F_{\calUdag },$$ où la limite est prise pour tous les ouverts $\calUdag $ du site qui relèvent $U$. On obtient ainsi un
{faisceau de Zariski}
défini localement qui définit un inverse au foncteur $P_X$. 
\enddemo


\begin{Rema}En fait, le théorème  précédent vaut  pour les catégories des ensembles.
\end{Rema} 

\begin{coro}\label{hom-spe}
Si $\calNdaginf $ et $\calMdaginf $ sont deux $ \calDdag_{\Xdaginf /R}$-modules à gauche \em{spéciaux}, l'action du faisceau de groupes $ \cal G_{\Xdaginf }$   sur  le faisceau: $$ \cHom_{\calDdag_{\Xdaginf /R},\Sp}(\calNdaginf ,\calMdaginf )=\cHom_{\calDdag_{\Xdaginf /R}}(\calNdaginf ,\calMdaginf )$$ est triviale, qui est alors un  {\bf faisceau de Zariski} de $R_X$-modules sur $X$.
\end{coro}

\demo
En effet, par construction un élément du groupe $g$ agit par définition \ref{fai-hom} par $\varphi\mapsto g\varphi g^{-1}$, mais comme $g$ et  $g^{-1}$ sont des opérateurs
différentiels,  en vertu du théorème \ref{inc-grp}, et que
$\varphi$ est $ \calDdag_{\Xdaginf /R}$-{\bf linéaire}, cette action est {\bf  triviale} et on applique le théorème précédent.
\enddemo

\begin{coro}Soit  $\calNdaginf $ un complexe appartenant à la catégorie $ \Dmoins ((\calDdag_{\Xdaginf /R},\Sp)\Mod)$. Le foncteur 
$$ \calMdaginf \mapsto \bfR \cHom_{\calDdag_{\Xdaginf /R},\Sp}(
\calNdaginf ,\calMdaginf )$$ est un foncteur exact de catégories triangulées de $ \Dplus ((\calDdag_{\Xdaginf /R}, \Sp)\Mod)$ vers   $ \Dplus (R_X)$.
\end{coro}

\demo
En effet, le complexe $ \bfR \cHom_{\calDdag_{\Xdaginf /R},\Sp}(
\calNdaginf ,\calMdaginf )$ se représente par le complexe: $$ \cHom^\bullet_{\calDdag_{\Xdaginf /R}, \Sp}(
\calNdaginf ,\calIdaginf  )=\cHom^\bullet_{\calDdag_{\Xdaginf /R}}(
\calNdaginf ,\calIdaginf  )\,,$$ où $\calIdaginf  $ 
est une résolution injective de $\calMdaginf $
dans  la catégorie des {\bf modules spéciaux}, qui est un complexe pour la topologie de Zariski de $X$ en vertu du corollaire précédent. 
\enddemo

\begin{prop}Soient $\calNdaginf $ et $\calMdaginf $  deux $ \calDdag_{\Xdaginf /R}$-modules à gauche \em{\bf spéciaux}. Il
existe alors un isomorphisme canonique 
$$ \hom_{\calDdag_{\Xdaginf /R},\Sp}(
\calNdaginf ,\calMdaginf )\simeq \Gamma(X,\cHom_{\calDdag_{\Xdaginf /R},\Sp}(
\calNdaginf ,\calMdaginf ))\,.$$
\end{prop}

\demo  Comme :
$$ \hom_{\calDdag_{\Xdaginf /R},\Sp}(\calNdaginf , \calMdaginf )= \hom_{\calDdag_{\Xdaginf /R}}(\calNdaginf , \calMdaginf )\,,$$ un élément $\phi$ de 
$ \hom_{\calDdag_{\Xdaginf /R},\Sp}(\calNdaginf , \calMdaginf )$ est la donnée pour tout ouvert $\calUdag $ d'un morphisme $\calDdag_{\cal
U\dag/R}$-linéaire $\varphi_{\calUdag }: \cal N\dag_{\calUdag }\to\cal M\dag_{\calUdag }$ tel que pour tout morphisme $r\dag: \calWdag \to \cal
U\dag$  on ait l'égalité:
$$\varphi_{\calWdag }= (\sharp)_{r\dag}\circ \varphi_{\calUdag |W}\circ
(\sharp)^{-1}_{r\dag}.$$ 
En particulier,  pour tout ouvert $\calUdag $, le morphisme $\phi$  définit une section de  $$\Gamma(U, \cHom_{\calDdag_{\cal
U\dag/R}\Sp}(
\cal N\dag_{\calUdag },\cal M\dag_{\calUdag }))= \Gamma(U,\cHom_{\calDdag_{\Xdaginf /R}}(
\calNdaginf ,\calMdaginf )).$$ En prenant un recouvrement de $X$ par des ouverts affines, on trouve que $\phi$ définit 
des sections locales qui se recollent, donc une section globale du faisceau $\cHom_{\calDdag_{\Xdaginf /R},\Sp}(
\calNdaginf ,\calMdaginf )$. Réciproquement, une section globale de ce dernier faisceau  détermine un morphisme $\calNdaginf \to
\calMdaginf $.
\enddemo

\begin{Rema}
Dans le cas des modules {\bf spéciaux}, c'est là un point tout à fait remarquable qui est à la base du succès des modules spéciaux contrairement à la situation générale \ref{fai-hom}, une section globale du faisceau de Zariski:
$$\cHom_{\calDdag_{\Xdaginf /R},\Sp}(
\calNdaginf ,\calMdaginf )=\cHom_{\calDdag_{\Xdaginf /R}}(
\calNdaginf ,\calMdaginf )$$ est un {\bf morphisme global}:  $\calNdaginf \to\calMdaginf $.
\end{Rema}

\begin{coro}\label{loc-deR}
En particulier,  pour des complexes spéciaux $\calNdaginf $ et $\calMdaginf $,  la cohomologie du complexe $ \bfR \hom_{\calDdag_{\Xdaginf /R},\Sp}(
\calNdaginf ,\calMdaginf )$ n'est rien d'autre que l'hypercohomologie au-dessus de $X$ du complexe de {\bf Zariski}: $$ \bfR \cHom_{\calDdag_{\Xdaginf /R},\Sp}(
\calNdaginf ,\calMdaginf ),$$ exactement comme dans la théorie usuelle.
\end{coro}

Cela permet de construire une fort utile suite spectrale du passage du local au global. 
Notons $\widetilde{\cal Ext^i}_{\calDdag_{\Xdaginf /R},\Sp}(
\calNdaginf ,\calMdaginf )$ le préfaisceau sur $X$, qui à un ouvert $U$ associe l'hypercohomologie de degré $i$ au-dessus de $U$ du complexe
de Zariski précédent.

\begin{theo}\label{glo-loc}Soit $\cal B$ un recouvrement ouvert de $X$,  alors il existe une suite spectrale dont le terme $\EEE_2$ est donné par la cohomologie de
\v Cech du préfaisceau précédent $ H^{j}(\cal B,
\widetilde{\cal Ext^i}_{\calDdag_{\Xdaginf /R},\Sp}(
\calNdaginf ,\calMdaginf ))$ et dont le terme $\EEE_\infty$ est le $R$-module  bigradué associé à une filtration convenable du $R$-module 
gradué
$ \ext^{i+j}_{\calDdag_{\Xdaginf /R},\Sp}(
\calNdaginf ,\calMdaginf )$.
\end{theo}

\demo En effet, si l'on prend une résolution $\calIdaginf  $ de $\calMdaginf $ par des modules spéciaux injectifs,   le complexe $ \cHom^{\bullet}_{\calDdag_{\Xdaginf /R},\Sp}(
\calNdaginf ,\calIdaginf  )$ est un complexe de faisceaux flasques, puisque la propriété
d'être flasque est locale sur $X$. La suite spectrale du théorème  est la suite spectrale du complexe double de \v Cech 
$ C^{\bullet}(\cal B,\cHom^{\bullet}_{\calDdag_{\Xdaginf /R},\Sp}(
\calNdaginf ,\calIdaginf  )).$
\enddemo

\begin{prop}\label{may-vie} Soient $X_1\cup X_2$ un recouvrement ouvert de $X$, 
$\calMdaginf $ un complexe spécial borné à gauche et $\calNdaginf $ un complexe spécial borné à droite. Il existe alors un triangle distingué de Mayer-Vietoris 
$$\def\hom_#1{\mathop{\rm Hom}\nolimits\goodSub{0pt}{4pt}{-8mm}{#1}}
\DeuxLignes  
 \bfR \hom_{\calDdag_{\Xdaginf /R},\Sp}(
\calNdaginf ,\calMdaginf )
\to
 \bfR \hom_{\calDdag_{(X_1)\daginf/R},\Sp}(
\cal N\dag_{1\inf},\cal M\dag_{1\inf})\oplus \bfR \hom_{\calDdag_{(X_2)\daginf/R},\Sp}(
\cal N\dag_{2\inf},\cal M\dag_{2\inf})\to
\\\to
\bfR \hom_{\calDdag_{(X_{12})\daginf/R},\Sp}(
\cal N\dag_{12\inf},\cal M\dag_{12\inf})\to\ ,
\endlignes$$
où $X_{12}:=X_1\cap X_2$ et 
$\cal M\dag_{1\inf},\cal M\dag_{2\inf},\cal M\dag_{12\inf}$
désignent les restrictions de $\calMdaginf $ à $X_1, X_2$ et $ X_1\cap X_2$. 
\end{prop}

\demo En effet, si $\calIdaginf  $ est une résolution injective de $\calMdaginf $ par des modules spéciaux, le complexe de Zariski 
$\cHom^\bullet_{\calDdag_{\Xdaginf /R},\Sp}(
\calNdaginf ,\calIdaginf  )$ est à {\bf termes flasques} et donne donc naissance au triangle distingué de Mayer-Vietoris.
\enddemo

\begin{defi} Soit $\calMdaginf $ un complexe spécial, on définit son complexe de de Rham local $\dagger$-adique:
$$ \dR(\calMdaginf ):= \bfR \cHom_{\calDdag_{\Xdaginf /R},\Sp}(
\cal O_{\Xdaginf /R},\calMdaginf ).$$ On obtient ainsi un foncteur covariant exact de catégories triangulées
$$  \Dplus ((\calDdag_{\Xdaginf /R}, \Sp)\Mod)\to  \Dplus (R_X)\,.\leqno \dR(X/R,-):$$
\end{defi}
En particulier, si $Z\subset X$ est un fermé, le complexe de cohomologie locale $ \bfR \Gamma_Z(\dR(\calMdaginf ))$ est un complexe de Zariski
parfaitement défini dont l'hypercohomologie $ \bfR \Gamma_Z(X, \dR(\calMdaginf ))$ fournit la cohomologie de de Rham $\dagger$-adique à
support dans $Z$, donnant lieu aux triangles de cohomologie locale usuels.

\begin{prop} Soit
$\calXdag $ un relèvement plat  de $X$, et soient $\cal N\dag_{\calXdag }$ un complexe de $ \Dmoins(\calDdag_{\calXdag /R})$ et  $\cal M\dag_{\calXdag }$ un complexe de $ \Dplus (\calDdag_{\calXdag /R})$. Alors, il existe un isomorphisme canonique de la catégorie $ \rmD(R_X)$:
$$ \bfR \cHom_{\calDdag_{\calXdag /R}}(
\cal N\dag_{\calXdag },\cal M\dag_{\calXdag })\simeq\bfR \cHom_{\calDdag_{\Xdaginf /R},\Sp}(
P_{\calXdag }(\cal N\dag_{\calXdag }),P_{\calXdag }(\cal M\dag_{\calXdag })).$$
\end{prop}

\demo En effet, pour deux $\calDdag_{X\dag/R}$-modules  à gauche on dispose d'un isomorphisme canonique 
$$ \cHom_{\calDdag_{\calXdag /R}}(
\cal N\dag_{\calXdag },\cal M\dag_{\calXdag })\simeq\cHom_{\calDdag_{\Xdaginf /R}}(
P_{\calXdag }(\cal N\dag_{\calXdag }),P_{\calXdag }(\cal M\dag_{\calXdag })),$$ qui se dérive naturellement.
\enddemo

\begin{Rema}
Si $R\to S$ est un morphisme d'anneaux,  on peut considérer le    faisceau d'anneaux sur le site $\calDdag_{X\daginf/S}:= \calDdag_{X\daginf/R}\otimes_RS$, et la théorie précédente a lieu sur  ce faisceau.
\end{Rema}

\subsection{Le théorème de finitude des nombres de Betti $p$-adiques d'une variété algébrique lisse}
Jusqu'à présent, nous avons construit la théorie sur un anneau noethérien   $R$ muni de la topologie définie par un idéal $I$, ce qui a un grand intérêt, mais les nombres de Betti \ps
d'une variété algébrique sont des dimensions d'espaces vectoriels sur le corps $K$ des fractions d'un anneau de valuation discrète complet $V$ d'inégales
caractéristiques. Aussi, il nous faut considérer les faisceaux d'anneaux:
$$ K_{\Xdaginf }\subset\cal O_{\Xdaginf /K}:= \cal O_{\Xdaginf /V}\otimes_VK\subset\calDdag_{\Xdaginf /K}:= \calDdag_{\Xdaginf /V}\otimes_VK\,,$$
 sur le site
$\Xdaginf $ infinitésimal pour la topologie $\goth m$-adique de $V$.\glossary{$\cal O_{\Xdaginf /K}$}\glossary{$\calDdag_{\Xdaginf /K}$}
Le faisceau de groupes $\cal G_{\Xdaginf }$ agit sur tous ces faisceaux.  On peut alors considérer la catégorie  $\RmMod(K_{\Xdaginf})$ de faisceaux d'espaces vectoriels sur le site infinitésimal, la catégorie $ \calDdag_{\Xdaginf /K}\Mod$
des $\Rm\calDdag_{\Xdaginf /K}$-modules à gauche et sa sous-catégorie pleine $ (\calDdag_{\Xdaginf /K}, \Sp)\Mod$ des
$ \calDdag_{\Xdaginf /K}$-modules à gauche spéciaux qui est une catégorie abélienne
admettant suffisamment d'objets injectifs.

\begin{defi}Si $X$ est un schéma lisse sur $k$, on définit ses nombres de Betti $p$-adiques comme les dimensions $B_{p,i}(X)$ de ses
$K$-espaces
de cohomologie de de Rham $p$-adique:
$$ H^i_{\DR}(X/K):= H^i_{\DR}(X/K,\cal O_{\Xdaginf /K}):= \ext^{i}_{\calDdag_{\Xdaginf /K},\Sp}(
\cal O_{\Xdaginf /K},\cal O_{\Xdaginf /K}).$$
\end{defi}
Nous rappelons qu'une variété algébrique sur un corps est un schéma séparé et de type fini sur ce corps [EGA I].
\begin{theo}\label{fin}Si $X$ est une variété algébrique lisse sur $k$, ses nombres de Betti $p$-adiques $B_{p,i}(X)$ sont \em{finis} pour tout $i\geq0$.
\end{theo}

\demo
Considérons un recouvrement fini $\cal B$ de $X$ par des ouverts affines. La filtration de la suite spectrale local-global \ref{glo-loc} sur 
le groupe gradué $ H^\bullet_{\DR}(X/K)$ est finie. Les termes $\EEE_2$ sont des produits finis de cohomologie de de Rham $p$-adiques de variétés 
affines lisses qui, en vertu du théorème de finitude [Me$_2$], sont de dimension finie sur $K$. Cela entraîne que les nombres $B_{p,i}(X)$
sont finis.
\enddemo

\subsection{La fonctorialité de la cohomologie de de Rham $p$-adiques pour les schémas affines}
Si $X$ est un schéma lisse sur $k$ et admet un relèvement  plat  $\calXdag $ sur $V$, la catégorie des $\calDdag_{\calXdag /K}$-modules à gauche  est canoniquement
équivalente à la catégorie  des modules spéciaux 
$ (\calDdag_{\Xdaginf /K}, \Sp)\Mod$  en vertu du théorème fondamental \ref{can}. En particulier, 
les espaces de cohomologie $\ext^{i}_{\calDdag_{\calXdag /K}}(
\cal O_{\calXdag /K},\cal M\dag_{\calXdag })$ sont canoniquement isomorphes aux espaces de cohomologie de de Rham du site infinitésimal: 
$$ \ext^{i}_{\calDdag_{\Xdaginf /K},\Sp}(
\cal O_{\Xdaginf /K},P_{\calXdag }(\cal M\dag_{\calXdag }))$$ pour tout complexe $\cal M\dag_{\calXdag }$ de $ \Dplus (\calDdag_{\calXdag /K}\Mod)$. 

Si $X$ est affine et lisse  et si $\calXdag $ est un relèvement  plat d'algèbre $A\dag$, il est montré  ([M-N$_1$], [M-N$_2$]) que 
les espaces $\ext^{i}_{\calDdag_{\calXdag /K}}(
\cal O_{\calXdag /K},\cal O_{\calXdag /K})$ coïncident canoniquement avec la cohomologie du complexe de de Rham $\Omega^{\bullet}_{A\dag/K}$ des formes
différentielles séparées de l'algèbre $A\dag$ ([M-W]). Nous allons reprendre le principe de la démonstration de ce résultat pour faire le lien avec le foncteur
construit dans [M-W], entre la catégorie des variétés algébriques affines lisses sur $k$ et la catégorie homotopique $ K^b(K)$ des complexes d'espaces
vectoriels sur $K$.

Supposons plus généralement que $\calXdag $ est un relèvement  plat d'un $k$-schéma lisse $X$ et considérons le complexe de Spencer: 
$$\DeuxLignes
\Sp^{\bullet}(\cal O_{\calXdag/K}):= 0\rightarrow \calDdag_{\calXdag/K}\otimes_{\cal O_{\calXdag/K}}
\smash{\textstyle\bigwedge\nolimits^n} \cal T_{\calXdag/K}\rightarrow\cdots
\\
\cdots\rightarrow\calDdag_{\calXdag/K}\rightarrow \cal O_{\calUdag /K}\rightarrow0
\endlignes
$$ 
où $\bigwedge^i\cal T_{\calXdag /K}$ est l'algèbre extérieure du fibré tangent et la
différentielle est définie par:
$$\ecartlignes15pt\DeuxLignes  
\delta(P\otimes v_{1}\wedge\cdots\wedge v_{i}):= 
\smash{\sum_{j=1,\dots,i} }(-1)^{i-1}Pv_j\otimes(v_{1}\wedge\cdots\wedge\hat{v}_j\wedge\cdots\wedge v_{i})-
\\
\sum_{1\leq j< l\leq i}
(-1)^{i+l}P\otimes([v_j,v_l]\wedge v_{1}\wedge\cdots\wedge\hat{v}_j\wedge\cdots\wedge\hat{v}_l\wedge\cdots\wedge v_{i})\,,
\endlignes$$ où $P$ est un opérateur différentiel et les $v_j$
sont des vecteurs tangents. 

On rappelle le résultat suivant ([Me$_3$], 7.2.5), qui  ne se démontre pas uniquement à partir du complexe de Koszul vu que les ordres des opérateurs différentiels sont infinis.

\begin{lemm}\label{Spen} Le complexe de Spencer $ \Sp^{\bullet}(\cal O_{\cal 
X\dag/K})$ est une résolution du faisceau structural $\cal O_{\cal 
X\dag/K}$ par des $ \calDdag_{
\calXdag /K}$-modules à gauche localement libres de type fini. En particulier, $\cal O_{\cal 
X\dag/K}$ est un  $ \calDdag_{
\calXdag /K}$-module à gauche parfait.
\end{lemm}

Le complexe de de Rham $\Omega^{\bullet}_{\calXdag /K}$ est canoniquement isomorphe au complexe:
$$\cHom_{\calDdag_{\cal 
X\dag/K}}(\Sp^{\bullet}(\cal O_{\cal 
X\dag/K}), \cal O_{\cal 
X\dag/K})\,.$$
Soit $\cal O_{\cal 
X\dag/K}\to\cal I^{\bullet}_{\calXdag }$ une résolution injective de $\cal O_{\cal 
X\dag/K}$ par des $\calDdag_{\cal 
X\dag/K}$-modules à gauche. Alors, les morphismes canoniques de  complexes  :
$$\DeuxLignes
\Omega^{\bullet}_{\calXdag /K}\simeq\cHom_{\calDdag_{\cal 
X\dag/K}}(\Sp^{\bullet}(\cal O_{\cal 
X\dag/K}), \cal O_{\cal 
X\dag/K})\rightarrow
\\
\rightarrow \cHom^{\bullet}_{\calDdag_{\cal 
X\dag/K}}(\Sp^{\bullet}(\cal O_{\cal 
X\dag/K}), \cal I^{\bullet}_{\cal 
X\dag/K})\leftarrow\cHom_{\calDdag_{\cal 
X\dag/K}}(\cal O_{\cal 
X\dag/K}, \cal I^{\bullet}_{\cal 
X\dag/K})\endlignes$$ 
induisent des isomorphismes en cohomologie. Ce sont donc des isomorphismes canoniques de $\Rm
\Db (K_X)$. Cela entraîne que les hypercohomologies des complexes $\Omega^{\bullet}_{\calXdag /K}$ et $\cHom_{\calDdag_{\cal 
X\dag/K}}(\cal O_{\cal 
X\dag/K}, \cal I^{\bullet}_{\cal 
X\dag/K})$ sont canoniquement isomorphes. Mais si $X$  est affine, l'hypercohomologie du complexe $\Omega^{\bullet}_{\calXdag /K}$ se réduit à la cohomologie du complexe des sections globales, en vertu du théorème d'acyclicité de Meredith \ref{acy}.

\begin{theo}\label{ind-rel}
Soient $X$ un  schéma lisse sur $k$ et $\cal X\dag_1, \cal X\dag_2$ deux  relèvements  plats de
$X$ et
$f\dag_1$, $f\dag_2$ des relèvements de l'identité. Alors, les morphismes de complexes: $$\Omega^{\bullet}_{\cal X\dag_2/K}\to\Omega^{\bullet}_{\cal X\dag_1/K}\leqno f^{ *}_1, f^{ *}_2:$$ définissent le {\bf même élément}  du groupe $
 \hom_{\Db (K_X\Mod)}(\Omega^{\bullet}_{\cal X\dag_2/K},\Omega^{\bullet}_{\cal X\dag_1/K})$.
\end{theo}

\demo Rappelons,  en vertu de la prop. \ref{fai-hom}, que si $\calMdaginf $ et $\calNdaginf $ sont deux $\calDdag_{ 
\Xdaginf /K}$-modules à gauche, le faisceau $\cHom_{\calDdag_{ 
\Xdaginf /K}}(\calNdaginf , \calMdaginf )$ est un faisceau infinitésimal de $K$-espaces vectoriels.
Si $f\dag$ est un morphisme du site $\cal X\dag_1\to\cal X\dag_2,$ le morphisme de restriction 
$$\cHom_{\calDdag_{ 
\cal X\dag_{2}/K}}(\cal N\dag_{\cal X_2}, \cal M\dag_{\cal X_2})\to\cHom_{\calDdag_{\cal 
X\dag_{1}/K}}(\cal N\dag_{\cal X_1}, \cal M\dag_{\cal X_1})\leqno (\sharp)_{f\dag}:$$ est donné par la définition \ref{fai-hom}
par $\varphi\mapsto (\sharp)_{f\dag}\circ\varphi\circ((\sharp)_{f\dag})^{-1}
$. Mais si $\cal X\dag_1=\cal X\dag_2$ et que
$\calMdaginf $ et $\calNdaginf $ sont deux $\calDdag_{ 
\Xdaginf /K}$-modules à gauche {\bf spéciaux}, alors $\varphi$ et les restrictions $(\sharp)_{f\dag}$ sur ces modules  {\bf commutent}, et les restrictions $(\sharp)_{f\dag}$  
sur le faisceau \relax{$ \cHom_{\calDdag_{ 
\cal X\dag_{1}/K}}(\cal N\dag_{\cal X_1}, \cal M\dag_{\cal X_1})$} sont {\bf triviales}.

Soit $\cal O_{ 
\Xdaginf /K}\to\cal I^{\bullet}_{ \Xdaginf }$ une résolution injective de $\cal O_{ 
\Xdaginf /K}$ par des $\calDdag_{ 
\Xdaginf /K}$-modules à gauche spéciaux, par exemple $P_{\cal X\dag_1}(\cal I^{\bullet}_{\cal X\dag_1})$. Alors,
pour chaque relèvement $f\dag_i$ on a un diagramme commutatif {\bf de morphismes de complexes} qui sont des {\bf isomorphismes dans la catégorie dérivée}:
$$\hss\scriptspace0.pt\mathrigid2mu
\def\cHom_#1{\mathop{\cal H\!\it om}\nolimits\sub{0pt}{#1}}
\def\Un#1{\Omega^{\bullet}_{\cal X_{#1}\dag/K}}
\def\Deux#1{\cHom_{\calDdag_{\cal 
X\dag_{#1}/K}}(\Sp\sp{\bullet}(\cal O_{\cal 
X\dag_{#1}/K}), \cal O_{\cal 
X\dag_{#1}/K})}
\def\Trois#1{\cHom_{\calDdag_{\cal 
X\dag_{#1}/K}}\sp{\bullet}(\Sp\sp{\bullet}(\cal O_{\cal 
X\dag_{#1}/K}), R_{\cal X\dag_{#1}}(\cal I^{\bullet}_{
\Xdaginf /K}))}
\def\Quatre#1{\cHom_{\calDdag_{\cal 
X\dag_{#1}/K}}(\cal O_{\cal 
X\dag_{#1}/K}, R_{\cal X\dag_{#1}}(\cal I^{\bullet}_{
\Xdaginf /K}))}
\let\quad\relax
\matrix{
\Un2&\smash{\varto{0.5cm}{4.3cm}_{(\sharp)_{f\dag_i}=f\sb{i}\sp{*}}}&\Un1\cr\noalign{\kern0.5mm}
\vegal{}{\simeq}{0.4cm}&&\vegal{}{\simeq}{0.4cm}\cr
\Deux2&\smash{\varto{0cm}{1.8cm}_{(\sharp)_{f\dag_i}}}&\Deux1\cr
\downarrow&&\downarrow\cr
\Trois2&\smash{\varto{0.9cm}{0.8cm}_{(\sharp)_{f\dag_i}}}&\Trois1\cr
\uparrow&&\uparrow\cr
\Quatre2&\smash{\varto{0cm}{1.6cm}_{(\sharp)_{f\dag_i}}}&\Quatre1\cr\noalign{\kern10pt}
}\hss$$
 Mais en vertu du théorème \ref{iso-pla}, les morphismes $f^{\dagger}_1$ et $ f^{\dagger}_2$ sont inversibles et $f^{\dagger }_2 =
g\circ f^{\dagger }_1$, où $g$ est un automorphisme d'algèbre de $\cal O_{\cal X\dag_1/V}$ qui se réduit à l'identité modulo $\goth m$, 
c'est donc un opérateur
différentiel en vertu du théorème \ref{inc-grp}, et cela montre 
$f^{\dagger }_2$ et $f^{\dagger }_1$ induisent le {\bf même} morphisme sur la dernière ligne du diagramme précédent, soit:
$$\vrule height13pt width0pt\cHom\SubX{\calDdag_{\cal 
X\dag_2/K}}\mBig(\cal O_{\cal 
X\dag_2/K}, R_{\cal X\dag_2}(\cal I^{\bullet}_{
\Xdaginf /K})\mBig)\to \cHom\SubX{\calDdag_{\cal 
X\dag_1/K}}\mBig(\cal O_{\cal 
X\dag_1/K}, R_{\cal X\dag_1}(\cal I^{\bullet}_{
\Xdaginf /K})\mBig)\,,\leqno (\sharp)_{f\dag_1}=(\sharp)_{g}\circ(\sharp)_{f\dag_1}=(\sharp)_{f\dag_2}:$$ puisque la restriction $(\sharp)_{g}$ 
sur
$\cHom_{\calDdag_{\cal 
X\dag_1/K}}\big(\cal O_{\cal 
X\dag_1/K}, R_{\cal X\dag_1}(\cal I^{\bullet}_{
\Xdaginf /K})\big)$
est triviale. Les morphismes $f_1\sp{*}$ et $f_2\sp{*}$
définissent donc le {\bf même élément} du groupe $ \hom_{\Db (K_X\Mod)}(\Omega^{\bullet}_{\cal X\dag_2/K},\Omega^{\bullet}_{\cal X\dag_1/K})$.
\enddemo

\begin{coro}
Si $ A\dag_1$ et  $ A\dag_2$ sont deux relèvements  plats d'une algèbre $A$ affine lisse sur $k$, alors  tous les morphismes entre  les complexes
$\Omega^{\bullet}_{A\dag_1/K}$ et  $\Omega^{\bullet}_{A\dag_2/K}$ induits par des relèvements $ A\dag_2\to A\dag_1$ de l'identité sont
{\bf homotopes}.
\end{coro}

\demo En effet, en vertu du théorème précédent tous les relèvements induisent le même morphisme de $\Rm
\hom_{\Db (K)}(\Omega^{\bullet}_{ A\dag_2/K},\Omega^{\bullet}_{A\dag_1/K})$, mais: $$ \hom_{K^b(K)}(\Omega^{\bullet}_{
A\dag_2/K},\Omega^{\bullet}_{A\dag_1/K})=\hom_{\Db (K)}(\Omega^{\bullet}_{
A\dag_2/K},\Omega^{\bullet}_{A\dag_1/K}),$$ où 
$K^b(K)$ désigne comme d'habitude la catégorie des complexes bornés de $K$-espaces vectoriels à homotopie près.
\enddemo

\begin{Rema}
La méthode de Monsky-Washnitzer ([M-W], remark p. 203) en construisant une homotopie explicite, montre  en fait pour $X$ affine et lisse sur $k$, que tous les
morphismes induits par des relèvements définissent le {\bf même} élément dans le groupe: $$ \hom_{K^b(K_X\Mod)}(\Omega^{\bullet}_{\cal X\dag_2/K},\Omega^{\bullet}_{\cal X\dag_1/K}).$$ La méthode précédente montre seulement que c'est le même élément dans le groupe: $$\Rm
\hom_{\Db (K_X\Mod)}(\Omega^{\bullet}_{\cal X\dag_2/K},\Omega^{\bullet}_{\cal X\dag_1/K}).$$ Cette subtilité disparaît
au niveau des complexes globaux puisque la catégorie homotopique de la catégorie des espaces vectoriels sur un corps coïncide avec sa catégorie dérivée. Le passage de la catégorie homotopique à la catégorie
dérivée  au niveau des complexes de faisceaux est le prix à payer pour passer des variétés
affines aux variétés non affines, et surtout pour passer du cas du fibré trivial au cas d'un module spécial. Mais en même temps, 
c'est une simplification considérable de la construction originale de Monsky-Washnitzer puisqu'elle 
ne nécessite pas une homotopie explicite, qui est le point délicat de la construction originelle et qui ne semble pas se recoller.
\end{Rema}

\begin{coro}Soient $f: Y\to X$ un morphisme de schémas  {\bf affines} lisses sur $k$, et $\calYdag $  et $\calXdag $ des relèvements 
plats. Alors, tous les morphismes $$f^{-1}\Omega^{\bullet}_{\calXdag/K}\to\Omega^{\bullet}_{\calYdag /K}$$
induits par des relèvements $f\dag$ définissent le {\bf même} élément du groupe: 
$$\Rm
\hom_{\Db (K_Y\Mod)}(f^{-1}\Omega^{\bullet}_{\calXdag /K},\Omega^{\bullet}_{\cal
Y\dag/K})\,.$$
\end{coro}

\demo 
En effet, si $f\dag_1, f\dag_2$ sont deux relèvements de $f$, on a les factorisations canoniques:
$f\dag_i= i_{f\dag_i}\circ p\dag$, $i=1,2$, où $i_{f\dag_i}$ est l'immersion graphe $\calYdag \to \calYdag \times_{V}\calXdag $ de $f\dag_i$
et $p\dag$ est la projection canonique $\calYdag \times_{V}\calXdag \to \calXdag $. Cela ramène à supposer que $f$ est une immersion fermée.
Dans ce cas   affine, et en vertu de la propriété 3) de la lissité, il existe un automorphisme 
$g$ de $\calXdag $ tel que $f\dag_2= g\circ f\dag_1$. On est réduit au théorème précédent.
\enddemo

\begin{coro}\label{fon-aff} La correspondance $X\mapsto \Omega^{\bullet}_{
A\dag/K}$ définit un foncteur contravariant $$ \SmAff(k)\to K^+(K)\leqno \DR(-/K):$$ de la catégorie des schémas affines lisses $ \SmAff(k)$ sur $k$
dans la catégorie homotopique
$ K^+(K)$.
\end{coro}
\glossary{$\DR(X/K)$}

\demo En effet, dans le cas affine un morphisme $f$ admet toujours des relèvements (\ref{rel-mor}) et on est réduit au corollaire précédent.
\enddemo

Le foncteur ainsi construit coïncide,  par construction, avec celui de Monsky-Washnitzer [M-W], mais on ne peut pas étendre à ce stade la fonctorialité au cas non
affine
parce qu'en général un morphisme ne se relève pas. L'un des buts, qui est aussi un test,  des chapitres
suivants, est d'étendre cette fonctorialité à la catégorie de tous les  schémas lisses séparés sur
$k$. On est encouragé de retrouver par cette méthode la fonctorialité du cas affine de Monsky-Washnitzer ([M-W]), sans quoi rien n'est possible.

\begin{coro}Soient $X$ un  schéma lisse sur $k$ et  $\calMdaginf $  un $ \calDdag_{\Xdaginf /K}$-module à gauche spécial. Les faisceaux de Zariski: 
$$\cal Ext^i_{\calDdag_{\Xdaginf /K},\Sp}(
\cal O_{\Xdaginf /K},\calMdaginf )$$ sont nuls si $i>\dim X$ et les espaces: 
$$\ext^i_{\calDdag_{\Xdaginf /K},\Sp}(
\cal O_{\Xdaginf /K},\calMdaginf )$$ sont nuls si $i>2\dim X$. En particulier, les nombres de Betti $B_{p,i}(X)$ sont nuls si $i>2\dim X$.
\end{coro}

\demo  En effet, si $\calXdag $ est un relèvement plat, le complexe de faisceaux de Zariski $\bfR \cHom_{\calDdag_{\Xdaginf /K},\Sp}(
\cal O_{\Xdaginf /K},\calMdaginf )$ est isomorphe au-dessus de $X$ au complexe 
de faisceaux $\cHom_{\calDdag_{\cal 
X\dag/K}}(\Sp^{\bullet}(\cal O_{\cal 
X\dag/K}), \cal M\dag_{\cal 
X\dag})$, concentré en degrés $[0, \dim X]$. La seconde assertion résulte du fait que la dimension topologique de $X$ est bornée par la dimension 
algébrique $\dim X$ en vertu du théorème de Grothendieck.
\enddemo
\begin{Rema} On trouve par construction la suite spectrale conjuguée: $$\rm H^{j}(X,\cal Ext^i_{\calDdag_{\Xdaginf /K},\Sp}(
\cal O_{\Xdaginf /K},\calMdaginf ))\Longrightarrow Ext^{i+j}_{\calDdag_{\Xdaginf /K},\Sp}(
\cal O_{\Xdaginf /K},\calMdaginf ),$$ et l'on s'attend, bien sûr, dans le cas géométrique $(\calMdaginf =\cal O_{\Xdaginf /K})$ à ce que la filtration de cette suite spectrale ait une interprétation géométrique.\end{Rema}
\subsection{Le théorème de comparaison dans le cas d'un relèvement propre et lisse}
On trouve ou on retrouve comme corollaire l'invariance  de la cohomologie de de Rham {\bf algébrique} 
d'un relèvement  propre et lisse,  sans invoquer la cohomologie cristalline qui  a fourni  historiquement la première démonstration de ce résultat [B-O].

\begin{coro} Soit ${\mathbf X}/V$ un relèvement \em {algébrique} propre et lisse de $X/k$.  La cohomologie de de Rham algébrique usuelle $\Rm
H^{\bullet}_{\DR}({\mathbf{X}}/K):= H^{\bullet}({\mathbf{X}}, \Omega^{\bullet}_{{\mathbf{X}}/K})$ est canoniquement isomorphe à la 
cohomologie de de Rham $p$-adique  $ H^\bullet_{\DR}(X /K, \cal O_{\Xdaginf /K})$ et 
est indépendante à isomorphisme canonique près du relèvement ${\bf
X}/V$.
\end{coro}

\demo  Le relèvement  algébrique ${\mathbf{X}}/V$ induit un relèvement  plat $\calXdag $ ([Mr]). Définissons:  
$$ H^i_{\DR}(\calXdag /K):= H^i_{\DR}(X/K,\cal O_{\calXdag  /K}):= \ext^{i}_{\calDdag_{\calXdag  /K}}(
\cal O_{\calXdag  /K},\cal O_{\calXdag  /K}).$$
Au vu des  résultats précédents,  on a alors l'isomorphisme canonique 
$$ H^i_{\DR}(\calXdag /K)\simeq H^i_{\DR}(X/K, \cal O_{\Xdaginf /K}).$$
En vertu du théorème de comparaison de Meredith [Mr]  pour les schémas propres, le morphisme canonique  
$$ H^{i}_{\DR}({\mathbf{X}}/K)\to   H^i_{\DR}(\calXdag /K)$$ est un isomorphisme de $K$-espaces vectoriels, ce qui montre le corollaire. \enddemo

\begin{Rema}Le morphisme naturel: $$ \cHom_{\calDdag_{\Xdaginf /V}}(
\cal O_{\Xdaginf /V},\calMdaginf )\otimes_VK\to \cHom_{\calDdag_{\Xdaginf /K}}(
\cal O_{\Xdaginf /K},\calMdaginf \otimes_VK)$$ \em{n'est pas} un isomorphisme en général. En particulier, la cohomologie de de Rham \em{ne commute
pas} au changement de base
$V\to K$. C'est une des raisons pour laquelle il faut modifier la définition de la cohomologie de de Rham $\dagger$-adique sur $R$, mais nous n'insisterons pas dans cet article sur ce point.
\end{Rema}

\section{La catégorie des modules à gauche spéciaux et la cohomologie de de Rham $p$-adique  $ H^\bullet_{\DR,h}(\Xdaginf /V)$ d'échelon $h\geq{}$0}

\subsection{La filtration par les échelons du faisceau $ \calDdag_{\Xdaginf /R}$}
Dans le cas où $R$ est une $\Bbb Z_{(p)}$-algèbre où $\Bbb Z_{(p)}$ est le localisé de $\Bbb Z$ en l'idéal $(p)$
pour un nombre premier
$p>0$, et si $\calXdag $ est un relèvement  plat  d'un $R_1$-schéma lisse $X$,
nous avons défini dans les articles ([M-N$_1$], [Me$_3$]) une filtration croissante et
exhaustive du faisceau \smashbot{$\calDdag_{\calXdag /R}$} par les faisceaux d'anneaux $\cal
D^{\dagger,h}_{\calXdag /R}$\glossary{$\cal
D^{\dagger,h}_{\calXdag /R}$} des opérateurs différentiels d'échelon $h\geq0$,  que nous appelons \og la $p$-filtration\fg\
ou \og la filtration par les échelons\fg\ qui est similaire à la $p$-filtration du faisceau $\calD_{X/k}$ des opérateurs différentiels sur un schéma lisse
sur un corps de caractéristique $p>0$. Rappelons-en la définition ([M-N$_1$], [Me$_3$], 2.6.1):

\begin{defi}Soit $h\geq0$ un entier. On définit le faisceau $\calD^{<\infty,h}_{\calXdag /R}$ des opérateurs différentiels d'ordre fini et d'échelon $h$,
comme le faisceau de 
sous-$\cal O_{\calXdag /R}$-algèbres du faisceau  $\calD_{\calXdag /R}$,
localement engendré par les opérateurs d'ordre $\leq p^h$.
\end{defi}
Le faisceau  $\calD^{<\infty,h}_{\calXdag /R}$ est un sous-faisceau de $\cal O_{\calXdag /R}$-algèbres du faisceau
 $\calD_{\calXdag /R}$, qui est aussi localement libre comme $\cal O_{\calXdag /R}$-module [M-N$_1$].
 
\bigskip
Pour tout $h\geq 0$, nous pouvons considérer le complété $I$-adique $\calD^{\wedge,h}_{\cal X^\wedge/R}$ du faisceau $\calD^{<\infty,h}_{\calXdag /R}$.
Nous pouvons aussi considérer le complété $I$-adique $(X,\cal O_{\cal X^\wedge/R})$ du schéma $\dagger$-adique $(X,\cal O_{\calXdag /R})$ et le faisceau  $\calD^{\wedge,h}_{\cal
X^\wedge/R}$, qui opère de façon évidente sur le faisceau $\cal O_{\cal X^\wedge/R}$,  apparaît comme un sous-faisceau de $R$-algèbres du faisceau des endomorphismes $\cal
End_{R}(\cal O_{\cal X^ \wedge/R})$. Pour tout entier $s\geq 1$, nous pouvons considérer la réduction modulo l'idéal $I^s$
 du faisceau $\calD^{\wedge,h}_{\cal X^\wedge/R}$:
$$\calD^{\wedge,h}_{\cal X^\wedge/R}
I^s= \calD^{<\infty,h}_{\calXdag /R}/I^s.$$ C'est de façon naturelle un faisceau de $\cal
O_{X_s}$-algèbres filtré sur le schéma $X_s$,
réduction modulo  l'idéal $I^s$ 
du schéma $\dagger$-adique
$(X,\cal O_{\calXdag /R})$. 

\medskip
L'inclusion $\cal
D^{<\infty,h}_{\calXdag /R}\subset \cal
D_{\calXdag /R}$ induit,  pour tout $s\geq1$, un morphisme filtré: 
$$\cal
D^{<\infty,h}_{\calXdag /R}/I^s\rightarrow  \cal
D_{X_s
/R_s},$$  qui, bien entendu,  \em{n'est pas injectif}.

\bigskip
Le faisceau $\calDdag_{\calXdag/R}$, qui opère de façon évidente sur le faisceau $\cal O_{\cal X^\wedge/R}$,  apparaît alors aussi comme un sous-faisceau de $R$-algèbres du faisceau des
endomorphismes $\cal End_{R}(\cal O_{\cal X^\wedge/R})$. On peut donc pour tout $h\geq0$ considérer l'intersection $\calDdag_{\calXdag/R}\cap\calD^{\wedge,h}_{\cal X^\wedge/R}$ prise dans le faisceau $\cal End_{R}(\cal
O_{\cal X^\wedge/R})$. Rappelons la définition ([Me$_3$], 2.6.5):

\begin{defi}\label{def-ech}Sous les conditions précédentes, soit un entier $h\geq0$ et soit $\calXdag  =(X,\cal O_{\calXdag /R})$ un $R$-schéma $\dagger$-adique lisse, on définit le faisceau des opérateurs
d'ordre infini et d'échelon $h$ noté $\calD^{\dagger,h}_{\calXdag /R}$, comme le sous-faisceau du faisceau $\calDdag_{\calXdag/R}\cap\calD^{\wedge,h}_{\cal X/R}$ des opérateurs dont la réduction modulo $I^s$ dans $\cal
D^{\wedge,h}_{\calXdag /V}/I^s$ est d'ordre localement borné par une fonction linéaire en $s$.
\end{defi}

\begin{Rema}Le lecteur prendra garde que l'inclusion $\calD^{\dagger,h}_{\calXdag /R}\subset \calDdag_{\calXdag/R}\cap\calD^{\wedge,h}_{\cal X^\wedge/R}$ est {\bf stricte}.
\end{Rema}

\begin{notation}Si $\beta =(\beta_1,\dots,\beta_n)$ est un $n$-uplet et $h\geq 0$, on  note  dans l'énoncé suivant 
 $\beta\in \Bbb N^n_{\leq p-1}$ si $\beta_i\leq p-1$ et: $$(\Delta^{p^h})^\beta:=
(\Delta_1^{p^h})^{\beta_1}\cdots(\Delta_n^{p^h})^{\beta_n}.$$ 
\end{notation}

\begin{theo}\label{sym-toth}Soient   $\calXdag $  un schéma $\dagger$-adique 
lisse sur $V$ et 
$x_1,\dots,x_n$  une famille de sections locales du faisceau $\cal O_{\calXdag /V}$ au-dessus d'un ouvert affine $U$ d'algèbre $A\dag
:= \Gamma (U,\cal O_{\calXdag /V})$ telles que leurs différentielles
$dx_1,\dots,dx_n$ forment une  base locale du $\cal O_{\calXdag /V}$-module des formes différentielles séparées ${\Omega}_{\calXdag /V}.$ Alors, 
pour 
$h\geq0,$  un  opérateur différentiel $P$ d'échelon  $h$ s'écrit de {\bf manière unique}: 
\glossary{$P(a,\Delta,\dots,\Delta^{p^h})$}
$$P(a,\Delta,\dots,\Delta^{p^h})= 
\kern-1.2cm\sum_{\vrule height10pt width0pt\beta^0,\dots,\beta^{h-1}\in (\Bbb N^n_{\leq p-1})^h, \beta^h\in\Bbb
N^n}\kern-1.2cm
a_{\beta^0,\dots,\beta^h}(\Delta)^{\beta^0}\cdots(\Delta^{p^h})^{\beta^h}$$ 
où $a_{\beta^0,\dots,\beta^h}$ est une suite d'éléments de l'algèbre $A\dag$ telle que la série \og symbole total\fg:
\glossary
{$\sigma_{P, h}(x,\xi^0,\dots,\xi^h):=P(a,\xi^0,\dots,\xi^h)$ \hbox to4.1cm{}}
$$\sigma_{P, h}(x,\xi^0,\dots,\xi^h):=P(a,\xi^0,\dots,\xi^h):= 
\kern-1.5cm\sum_{\vrule height10pt width0pt\beta^0,\dots,\beta^{h-1}\in (\Bbb N^n_{\leq p-1})^h,\beta^h\in \Bbb
N^n}\kern-1.5cm
a_{\beta^0,\dots,\beta^h}(\xi^0)^{\beta^0}\cdots(\xi^{h})^{\beta^h}$$ est un élément de l'algèbre
$(A\dag[\xi^0,\dots,\xi^h])\dag$,
où  la $A\dag$-algèbre $A\dag[\xi^0,\dots,\xi^h]$ est engendrée par les générateurs $\xi^j=(\xi^j_1,\dots,\xi^j_n)$, $0\leq j\leq h, $  soumis aux relations
$(\xi^j_i)^p= u_{j}\xi_i^{j+1}$ $0\leq j\leq h-1$,  $ 1\leq
i\leq n$ et  $u_j:=\frac{p^{j+1}!}{(p^{j}!)^p}$.  En particulier, on a la majoration $v_{\goth m}(a_{\beta^0,\dots,\beta^h})\geq \lambda(|\beta^0|+\cdots+|\beta^h|p^h),$ $\lambda>0$.

De plus, l'application symbole total qui à $P$ associe son symbole total,
$$P\mapsto \sigma_{P, h}(x,\xi^0,\dots,\xi^h),$$ est un {\bf isomorphisme} de $A\dag$-modules à gauche entre l'anneau $D^{\dagger,h}_{A\dag/V}:= 
\Gamma(U, \cal
D^{\dagger,h}_{\calXdag /V} )$ des sections globales des
opérateurs d'échelon $h$ et la
$A\dag$-algèbre $(A\dag[\xi^0,\dots,\xi^h])\dag$.
\end{theo}

Le théorème \ref{sym-toth} est démontré en deux parties ([Me$_3$], Thm. 3.1.2 et Thm. 3.1.4). La démonstration utilise le théorème \ref{sym-tot}.
Là aussi les éléments de l'algèbre $(A\dag[\xi^0,\dots,\xi^h])\dag$ sont faciles à décrire ce qui produit beaucoup d'opérateurs d'échelon $h$. On déduit du théorème précédent  que le faisceau $\calD^{\dagger,h}_{\calXdag /V}$ est aussi acyclique sur les ouverts affines assez petits ([Me$_3$], 3.2.3), que l'idéal $\goth m$ est contenu dans la radical de $\Gamma(U,\calD^{\dagger,h}_{\calXdag /V})$ pour un ouvert affine assez petit $U$ ([Me$_3$], 4.1.1) et que la transposition  est une anti-involution de cet anneau ([Me$_3$], 4.2.1).

L'intérêt de la $p$-filtration est le théorème:

\begin{theo}\label{noe}Pour tout  ouvert assez petit affine $U$ d'algèbre $A\dag$ et  tout $h\geq0$ l'anneau
$D^{\dagger,h}_{A\dag/V}:= 
\Gamma(U, \cal
D^{\dagger,h}_{\calXdag /V} )$ \glossary{$\cal
D^{\dagger,h}_{\calXdag /V}$}est \em{noethérien}.
\end{theo}

Le théorème \ref{noe} est démontré dans l'article ([Me$_3$], Thm. 5.3.6). 
Ce théorème ne se déduit pas de son analogue formel, il est beaucoup plus profond et
sa démonstration utilise en fait le théorème précédent. 

\begin{Rema}\label{pfi-pla}La filtration par les échelons permet d'appliquer le critère de platitude locale \ref{cri-pla}
pour  montrer la platitude des extensions de faisceaux d'opérateurs différentiels que nous rencontrerons dans cet article. En particulier, les extensions $\cal O_{
\calXdag /V}\rightarrow\calD_{\calXdag /V}\rightarrow\calDdag_{\calXdag /V}$ sont   plates ([Me$_3$], Coro. 6.1.2).
\end{Rema}

Nous allons étudier l'extension de cette filtration
au faisceau
\smash{$ \calD^{\dagger}_{\Xdaginf /R}$} des opérateurs différentiels sur le site $ \Xdaginf $. 

\begin{prop}Soit $R$ une $\Bbb Z_{(p)}$-algèbre. Alors, pour $h\geq0$  et pour tout morphisme du site $r\dag: \calWdag \to \calUdag $, le morphisme
de restriction $(\sharp)_{r\dag}: r^{-1}\calDdag_{\calUdag /R}\to\calDdag_{\calWdag /R}$ induit un morphisme de restriction
$$r^{-1}\cal
D^{\dagger,h}_{\calUdag /R}\to\calD^{\dagger,h}_{\calWdag /R}\leqno (\sharp)_{r\dag,h}:$$ compatible à la composition des morphismes.
\end{prop}

\demo La question est locale. Il s'agit de montrer que si $P$ est un opérateur différentiel d'échelon
$h$ l'opérateur différentiel
$r^*Pr^{*-1}$ est d'échelon h. Si $x_1,\dots,x_n$ sont des coordonnées locales sur $\calUdag |W$ il suffit de montrer, par
définition des opérateurs différentiels d'échelon $h$ ([M-N$_1$], [Me$_3$]), que les opérateurs différentiels:
$$r^*\Delta_{x_j}r^{*-1},\dots, r^*\Delta^{p^h}_{x_j}r^{*-1}$$ sont d'ordre $1,\dots,p^h$ respectivement, ce qui est conséquence  du fait $r^*\circ
P\circ r^{*-1}$ est un isomorphisme d'algèbres  respectant l'ordre des opérateurs différentiels.
\enddemo

\begin{defi} Soit $R$ une $\Bbb Z_{(p)}$-algèbre. Alors, la famille $\cal
D^{\dagger,h}_{\calUdag /R}$, munie des morphismes $(\sharp)_{r\dag,h}$, définit un faisceau d'algèbres $ \cal
D^{\dagger,h}_{ \Xdaginf /R} $ \glossary{$ \cal
D^{\dagger,h}_{ \Xdaginf /R} $}sur le site $ \Xdaginf $, le faisceau des opérateurs différentiels d'échelon $h$. La filtration $\Rm
\calD^{\dagger,h}_{ \Xdaginf /R} $ est une filtration exhaustive du faisceaux $ \cal
D^{\dagger}_{ \Xdaginf /R} $ par des sous-faisceaux d'algèbres.
\end{defi}

Rappelons que dans le cas d'un anneau de valuation discrète l'on définit l'indice de ramification absolu $e$ comme
la valuation
$\goth m$-adique $v_{\goth m}(p)$
de $p$.

\begin{prop}Soit $\calXdag $ un schéma $\dagger$-adique  sur $V$ lisse. Si $h\geq0$ est un entier tel que $e<p^h(p-1)$, alors  le
faisceau 
$\cal G_{\calXdag }$ est un \em{sous-faisceau de groupes pour la structure multiplicative}  du  faisceau $\cal
D^{\dagger,h}_{\calXdag /V}$ des opérateurs différentiels d'échelon $h$, et  le
faisceau 
$\cal G_{\Xdaginf}$ est un \em{sous-faisceau de groupes pour la structure multiplicative}  du  faisceau $\cal
D^{\dagger,h}_{\Xdaginf/V}$ des opérateurs différentiels d'échelon $h$ pour un schéma $X$ lisse sur $k$.
\end{prop}

\demo
Les faisceaux $\cal G_{\calXdag}$ et $\cal
D^{\dagger,h}_{\calXdag /V}$ sont par définition des sous-faisceaux du faisceau $\cal End_{V}(\cal O_{\calXdag /V})$. La question est donc locale.
Soit $(x_1,\dots,x_n)$ une famille de sections locales du faisceau $\cal O_{\calXdag /V}$ au-dessus d'un ouvert affine $U$ 
d'algèbre $A\dag := \Gamma (U,\cal O_{\calXdag /V})$ telles que leurs différentielles
$(dx_1,\dots,dx_n)$ forment une  base locale du $\cal O_{\calXdag /V}$-module des formes différentielles séparées ${\Omega}_{\cal 
X\dag/V}$. Si
$g$  est une section de $\Gamma(U,\cal G_{\calXdag })$, c'est un automorphisme de l'algèbre $ A\dag$ qui se réduit à l'identité modulo $\goth
m$ et, en vertu du  corollaire \ref{dev-grp},
c'est donc un opérateur différentiel  qui admet comme développement en série:
$$g =\sum_{\beta\in \Bbb N^n}(g(x_1)-x_1)^{\beta_1}\cdots(g(x_n)-x_n)^{\beta_n}\Delta_1^{\beta_1}\cdots\Delta_n^{\beta_n}.$$ La valuation $\goth m$-adique du coefficient
$a_\beta:= (g(x_1)-x_1)^{\beta_1}\cdots(g(x_n)-x_n)^{\beta_n}$ est supérieure à $|\beta|:= \beta_1+\cdots+\beta_n$. 
Soit $$\beta_i := q_{\beta,i}p^h+a_{\beta,i,h-1}p^{h-1}+\cdots+a_{\beta,i,0}, \,\,\,0\leq a_{\beta,i,j}\leq p-1\,,$$ la division euclidienne de $\beta_i$
par
$p^h$. Posons 
$$b_\beta := a_\beta\frac{(p^h!)^{|q_\beta|}({p^{h-1}!})^{|a_{\beta,h-1}|}\cdots (p!)^{|a_{\beta,1}|}}{\beta!}.$$ La valuation $\goth m$-adique de $b_\beta$
est supérieure à:
$$\sum_{1\leq i\leq n} \mBig(\beta_i +q_{\beta,i}\frac{(p^h-1)e}{p-1}+\cdots+a_{\beta,i,1}\frac{(p-1)e}{p-1}-\frac{(\beta_i-s_{\beta_i})e}{p-1}\mBig),$$ où $s_{\beta_i}$ est la somme
des coefficients du développement $p$-adique de $\beta_i$. Si $p^h(p-1) -e>0$, la valuation $\goth m$-adique de $b_\beta$ est supérieure à $\lambda|\beta|$ pour une
constante $\lambda:= \frac{(p^h(p-1)-e)}{(p-1)p^h}>0$, et la série 
$$\sum_{\beta}b_\beta (\Delta)^{a_{\beta,0}}\cdots(\Delta^{p^h})^{q_\beta}$$ est un opérateur différentiel d'échelon $h$ en vertu de la définition du faisceau $\cal
D^{\dagger,h}_{\calXdag /V}$.
\enddemo

\subsection{La cohomologie de de Rham $p$-adique  $ H^\bullet_{\DR,h}(\Xdaginf /V)$ d'échelon $h\geq{}$0}
En vertu de la proposition précédente et sous la condition $e<p^h(h-1)$, le faisceau $\cal G_{\Xdaginf}$ est, pour la structure multiplicative,  un sous-faisceau  de groupes du faisceau d'anneaux   $\cal
D^{\dagger,h}_{\Xdaginf/V}$.

\begin{defi}Sous la condition $e<p^h(h-1)$, on dit qu'un $\cal
D^{\dagger,h}_{ \Xdaginf /V}$-module à gauche $\cal M^{\dagger,h}_{\inf}$ est  spécial si l'action $(\sharp)$ de $\cal G_{
\Xdaginf }$ sur $\cal M^{\dagger,h}_{\inf}$  se fait à travers $\cal
D^{\dagger,h}_{ \Xdaginf /V}$.
\end{defi}
 On note $ (\cal
D^{\dagger,h}_{ \Xdaginf /V},  \Sp)\Mod$ la catégorie des $\cal
D^{\dagger,h}_{\Xdaginf/V}$-modules à gauche spéciaux.
Tous les raisonnements que nous avons faits
pour le couple $\cal
G_{ \Xdaginf }\hookrightarrow\cal
D^{\dagger}_{ \Xdaginf /V}$ se transposent mutatis mutandis au couple $\cal
G_{ \Xdaginf }\hookrightarrow\cal
D^{\dagger,h}_{ \Xdaginf /V}$ sous la condition $e<p^h(p-1)$. 
\begin{liste}\item 1)
On trouve que la donnée d'un module à gauche spécial d'échelon $h$ est équivalente 
la donnée d'un famille de $\calD^{\dagger,h}_{\calUdag /V}$-modules à gauche munie d'isomorphismes $(\diamond)_h$ définis à l'aide des modules de
transfert
$\cal
D^{\dagger,h}_{ \calWdag \to\calUdag /V}$.
\item 2)
On trouve que si $\calXdag $ est un relèvement
plat de $X$, la catégorie $ \cal
D^{\dagger,h}_{ \calXdag /V}\Mod$ est canoniquement équivalente à la catégorie des modules spéciaux $ (\cal
D^{\dagger,h}_{ \Xdaginf /V},  \Sp)\Mod$ par le foncteur prolongement canonique $ P_{\calXdag , h}$. 
\item 3)
On trouve que
la catégorie $ (\calD^{\dagger,h}_{ \Xdaginf /V},  \Sp)\Mod$
des modules spéciaux à gauche 
  est une catégorie abélienne qui a suffisamment d'injectifs et est un champ. 
\item 4)
On trouve
que  si $\cal N^{\dagger,h}_{\inf}, \cal M^{\dagger,h}_{\inf}$ sont deux modules à gauche spéciaux d'échelon 
$h$, le complexe $\bfR \cHom_{\calD^{\dagger,h}_{\Xdaginf /V},\Sp} (\cal N^{\dagger,h}_{\inf}, \cal M^{\dagger,h}_{\inf})$
est un complexe de Zariski. On définit 
en particulier la cohomologie de de Rham de $X$
$p$-adique d'échelon $h$  à coefficients dans $\cal M^{\dagger,h}_{\inf}$ par :
\glossary{$ H^\bullet_{\DR,h}(\Xdaginf /V, \cal M^{\dagger,h}_{\inf})$}\glossary{$\ext^{\bullet}_{\cal
D^{\dagger,h}_{\Xdaginf /V},\Sp}(\cal O_{\Xdaginf /V},\cal M^{\dagger,h}_{\inf})$}
$$ H^\bullet_{\DR,h}(\Xdaginf /V, \cal M^{\dagger,h}_{\inf}):= \ext^{\bullet}_{\calD^{\dagger,h}_{\Xdaginf /V},\Sp}(\cal O_{\Xdaginf /V},\cal
M^{\dagger,h}_{\inf}).$$ 
\item 5)
On trouve que si $\calXdag $ est un relèvement  plat, on a des
isomorphismes canoniques:
$$ \ext^{\bullet}_{\cal
D^{\dagger,h}_{ \calXdag /V}}(\cal O_{\calXdag /V},R_{\calXdag }(\cal
M^{\dagger,h}_{\inf}))\simeq \ext^{\bullet}_{\calD^{\dagger,h}_{\Xdaginf /V},\Sp}(\cal
O_{\Xdaginf /V},\cal
M^{\dagger,h}_{\inf}).$$ 
\end{liste}Le cas
$h=0$  est particulièrement intéressant et donne:

\begin{prop} Si  $\calXdag $ est une relèvement $\dagger$-adique 
lisse de $X$ et si $e<p-1$, la cohomologie de de Rham  usuelle  $ H^\bullet_{\DR}(\calXdag /V):=  H^\bullet(X, \Omega^{\bullet}_{{\calXdag }/V})$ 
est canoniquement isomorphe à la cohomologie d'échelon zéro $ H^\bullet_{\DR,0}(\Xdaginf /V, \cal O_{\Xdaginf /V})$ et est 
indépendante du relèvement.
\end{prop}

En particulier, si $X$ est affine et si $e<p-1$, la cohomologie de de Rham des formes différentielles séparées d'un relèvement, qui est canoniquement isomorphe à
la cohomologie de de Rham d'échelon  $0$ du site $\Xdaginf $, est indépendante du relèvement comme $V$-module gradué [M-W].

On trouve ou on retrouve comme corollaire l'invariance  de la cohomologie de de Rham {\bf algébrique} 
d'un morphisme propre et lisse sur
$V$ sans invoquer la cohomologie cristalline, qui  a fourni  historiquement la première démonstration de ce résultat [B$_1$].

\begin{coro}\label{alg-ada} Soit ${\mathbf{X}}/V$ un relèvement \em {algébrique} propre et lisse de $X/k$. Si $e<p-1$, la cohomologie de de Rham algébrique usuelle $\Rm
H^{\bullet}_{\DR}({\mathbf{X}}/V):= H^{\bullet}({\mathbf{X}}, \Omega^{\bullet}_{{\mathbf{X}}/V})$ est canoniquement isomorphe à la 
cohomologie de de Rham $p$-adique d'échelon zéro $ H^\bullet_{\DR,0}(\Xdaginf /V, \cal O_{\Xdaginf /V})$ et 
est indépendante à isomorphisme canonique près du relèvement ${\mathbf
X}/V$.
\end{coro}

\demo En effet, en vertu du théorème de comparaison $GAGA\dag$ de Meredith [Mr] pour les $V$-schémas propres, le morphisme canonique:
$$ H^{\bullet}_{\DR}({\mathbf{X}}/V)\to   H^\bullet_{\DR}(\calXdag /V)$$ est un isomorphisme de $V$-modules gradués. On est alors réduit à la proposition précédente.
\enddemo

\begin{coro} Soit ${\mathbf{X}}/V$ un relèvement \em {algébrique} propre et lisse de $X/k$. Si $e<p-1$, la cohomologie de de Rham algébrique usuelle $ H^{\bullet}_{\DR}({\mathbf{X}}/V):= H^{\bullet}({\mathbf{X}}, \Omega^{\bullet}_{{\mathbf{X}}/V})$ est canoniquement isomorphe à la cohomologie 
de de Rham $p$-adique d'échelon zéro  du site infinitésimal formel
$ H^\bullet_{\DR,0}( X_{\inf}^\wedge/V)$. 
\end{coro}

\demo En vertu du théorème de comparaison $GAGF$ de Grothendieck  [EGA III$_1$] pour les $V$-schémas propres, la cohomologie de de Rham du schéma $\bf X$, \em{i.e.} $ H^{\bullet}_{\DR}({\mathbf{X}}/V):= H^{\bullet}({\mathbf{X}}, \Omega^{\bullet}_{{\mathbf{X}}/V})$, est canoniquement isomorphe à la cohomologie de
de Rham formelle $ H^\bullet_{\DR}(
\cal X^\wedge/V):= H^{\bullet}({ X}, \Omega^{\bullet}_{{\cal  X^\wedge}/V})$ ({\it cf. }\ref{site-infinitesimal-formel})
et donc canoniquement à celle du site infinitésimal formel d'échelon nul
$ H^\bullet_{\DR,0}( X_{\inf}^\wedge/V)$.
\enddemo

Ainsi  toutes les cohomologies de type de de Rham qu'on peut définir pour une variété propre et lisse sur $k$ qui se {\bf relève} sur $V$ en schéma propre et lisse, 
sont toutes {\bf canoniquement} isomorphes.

\begin{Rema}
Un peu plus généralement et  comme nous l'avons signalé dans l'introduction (\ref{site-infinitesimal-formel}), la cohomologie de de Rham $p$-adique
formelle  $ H^\bullet_{\DR}(
X^\wedge_{\inf}/K)$  d'une variété algébrique $X$ lisse sur $k,$ {\bf non nécessairement propre}, doit être canoniquement isomorphe à la cohomologie cristalline de
$X$ sur $K$ et, si $e<p-1$, la cohomologie $p$-adique $ H^\bullet_{\DR,0}(
X^\wedge_{\inf}/V)$ d'échelon $0$  doit être canoniquement isomorphe à la cohomologie cristalline de
$X$ sur $V$. C'est le cas des variétés qui se relèvent en vertu du théorème de comparaison de Berthelot-Ogus [B-O]. Il faut donc  construire un morphisme de comparaison et passer  du local au global. 
\end{Rema}

\section{La cohomologie infinitésimale et la cohomologie de de Rham $p$-adiques}
Soit $X/k$ une variété algébrique lisse sur un corps $k$. Lorsque la caractéristique du corps de base $k$ est nulle, Grothendieck  a montré ([G$_3$]) que la cohomologie
$ H^\bullet_{\inf}(X/k,\cal O_{\Xdaginf /k})$ du site infinitésimal de $X$ à valeurs dans le faisceau structural est canoniquement isomorphe à sa
cohomologie de de Rham $ H^\bullet_{\DR}(X/k):= {H}^\bullet(X,\Omega^\bullet_{X/k})$, et fournit donc les 
bons nombres de Betti en vertu du théorème de comparaison de Grothendieck ([G$_1$]) entre la cohomologie de de Rham algébrique et la cohomologie transcendante, quand $k$ est plongé dans le corps des nombres complexes $\Bbb C$. Nous allons voir que {\bf ce n'est pas} le cas de la cohomologie du site infinitésimal
$p$-adique à valeurs dans le faisceau structural.

\subsection{La cohomologie infinitésimale $\dagger$-adique des complexes de la catégorie $ \Dplus (R_{\Xdaginf })$}
Soient $X$  un schéma lisse sur $R_1$ et $ \cal A_{\Xdaginf }$ un faisceau d'anneaux sur le site $ \Xdaginf $.

\begin{prop}La catégorie $ \cal A_{\Xdaginf }\Mod$ des $ \cal A_{\Xdaginf }$-modules à gauche sur le site $ \Xdaginf $ admet suffisamment
d'objets injectifs.
\end{prop}

\demo C'est un cas particulier de l'existence de suffisamment d'injectifs dans la catégorie de modules sur un faisceau d'anneaux sur un site. Mais ici
on peut explicitement construire des résolutions injectives à l'aide des foncteurs $P_{\cal A_{\Xdaginf }, \calUdag }$ exactement comme on a construit des
résolutions injectives dans la catégorie des modules spéciaux $ (\calDdag_{\Xdaginf /R}, \Sp)\Mod$ dans le théorème \ref{obj-inj}.
\enddemo

\bigskip
Le foncteur covariant $ \cal F_{\inf}\fonct \hom_{R_{\Xdaginf }}(R_{\Xdaginf }, \cal F_{\inf})$  de la catégorie $\RmMod(R_{\Xdaginf })$ dans
la catégorie
$\RmMod(R)$ est exact à gauche et se dérive à droite $ \cal F_{\inf}\fonct \bfR \hom_{R_{\Xdaginf }}(R_{\Xdaginf }, \cal F_{\inf})$ comme foncteur exact
de catégorie triangulées, de la catégorie $ \Dplus (R_{\Xdaginf })$  vers la catégorie $ \Dplus (R)$. 

\begin{defi}Si $\cal F_{\inf}$ est un complexe appartenant à la catégorie $ \Dplus (R_{\Xdaginf })$, on définit la cohomologie infinitésimale $\dagger$-adique: 
$$ H^{\bullet}_{\inf}(X/R, \cal F_{\inf}),$$\glossary{$ H^{\bullet}_{\inf}(X/R, \cal F_{\inf})$}de $X$ à valeurs dans $\cal F_{\inf}$, comme la
cohomologie du complexe: $$\Inf(X/R,\cal F_{\inf}):= \bfR \hom_{R_{\Xdaginf }}(R_{\Xdaginf },
\cal F_{\inf}).$$
\end{defi}

\begin{theo}\label{coh-inf}Soit $\calXdag $ un relèvement  plat de $X$.  
Si $\cal F_{\inf}$ est  un complexe de la catégorie $ \Dplus (R_{\Xdaginf })$, alors il existe
un  isomorphisme canonique :
$$\bfR\hom_{R[\cal G_{\calXdag }]}(R_{X},
R_{\calXdag }(\cal F_{\inf}))\simeq  \bfR\hom_{R_{\Xdaginf }}(R_{\Xdaginf },
\cal F_{\inf}),$$ et donc des isomorphismes canoniques:
$$ \ext^{\bullet}_{R[\cal G_{\calXdag }]}(R_{X}, R_{\calXdag }(\cal F_{\inf}))\simeq  H^{\bullet}_{\inf}(X/R, \cal F_{\inf}).$$
\end{theo}

\demo C'est une conséquence directe de l'équivalence du théorème \ref{equ-god}.
\enddemo

\begin{Rema}
On peut de façon évidente considérer la théorie sur le corps des fractions $K$ d'un anneau $V$ de valuation discrète complet, où se pose la question
redoutable  du calcul des espaces: 
$$ H^i_{\inf}(X/K, K_{\Xdaginf }),\ H^i_{\inf}(X/K, \cal O_{\Xdaginf /K}),
$$$$H^i_{\inf}(X/K,
\Omega^j_{\Xdaginf /K}),\  H^i_{\inf}(X/K,
\Omega^{\bullet}_{\Xdaginf /K}),
$$ même pour la droite affine ou la droite projective. Nous allons montrer,  dans le paragraphe qui suit, des résultats très partiels dans cette direction. \end{Rema}

\subsection{La cohomologie infinitésimale en degré zéro  à valeurs dans le faisceau structural fournit le bon nombre de Betti}
Notre première approche a été de montrer, de même qu'en caractéristique nulle, qu'en caractéristique  positive  la cohomologie $ H^\bullet_{\inf}(X/K,\cal
O_{\Xdaginf /K})$ du site infinitésimale $p$-adique de $X$, supposée affine pour commencer, à valeurs dans le faisceau structural est
canoniquement isomorphe à sa cohomologie de de Rham
$p$-adique, ce qui éventuellement établit l'indépendance du relèvement de la cohomologie de de Rham. Pendant longtemps, plusieurs années, nous avons cru que tel était  bien le
cas, mais nos efforts nous ont conduits à  montrer qu'il n'en était rien, de façon non triviale,  ce qui soulève des problèmes fort intéressants en liaison avec la
cohomologie  de groupes et  l'algèbre (non)-commutative.

Si $\cal F_{\inf}$ est un $ K_{\Xdaginf }$-module, son $K$-espace des sections globales infinitésimales est par définition le $K$-espace  vectoriel
$ \hom_{K_{\Xdaginf }}(K_{\Xdaginf },\cal F_{\inf})$, et   sa cohomologie infinitésimale est donc :
$$ H^{\bullet}_{\inf}(X/K,\cal F_{\inf}):= \ext^\bullet_{K_{\Xdaginf }}(K_{\Xdaginf },\cal F_{\inf}).$$ 
On peut faire le changement de base: $$ K_{\Xdaginf }\rightarrow 
\cal O_{\Xdaginf /K},$$ qui est plat  et trouver,
pour un $ \cal O_{X\dag_{\inf/K}}$-module $\cal M_{\inf},$ les  isomorphismes canoniques:
$$ H^{\bullet}_{\inf}(X/K, \cal M_{\inf})\simeq \ext^\bullet_{\cal O_{\Xdaginf /K}}(\cal O_{X\dag_{\inf/K}},\cal
M_{\inf}).$$ 

\begin{prop} Si $\calMdaginf $ est un $\calDdag_{\Xdaginf /K}$-module spécial, il existe un isomorphisme canonique:
$$ H^{0}_{\DR}(X/K,\calMdaginf )\simeq H^{0}_{\inf}(X/K, \calMdaginf ).$$
\end{prop}

\demo Comme $\cal O_{\Xdaginf /K}$ est un sous-faisceau d'anneaux du faisceau $\calDdag_{\Xdaginf /K}$, on a un morphisme canonique:
$$ \hom_{\calDdag_{\Xdaginf /K}}(\cal O_{\Xdaginf /K},\calMdaginf )\to\hom_{\cal O_{\Xdaginf /K}}(\cal
O_{\Xdaginf /K},\calMdaginf ).$$ 
En vertu de \ref{rec-fai},  le fait que ce soit un isomorphisme est une question de nature locale pour la topologie de Zariski.  Aussi,  on peut supposer que $X$ est affine. 
Soit $\calXdag $ un relèvement plat. Il suffit de
montrer que le  morphisme canonique:
$$\cHom_{\calDdag_{\calXdag /K}}(\cal O_{\calXdag /K},\cal M\dag_{\calXdag })\to\cHom_{\cal O_{\calXdag /K}[\cal G_{\calXdag}]}(\cal O_{\calXdag /K},\cal M\dag_{\calXdag })$$ est un isomorphisme.
En vertu du théorème
\ref{inc-grp}, on a morphisme canonique de
$\cal O_{\calXdag /K}$-algèbres:
$$\cal A_{\calXdag /K}:= \cal O_{\calXdag /K}[\cal G_{\calXdag }]\to \calDdag_{\calXdag /K}.$$   Il suffit de montrer  que le morphisme
canonique: 
$$\calDdag_{\calXdag /K}\otimes_{\cal A_{\calXdag /K}}\cal O_{\calXdag /K}\rightarrow \cal O_{\calXdag /K}$$ est un isomorphisme de 
$\calDdag_{\calXdag /K}$-modules à gauche. La question est encore locale, et on peut supposer que $x_1,\dots,x_n$ sont des sections locales  du faisceau
structural dont les  différentielles forment une  base de l'espace des formes différentielles séparées. Il suffit enfin de montrer en vertu de
([Me$_3$], Prop. 6.2.3) que l'idéal
d'augmentation
$\cal A_{\calXdag /K}\rightarrow\cal O_{\calXdag /K}$ engendre sur
$\calDdag_{\calXdag /K}$ l'idéal $\Delta_1,\dots,\Delta_n$. Mais en vertu du corollaire  \ref{dev-grp}, les générateurs de l'idéal 
d'augmentation sont de la forme:
$$g -1 =\sum_{\beta\in \Bbb N^n,\beta\neq0}a_1^{\beta_1}\cdots a_n^{\beta_n}\Delta_1^{\beta_1}\cdots\Delta_n^{\beta_n},$$ pour $a_1,\dots,a_n$ dans l'idéal engendré par
$\goth m$. L'idéal sur $\calDdag_{\calXdag /K}$ engendré par l'idéal d'augmentation contient des éléments de la forme $(1+P_i)\Delta_i$,
qu'on obtient, par exemple, pour 
$g =(\sum_{0\leq \beta_i\leq \infty}((p^2)^{\beta_i})\Delta_i^{\beta_i})$, où $P_i$
appartient à l'idéal engendré par $\goth m$. En vertu du théorème  ([Me$_3$], Thm. 4.2.1),  l'opérateur $1+P_i$ est inversible et donc $\Delta_i$
appartient à l'idéal  engendré par l'idéal d'augmentation pour $i=1,\dots,n$. En fait, on obtient
des opérateurs  $1+P_i$ comme dans l'exemple précédent qui sont de façon évidente inversibles.
\enddemo

La cohomologie infinitésimale $p$-adique d'une variété algébrique lisse sur $k$  fournit le bon nombre de Betti en degré zéro.

\subsection{La cohomologie infinitésimale en degré un à valeurs dans le faisceau structural \em{ne fournit pas} le bon nombre de Betti}

\begin{theo}\label{lem-poi}Soit $\calXdag $ une relèvement  plat d'un schéma affine lisse $X$ muni de sections $x_1,\dots,x_n$   du faisceau
structural telles que leurs différentielles forment  une base de l'espace des formes différentielles séparées. Alors, le $K$-espace $ H^1_{\inf}(X/K, \cal
O_{ \Xdaginf /K})$ est
{\bf non nul}.
\end{theo}

\demo En vertu de l'équivalence de catégories \ref{equ-god}, on a un isomorphisme canonique:
$$ \ext^1_{\cal A_{\calXdag /K}}(\cal O_{\calXdag /K}, \cal O_{\calXdag /K})\simeq H^1_{\inf}(X/K,\cal O_{
\Xdaginf /K}).$$ Soit la suite exacte:
$$0\rightarrow \cal I_{\calXdag }\rightarrow\cal
A_{\calXdag /K}\rightarrow\cal O_{\calXdag /K}\rightarrow0\,,$$ où $\cal I_{\calXdag }$ est l'idéal d'augmentation. 

En appliquant le foncteur
$\hom_{\cal A_{\calXdag /K}}(?,\cal O_{\calXdag /K})$, on trouve une suite exacte:
$$0\rightarrow\hom\Sub{6pt}{\cal A_{\calXdag /K}}(\cal O_{\calXdag /K},\cal O_{\calXdag /K})\rightarrow\Gamma(X, \cal O_{\calXdag/K})\rightarrow\hom\Sub{6pt}{\cal A_{\calXdag /K}}(\cal I_{\calXdag },\cal O_{\calXdag /K})
\newdisplayline{0pt}{0pt}{-1ex}
\rightarrow\ext^1\Sub{6pt}{\cal
A_{\calXdag /K}}(\cal O_{\calXdag /K},\cal O_{\calXdag /K}).$$
Il suffit de montrer que le morphisme  induit:
$$\Gamma(X,\cal O_{\calXdag/K})\rightarrow\hom_{\cal A\dag_{\calXdag /K}}(\cal I_{\calXdag },\cal O_{\calXdag /K})$$ \em{n'est pas surjectif}. Un élément
de l'espace figurant à  droite est un morphisme de faisceaux, et est  défini par ses restrictions  sur les ouverts affines qui forment une base de la topologie. Soit
$U$ un ouvert affine d'algèbre $\dagger$-adique 
$A\dag_U:= \Gamma(U,\cal O_{\calXdag /V})$. L'idéal $\Gamma(U,\cal I_{\calXdag })$ est engendré par les éléments $g-1$, où $g$ est un automorphisme de
l'algèbre
$A\dag_U$ qui se réduit à l'identité modulo
$\goth m$. En vertu du corollaire  \ref{dev-grp}, il s'agit   d'opérateurs différentiels qui admettent des développements de la forme:
$$g =\sum_{\beta\in \Bbb N^n}a_1^{\beta_1}\cdots a_n^{\beta_n}\Delta_1^{\beta_1}\cdots\Delta_n^{\beta_n},$$ pour $a_1,\dots,a_n$ dans l'idéal engendré par
$\goth m$.

\begin{lemm} L'opérateur différentiel 
$$\trans g:= \sum\nolimits_{\beta\in \Bbb N^n}(-1)^{|\beta|}\Delta_1^{\beta_1}\cdots\Delta_n^{\beta_n}a_1^{\beta_1}\cdots
a_n^{\beta_n}$$  est égal à $\trans g(1)g^{-1}$.
\end{lemm}

\demo
La transposition $P\mapsto \trans P$ est une anti-involution de l'algèbre des opérateurs différentiels ([Me$_3$], Coro. 4.2.2). Si $g$ est un
automorphisme de
$\Gamma(U,\cal G_{\calXdag })$ et si $f$ une fonction de $A\dag_U$, on a $gf= g(f)g$ et $f\trans g=\trans gg(f)$ soit $f\trans g(1):=\trans g(g(f))$ ou encore 
$\trans g(1)g^{-1}(f)=\trans g(f)$.
\endsubdemo

\begin{lemm}\label{coc-nnu} Si $g$ est un élément du groupe $\Gamma(U, \cal G_{\calXdag })$, l'élément: $$\log(\trans g(1)):= \sum_{m=0,\infty}(-1)^m
\frac{(\trans g(1)-1)^{m+1}}{m+1}$$ est un élément bien défini de l'algèbre
$A\dag_{U}\otimes_VK$.
\end{lemm}

\demo
En effet, si l'on choisit une présentation de l'algèbre $A\dag_U$ comme quotient du complété $\dagger$-adique d'une algèbre de polynômes,
l'élément  $\trans g(1)-1$ appartient à l'algèbre de Banach $A\dag_{U,\rho}$, image  de l'algèbre des séries qui convergent dans un domaine $|x|\leq \rho$ pour un nombre
$\rho>1$, et sa norme quotient vérifie $||\trans g(1)-1||_\rho<1$, ce qui entraîne que la série converge dans l'algèbre $A\dag_{U,\rho}$.
\endsubdemo

\begin{lemm}\label{coc}L'application $g\mapsto \delta(g):=\log(\trans g^{-1}(1))$ \glossary{$g\mapsto \delta(g):=\log(\trans g^{-1}(1))$}est un cocycle de 
$\Gamma(U,\cal
G_{\calXdag })$ à valeurs dans
$A\dag_U\otimes_VK$ qui n'est pas un cobord.
\end{lemm}

\demo
Il faut vérifier la propriété de cocycle: $$\delta(g_1g_2)= \delta(g_1)+g_1\delta(g_2).$$ 
Par définition,  on a l'égalité: $$\delta(g_1g_2)= \log(\trans (g_1g_2)^{-1}(1))=\log\trans g_1^{-1}(1)+\log(g_1(\trans g^{-1}_2(1))),$$  car $\trans g_1^{-1}(1)g_1=\trans g_1^{-1}$. Mais comme $g_1$ est un morphisme 
de l'algèbre
$A\dag_{U,\rho}$ vers l'algèbre $A\dag_{U, \rho '}$ pour
$\rho>1, \rho '>1$ assez près de $1$, on trouve que $\log(g_1(\trans g^{-1}_2(1)))= g_1(\log(\trans g^{-1}_2(1)))$, ce qui montre bien la propriété de cocycle.
\endsubdemo

\medskip
Le cocycle ainsi défini est nul sur les éléments du groupe à coefficients constants,  c'est-à-dire les éléments $g$ tels que
$g= \sum_{\beta} a^{\beta}\Delta_i^{\beta_i}$  où $a$ est un élément de $\goth m\subset V$.  Cela implique que ce cocycle n'est pas un cobord. En effet, si ce
cocycle était un cobord, il existerait un élément $f$ 
tel que $\delta(g)= (g-1)f$, d'où $\Delta_i(f)=0, i=1,\dots,n$
et donc $\delta(g)=0$, ce qui n'est pas le cas.
D'autre part, l'application
$\delta$ commute aux restrictions et définit un morphisme de
$\hom_{\cal A_{\calXdag /K}}(\cal I_{\calXdag },\cal O_{\calXdag /K})$ dont l'image dans $\ext^1_{\cal
A_{\calXdag /K}}(\cal O_{\calXdag /K},\cal O_{\calXdag /K})$ n'est pas nulle.
\enddemo

\begin{Rema}Sous la condition précédente, la cohomologie de groupes $H^1(G_{A\dag}, A\dag_U\otimes_VK)$ de degré $1$ est non nulle.
\end{Rema}

\begin{coro} Le lemme de Poincaré pour la cohomologie infinitésimale n'a pas lieu: le $ H^{1}_{\inf}(A^n_k/K, \cal O_{({A^n_k})\daginf /K})$ de l'espace affine 
$\Spec(k[x_1,\dots,x_n])$ n'est
pas nul.
\end{coro}

\demo
En effet, le faisceau des formes différentielles séparées est libre sur l'espace affine. 
\enddemo

En particulier, la comparaison entre la cohomologie infinitésimale et la
cohomologie de de Rham n'a pas lieu.

\begin{coro}Sous les condition du théorème \ref{lem-poi}, si l'espace $ H^1_{\DR}(X/K)$ est  nul,  l'extension $\cal A_{\calXdag /K}\rightarrow
\calDdag_{\calXdag /K}$
\em{\bf n'est pas plate}.
\end{coro}

\demo En effet, si l'extension $\cal A_{\calXdag /K}\rightarrow \calDdag_{\calXdag /K}$  était  plate,
il existerait  un isomorphisme canonique: 
$$\bfR \cHom\Sub{4pt}{\cal A_{\calXdag /K}}(\cal O_{\calXdag /K},\cal O_{\calXdag /K})\simeq\bfR \cHom\Sub{4pt}{\calDdag_{\calXdag /K}}(\calDdag_{\calXdag /K}\otimes_{\cal A\dag_{\calXdag /K}}\cal O_{\calXdag /K},\cal O_{\calXdag /K})$$ et la cohomologie
infinitésimale serait isomorphe à la cohomologie de de Rham $p$-adique. 
\enddemo

\medskip 

On rappelle que l'on a l'isomorphisme  de changement de base:
$$\bfR \cHom\Sub{7pt}{\calDdag_{\calXdag /K}}(\calDdag_{\calXdag /K}\LOtimes_{\cal A\dag_{\calXdag /K}}\cal O_{\calXdag /K},\cal O_{\calXdag /K})\simeq\bfR\cHom\Sub{7pt}{\cal A_{\calXdag /K}}(\cal O_{\calXdag /K},\cal O_{\calXdag /K})\,,$$ qui a lieu dans le contexte topologique très général.

En fait pendant longtemps nous avons cru et essayé   de montrer que cette extension est plate.

\begin{coro}Sous les conditions du théorème \ref{lem-poi},  si l'espace $ H^1_{\DR}(X/K)$ est  nul, le faisceau: $$\cTor_1^{\cal A_{\calXdag /K}}(\calDdag_{\calXdag /K},\cal O_{\calXdag/K})$$ n'est pas nul.
\end{coro}

\demo 
En effet, en vertu de l'isomorphisme précédent, la cohomologie  infinitésimale $ H^{1}_{\inf}(X/K,\cal O_{\Xdaginf /K})$ est isomorphe à l'espace:
$$\ext^1_{\calDdag_{\calXdag /K}}(\calDdag_{\calXdag /K}\Lotimes_{\cal A\dag_{\calXdag /K}}\cal O_{\calXdag/K},\cal O_{\calXdag /K}).$$

Mais comme le morphisme $\calDdag_{X\dag/K}\otimes_{\cal A_{\calXdag /K}}\cal O_{\calXdag/K}\rightarrow \cal O_{\calXdag /K}$ est un isomorphisme et que l'espace $ H_{\DR}^{1}(X/K)$  est nul,
en considérant la suite longue de cohomologie
du triangle construit sur le morphisme de comparaison entre la cohomologie de de Rham $p$-adique et la cohomologie infinitésimale:
$$\DeuxLignes
\DR(X/K):= \bfR\hom_{\calDdag_{\calXdag /K}}(\cal O_{\calXdag /K},\cal O_{\calXdag /K})\to\\\to \Inf(X/K)\simeq\bfR \hom\Sub{0pt}{\calDdag_{\calXdag /K}}(\calDdag_{\calXdag /K}\LOtimes_{\cal A\dag_{\calXdag /K}}\cal O_{\calXdag /K},\cal O_{\calXdag /K})\,,
\endlignes
$$ 
on trouve que si $\cTor_1^{\cal A_{\calXdag /K}}(\calDdag_{\calXdag /K},\cal O_{\calXdag/K})$ est nul, l'espace: 
$$\ext^1_{\calDdag_{\calXdag /K}}(\calDdag_{\calXdag /K}\Lotimes_{\cal A\dag_{\calXdag /K}}\cal O_{\calXdag/K},\cal O_{\calXdag /K})$$ est nul, et le corollaire est conséquence du théorème précédent.
\enddemo

En fait, on trouve 
une injection:
$$\deuxlignes{0cm}{0cm}
 0\to \ext^1_{\calDdag_{\calXdag /K}}(\calDdag_{\calXdag /K}\Lotimes_{\cal A\dag_{\calXdag /K}}\cal O_{\calXdag/K},\cal O_{\calXdag /K})\to
\\
\to\hom_{\calDdag_{\calXdag /K}}(\cTor_1^{\cal A_{\calXdag /K}}(\calDdag_{\calXdag /K},\cal O_{\calXdag/K}),\cal O_{\calXdag /K}).
\endlignes$$ 

\begin{Remas} \begin{liste}\item 1) La non-platitude précédente est due à la non-noethérianité de l'algèbre du groupe. 
De façon plus précise, supposons que $e< p^h(p-1)$. 
On a alors  un morphisme
$\cal A_{\calXdag /V}\to
\calD^{\dagger,h}_{\calXdag /V}$ d'image $Im(\cal A_{\calXdag /V})$,  et  on peut montrer que ce morphisme est injectif. Définissons le
faisceau  d'algèbres
$$\cal B_{\calXdag /V}:=
(\Im(\cal A_{\calXdag })\otimes_VK)\cap
\calD^{\dagger,h}_{\calXdag /V}\subset\cal
D^{\dagger,h}_{\calXdag /K}.$$ 
Par construction, $\Im(\cal A_{\calXdag /V})\otimes_VK\simeq \cal B_{\calXdag /V}\otimes_VK$. On peut montrer que
l'inclusion
$\cal B_{\calXdag /V}\subset \cal
D^{\dagger,h}_{\calXdag /V}$ induit un isomorphisme modulo $\goth m^s$ pour tout $s\geq 1$. Cela entraîne que l'algèbre $\cal B_{\calXdag /V}$ 
\em{n'est pas noethérienne}. Parce que si elle l'était,
l'extension précédente serait plate en vertu du critère de platitude local
\ref{cri-pla}, et la cohomologie infinitésimale fournirait les bons nombres de Betti.
\item 2) Il n'est par contre pas impossible que les  idéaux de type fini
au-dessus d'un ouvert affine assez petit de l'algèbre $\cal B_{\calXdag /V}$ soient de présentation finie. Dans ce cas les faisceaux
$$\cTor_i^{\cal A_{\calXdag /K}}(\calDdag_{\calXdag /K},\cal O_{\calXdag/K})$$pour $i\geq2$, seraient nuls, et l'obstruction au théorème de comparaison entre  la cohomologie infinitésimale et la cohomologie de de Rham
$p$-adique se trouverait dans les foncteurs dérivés à droite du foncteur 
covariant exact à gauche des morphismes 
du $\calDdag_{\calXdag /K}$-module $\cTor_1^{\cal A_{\calXdag /K}}(\calDdag_{\calXdag /K},\cal O_{\calXdag/K})$ à valeurs dans le faisceau structural. En effet, il est difficile de construire des relations selon l'expérience que l'on a de cette
algèbre.
\item 3) Le faisceau $\cTor_1^{\cal A_{\calXdag /K}}(\calDdag_{\calXdag /K},\cal O_{\calXdag/K})$ est non nul en général. Aussi, il serait intéressant de construire un élément non trivial dans le cas de la droite affine. Nous allons proposer un candidat. Supposons que $X$ est la
droite affine sur un corps fini, soit  $X= \Spec(\Bbb F_p[x])$, et soit  $(\Bbb Z_p[x])\dag$ son relèvement canonique. Notons $\theta_p$ l'automorphisme $x\mapsto x+p$, et
$\theta_{px^2}$ l'automorphisme
$x\mapsto x+px^2$ de l'algèbre $(\Bbb Z_p[x])\dag$. On peut montrer de façon non triviale que le groupe $G(\theta_p,\theta_{px^2})$ engendré par les éléments $\theta_{p}, \theta_{px^2}$ est 
\em{libre}, de sorte que les éléments
$\theta_{p}-1$ et $\theta_{px^2}-1$ n'admettent aucune relation non triviale dans l'algèbre $(\Bbb Z_p[x])\dag[G(\theta_p,\theta_{px^2})]$. Pourtant,  l'opérateur 
$\theta_{px^2}-1$ appartient à l'idéal engendré par l'opérateur 
 $\theta_{p}-1$
 dans l'anneau $D\dag_{(\Bbb Z_p[x])\dag/\Bbb Q_p}$, ce qui fournit une relation non triviale dans cet anneau. Il est possible que
$\theta_{p}-1$ et $
\theta_{px^2}-1$ n'admettent aucune relation non triviale dans l'algèbre $(\Bbb Z_p[x])\dag[G_{(\Bbb Z_p[x])\dag}]$, et la relation précédente fournit un élément non trivial du faisceau $\cTor_1^{\cal A_{\calXdag /K}}(\calDdag_{\calXdag /K},\cal O_{\calXdag/K})$, pour $\calXdag $ le schéma  $\dagger$-adique associé à $(\Bbb Z_p[x])\dag$.

\item 4) Le lecteur notera que la découverte du
cocycle \ref{coc}
 a été une chance inespérée pour le fondement de la théorie de la cohomologie de de Rham $p$-adique. Il a permis 
le passage  de la cohomologie 
infinitésimale $p$-adique à la cohomologie de de Rham $p$-adique une fois pour toutes.
\glossary{$\cal A_{\calXdag /V}$}\glossary{$\cal B_{\calXdag /V}$}
\end{liste}
\end{Remas}

\subsection{Le foncteur  de de Rham infinitésimal local}
Si $\calMdaginf $ et $\calNdaginf $ sont deux $ \calDdag_{\Xdaginf /R}$-modules à gauche, le faisceau: 
$$\cHom_{\calDdag_{\Xdaginf /R}}(
\calNdaginf ,\calMdaginf )$$ est un faisceau infinitésimal de $ R_{\Xdaginf }$-modules par la construction \ref{fai-hom}, mais dont l'action de $G_{\Xdaginf}$ est triviale si $\calNdaginf ,\calMdaginf $ sont {\bf spéciaux}. Nous allons voir que ce phénomène survit en cohomologie. Cela permet de définir le complexe de
de Rham infinitésimal local.

\begin{defi}Soit $\calMdaginf $ un complexe appartenant à la catégorie $ \Dplus (\calDdag_{\Xdaginf /R}\Mod)$,
on définit son complexe de de Rham {\bf infinitésimal local}:
$$ \dR_{\inf}(\calMdaginf ):= \bfR \cHom_{\calDdag_{\Xdaginf /R}}(
\cal O_{\Xdaginf /R},\calMdaginf )$$ comme foncteur dérivé du foncteur covariant exact à 
gauche: $$\calMdaginf \fonct \cHom_{\calDdag_{\Xdaginf /R}}(
\cal O_{\Xdaginf /R},\calMdaginf )\,.$$
On obtient un foncteur
covariant exact de catégories triangulées\glossary{$\dR, \dR_{\inf}(X/R,-)$}
$$  \Dplus (\calDdag_{\Xdaginf /R}\Mod)\to  \Dplus (R_{\Xdaginf }).\leqno  \dR_{\inf}(X/R,-) :$$
\end{defi}

\begin{theo}Soit $\calMdaginf $ un complexe de $\calDdag_{\Xdaginf /R}$-modules à gauche {\bf spéciaux} borné inférieurement,  alors il existe un isomorphisme canonique:
$$ \bfR \cHom_{\calDdag_{\Xdaginf /R},\Sp}(
\cal O_{\Xdaginf /R},\calMdaginf )\simeq\bfR \cHom_{\calDdag_{\Xdaginf /R}}(
\cal O_{\Xdaginf /R},\calMdaginf )$$ dans la catégorie $ \Dplus (R_{\Xdaginf })$.
\end{theo}

\demo Soit $\calIdaginf  $ une résolution de $\calMdaginf $ par des $\calDdag_{\Xdaginf /R}$-mo\-dules à gauche spéciaux et injectifs et $\cal J\daginf $ une résolution de $\calIdaginf  $ par des $\calDdag_{\Xdaginf /R}$-modules à gauche  injectifs. Alors, 
le morphisme naturel:
$$ \cHom_{\calDdag_{\Xdaginf /R}}(
\cal O_{\Xdaginf /R},\calIdaginf  )\to \cHom_{\calDdag_{\Xdaginf /R}}(
\cal O_{\Xdaginf /R},\cal J\daginf )$$ représente le morphisme du théorème précédent. 
Pour montrer que c'est un isomorphisme la question est locale. Soit $\calUdag $ un ouvert du site il suffit de montrer que
le morphisme  induit:
$$ \cHom_{\calDdag_{\calUdag /R}}(
\cal O_{\calUdag /R},\cal I\dag_{\calUdag })\to \cHom_{\calDdag_{\calUdag /R}}(
\cal O_{\calUdag /R},\cal J\dag_{\calUdag })$$ est  un quasi-isomorphisme. Mais en vertu de l'équivalence de catégories \ref{can}, le complexe
$\cal I\dag_{\calUdag }:= R_{\calUdag }(\calIdaginf  )$ est un complexe de $\calDdag_{\calUdag /R}$-modules à gauche injectifs et en
vertu de l'équivalence de catégories \ref{equ-god} le complexe
$\cal J\dag_{\calUdag }:= R_{\calUdag }(\cal J\daginf )$ est un complexe de $\calDdag_{\calUdag /R}[\cal G_{\calUdag }]$-modules à
gauche injectifs. Mais l'extension $\calDdag_{\calUdag /R}\to\calDdag_{\calUdag /R}[\cal G_{\calUdag }]$ est plate,  donc
le complexe $\cal J\dag_{\calUdag }$ est une résolution  par des $\calDdag_{\calUdag /R}$-modules à gauche injectifs de $\cal I\dag_{\cal
U\dag}$ et donc de $\calMdaginf $, ce qui entraîne que le morphisme de complexes précédent est un quasi-isomorphisme.
\enddemo

\begin{coro}Soit  $\calMdaginf $ un complexe appartenant à la catégorie $ \Dplus ((\calDdag_{\Xdaginf /R}, \Sp)\Mod)$ alors le $\cal G_{\Xdaginf }$-module: $$ \cal Ext^{j}_{\calDdag_{\Xdaginf /R}}(
\cal O_{\Xdaginf /R},\calMdaginf )$$ est {\bf trivial} pour tout $j\geq 0$.
\end{coro}

\demo
En effet, en vertu du théorème précédent on a un isomorphisme canonique $\cal G_{\Xdaginf }$-modules:
$$ \cal Ext^{j}_{\calDdag_{\Xdaginf /R},\Sp}(
\cal O_{\Xdaginf /R},\calMdaginf )\simeq\cal Ext^{j}_{\calDdag_{\Xdaginf /R}}(
\cal O_{\Xdaginf /R},\calMdaginf )$$ et le module de gauche est trivial.\glossary{$\cal Ext^{j}_{\calDdag_{\Xdaginf /R},\Sp}(
\cal O_{\Xdaginf /R},\calMdaginf )$}\glossary{$\cal Ext^{j}_{\calDdag_{\Xdaginf /R}}(
\cal O_{\Xdaginf /R},\calMdaginf )$}
\enddemo

\begin{coro} Les {\bf faisceaux de cohomologie} du complexe de de Rham du site $ \Omega^{\bullet}_{\Xdaginf /K}$ sont des $ \cal G_{\Xdaginf }$-modules triviaux.
\end{coro}

\demo Il suffit de montrer en vertu du corollaire précédent qu'il existe un isomorphisme canonique de la catégorie $ \Dplus (K_{\Xdaginf })$:
$$ \Omega^{\bullet}_{\Xdaginf /K}\simeq \bfR \cHom_{\calDdag_{\Xdaginf /K}}(
\cal O_{\Xdaginf /K},\cal O_{\Xdaginf /K}).$$ 
Si $r\dag : \calWdag \to \calUdag $ est un morphisme du site on a un morphisme de restriction
$$ r^{-1}\Sp^{\bullet}(\cal O_{\cal
U\dag/K})\to \Sp^{\bullet}(\cal O_{\cal
W\dag/K})\leqno(\sharp)_{r\dag}:$$ compatible à la composition ce qui définit le complexe de Spencer sur le site $ \Xdaginf $, noté $ \Sp^{\bullet}(\cal O_{
\Xdaginf /K})$ \glossary{$ \Sp^{\bullet}(\cal O_{
\Xdaginf /K})$}.

\begin{lemm} Le complexe de Spencer $ \Sp^{\bullet}(\cal O_{
\Xdaginf /K})$ est une résolution du faisceau structural par des $ \calDdag_{
\Xdaginf /K}$-modules à gauche localement libres de type fini, en particulier le module $\cal O_{
\Xdaginf /K}$ est un  $ \calDdag_{
\Xdaginf /K}$-module à gauche parfait.
\end{lemm}

\demo
C'est une conséquence du fait que le complexe de Spencer $ \Sp^{\bullet}(\cal O_{\cal 
U\dag/K})$ est une résolution du faisceau structural en vertu du lemme \ref{Spen}.
\endsubdemo

En vertu du lemme précédent, on a un isomorphisme canonique de $ \Dplus (K_{\Xdaginf })$:
$$ \cHom\SubX{\calDdag_{\Xdaginf /K}}(
\Sp^{\bullet}(\cal O_{
\Xdaginf /K}),\cal O_{\Xdaginf /K})\simeq\bfR \cHom\SubX{\calDdag_{\Xdaginf /K}}(
\cal O_{\Xdaginf /K},\cal O_{\Xdaginf /K}).$$ Mais, par construction, le complexe de de Rham $ \Omega^{\bullet}_{\Xdaginf /K}$ est canoniquement
isomorphe au complexe $ \cHom_{\calDdag_{\Xdaginf /K}}( \Sp^{\bullet}(\cal O_{
\Xdaginf /K}),\cal O_{\Xdaginf /K})$.\enddemo

\bigskip
Nous allons maintenant étudier le rapport entre la cohomologie infinitésimale du complexe de de Rham du site et la cohomologie de de Rham $\dagger$-adique
du site.

\begin{prop}Soient $\calMdaginf $ et $\calNdaginf $ deux $\calDdag_{\Xdaginf /R}$-modules à gauche. 
Alors, on a l'égalité
de
$R$-modules:
$$ \hom_{\calDdag_{\Xdaginf /R}}(\calNdaginf , \calMdaginf )= \hom_{R_{\Xdaginf }}(R_{\Xdaginf },\cHom_{\calDdag_{\Xdaginf /R}}(\calNdaginf ,
\calMdaginf )).$$
\end{prop}

\demo Un élément de $ \hom_{\calDdag_{\Xdaginf /R}}(\calNdaginf , \calMdaginf )$ est la donnée
pour tout ouvert $\calUdag $ d'un morphisme $\calDdag_{\calUdag /R}$-linéaire $\varphi_{\calUdag }: \cal N\dag_{\calUdag }\to \cal M\dag_{\calUdag }$
tel que pour tout morphisme $r\dag: \calWdag \to \calUdag $  on ait l'égalité:
$$\varphi_{\calWdag }= (\sharp)_{r\dag}\circ \varphi_{\calUdag |W}\circ
(\sharp)^{-1}_{r\dag}.$$ Un élément de $ \hom_{R_{\Xdaginf }}(R_{\Xdaginf },\cHom_{\calDdag_{\Xdaginf /R}}(\calNdaginf ,
\calMdaginf ))$  est la donnée pour tout ouvert $\calUdag $ de l'image de $1$, qui est une section globale  du faisceau $\cHom_{\calDdag_{\calUdag /R}}(\cal N\dag_{\calUdag },
\cal M\dag_{\calUdag }))$, c'est-à-dire d'un
morphisme
$\calDdag_{\calUdag /R}$-linéaire
$\varphi'_{\calUdag }:
\cal N\dag_{\calUdag }\to \cal M\dag_{\calUdag }$ tel que pour tout morphisme $r\dag: \calWdag \to \calUdag $  on ait l'égalité:
$$\postskip1em
\varphi'_{\calWdag }= (\sharp)_{r\dag}\circ \varphi'_{\calUdag |W}\circ
(\sharp)^{-1}_{r\dag}\,.$$ Cela montre bien l'égalité de la proposition.
\enddemo

\begin{Rema}\em {\bf Attention!} Cette égalité {\bf ne se dérive pas} en un isomorphisme de foncteurs:
$$ \bfR \hom\SubX {\calDdag_{\Xdaginf /R}}(\calNdaginf , \calMdaginf )\simeq \bfR \hom\goodSub{10pt}{4pt}{-2mm}{R_{\Xdaginf }}(R_{\Xdaginf },\bfR \cHom\SubX {\calDdag_{\Xdaginf /R}}(\calNdaginf ,
\calMdaginf ))$$ parce que le foncteur $ \calMdaginf \fonct \cHom_{\calDdag_{\Xdaginf /R}}(\calNdaginf ,
\calMdaginf ))$ {\bf ne transforme pas} un $\calDdag_{\Xdaginf /R}$-module injectif en module acyclique pour le foncteur
$?\fonct \hom_{R_{\Xdaginf }}(R_{\Xdaginf }, ?)$, voir \ref{con-exe},
et tout le sel de cette théorie est là.
\end{Rema}
Cependant, ceci nous conduit à faire la conjecture d'annulation suivante.

\begin{conj}\label{conj}Soit $\calIdaginf  $ un $\calDdag_{\Xdaginf /K}$-module à gauche {\bf spécial} et injectif, alors
le faisceau de Zariski $\cHom_{\calDdag_{\Xdaginf /K}}(\cal O_{\Xdaginf /K},
\calIdaginf  )$ est acyclique pour le foncteur $?\fonct \hom_{K_{\Xdaginf }}(K_{\Xdaginf }, ?)$:
$$\forall  j>0, \quad  \ext^j_{K_{\Xdaginf }}(K_{\Xdaginf },\cHom_{\calDdag_{\Xdaginf /K}}(\cal O_{\Xdaginf /K}, \calIdaginf  ))=0.$$
\end{conj}
La conjecture précédente entraîne
l'isomorphisme canonique:
$$ \bfR \hom\SubX {\calDdag_{\Xdaginf /K},\Sp}(\cal O_{\Xdaginf /K}, \calMdaginf )\simeq 
\bfR \hom\goodSub{10pt}{4pt}{-2mm}{K_{\Xdaginf }}(K_{\Xdaginf },\bfR \cHom\SubX {\calDdag_{\Xdaginf /R}}(\cal O_{\Xdaginf /K},
\calMdaginf ))$$ pour tout complexe spécial $\calMdaginf $. En particulier,
cette conjecture fournirait éventuellement l'isomorphisme canonique:
$$ \bfR \hom_{\calDdag_{\Xdaginf /K},\Sp}(\cal O_{\Xdaginf /K}, \cal O_{\Xdaginf /K})\simeq 
\bfR \hom_{K_{\Xdaginf }}(K_{\Xdaginf },\Omega^{\bullet}_{\Xdaginf /K})$$ montrant que la cohomologie infinitésimale du complexe de Rham du site infinitésimal 
$ \Omega^{\bullet}_{\Xdaginf /K}$ fournit les bons nombres de Betti $p$-adiques et donne lieu à la suite spectrale Hodge-de Rham: $$\rm H^i_{\inf}(X,\Omega^{j}_{\Xdaginf /K})\Longrightarrow H^{i+j}_{dR}(X/K)$$ dans le contexte de cette  théorie, c'est dire l'importance de cette question.

Nous allons indiquer quelques points positifs en direction de la conjecture précédente qui semble assez profonde selon les calculs que nous avons faits pour essayer de démontrer cette
annulation.

\begin{Remas}\label{con-exe}

\begin{liste}\item 1) Dans le cas d'une variété affine $X$, on peut montrer que la non nullité de $ H^1_{\inf}(X/K, \cal O_{\Xdaginf /K})$
entraîne la non nullité de $ \ext^1_{\calDdag_{\Xdaginf /K}}(\cal O_{\Xdaginf /K}, \cal O_{\Xdaginf /K})$. Cela montre qu'un objet 
injectif de la catégorie $ (\calDdag_{\Xdaginf /K}, \Sp)\Mod$  n'est pas en général un  objet injectif de la catégorie $ \calDdag_{\Xdaginf /K}\Mod $ et la catégorie $ (\calDdag_{\Xdaginf /K}, \Sp)\Mod$ n'est pas stable par extensions comme sous-catégorie de la catégorie $ \calDdag_{\Xdaginf /K}\Mod$.

\item 2) On peut montrer que pour la droite affine $X$, l'espace: $$ \ext^1_{K_{\Xdaginf }}(K_{\Xdaginf },
\bfR \cHom_{\calDdag_{\Xdaginf /K}}(\cal O_{\Xdaginf /K},
\cal O_{\Xdaginf /K}))$$ est nul, ce qui  effectivement montre que le foncteur: 
$$ \calMdaginf \fonct \cHom_{\calDdag_{\Xdaginf /K}}(\cal O_{\Xdaginf /K},
\calMdaginf )$$ ne transforme pas un $\calDdag_{\Xdaginf /R}$-module à gauche injectif en objet acyclique pour le foncteur
$ \hom_{R_{\Xdaginf }}(R_{\Xdaginf }, ?)$.

\item 3) 
Plus généralement, on peut montrer que pour tout ouvert $U$ de la droite
affine, l'espace:
$$
\deuxlignes{0cm}{0cm}
 \ext^1_{K_{\Udaginf }}(K_{\Udaginf },
\bfR \cHom_{\calDdag_{\Udaginf /K},\Sp}(\cal O_{\Udaginf /K},
\cal O_{\Udaginf /K}))
\\
= \ext^1_{K_{\Udaginf }}(K_{\Udaginf },
\bfR \cHom_{\calDdag_{\Udaginf /K}}(\cal O_{\Udaginf /K},
\cal O_{\Udaginf /K}))
\endlignes$$
fournit le bon nombre de Betti ce qui donne un support crédible à la conjecture \ref{conj}.

\item 4) On peut montrer que pour tout ouvert $U$ de la droite affine, on a:
$$ H^1_{\inf}(U/K, K_{\Udaginf }) = \ext^1_{K_{\Udaginf }}(K_{\Udaginf }, K_{\Udaginf })=0.$$

\item 5) Tout ce qui précède montre
tout l'intérêt,    du passage de la catégorie $ \calDdag_{\Xdaginf /K}\Mod$ à la catégorie $ (\calDdag_{\Xdaginf /K}, \Sp)\Mod$. En fait tout notre problème pendant de longues années  venait clairement  de là. Mais comme toute idée féconde
il faut plusieurs années d'efforts et d'expériences avant de tomber dessus. Nous avons maintenant une idée très claire de la situation et il ne fait aucun doute pour nous que c'est le bon point de vue de la théorie de de Rham $p$-adique qui nous conduira à la bonne théorie que nous cherchons.
\end{liste}
\end{Remas}

\section{La catégorie $ \protect\Modd(\calDdag_{\Xdaginf /R},\Sp)$ des Modules à droite spéciaux sur le site infinitésimal}
Pour définir le foncteur image directe et le foncteur de dualité pour les modules à gauche spéciaux, il nous faut comprendre  la structure des modules à {\bf droite}  spéciaux sur le site
infinitésimal, qui a en soi-même un grand intérêt.

\subsection{La structure de $ \calDdag_{ \Xdaginf /V}$-module  à droite du fibré des formes différentielles de degré maximum}
Rappelons d'abord que si $X_s/R_s$ est un schéma lisse sur $R_s$, le faisceau des formes différentielles $\omega_{X_s/R_s}$ de degré maximum est un $\calD_{X_s/R_s}$-module à {\bf droite}. Si $x_1,\dots,x_n$ est un système de coordonnées locales, cette action
à droite est donnée par $(\omega)P:= (fdx)P= \trans P(f)dx$ où $\trans P := \sum_\alpha (-1)^{|\alpha|}\Delta_x^\alpha a_\alpha$ si 
$P = \sum_\alpha a_\alpha\Delta_x^\alpha$ ([EGA IV$_4$], $\S 16$).

Soit $\calXdag  =(X,\cal O_{\calXdag /R})$ un schéma $\dagger$-adique lisse sur $R$ et 
soit $\omega_{\calXdag /R}$ le fibré de rang $1$ des
formes  différentielles séparées de degré maximum.

\begin{defi} On définit le faisceau $\calDdag_{\calXdag /R}(\omega_{\calXdag /R})$ comme le sous-faisceau du faisceau des endomorphismes $\cal
End_R(\omega_{\calXdag /R})$  dont la réduction modulo $I^s$ est un opérateur différentiel de degré localement borné par une fonction linéaire en $s$.
\end{defi}

\begin{prop}\label{gau-dro} Si $R$ est un anneau de valuation discrète complet $V$, il existe un anti-morphisme canonique de faisceaux d'algèbres: $$\calDdag_{\calXdag /V}\rightarrow\calDdag_{\calXdag /V}(\omega_{\calXdag /V})$$ qui munit le faisceau $\omega_{\calXdag /V}$ 
d'une
structure de $\calDdag_{\calXdag /V}$-module à droite et qui est un isomorphisme.
\end{prop}

\demo 
Soient  $P$ et $\omega$ un opérateur différentiel et une forme différentielle séparée 
de degré maximum sur un ouvert $U$ muni d'un recouvrement par des ouverts
$U_i$ au-dessus desquels le faisceau des $1$-formes différentielles séparées est trivial. 
Si $x_1,\dots,x_n$ sont des fonctions sur $U_i$ dont les
différentielles forment une base du module des formes différentielles séparées, en vertu du théorème \ref{sym-tot},
l'opérateur $P$ admet sur $U_i$ un
développement
$\sum_{\alpha}a_\alpha \Delta_x^\alpha$ et en vertu du théorème
([Me$_3$], Coro. 4.2.2), la série $\trans P:= \sum_{\alpha}(-1)^{|\alpha|}\Delta^\alpha_x a_\alpha$ est un opérateur différentiel sur $U_i$. Si $\omega =
fdx_1\dots dx_n$ on pose:
$$(\omega) P:=\trans P(f)dx_1\dots dx_n.$$ Cela définit une action à droite qui ne dépend pas des coordonnées parce qu'il est en ainsi modulo $\goth m^s$ pour tout
$s\geq1$. Cette action se globalise et 
définit une structure de $\calDdag_{\calXdag /V}$-module à droite sur le faisceau $\omega_{\calXdag /V}$. La démonstration du théorème du symbole total
\ref{sym-tot} [M-N$_3$] montre  localement  que les sections du faisceau $\calDdag_{\calXdag /R}(\omega_{\calXdag /V})$ sont de
cette forme là et le morphisme de la proposition est donc un isomorphisme.
\enddemo


Un  morphisme  $r\dag: \calWdag \rightarrow \calUdag $  du site $ \Xdaginf $ induit par 
image inverse un morphisme de $\cHom_{R}(r^{-1}\omega_{\calUdag /V},
\omega_{\calWdag /V})$  par: $$r^*(fd(f_1)\wedge\cdots\wedge d(f_n)):= r^*fd(r^*f_1)\wedge\cdots\wedge d(r^*f_n)).\leqno(\sharp)_{r\dag}:$$

\begin{theo}\label{act-drt}Soit $\calUdag $ un relèvement d'un ouvert affine, alors l'action géométrique $(\sharp)_g$ d'un élément $g$  du groupe 
$\cal G_{\calUdag }$
sur
$\omega_{\calUdag /V}$ se fait à travers l'opérateur différentiel $g^{-1}$ autrement dit: $$(\sharp)_g(\omega)= (\omega) g^{-1}$$ où $\omega$ est une section du faisceau $\omega_{\calUdag /V}$. 
\end{theo}

\demo
La question est locale. Soit $x=(x_1,\dots,x_n)$ un système de coordonnées locales au-dessus d'un ouvert affine d'algèbre $A\dag$ et $f(x)dx$ une
forme différentielle, il faut montrer que: $$g(f)dg(x_1)\dots dg(x_n)=\trans g^{-1}(f)dx=\trans g^{-1}(1)g(f)dx.$$ 
En vertu du corollaire \ref{dev-grp}, l'application:
$$g\mapsto \delta(g):= (\delta_1(g),\dots,\delta_n(g)), \delta_i(g):= (g-1)(x_i)$$ est une bijection entre le groupe $G_{A\dag}$ et $(\goth
mA\dag)^n$.

\begin{lemm}Soit $g$ un élément du groupe $G_{A\dag}$ alors $g$ est égal au produit $g_n\dots g_1$ d'éléments du groupe tels que $\delta_i(g_j)= 0$ si
$i\neq j$.
\end{lemm}

\demo
L'application: 
$$(a_1,\dots,a_n)\mapsto \theta_{(a_1,\dots,a_n)}:= \sum_{\alpha\in \Bbb N^n}a^\alpha\Delta^\alpha$$ est l'inverse de l'application $\delta$. 
Posons $g_1:= (\theta_{(-g^{-1}(\delta_1(g)),0,\dots,0)})^{-1}$ alors $g= gg_1^{-1}g_1$ et $\delta_1(gg_1^{-1}) := \delta_1(g) +
g\delta_1(g_1^{-1})=0$. On a
$\delta_1(gg_1^{-1})=0$ et
$\delta_i(g_1)=0$ si $i\neq 1$. De proche en proche on construit la suite $g_1,\dots,g_n$ qui a la propriété du lemme.
\endsubdemo
\medskip
En vertu du lemme, pour démontrer le théorème on peut
supposer 
que $g:=g_i= \sum_{\alpha\in \Bbb N}a_i^\alpha\Delta^\alpha_i$,
de sorte que $\delta_j(g)=0$ si $ j\neq i$, et  $\delta_i(g)= a_i$. 
Pour simplifier les notations, posons aussi $a:= a_i$, et $\Delta:= \Delta_i$.
La fonction
$a$ appartient par définition à l'idéal engendré par $\goth m$ et il s'agit d'établir l'égalité:
$$\trans g^{-1}(1)=(1+ \Delta(a)).$$
Mais $g^{-1}= \sum_{k}(-1)^k(g-1)^k$ et
$$\trans g^{-1}=1+ 
\sum_{k>0}(-1)^k(\trans g-1)^k=\sum_{k}(-1)^k\Big(\sum_{\alpha}(-1)^{\alpha}\Delta^\alpha a^\alpha\Big)^k,$$ 
de sorte que l'on est ramené à montrer  l'égalité : 
$$\sum_{k\geq1}(-1)^k\Big(\sum_{\alpha\geq1}(-1)^{\alpha}\Delta^\alpha a^\alpha\Big)^k(1)= \Delta(a),$$
qu'il suffit de montrer modulo $\goth m^{s+1}$ pour tout $s\geq 1$, 
soit:
$$-\Delta(a)+\sum_{1\leq k\leq s}(-1)^k\Big(\sum_{1\leq\alpha\leq s}(-1)^{\alpha}\Delta^\alpha a^\alpha\Big)^k(1)\equiv 0.$$ 

En raisonnant par récurrence sur $s$, en appliquant l'opérateur différentiel
$-\sum_{1\leq\alpha\leq s}(-1)^{\alpha}\Delta^\alpha a^\alpha$ à la relation modulo $\goth m^{s+1}$,
on trouve la relation modulo $\goth m^{s+2}$, en tenant compte de l'égalité:
$$-\Big(\sum_{1\leq\alpha\leq s}(-1)^{\alpha}\Delta^\alpha a^\alpha\Big)\,(\Delta a)=\sum_{2\leq\alpha\leq s+1}(-1)^{\alpha}\Delta^\alpha (a^\alpha)$$ qui résulte
de l'identité:
$$\Delta^\alpha(a^\alpha\Delta(a))= \Delta^{\alpha+1}(a^{\alpha+1})$$ puisque l'algèbre $A\dag$ se plonge dans l'algèbre $A\dag_K$ et qu'on peut intégrer dans cette dernière algèbre.
\enddemo

\begin{coro}Soient $x_1,\dots,x_n$ des coordonnées locales sur un ouvert $U$ affine d'algèbre $A\dag$. Alors, l'application
$g\mapsto \log(\det(\Delta_j(g(x_i))_{ij}))$ est un cocycle de $G_{A\dag}$ à valeurs dans $A\dag_K$, qui n'est pas un cobord.
\end{coro}

\demo En effet, $\trans g^{-1}(1)=\det(\Delta_j(g(x_i))_{ij})$ en vertu du théorème précédent, et l'application $g\mapsto \log(\trans g^{-1}(1))$ est un cocycle qui n'est pas
un cobord en vertu du lemme
\ref{coc-nnu}.
\enddemo

\subsection{Le module de transfert $\calDdag_{\cal
U\dag\leftarrow \calWdag /V}$ pour une immersion ouverte}\glossary{$\calDdag_{\cal
U\dag\leftarrow \calWdag /V}$}

\`A cause de l'isomorphisme de la proposition \ref{gau-dro}, nous ne considérons que le cas du couple $(V,\goth m$) pour définir le module de transfert $\calDdag_{\cal
U\dag\leftarrow \calWdag /V}$.
\begin{defi}Soit $r\dag : \calWdag \to \calUdag $ un morphisme du site $ \Xdaginf $.  On définit le module de transfert
$\calDdag_{\cal
U\dag\leftarrow \calWdag /V}$ comme le module des morphismes $V$-linéaires  $r^{-1}\omega_{\calUdag /V}\to \omega_{\calWdag /V}$ 
dont la réduction  modulo $\goth m^s$ est un opérateur différentiel construit sur $r^*_s$ d'ordre localement borné par une fonction linéaire en $s$.
\end{defi}

\begin{prop} Le module de transfert $\calDdag_{\cal
U\dag\leftarrow \calWdag /V}$ ne dépend pas du relèvement  $r^*$, et est un $(r^{-1}\calDdag_{\calUdag /V}, \calDdag_{\cal
W\dag/V})$-bimodule engendré par toute section globale du faisceau  $\cal G_{\calWdag \to \calUdag }$.
\end{prop}

\demo La démonstration est la même que celle du cas des modules à gauche \ref{ouv-tra}.
\enddemo

\subsection{La catégorie des $\calDdag_{\Xdaginf /V}$-modules à \em{droite} spéciaux $\Modd(\calDdag_{\Xdaginf /V}, \Sp)$}
En vertu de la proposition précédente, le module de transfert $\calDdag_{\cal
U\dag\leftarrow \calWdag /V}$ est défini pour tout couple $(\calWdag , \calUdag)$ 
avec  $r: W\hookrightarrow U$. Par construction, on a un
plongement canonique
$\cal G_{\calWdag \rightarrow \calUdag }\hookrightarrow \calDdag_{\cal
U\dag\leftarrow \calWdag /V}$ et donc un morphisme canonique $V[\cal G_{\calWdag \rightarrow \calUdag }]\to \calDdag_{\cal
U\dag\leftarrow \calWdag /V}$.

\begin{defi}\label{mod-dro}
Soit  $X$ un  $V$-schéma lisse.
Un $\calDdag_{\Xdaginf /V}$-module à droite {\rm spécial} $\calMdaginf$ sur le site
$\Xdaginf $ est la donnée, pour tout ouvert $\calUdag $ du site
$\Xdaginf $, d'un $\calDdag_{\calUdag /V}$-module à droite $\cal M\dag_{\cal
U\dag}$, 
et la donnée, pour tout  couple d'objets $(\calWdag , \calUdag )$ du site
$\Xdaginf $, avec $r:
W\hookrightarrow U$,
d'un morphisme de $\calDdag_{\calWdag /V}$-modules à droite:
$$r^{-1}\cal N\dag_{\calUdag }\otimes_{r^{-1}\calDdag_{\calUdag /V}}\calDdag_{\cal
U\dag\leftarrow \calWdag /V}\rightarrow \cal N\dag_{\calWdag }.\leqno (\diamond):$$

En outre, ces données doivent satisfaire les conditions:
\begin{liste}\item 1) pour un couple $(\calUdag |W , \calUdag )$,  avec $r:
W\hookrightarrow U$, on a $\cal M\dag_{\calUdag |W}=\cal M\dag_{\calUdag }|W$ et le morphisme $(\diamond)$ coïncide avec le morphisme canonique,
\item 2) pour un triplet  $(\cal W'{}\dag ,\calWdag , \calUdag )$, avec $r': W'\hookrightarrow
W$ et $r: W\hookrightarrow U$,  le diagramme suivant est commutatif: 
$$\resetdisplay
\def\smt{\let\bigotimes\otimes\scriptwd1.3em\scriptdp=6pt
}\hss\mathrigid0mu  \scriptspace0,5pt\def\quad{\hskip0.8ex}
\matrix{\smt{r'^{-1}\mBig(r^{-1}\cal N\dag_{\calUdag }\Otimes_{r^{-1}\calDdag_{\calUdag /V}}\calDdag_{\cal
U\dag\leftarrow \calWdag /V}\mBig)}\Otimes_{r'^{-1}\calDdag_{\calWdag /V}}\calDdag_{\cal
W\dag\leftarrow \cal W'{}\dag /V}&\longrightarrow &\smt r'^{-1}\cal N\dag_{\calWdag }\Otimes_{r'^{-1}\calDdag_{\calWdag /V}}\calDdag_{\calWdag \leftarrow \cal W'{}\dag /V}\cr
\noalign{\kern-2pt}\downarrow&&\downarrow\cr
{(r\circ r')^{-1}\cal N\dag_{\calUdag }}\fbigotimes_{(r\circ r')^{-1}\calDdag_{{\calUdag /V}}}\calDdag_{\cal
U\dag\leftarrow \cal W'{}\dag /V}&\decale{-0.5cm}{\varto{0pt}{15ex} }&{\cal N\dag_{\cal W'{}\dag }}\,.}$$ 
\end{liste}
\end{defi}


Toutes les propriétés de la catégorie des modules spéciaux à gauche se transposent à la catégorie des modules à droite spéciaux $ \Modd(\calDdag_{\Xdaginf /V}, \Sp)$ avec les mêmes  arguments. En particulier,
 les morphismes de restriction $(\diamond)$ sont des isomorphismes et le module $\omega_{{\Xdaginf /V}}$ des formes différentielles séparées de degré maximum sur $ \Xdaginf $ est un $\calDdag_{\Xdaginf /V}$-module à droite  spécial en vertu de la proposition \ref{gau-dro}.

\medskip
Cependant, il y a un différence  saillante dans  l'action du
groupe $\cal G_{\Xdaginf}$. Considérons le morphisme $ \Inv: V[\cal G_{\calUdag }]\rightarrow \calDdag_{\calUdag /V}$ qui envoie un élément $g$ du groupe
sur l'opérateur différentiel $g^{-1}$.

\begin{prop}Supposons que $X$ est lisse sur $k$, alors pour un  $\calDdag_{\Xdaginf /V}$-module à droite $\calNdaginf $ spécial
et  pour tout relèvement
$\calUdag $ d'un ouvert affine $U$, l'action du groupe $\cal G_{\calUdag }$ sur $\cal N\dag_{\calUdag }$ se fait à travers le morphisme 
$\Inv$.
\end{prop}

\demo C'est une conséquence directe du théorème \ref{act-drt}.
\enddemo

\begin{coro}Soient $\calNdaginf $ un $\calDdag_{\Xdaginf /V}$-module à droite spécial et  $\calMdaginf $ un $\calDdag_{\Xdaginf /V}$-module à gauche spécial. Alors, le produit tensoriel $\cal
N\daginf \otimes_{\calDdag_{\Xdaginf /V}}\calMdaginf $ est un {\bf faisceau de Zariski} 
sur $X$ de $V$-modules.
\end{coro}

\demo En effet, l'action de groupe 
$g(n\otimes m):= gn\otimes g m= ng^{-1}\otimes g m= n\otimes g^{-1}gm=n\otimes m$ est triviale, et l'on applique le théorème 
\ref{equ-fai}. 
\enddemo

\section{Le foncteur image inverse dans la catégorie $ \Dmoins ((\calDdag_{\Xdaginf /V},\Sp)\protect\Mod )$}
Pour un morphisme $f$ de schémas  lisses sur $k$, on peut, bien sûr,  définir avec Grothendieck  un morphisme de topos entre les topos associés aux sites 
infinitésimaux $\dagger$-adiques des schémas. Mais, le paragraphe $8$ montre que ce sont des opérations insuffisantes et pathologiques, tout au moins pour la théorie de de Rham.

 Aussi, nous allons développer les opérateurs
cohomologiques pour les catégories des modules {\bf spéciaux} sur les sites infinitésimaux, et nous  allons voir que ce sont les bonnes  opérations fournissant à terme la bonne théorie que nous cherchons.

Le lecteur prendra garde tout de suite que les foncteurs cohomologiques que nous définissons
{\bf ne se prolongent pas} à la catégorie de tous les modules sur le site. Cela indique
que l'on ne pouvait espérer, 
 sans la notion essentielle de {\bf module spécial},
 développer les opérations cohomologiques pour la cohomologie de de Rham $p$-adique, ni obtenir en particulier  la fonctorialité de cette cohomologie et l'expression cohomologique $p$-adique de la fonction Zêta d'une variété algébrique lisse sur un corps fini.

\bigskip
Nous rappelons que l'on suppose le schéma $X$ lisse sur $R_1$.

\subsection{La catégorie $ (\calDdag_{ \Xdaginf /R}, \Sp)\Mod$ a suffisamment d'objets plats}
Nous allons d'abord montrer que la catégorie des modules spéciaux  a suffisamment d'objets plats, ce qui nous permettra de définir le foncteur image
inverse.


\begin{lemm}Soient  $\calUdag $ un ouvert du site $ \Xdaginf $ et $\cal P_{\calUdag }$ un $\calDdag_{\calUdag /R}$-module à gauche plat. Le
prolongement $ P_{\calUdag !}(\cal P_{\calUdag })$ est un $\calDdag_{\Xdaginf /R}$-module à gauche spécial, et plat comme objet de la catégorie
$ (\calDdag_{\Xdaginf /R}, \Sp)\Mod$.
\end{lemm}

\demo
Si $\calWdag $ est un ouvert, notons $r : U\cap W\hookrightarrow W$, $r' : U\cap W\hookrightarrow U$  les inclusions canoniques. Par définition:
$$P_{\calUdag !}(\cal P_{\calUdag })(\calWdag ):= r_!\mBig(\calDdag_{\calWdag |W\cap U\to\cal
U\dag/R}\otimes_{r'^{-1}\calDdag_{\calUdag /R}}r'^{-1}\cal P_{\calUdag }\mBig).$$ Mais $ \calDdag_{\calWdag |W\cap U\to\cal
U\dag/R}\otimes_{r'^{-1}\calDdag_{\calUdag /R}}r'^{-1}\cal P_{\calUdag }$ est un $\calDdag_{\calWdag |_{U\cap W}/R}$-module à gauche plat,  ce
qui entraîne que  $r_!\mBig(\calDdag_{\calWdag |W\cap U\to\cal
U\dag/R}\otimes_{r'^{-1}\calDdag_{\calUdag /R}}r'^{-1}\cal P_{\calUdag }\mBig)$ est aussi un $\calDdag_{\cal
W\dag/R}$-module à gauche plat et donc  $ P_{\calUdag !}(\cal P_{\calUdag })$ est plat comme objet
de la catégorie
$ (\calDdag_{\Xdaginf /R}, \Sp)\Mod$.
\enddemo

\begin{theo}\label{exi-pla} La catégorie $ (\calDdag_{\Xdaginf /R}, \Sp)\Mod$ a suffisamment d'objets plats.
\end{theo}
\demo
Soit  
$\calMdaginf $  un $\calDdag_{\Xdaginf /R}$-module à gauche spécial. 
Pour chaque ouvert $\calUdag $ de $\Xdaginf$, 
soit $\cal P_{\calUdag }\to \cal M\dag_{\calUdag }
\to0 $ une surjection d'un module plat. Alors: 
$$ P_{\calUdag !}(\cal P_{\calUdag })\to  P_{\calUdag !}(\cal M\dag_{\calUdag })
\to0 $$ est une surjection  d'un $\calDdag_{\Xdaginf /R}$-module à gauche spécial et plat. 
Si maintenant $\set\cal U\dag_i, i\in I/$ est une famille d'ouverts  du site 
tels que $\set U_i, i\in I/$ est un recouvrement de $X$, le morphisme
somme, composé des surjections: 
$$\oplus_i P_{\cal U\dag_i!}(\cal P_{\cal U\dag_i})\to\mkern-14mu\to \oplus_i P_{\cal U\dag_i!}(\cal M\dag_{\cal U\dag_i})\to\mkern-14mu\to  \calMdaginf 
\to0, $$ est une surjection  d'un $\calDdag_{\Xdaginf /R}$-module à gauche spécial,  plat  comme objet de la catégorie $ (\calDdag_{\Xdaginf /R}, \Sp)\Mod$. 
\enddemo

On pourra donc dériver à gauche dans la catégorie $ (\calDdag_{\Xdaginf /V}, \Sp)\Mod$ 
le foncteur produit tensoriel, et on note: $$\cal
N\daginf \Lotimes_{\calDdag_{\Xdaginf /V}, \Sp}\calMdaginf $$\glossary{$\cal
N\daginf \Lotimes_{\calDdag_{\Xdaginf /V}, \Sp}\calMdaginf $}le foncteur dérivé.

Les raisonnements  que nous avons faits pour les modules à gauche spéciaux  peuvent se faire pour les modules à droite spéciaux  et on pourra donc dériver à gauche le foncteur produit tensoriel dans la catégorie des modules à droite spéciaux.

\subsection{Le module de transfert $\calDdag_{\calYdag \rightarrow \calXdag /R}$ pour une immersion fermée}
Soit $f\dag: (Y,\cal O_{\calYdag /R})\rightarrow  (X,\cal O_{\calXdag /R})$ un morphisme de schémas $\dagger$-adiques sur $R$. Pour tout
$s\geq 1$, la réduction modulo $I^s$ est un morphisme de $R_s$-schémas $Y_s\rightarrow X_s$. Le module de transfert  $\calD_{Y_s\rightarrow 
X_s/R_s}$ est défini comme l'image inverse:
$$\calD_{Y_s\rightarrow 
X_s/R_s}:= \cal O_{Y_s/R_s}\otimes_{f^{-1}_s\cal O_{X_s/R_s}}f_s^{-1}\calD_{X_s/R_s}.$$ Il s'agit  du
faisceau des opérateurs différentiels $R_s$-linéaires de
$f^{-1}_s\cal O_{X_s/R_s}$ dans $\cal O_{Y_s/R_s}$, noté   $\cDiff_{R_s}(f^{-1}_s\cal O_{ X_s/R_s}, \cal O_{ Y_s/R_s})$  dans  ([EGA IV$_4$], $\S 16$).
C'est de façon naturelle un sous-bimodule filtré du
$(\calD_{Y_s/R_s}, f^{-1}_s\calD_{X_s/R_s})$-bimodule: $$\cHom_{R_s}(f^{-1}_s\cal O_{ X_s/R_s}, \cal O_{ Y_s/R_s}).$$

\begin{defi}\label{tra-rel}Le module de transfert $\calDdag_{\calYdag \rightarrow \calXdag /R}$,
\glossary{$\calDdag_{\cal
Y\dag\rightarrow \calXdag /R}$},construit sur le morphisme
$f\dag$, est le sous-bimodule  du $(\calD_{\calYdag /R}, f^{-1}\calD_{\calXdag /R})$-bimodule $\cHom_{R}(f^{-1}\cal O_{\calXdag /R}, \cal O_{\cal
Y\dag/R})$ des $R$-homomorphismes $P$ dont la réduction modulo $I^s$ est un élément de $\calD_{Y_s\rightarrow 
X_s/R_s}$ de degré localement borné par une fonction linéaire en $s$.
\end{defi}
Nous démontrons d'abord l'indépendance locale du module de transfert 
relativement aux relèvements du morphisme.

\begin{theo}\label{ind-tra}Soit $i: Y\rightarrow X$ une immersion \em{fermée}
de schémas  affines lisses sur $R_1$. Soient $\calYdag $ et $\calXdag $ des relèvements   
plats sur $R$, et soit $i\dag:\calYdag\to\calXdag$  un relèvement de l'immersion.
Alors, le module de transfert
$\calDdag_{\calYdag \rightarrow \calXdag /R}$ {\bf ne dépend pas} du relèvement $i\dag$ choisi.
\end{theo}

\demo
Une section $P$ du module de transfert construit sur $i\dag$ est, par construction,   une section du faisceau $\cHom_{R}(i^{-1}\cal O_{\calXdag/R}, \cal O_{\calYdag /R})$. Soit $i\dag_1$ un autre relèvement de $i$. Comme le relèvement $\calXdag $ est plat, il existe
un morphisme d'algèbres $g$ de $\cal O_{\calXdag /R}$ tel que l'on a la factorisation  $i^*= i_1^*\circ g$. En vertu du théorème \ref{inc-grp}, $g$ est un
opérateur différentiel et sa réduction $g_s$ modulo $\goth m^s$ est un opérateur différentiel d'ordre $s-1$.
Pour tout $s\geq1$, la réduction $P_s$ modulo $\goth m^s$
est par définition  un opérateur différentiel de $\calD_{Y_s\rightarrow 
X_s/R_s}:= \cal O_{Y_s/R_s}\otimes_{f^{-1}_s\cal O_{X_s/R_s}}f_s^{-1}\calD_{X_s/R_s}$, d'ordre
localement borné par une fonction linéaire en $s$. Mais 
le module $\cal O_{Y_s/R_s}\otimes_{f^{-1}_s\cal O_{X_s/R_s}}f_s^{-1}\calD_{X_s/R_s}$ est un $\cal
D_{X_s/R_s}$-module à droite engendré par
$i^*_s$. Autrement dit, $P$ admet la factorisation: 
$$i^{-1}\cal O_{X_s/R_s}\stackrel{i^{-1}Q_s}{\hbox to1cm{\rightarrowfill}} i^{-1}\cal O_{X_s/R_s}\stackrel{i^*_s}{\too}\cal O_{Y_s/R_s},$$ où $Q_s$ est un opérateur différentiel de
$\calD_{X_s/R_s}$ d'ordre borné par une fonction linéaire en $s$. Mais  $i^*_s= i_{1s}^*\circ g_s$ et  $P_s$ admet aussi la factorisation
$P_s= i_{1s}^*\circ g_s\circ Q_s$. Comme  $g_s\circ Q_s$ est un opérateur différentiel d'ordre localement borné par une fonction linéaire en $s$, cela montre que
$P$ est une section du module de transfert construit sur $i\dag_1$.
\enddemo

\begin{theo}\label{rg1}  Sous les conditions
précédentes, 
le module de transfert
$\calDdag_{\calYdag \rightarrow
\calXdag /V}$ est un
$i^{-1}\calDdag_{\calXdag /V}$-module à droite engendré par un relèvement
$i^*$. En particulier, c'est un $i^{-1}\calDdag_{\calXdag /V}$-module de type fini.
\end{theo}

\demo 
Il faut montrer que le morphisme:
$$i^{-1}\calDdag_{\calXdag /V}\to\calDdag_{\calYdag \rightarrow
\calXdag /V}$$ défini par la composition avec $i^*$ est {\bf surjectif}. 
La question est donc locale. Soit $U$ un ouvert affine de $X$ d'algèbre $A\dag$ munie 
d'un système de coordonnées 
$x_1,\dots,x_n$ tel que les différentielles $dx_1,\dots,dx_n$ 
forment une base du module des formes différentielles séparées,  soit $B\dag$ l'algèbre de
l'ouvert affine $U\cap Y$ de $Y$, et soit $i^*: A\dag\rightarrow B\dag$ le morphisme surjectif d'algèbres induit par  $i\dag$.  Pour $P :
A\dag\rightarrow B\dag$, section globale au-dessus de $U$ du faisceau  $\calDdag_{\calYdag \rightarrow \calXdag /V}$, définissons les
éléments $b_\alpha$ de l'algèbre
$B\dag$  par la formule:
$$b_\alpha := \sum_{0\leq \beta\leq\alpha}\comb\alpha\beta i^*((-x)^{\beta})P(x^{\alpha-\beta}).$$ Pour tout élément $a$ de l'algèbre $A\dag$, considérons la série
$$\phi(a):= \sum_{\alpha}b_\alpha i^*(\Delta^\alpha_x(a)),$$ qui converge formellement. Par hypothèse, $P(a)$ est égal à $\phi(a)$ modulo $\goth m^s$
pout tout
$s\geq1$ et donc: $$P(a)= \phi(a).$$

\medskip
Soit $(V[Y_1,\dots,Y_m])\dag\rightarrow A\dag$ une présentation de l'algèbre $A\dag$. En composant avec le morphisme d'algèbres $i^*$, on trouve une
présentation de l'algèbre $B\dag$. En vertu du théorème de ([M-N$_3$], Thm. 2.4.8), les topologies quotient sur les espaces $A\dag_K,B\dag_K$ sont {\bf séparées}, ce sont donc des topologies  $\cal L\cal F$, limite inductive dénombrable séparée d'espaces métriques complets.
Le théorème
du graphe fermé   
de Grothendieck pour les espaces vectoriels topologiques de type limite inductive dénombrable d'espaces métriques complets sur un corps valué complet
(cf. [M-N$_3$], Thm. 3.0.8) montre que le morphisme induit:
$$A\dag_K\rightarrow B\dag_K$$ est {\bf continu} pour les topologies de type $\cal L\cal F$ et le théorème de factorisation  
de Grothendieck pour les espaces vectoriels topologiques de type limite inductive dénombrable d'espaces métriques complets sur un corps valué complet (cf. [M-N$_3$], Thm. 3.0.9)
montre que les fonctions
$b_\alpha$ appartiennent à {\bf un cran $B\dag_{\rho}$}, image par la présentation induite du module  $ (V\langle Y_1,\dots,Y_m\rangle)_{\rho}$ des séries qui convergent  
dans un domaine $|Y_i|\leq \rho$, avec $\rho>1,$  pour une valeur absolue $|-|$ de $K$. Considérons la division  dans l'algèbre $(V[Y_1,\dots,Y_m])\dag$
par une base de division $J$ ([M-N$_3$], Thm. 2.3.4) de 
l'idéal $N$, noyau du morphisme d'algèbres
$(V[Y_1,\dots,Y_m])\dag\rightarrow B\dag$. 
La rétraction reste  de la division
$  {\Retr_J}: B\dag_{\rho} \rightarrow (V\langle Y_1,\dots,Y_m\rangle)_{\rho'}$, avec $\rho'\leq \rho$, permet de définir les éléments $a_\alpha$ de l'algèbre
$A\dag_{\rho'}$ comme les images des $ \Retr_J(b_\alpha)$ par la présentation $(V\langle Y_1,\dots,Y_m\rangle)_{\rho'}\rightarrow A\dag_{\rho'}$. Le
raisonnement, dans l'article [M-N$_3$], de la démonstration du théorème  du symbole total \ref{sym-tot} montre qu'il existe deux constantes $C_1>0$ et $0<\lambda<1$ telles que: 
$$||b_\alpha||_\rho\leq C_1\lambda^{|\alpha|}\,,$$ pour un nombre $\rho>1$ assez près de $1$, où $||-||_\rho$ est la norme quotient de la norme en $\rho$ de l'algèbre $(V\langle Y_1,\dots,Y_m\rangle)_{\rho}$. En vertu de la majoration  du reste du corollaire   
([M-N$_3$], Coro. 2.3.5), il existe une constante
$C_2>0$ telle que: $$|| \Retr_J(b_\alpha)||_\rho\leq C_2\lambda^{|\alpha|}\,,$$ pour $\rho>1$ assez près de $1$ et donc :$$||a_\alpha||_\rho\leq
C_2\lambda^{|\alpha|}\,,$$ pour $\rho>1$ assez près de $1$. En vertu du théorème
\ref{sym-tot}, la série: $$\tilde P(x,\Delta):= \sum_{\alpha}a_\alpha \Delta_x^\alpha$$ opère sur l'algèbre $A\dag$ et définit un opérateur différentiel.
Autrement dit, l'opérateur:
$$P: A\dag\rightarrow B\dag$$ admet une factorisation: $$A\dag\stackrel{\tilde P}{\rightarrow} A\dag\stackrel{i^*}{\rightarrow}B\dag\,,$$ où $\tilde
P$ est un opérateur différentiel de $A\dag$. D'où le théorème \ref{rg1}.
\enddemo

La même démonstration montre:

\begin{theo}\label{rg1h}Soit un entier $h\geq0$. 
Sous les conditions précédentes, le module de transfert
$\calD^{\dagger,h}_{\calYdag \rightarrow
\calXdag /V}$ est un
$i^{-1}\calD^{\dagger,h}_{\calXdag /V}$-module à droite engendré par un relèvement
$i^*$. En particulier, c'est un $i^{-1}\calD^{\dagger,h}_{\calXdag /V}$-module de type fini.\end{theo} 
\begin{Rema}Dans le cas d'un anneau de valuation discrète complet, la factorisation de $P_s$ du théorème \ref{ind-tra} est localement uniforme en $s$ et provient d'une
factorisation de $P$. En fait, le théorème précédent qui est un théorème de finitude algébro-topologique remarquable, a des conséquences importantes pour le calcul de l'image inverse, en particulier pour les fibrés $p$-adiques \ref{fib-pad}. Par ailleurs, cette propriété de finitude jouera un rôle essentiel dans la suite de cet article.
\end{Rema}

\begin{lemm}\label{ima-invd} 
Soit $i: Y\rightarrow X$ une immersion \em{fermée}
de schémas  affines lisses sur $k$. Soient
$\calYdag $ et $\calXdag $ des relèvements   plats sur $V$, 
et $i\dag$  un relèvement de $i$ d'idéal d'augmentation $\cal I_{\calYdag /V}$. Alors, 
on a une suite exacte de $i^{-1}\calDdag_{\calXdag /V}$-modules à droite:
$$0\rightarrow\cal I_{\calYdag /V}\otimes_{i^{-1}\cal O_{\calXdag /V} }i^{-1}\calDdag_{\calXdag /V}\rightarrow i^{-1}\calDdag_{\calXdag /V}\stackrel{i^*}{\rightarrow}\calDdag_{\calYdag \rightarrow\calXdag /V}\rightarrow0\,,$$
où le dernier morphisme est défini par le générateur $i^*$.
\end{lemm}

\demo La question est locale. 
Soit $U$ un ouvert affine de $ X$ d'algèbre $\dagger$-adique $A\dag$, munie d'éléments  $z_1,\dots,z_l$ tels que leur
différentielles forment une base du module des formes différentielles séparées. L'ouvert  $W=U\cap Y$ est un ouvert affine d'algèbre $B\dag$. 
Si $P$ est un opérateur différentiel au-dessus de $U$,  il admet en vertu du théorème \ref{sym-tot} un développement 
$P=\sum_{\alpha}a_\alpha\Delta_z^\alpha$, où $a_\alpha$ est une suite d'éléments de $A\dag$. Si $i^*\circ  P$ est nul dans $\calDdag_{\cal
W\dag\rightarrow\calUdag /V}$, nécessairement $i^*(a_\alpha)$ est nul dans $B\dag$ pour tout $\alpha$. Soit
$(V[Y_1,\dots,Y_m])\dag\rightarrow A\dag$ une présentation de l'algèbre $A\dag$, et soit $(V[Y_1,\dots,Y_m])\dag\rightarrow B\dag$ la
présentation de $B\dag$ induite, dont on notera $N$ le noyau. En vertu du théorème
\ref{sym-tot}, les fonctions $a_\alpha$ se relèvent en de fonctions 
$\tilde a_\alpha$ définies dans un même domaine. On peut effectuer la division des éléments $\tilde a_\alpha$ par une base de division de l'idéal $N$ ([M-N$_3$], Thm. 2.3.4).
Le reste de cette division est nul, et le théorème de continuité de la division ([M-N$_3$], Thm. 2.3.4) montre que les quotients de cette division sont
définis dans un même domaine, dont les normes admettent de majorations par des constantes du type $C\lambda^{|\alpha|}$. Cela montre que le noyau du
morphisme défini par $i^*$ est l'idéal à droite de
$\calDdag_{\calXdag /V}$ engendré par $\cal I_{\calYdag /V}$. Pour terminer, le morphisme naturel $\cal I_{\calYdag /V}\otimes_{\cal O_{\calXdag /V}}\calDdag_{\calXdag /V}\rightarrow
\cal I_{\calXdag /V}\calDdag_{\calXdag /V}$ est un isomorphisme parce que l'extension $\cal O_{
\calXdag /V}\rightarrow\calDdag_{\calXdag /V}$ est plate ([Me$_3$], Coro. 6.1.2). D'où le lemme.
\enddemo

\begin{Rema}Le lemme précédent dit que le morphisme défini par $i\dag$ induit un isomorphisme:
$$0\rightarrow\cal O_{\calYdag /V}\otimes_{i^{-1}\cal O_{\calXdag /V} }i^{-1}\calDdag_{\calXdag /V}\stackrel{i^*}{\rightarrow}
\calDdag_{\calYdag \rightarrow\calXdag /V}\rightarrow0$$
de $i^{-1}\calDdag_{\calXdag /V}$-modules à droite.  Autrement dit, de façon tout à fait remarquable, si on oublie un instant la structure de $\calDdag_{\calYdag /V}$-module à gauche, l'image inverse pour une immersion fermée coïncide avec l'image inverse  en tant que $\cal O_{\calXdag /V}$-module, ce qui est très utile, en particulier, pour étudier l'image inverse des fibrés $p$-adiques. 
\end{Rema}

\begin{defi} Si $i: Y\rightarrow X$ est une immersion fermée de schémas lisses sur $R_1$ admettant des relèvements $\calYdag $ et $\calXdag$  plats  sur $R$, on définit le module de transfert $\calDdag_{\calYdag \rightarrow \calXdag /R}$ comme le
sous-bimodule  du $(\calDdag_{\calYdag /R}, i^{-1}\calDdag_{\calXdag /R})$-bimodule $\cHom_{R}(i^{-1}\cal O_{\calXdag /R}, \cal O_{\cal
Y\dag/R})$  des homomorphismes qui sont localement des homomorphismes de transfert pour un relèvement local de $i$. 
\end{defi}
En vertu du théorème \ref{inc-grp'}, le module de transfert  $R[\cal G_{\calYdag \rightarrow \calXdag }]$ est alors un sous $R$-module du module de
transfert
$\calDdag_{\calYdag \rightarrow \calXdag /R}$.  

\begin{coro}Sous les conditions de la définition précédente, le module de transfert $\calDdag_{\calYdag \rightarrow \calXdag /V}$ est un $i^{-1}\calDdag_{\calXdag /V }$-module à droite engendré localement par un relèvement local de $i$.
\end{coro}

On définit alors le foncteur image inverse.
\begin{defi} Si $i: Y\rightarrow X$ est une immersion fermée de schémas lisses  sur $R_1$ admettant des relèvements $\calYdag $ et $\calXdag$  plats sur $R$, on définit le foncteur image inverse :
$$\calDdag_{\calXdag /R}\Mod \to \calDdag_{\calYdag /R}\Mod
\,,\qquad\cal M\dag_{\calXdag }\mapsto \calDdag_{\calYdag \rightarrow \calXdag /R}\smash{\Otimes_{i^{-1}\calDdag_{\calXdag/R}}}i^{-1}\cal M\dag_{\calXdag }\leqno i_{\diff}^{*,0} :$$ et le foncteur image inverse: 
$$
\left\{\matrix{\Dmoins((\calDdag_{\calXdag /R}))\Mod &\too& \Dmoins((\calDdag_{\calYdag /R}\Mod))\,,\hfill\cr
\cal M\dag_{\calXdag }&\longmapsto&\calDdag_{\calYdag \rightarrow \calXdag /R}\LOtimes_{i^{-1}\calDdag_{\calXdag/R}}i^{-1}\cal M\dag_{\calXdag },}\right.
\leqno{ i_{\diff}^{*} :=\bfL  i_{\diff}^{*,0}:}$$ 
$$\preskip-1ex
i_{\diff}^*\cal M\dag_{\calXdag }:= \calDdag_{\calYdag \rightarrow \calXdag /R}\Lotimes_{i^{-1}\calDdag_{\calXdag/R}}i^{-1}\cal M\dag_{\calXdag }.$$
\end{defi}
Ainsi, pour éviter toute confusion, le foncteur $i_{\diff}^{*,0}$  opère sur la catégorie des $\calDdag_{\calXdag /R}$-modules alors que le foncteur $i_{\diff}^{*}$ opère sur la catégorie des complexes de $\calDdag_{\calXdag /R}$-modules. Le foncteur image inverse $i_{\diff}^*$ est, par construction, un foncteur exact de catégories triangulées, de la catégorie $ \Dmoins ((\calDdag_{\calXdag/R}\Mod ))$ vers la catégorie $ \Dmoins ((\calDdag_{\calYdag /R})\Mod )$, 
et
en vertu de l'équivalence fondamentale \ref{can},
il
définit aussi
un foncteur exact de catégories triangulées de la catégorie $ \Dmoins ((\calDdag_{\Xdaginf /R},\Sp)\Mod)$ vers la catégorie $\Rm
\Dmoins ((\calDdag_{\Ydaginf /R},\Sp)\Mod)$.

\subsection{Le foncteur image inverse  dans la catégorie $ \Dmoins ((\calDdag_{\Xdaginf /R},\Sp)\protect\Mod )$ pour une immersion}
Soit $i: Y\rightarrow X$  une immersion  fermée de schémas lisses  sur $R_1$. On se propose de définir le foncteur image inverse:
$$i_{\diff}^* :   \Dmoins ((\calDdag_{\Xdaginf /R},\Sp)\Mod)\rightarrow \Dmoins ((\calDdag_{ \Ydaginf /R},\Sp)\Mod).$$
Soit $\set X_j,j\in J/$ un recouvrement de $X$ par des ouverts affines  et soit
$\set Y_j,j\in J/$ le recouvrement par des ouverts affines de $Y$ induit
sur $Y$.  Soient $\cal X\dag_j$ et $\cal Y\dag_j$ des relèvements plats de $X_j$ et de $Y_j$ pour $j\in J$, et soit $\calMdaginf $ un $\calDdag_{\Xdaginf /R}$-module
à gauche {\bf spécial} et $\cal M\dag_{\cal X_j}$ sa restriction à $\cal X_j$.
L'inclusion $i_j: Y_j\rightarrow X_j$ définit le module de transfert $\calDdag_{\cal Y\dag_j\rightarrow\cal X\dag_j/R}$ et on pose:
$$\cal N\dag_{\cal Y\dag_j}:= \calDdag_{\cal Y\dag_j\rightarrow\cal X\dag_j/R}\otimes_{i_j^{-1}\calDdag_{\cal X\dag_j/R}}i_j^{-1}\cal M\dag_{\cal X\dag_j},$$ qui est alors un $\calDdag_{\cal Y\dag_j/R}$-module à gauche.

\begin{prop} Les modules $\cal N\dag_{\cal Y\dag_j}$ admettent des isomorphis\-mes canoniques de recollement qui satisfont aux conditions de cocycle et
définissent un $\calDdag_{\Ydaginf /R}$-module à gauche 
spécial  $\calNdaginf $ qui ne dépend pas du recouvrement $\set X_j, \cal X\dag_j,j\in J/$ et qui dépend fonctoriellement de $\calMdaginf $.
\end{prop}

\demo 
En vertu de l'équivalence de catégories fondamentale \ref{can},
les prolongements  $P_{\cal Y\dag_j}(\cal N\dag_{\cal Y_j})$ des modules $\cal
N\dag_{\cal Y_j}$  définissent des modules spéciaux sur les sites de
$Y_j$. Ces modules sont munis d'isomorphismes de recollement sur les intersections qui proviennent des isomorphismes de recollement des restriction de $\calMdaginf $
qui satisfont aux conditions de cocycles et définissent un module donné localement. En vertu du théorème de recollement \ref{rec-fai}, on obtient un module spécial $\calNdaginf $ sur le site de $Y$ qui dépend 
fonctoriellement de $\calMdaginf $. Il est facile de voir que $\calNdaginf $ ne dépend pas du recouvrement  $\set\cal X_j,j\in J/$.
\enddemo
\begin{Rema}Les modules $\calDdag_{\cal Y\dag_j\rightarrow\cal X\dag_j/R}$ ne sont pas canoniquement isomorphes et {\bf on ne peut pas} induire le faisceau {\bf non spécial} $\calDdag_{\Xdaginf /R}$ comme $\calDdag_{\Ydaginf /R}$-module à gauche.\end{Rema}

\begin{defi}\begin{liste} \item 1)Si $i: Y\to X$ est une immersion fermée,  on définit l'image inverse: $$i_{\diff}^{*,0}\calMdaginf :=\calNdaginf .$$

\item 2)Si $i$ est une immersion ouverte, on définit le foncteur image inverse $ i_{\diff}^{*,0}$ comme le foncteur de restriction naturel qui est un foncteur exact.

\item 3) Si $i$ est une immersion localement fermée, on définit le foncteur image inverse $i_{\diff}^{*,0}$
par factorisation en une immersion fermée suivie d'une immersion ouverte. Ce foncteur ne dépend pas de la factorisation choisie.
\end{liste}
Dans le cas d'une immersion, on obtient un foncteur
exact à droite:
$$ (\calDdag_{\Xdaginf /R},\Sp)\Mod \rightarrow (\calDdag_{ \Ydaginf /R},\Sp)\Mod .\leqno i_{\diff}^{*,0}: $$Comme la catégorie des
modules spéciaux a suffisamment d'objets plats, ce foncteur se dérive à gauche en un foncteur exact de catégories triangulées:
$$ \Dmoins ((\calDdag_{\Xdaginf /R},\Sp)\Mod)\rightarrow \Dmoins ((\calDdag_{
\Ydaginf /R},\Sp)\Mod).\leqno i_{\diff}^*:= \bfL  i_{\diff}^{*,0}:$$
\end{defi}

\begin{prop}Supposons que $R$ est un anneau de valuation discrète complet $V$, et soit $i: Y\to X$ une immersion fermée de schémas lisses sur le corps résiduel $k$. Alors, la dimension cohomologique du foncteur image inverse
$i^*_{\diff}$ est localement bornée par $\codim_XY$.
\end{prop}

\demo
La question est locale. 
On peut supposer que $X$ est affine. Soient $\calYdag $ et $\calXdag $ des relèvements plats  sur $V$. Il suffit
de  montrer que le module de transfert $\calDdag_{\calYdag \rightarrow\calXdag /V}$ admet une résolution par des $\calDdag_{\calXdag/V}$-modules à droite localement libres de type fini de longueur bornée par $\codim_XY$. Si $i\dag$ est un relèvement local de l'immersion,  on a  en vertu du
lemme \ref{ima-invd}  un isomorphisme canonique de $i^{-1}\calDdag_{\calXdag/V}$-modules à droite:
$$0\rightarrow\cal O_{\calYdag /V}\otimes_{i^{-1}\cal O_{\calXdag /V} }i^{-1}\calDdag_{\calXdag /V}\stackrel{i^*}{\rightarrow}
\calDdag_{\calYdag \rightarrow\calXdag /V}\rightarrow0.$$En vertu de l'équivalence entre la catégorie des  faisceaux
de $\cal O_{\calXdag /V}$-modules cohérents et la catégorie $\Gamma(X,
\cal O_{\calXdag /V})$-modules de type fini, le faisceau $\cal O_{\calYdag /V}$ admet une résolution finie par des $\cal O_{\calXdag /V}$-modules localement libres de rang fini $\cal L^{\bullet}$, de longueur bornée par
$\codim_XY$. Le complexe:
$$\cal L^{\bullet}\otimes_{i^{-1}\cal O_{\calXdag /V}}i^{-1}\calDdag_{\calXdag /V}$$ est alors
une résolution de $ \calDdag_{\calYdag \rightarrow\calXdag /V}$ par des
$i^{-1}\calDdag_{\calXdag /V}$-modules
à droite localement libres de type fini.
\enddemo

\subsection{Le foncteur image inverse  dans la catégorie 
$ \Dmoins ((\calDdag_{\Xdaginf /R},\Sp)\protect\Mod)$ pour une projection}
Soit $p: Y\times_{R_1} X\rightarrow X$ la projection de schémas lisses sur $R_1$ sur le second facteur. On se propose de définir le foncteur image inverse:
$$p^*_{\diff} :   \Dmoins ((\calDdag_{\Xdaginf /R},\Sp)\Mod)\rightarrow \Dmoins ((\calDdag_{ (Y\times_{R_1} X)\daginf /R},\Sp)\Mod).$$
Si $Y$ et $X$ admettent  des relèvements  plats  $\calYdag $ et $\calXdag $,  leur produit fibré  
$\calYdag \times_{R}\calXdag $ existe et est  un relèvement 
plat du produit $Y\times_{R_1} X$, muni d'une projection sur $\calXdag $
([A-M$_2$]). On définit le module de transfert
$\calDdag_{\calYdag \times_R\calXdag \rightarrow\calXdag /R}$ à l'aide de cette projection.


Si $Y$ et $X$ sont affines, on peut considérer des relèvements  plats $\calYdag $ et $\calXdag $ et si  $\calMdaginf $ est un $\calDdag_{\Xdaginf /R}$-module à gauche spécial il provient canoniquement d'un $\calDdag_{\calXdag /R}$-module $\cal M_{\calXdag }$. On
définit alors:
$$\cal N\dag_{\calYdag \times_R\calXdag }:= \calDdag_{\calYdag \times_{R}\calXdag \rightarrow\calXdag/R}\otimes_{p^{-1}\calDdag_{\calXdag }}p^{-1}\cal M\dag_{\calXdag }\,,$$ qui est un $\calDdag_{\calYdag \times_R\calXdag}$-module et qui donne naissance à un $ \calDdag_{(Y\times X)\daginf /R}$-module spécial qui est, par définition, l'image inverse
$p^{*,0}_{\diff}\calMdaginf $ de $\calMdaginf $.

\medskip
Dans le cas général, soient $\set Y_j,j\in J/$ et $\set X_k,k\in I/$ des recouvrements de $Y$ et de $X$ par des ouverts affines, soit $\calMdaginf $ un $\calDdag_{\Xdaginf /R}$-module spécial, $\cal M\dag_{\inf,k}$ sa restriction à $X_k$ et $\cal N\dag_{\inf,j,k}$ son image inverse par la projection
$Y_j\times_{R_1} X_k\rightarrow X_k$. Les modules $\cal N\dag_{\inf,j,k}$, avec $j,k\in J\times I$, admettent des isomorphismes de recollement qui satisfont à la
condition de cocycle. En vertu du théorème \ref{champ}, ils définissent un $\calDdag_{(Y\times_{R_1} X)\daginf /R}$-module spécial $\cal
N\daginf$ qui ne dépend pas des recouvrements et dépend fonctoriellement de $\calMdaginf $. 
C'est $p^{*,0}_{\diff}\calMdaginf$, par définition.

\begin{defi}Sous les conditions ci-dessus, on définit l'image inverse par:
$$p^{*,0}_{\diff}\calMdaginf :=\cal
N\daginf .$$ On obtient ainsi un foncteur  exact à droite
$$   (\calDdag_{\Xdaginf /R},\Sp)\Mod \rightarrow (\calDdag_{ (Y\times X)\daginf /R},\Sp)\Mod,\leqno p^{*,0}_{\diff} :  $$ qui se dérive
en un foncteur exact de catégories triangulées:
$$   \Dmoins ((\calDdag_{\Xdaginf /R},\Sp)\Mod)\rightarrow \Dmoins ((\calDdag_{ (Y\times X)\daginf /R},\Sp)\Mod).\leqno p^*_{\diff} :=\bfL  p^{*,0}_{\diff} :$$
\end{defi}

\begin{theo}Supposons
que $R$ est un anneau de valuation discrète complet $V$. 
Soit $p: Y\times_{k} X\rightarrow X$ la projection de schémas  affines lisses sur le corps résiduel $k$, et soient $\calYdag $ et $\calXdag $ des relèvements plats sur $V$. Alors, le module de transfert
$\calDdag_{\calYdag \times_V\calXdag \rightarrow\calXdag/V}$ est un ${p^{-1}\calDdag_{\calXdag }}$-module {\bf plat}.
\end{theo}

\demo
La question est locale sur $Y\times_k X$. Soit $U$ un  ouvert  affine voisinage d'un point de
 $Y\times_k X$, d'algèbre $B\dag$. 
Nous pouvons supposer que la projection de
$U$ est contenue dans un ouvert affine $W$ d'algèbre $A\dag$ au-dessus duquel il existe des éléments  $x_1,\dots,x_n$ de $A\dag$ dont les
différentielles forment une base du module des formes différentielles séparées. Notons $u$  le morphisme  d'algèbres $A\dag\rightarrow
B\dag$ induit par la projection.
Soit $P$ un opérateur différentiel de $A\dag\rightarrow B\dag$,
nous posons:
$$b_\alpha := \sum_{0\leq \beta\leq\alpha}\comb\alpha\beta 
u((-x)^{\beta})P(x^{\alpha-\beta}).$$  Les éléments $b_\alpha$ sont des éléments bien définis de l'algèbre $B\dag.$ Le raisonnement de la démonstration du théorème \ref{sym-tot}
montre que la série: $$\sum_{\alpha}b_\alpha\xi^\alpha$$ est un élément de l'algèbre $(B\dag[\xi_1,\dots,\xi_n])\dag$. Quitte à diminuer l'ouvert $U$, on
peut supposer que les différentielles des fonctions $x_1,\dots,x_n$ font partie d'une base du module des formes différentielles séparées au-dessus de
$U$, de sorte que les opérateurs 
$\Delta_x^\alpha$ opèrent sur les éléments de l'algèbre $B\dag$ et l'opérateur $\sum_{\alpha}b_\alpha\Delta_x^\alpha$ opère, en vertu du théorème
\ref{sym-tot}, sur l'algèbre $B\dag$. On trouve,  sous les conditions précédentes, que le module de transfert  $\calDdag_{\cal
U\dag\rightarrow
\calWdag /V}$,  défini en \ref{tra-rel}, est un   sous-faisceau   d'anneaux du faisceau  $\calDdag_{ \calUdag /V}$ et on peut considérer la
filtration
$\calD^{\dagger,h}_{\calUdag \rightarrow
\calWdag /V}$ par les échelons qui est en fait la filtration induite par la filtration de $\calDdag_{ \calUdag /V}$. La démonstration du
théorème
\ref{noe} dans  [Me$_3$] montre que l'anneau $\Gamma(U,\calD^{\dagger,h}_{\calUdag \rightarrow
\calWdag /V})$ est noethérien.

\medskip
Soit $N$ un idéal à gauche de type fini de $\Gamma(W,\calDdag_{\calWdag /V})$. Il faut montrer que 
le morphisme: $$\Gamma(U,\calD^{\dagger}_{\calUdag \rightarrow
\calWdag /V})\otimes_{\Gamma(W,\calDdag_{\calWdag /V})}N\rightarrow \Gamma(U,\calD^{\dagger}_{\calUdag \rightarrow
\calWdag /V})$$ est injectif. Soient $P_1,\dots,P_l$ un système de générateurs de $N$ et $\sum_{i}Q_iP_i=0$ une relation dans 
$\Gamma(U,\calD^{\dagger}_{\calUdag \rightarrow
\calWdag /V})$. Pour $h\geq0$ assez grand, c'est une relation dans $\Gamma(U,\calD^{\dagger,h}_{\calUdag \rightarrow
\calWdag /V})$ et la platitude de l'extension $\Gamma(W,\calD^{\dagger,h}_{\calWdag /V})\rightarrow\Gamma(U,\calD^{\dagger,h}_{\cal
U\dag\rightarrow
\calWdag /V})$ entraîne la  platitude de l'extension $\Gamma(W,\calD^{\dagger}_{\calWdag /V})\rightarrow\Gamma(U,\calD^{\dagger}_{\cal
U\dag\rightarrow
\calWdag /V})$. La réduction modulo $\goth m^s$ de l'extension $\Gamma(W,\calD^{\dagger,h}_{\calWdag /V})\rightarrow\Gamma(U,\calD^{\dagger,h}_{\cal
U\dag\rightarrow
\calWdag /V})$ est plate, parce que le morphisme $U\rightarrow W$ est lisse. On est dans les conditions d'application  du critère local de platitude \ref{cri-pla} : les anneaux d'opérateurs différentiels 
$\Gamma(W,\calD^{\dagger,h}_{\calWdag /V})$ et $\Gamma(U,\calD^{\dagger,h}_{\cal
U\dag\to 
\calWdag /V})$ sont noethériens et l'idéal $\goth m$ est contenu dans le radical de l'anneau $\Gamma(U,\calD^{\dagger,h}_{\calUdag \rightarrow
\calWdag /V})$ ([Me$_3$], Thm. 4.1.1).
\enddemo
\begin{coro}\label{pro-exa} Supposons que $R$ est un anneau de valuation discrète complet $V$. Le foncteur image inverse $$  (\calDdag_{\Xdaginf /V},\Sp)\Mod \rightarrow (\calDdag_{ (Y\times
X)\daginf /V},\Sp)\Mod\leqno p^{*,0}_{\diff} : $$ est {\bf exact}.
\end{coro}

\demo En effet, la question étant  locale, c'est donc une  conséquence du théorème précédent.
\enddemo

\subsection{Le foncteur image inverse dans le cas général}
Un morphisme de schémas  lisses sur $k$ se factorise en une immersion  suivie d'une projection. On peut définir le morphisme image
inverse par
composition, mais il nous faut montrer que l'image inverse par une projection transforme un module plat en un module acyclique pour le foncteur image inverse par  l'immersion, pour
que le foncteur image inverse soit indépendant de la factorisation. Nous notons,  pour simplifier, $Z:=Y\times_kX$ et  $\calZdag $ le produit fibré
$\calYdag \times_V\calXdag $. Si $f:Y\to X$ est un morphisme,  notons $f:= p\circ i_f$, où $i_f$ est l'immersion graphe et où $p$ est la projection de $Y\times_kX$ sur $X$.

\begin{theo}\label{acy-pro}  Soit $f: Y\rightarrow X$ un morphisme de schémas \em{affines} lisses sur le corps résiduel $k$ et  soient $\calYdag $ et
$\calXdag $ des relèvements  plats. Alors, le morphisme canonique:
$$\calDdag_{\calYdag \rightarrow\calZdag /V}\LOtimes_{i_f^{-1}\calDdag_{\cal
Z\dag/V}}i_f^{-1}\calDdag_{\calZdag \rightarrow\calXdag /V}
\ \longrightarrow\ \calDdag_{\calYdag \rightarrow\calZdag /V}\Otimes_{i_f^{-1}\calDdag_{\cal
Z\dag/V}}i_f^{-1}\calDdag_{\calZdag \rightarrow\calXdag /V}$$ est un isomorphisme. De plus, il existe un isomorphisme: 
$$\calDdag_{\calYdag \rightarrow\calZdag /V}\otimes_{i_f^{-1}\calDdag_{\cal
Z\dag/V}}i_f^{-1}\calDdag_{\calZdag \rightarrow\calXdag /V}\simeq\calDdag_{\calYdag \rightarrow \calXdag/V}$$ pour tout relèvement de $f$ sur lequel est construit le module de transfert $\calDdag_{\calYdag \rightarrow \calXdag/V}$.
\end{theo}

\demo
Soit $i_f^*: i_f^{-1}\cal O_{\calZdag /V}\rightarrow\cal O_{\calYdag /V}$ un relèvement du morphisme  graphe de $f$.
En vertu du lemme \ref{ima-invd}, le relèvement $i_f^*$ définit un isomorphisme de $\calDdag_{\calZdag /V}$-modules à droite:
$$\cal O_{\calYdag /V}\otimes_{i_f^{-1}\cal O_{\calZdag /V}}i_f^{-1}\calDdag_{\calZdag /V}\simeq \calDdag_{\calYdag \rightarrow\cal
Z\dag/V}.$$ En vertu de l'équivalence entre la catégorie des  faisceaux de $\cal O_{\calZdag /V}$-modules cohérents et la catégorie $\Gamma(Z,
\cal O_{\calZdag /V})$-modules de type fini, le faisceau $\cal O_{\calYdag /V}$ admet 
une résolution finie par des $i_f^{-1}\cal O_{\calZdag /V}$-modules localement libres de rang fini $\cal L^{\bullet}$, et  le complexe:
$$\cal L^{\bullet}\otimes_{i_f^{-1}\cal O_{\calZdag /V}}i_f^{-1}\calDdag_{\calZdag /V}$$ est
une résolution de $ \calDdag_{\calYdag \rightarrow\calZdag /V}$ par des
$\calDdag_{\calZdag /V}$-modules
à droite localement libres de type fini, en vertu de la platitude évoquée dans la remarque \ref{pfi-pla}. Le complexe 
$\calDdag_{\calYdag \rightarrow\calZdag /V}\Lotimes_{i_f^{-1}\calDdag_{\calZdag /V}}i_f^{-1}\calDdag_{\calZdag \rightarrow \calXdag /V}$ se représente donc par le complexe:
$$\cal L^{\bullet}\otimes_{i_f^{-1}\cal O_{\calZdag /V}}i_f^{-1}\calDdag_{\calZdag \rightarrow\calXdag /V}\,,$$ et il s'agit de montrer que ce complexe  n'a de la cohomologie qu'au dernier
cran. La question est locale sur
$Y\times X$ et l'on peut remplacer le produit par un ouvert affine $U$ où le module des formes différentielles séparées est libre. Dans ce cas, nous avons vu que 
$\calDdag_{\cal
U\dag\rightarrow\calXdag /V}$ est un  faisceau d'anneaux  d'opérateurs différentiels {\bf relatifs} filtré par $\calD^{\dagger,h}_{\cal
U\dag\rightarrow\calXdag /V}$, et on est ramené à montrer la même assertion pour ce dernier faisceau. On obtient un complexe de $i_f^{-1}\calD^{\dagger,h}_{\cal
U\dag\rightarrow\calXdag /V}$-modules {\bf à droite}, dont les termes  sont sans torsion sur
$\goth m$.  Le faisceau $\calD^{\dagger,h}_{\cal
U\dag\rightarrow\calXdag /V}$ étant  acyclique au-dessus d'un ouvert affine assez petit ([Me$_3$], Thm. 3.2.3), on est réduit à la même assertion
pour les sections globales. Mais la cohomologie est de type fini,  parce que  l'anneau des sections globales du faisceau $\calD^{\dagger,h}_{\cal
U\dag\rightarrow\calXdag /V}$ au-dessus d'un ouvert affine assez petit est noethérien  en vertu du théorème \ref{noe}.
Il suffit, en vertu du lemme de Nakayama, de montrer la même assertion après réduction modulo $\goth m$ qui est
contenu dans le radical de cet anneau ([Me$_3$], Thm. 4.1.1). Mais, modulo $\goth m$,  ce complexe  se réduit à la restriction à $i_f^{-1}U$ du faisceau $\cal O_{Y/k}\otimes_{f^{-1}\cal O_{ X/k}}f^{-1}(\cal
D^{\dagger,h}_{\calXdag /V}/\goth m)\simeq \cal O_{Y/k}\otimes_{i_f^{-1}\cal O_{ Z/k}}i_f^{-1}(\cal
D^{\dagger,h}_{\calZdag \to\calXdag /V}/\goth m)$. D'où la première partie du théorème \ref{acy-pro}.

\bigskip
Puisqu'on est dans le cas affine, soit $f\dag$ un relèvement de $f$, alors le morphisme:
$$\calDdag_{\calYdag \rightarrow\calZdag /V}\otimes_{i_f^{-1}\calDdag_{\cal
Z\dag/V}}i_f^{-1}\calDdag_{\calZdag \rightarrow\calXdag /V}\to\calDdag_{\calYdag \rightarrow\calXdag/V}\leqno (*):$$ provient du morphisme composé:
$$f^{-1}\cal O_{\calXdag /V}= i_f^{-1}p^{-1}\cal O_{\calXdag /V}\rightarrow i_f^{-1}\cal O_{\cal
Z\dag/V}\rightarrow \cal O_{\calYdag /V},$$ où $p$ est la projection $Y\times_kX\to X$. Pour montrer que c'est un isomorphisme la question est locale sur $Y$. On peut supposer
que le fibré des formes différentielles séparées sur  $X$ est trivial et alors une section globale $P$ de $\calDdag_{\calYdag \rightarrow \calXdag /V}$ se
représente par une série $\sum_{\alpha}b_\alpha\Delta_x^\alpha$, où $(b_\alpha)$ est une suite d'éléments  de $\Gamma(Y,\cal O_{\cal Y\dag /V})$. De plus, le théorème de
continuité de la division  ([M-N$_3$], Thm. 2.3.4) montre que la suite
$(b_\alpha)$ se relève en une suite d'éléments $(a_\alpha)$ de $\Gamma(Y\times_kX,\cal O_{\cal Y\dag\times_V\cal X\dag /V})$ tels
que la série
$\sum_{\alpha}a_\alpha\Delta_x^\alpha$ est un opérateur $\tilde P$ de $\calDdag_{\calZdag \rightarrow\calXdag /V}$.
Autrement dit,  si on note $i_f\dag$ le morphisme graphe de $f\dag$, l'opérateur $P$ est l'image de $i_f^*\otimes\tilde P$ et le morphisme $(*)$ est
surjectif. D'autre part, 
$\calDdag_{\calYdag \rightarrow\calZdag /V}$  est un ${i_f^{-1}\calDdag_{\calYdag \times\calXdag/V}}$-module à droite engendré par $i_f^*$ en vertu du théorème \ref{ind-tra}, mais si $\tilde P$ est un opérateur de $\calDdag_{\calZdag \rightarrow\calXdag /V}$ et si l'image par $(*)$ de $u\otimes\tilde P$ est nulle, 
nous avons vu, en vertu du théorème de division, que 
$\tilde P$ appartient à l'idéal à droite engendré par l'idéal d'augmentation de $i_f^*$ et donc $u\otimes\tilde P$ est nul. Le
morphisme $(*)$ est injectif. D'où la seconde partie du théorème \ref{acy-pro}.
\enddemo

\begin{coro}\label{ind-tra'}\looseness-1
Soit $f: Y\rightarrow X$ un morphisme de schémas \em{affines} lisses sur le corps résiduel, et soient $\calYdag $ et $\calXdag$ des relèvements  plats. Alors, le module de transfert $\calDdag_{\calYdag \rightarrow \calXdag /V}$ ne
dépend pas du relèvement de
$f$.
\end{coro}

\demo En effet, en vertu du théorème \ref{ind-tra} le module de transfert de l'immersion qui est fermée ne dépend pas du graphe du relèvement. Le corollaire est conséquence de la
factorisation précédente.
\enddemo

\begin{defi}\label{tra-gen}Soit $f: Y\rightarrow X$ un morphisme de schémas  lisses sur le corps résiduel, et soient $\calYdag $ et $\calXdag $ des relèvements
 plats. On définit le module de transfert $\calDdag_{\calYdag \rightarrow \calXdag /V}$ comme  le sous-faisceau du
faisceau $\cHom_{V}(f^{-1}\cal O_{\calXdag /V}, \cal O_{\calYdag /V})$ des morphismes qui sont  {\bf localement}
des sections du module de transfert  d'un  relèvement local du triplet  $(Y, X, f)$.

\end{defi}

En vertu du corollaire précédent, le module $\calDdag_{\calYdag \rightarrow \calXdag/V}$ est bien défini et c'est de façon naturelle un
$(\calDdag_{\calYdag /V}, f^{-1}\calDdag_{ \calXdag /V})$-bimodule. En outre, il coïncide avec le module de transfert construit sur un relèvement global de
$f$ quand il existe.

\begin{coro}Si  $f: Y\rightarrow X$ est un morphisme de schémas  lisses sur le corps résiduel, alors le foncteur  image inverse 
$$p^{*,0}_{\diff} :   (\calDdag_{\Xdaginf /V},\Sp)\Mod \rightarrow (\calDdag_{(Y\times X)\daginf /V},\Sp)\Mod$$ transforme un $\calDdag_{\Xdaginf /V}$-module spécial plat en un module spécial {\bf acyclique} pour le foncteur image inverse $i_{\diff}^{*,0}$.
\end{coro}

\demo C'est une conséquence du théorème \ref{acy-pro}, parce que l'image inverse se calcule localement à la source et à la base.
\enddemo

\begin{defi}\label{fon-inv}Soit $f: Y\rightarrow X$ un morphisme de schémas  lisses sur le corps résiduel. On définit l'image inverse 
$$   (\calDdag_{\Xdaginf /V},\Sp)\Mod \rightarrow (\calDdag_{ \Ydaginf /V},\Sp)\Mod\leqno f^{*,0}_{\diff} : $$ 
comme le composé $i_{\diff}^{*,0}\circ p^{*,0}_{\diff}$,
lequel  se dérive donc à gauche en un foncteur exact de catégories triangulées 
\glossary{$f_*^{\diff}$}$$   \Dmoins ((\calDdag_{\Xdaginf /V},\Sp)\Mod)\rightarrow \Dmoins ((\calDdag_{
\Ydaginf /V},\Sp)\Mod).\leqno f^*_{\diff}:= \bfL f^{*,0}_{\diff}:$$
\end{defi}

On obtient, comme corollaire des résultats précédents, ce qui suit.
\begin{theo}Soit $f: Y\rightarrow X$ un morphisme de schémas  lisses sur le corps résiduel, et soient $\calYdag $ et $\calXdag $ des relèvements
 plats. Alors, le diagramme suivant est commutatif:
$$\def\quad{\hskip0.75ex}
\matrix{f^*_{\diff}:=f^*_{\diff, \calYdag , \calXdag }:& \Dmoins (\calDdag_{\calXdag /V}\Mod )&\rightarrow& \Dmoins (\calDdag_{ \calYdag /V}\Mod )\cr
&\downarrow\rscript{ P_{\calXdag }}&&\downarrow\rscript{ P_{\calYdag }}
\cr f^*_{\diff}:& \Dmoins ((\calDdag_{\Xdaginf /V},\Sp)\Mod)&\rightarrow& \Dmoins ((\calDdag_{\Ydaginf /V},\Sp)\Mod),}$$ où le foncteur de
la première ligne est le foncteur défini par le module de transfert:
$$\cal M\dag_{\calXdag }\mapsto \calDdag_{\calYdag \rightarrow \calXdag/V}\Lotimes_{f^{-1}\calDdag_{\calXdag  /V}}f^{-1}\cal M\dag_{\calXdag }\,,$$ et où les flèches verticales sont les foncteurs prolongement canoniques.
\end{theo}

Le foncteur image inverse dans la catégorie des modules spéciaux sur le site infinitésimal prolonge  comme il se doit le foncteur image inverse lorsque les relèvements des objets 
existent sans que les morphismes se relèvent. 

\begin{theo}\label{inv-tri}Si $f: Y\rightarrow X$ est un morphisme de schémas   lisses sur le corps résiduel, alors le complexe $f^*_{\diff}\cal
O_{\Xdaginf /V}$ est
canoniquement isomorphe au fibré trivial $\cal O_{\Ydaginf /V}$ muni de sa structure de $\calDdag_{\Ydaginf /V}$-module à gauche
spécial canonique.
\end{theo}

\demo
Il suffit de traiter le cas d'une immersion et d'une projection. Si $Y$, $X$ sont affines, soient $\calYdag $, $\calXdag $ des relèvements
 plats. Dans le cas d'une projection, le foncteur image inverse est exact et  le morphisme canonique de $\calDdag_{\calYdag \times_V\calXdag /V}$-modules à gauche:
$$ \calDdag_{\calYdag \times_V\calXdag \rightarrow\calXdag/V}\otimes_{p^{-1}\calDdag_{\calXdag /V}}p^{-1}\cal O_{\calXdag /V}\to \cal O_{\calYdag \times_V \calXdag /V}\leqno (*):$$ est surjectif et un calcul local à l'aide du théorème  du symbole total \ref{sym-tot} montre qu'il est injectif.

\bigskip
Dans le cas d'une immersion, il suffit de considérer le cas d'une immersion  fermée $i: Y\to X$ de schémas affines lisses sur $k$. Si $\calYdag $ et $\calXdag $ 
sont des relèvements  plats,  
le complexe: $$ \calDdag_{\calYdag \rightarrow \calXdag /V}\Lotimes_{i^{-1}\calDdag_{\calXdag /V}}i^{-1}\cal
O_{\calXdag /V}$$ est en vertu du lemme \ref{ima-invd} concentré en degré zéro. Le morphisme canonique de $\calDdag_{\calYdag /V}$-modules à gauche:
$$\calDdag_{\calYdag \rightarrow \calXdag /V}\otimes_{i^{-1}\calDdag_{\calXdag /V}}i^{-1}\cal
O_{\calXdag /V}\to \cal O_{\calYdag /V}$$ est surjectif, et le lemme \ref{ima-invd} montre qu'il est injectif.
\enddemo
\begin{Rema}Tous les résultats précédents du foncteur image inverse valent également sur l'extension $V\to K$ et nous les utiliserons aussi dans ce contexte.\end{Rema}
\begin{defi} \label{fib-pad}On dit qu'un $\calDdag_{\Xdaginf /V}$-module, resp. $\calDdag_{\Xdaginf /K}$-module, à gauche spécial $\calMdaginf $  est un fibré
$p$-adique, si pour tout ouvert $\calUdag $ de $\Xdaginf$ sa restriction $\cal M\dag_{\calUdag }$ est un $\cal O_{\calUdag /V}$-module, resp. $\cal
O_{\calUdag /K}$-module, localement libre de rang fini. La notion de fibré $p$-adique sur le site infinitésimal $\Xdaginf $ est {\bf purement algébrique}.
\end{defi} 
La notion de fibré $p$-adique a été introduite dans l'article [Me$_2$] sous cette forme-là dans le cas d'un relèvement. Autrement dit,  un fibré sur $K$ à connexion intégrable
est un fibré $p$-adique si l'action à gauche du faisceau  $\cal
D_{\calXdag /K}$ se prolonge en une action à gauche du faisceau $\calDdag_{\calXdag /K}$. C'est là une propriété restrictive.
\begin{notation}On note $ \MLS (\cal O_{\Xdaginf /V})$, resp. $ \MLS (\cal O_{\Xdaginf /K})$, la catégorie des fibrés $p$-adiques.
\end{notation}

\begin{prop}Soit $i: Y\to X$ une immersion fermée de schémas lisses sur $k$. Alors, le foncteur  $i^{*,0}_{\diff}$ image inverse est un foncteur {\bf exact} de
la catégorie $ \MLS (\cal O_{\Xdaginf /V})$ dans la catégorie $ \MLS (\cal O_{\Ydaginf /V})$, resp. de
la catégorie $ \MLS (\cal O_{\Xdaginf /K})$ dans la catégorie $ \MLS (\cal O_{\Ydaginf /K})$.
\glossary{$ \MLS (\cal
O_{\Xdaginf /K})$}
\end{prop}

\demo La question est locale, et découle de l'isomorphisme:  $$0\rightarrow\cal O_{\calYdag /V}\otimes_{i^{-1}\cal O_{\calXdag /V} }i^{-1}\calDdag_{\calXdag /V}\stackrel{i^*}{\rightarrow}
\calDdag_{\calYdag \rightarrow\calXdag /V}\rightarrow0$$
de $i^{-1}\calDdag_{\calXdag /V}$-modules à droite induit par un relèvement de $i$ en vertu du corollaire \ref{ima-invd} et de la platitude de l'extension  $\cal O_{\calXdag /V}\to\calDdag_{\calXdag /V}$ ([Me$_3$], Coro. 6.1.2).
\enddemo

\begin{Rema}Le raisonnement précédent ne marche pas dans le cas d'une projection, mais on peut montrer que le foncteur $p^{*,0}_{\diff}$, qui est exact en vertu du corollaire \ref{pro-exa}, 
envoie la catégorie $ \MLS (\cal O_{\Xdaginf /V})$ dans la catégorie $ \MLS (\cal O_{Z\daginf /V})$, resp. 
la catégorie $ \MLS (\cal O_{\Xdaginf /K})$ dans la catégorie $ \MLS (\cal O_{Z\daginf /K})$, mais la démonstration n'est pas de même nature. Nous n'utiliserons pas ce résultat dans cet article, mais il permet de donner intrinsèquement la définition des fibrés $p$-adiques munis d'une structure de Frobenius.
\end{Rema}

\subsection{Le module de transfert $\calDdag_{\calXdag \leftarrow \calYdag /V}$}
Si  $f_s: Y_s\rightarrow X_s$ est un morphisme de $R_s$-schémas lisses,
on peut définir le module de transfert  $\calD_{X_s\leftarrow
Y_s/R_s}$  comme:
$$\calD_{X_s\leftarrow
Y_s/R_s}:= \omega_{Y_s/R_s}\otimes_{\cal O_{Y_s/R_s}}\calD_{Y_s/R_s\rightarrow  X_s/R_s}
\otimes_{f^{-1}\cal O_{ X_s/R_s}}f^{-1}\sb{s}\omega^{-1}_{ X_s/R_s}.$$ 
C'est de façon naturelle un sous-bimodule filtré du
$(f^{-1}_s\calD_{X_s/R_s}, \calD_{Y_s/R_s})$-bimodule: $$\cHom_{R}(f^{-1}_s\omega_{X_s/R_s}, \omega_{Y_s/R_s}).$$

Soit $f\dag: (Y,\cal O_{\cal
Y\dag/V})\rightarrow  (X,\cal O_{\calXdag /V})$ un morphisme de schémas $\dagger$-adiques sur $V$.
\begin{defi}
Le module de transfert $\calDdag_{\calXdag \leftarrow \calYdag /V}$ \glossary{$\calDdag_{\calXdag\leftarrow \calYdag /V}$}est le sous-$(f^{-1}\calD_{\calXdag /V}, \calD_{\calYdag /V})$-bimodule  du bimodule $\cHom_{V}(f^{-1}\omega_{\calXdag /V}, \omega_{\calYdag /V})$ des
$V$-homomorphismes $P$ tels que leur réduction modulo $\goth m^s$ est un élément de $\calD_{X_s\leftarrow 
Y_s/V}$ dont l'ordre  est localement borné par une fonction linéaire en $s$.
\end{defi}


Notons $\omega^{-1}_{\calXdag /V}:= \cHom_{\cal O_{\calXdag /V},}(\omega_{\calXdag /V},\cal O_{\calXdag /V})$ le fibré dual du fibré $\omega_{\calXdag /V}$.
\begin{lemm}  Il existe des morphismes canoniques:  
$$\omega_{\calYdag /V}\otimes_{\cal O_{\calYdag /V}}\calDdag_{\calYdag \rightarrow \calXdag /V}
\otimes_{f^{-1}\cal O_{\calXdag /V}}f^{-1}\omega^{-1}_{\calXdag /V}\rightarrow \calDdag_{\calXdag \leftarrow \calYdag /V}$$ 
et
$$f^{-1}\omega_{\calXdag /V}\otimes_{f^{-1}\cal O_{\calXdag /V}}
\calDdag_{\calXdag \leftarrow \calYdag /V}\otimes_{\cal O_{\calYdag /V}}\omega^{-1}_{\calYdag /V}\rightarrow 
\calDdag_{\calYdag \rightarrow \calXdag /V},$$
qui sont des isomorphismes de $(f^{-1}\cal O_{\calXdag /V}, \cal O_{\calYdag /V})$-bimodules. 
\end{lemm}

\demo On a un morphisme canonique:
$$\deuxlignes{0cm}{0cm}
\omega_{\calYdag /V}\Otimes_{\cal O_{\calYdag /V}}\cHom_{V}(f^{-1}\cal O_{\calXdag /V}, \cal O_{\cal
Y\dag/V})\Otimes_{f^{-1}\cal O_{\calXdag /V}}f^{-1}\omega^{-1}_{\calXdag /V}\rightarrow
\\
\to\cHom_{V}(f^{-1}\omega_{\calXdag /V}, \omega_{\cal
Y\dag/V}),
\endlignes$$ 
qui provient du morphisme d'évaluation: 
$$\deuxlignes{0cm}{0cm}
\cHom_{V}(f^{-1}\cal O_{\calXdag /V}, \cal O_{\calYdag /V})\otimes_{f^{-1}\cal O_{\calXdag /V}}f^{-1}\omega^{-1}_{\calXdag /V}\rightarrow
\\\to
\cHom_{V}\big(f^{-1}\cHom_{V}(\omega^{-1}_{\calXdag /V}, \cal O_{\calXdag /V}),
\cal O_{\calYdag /V}\big)
\endlignes$$ 
suivi du morphisme canonique du produit tensoriel:
$$
\troislignes{0cm}{0cm}
\omega_{\calYdag /V}\otimes_{\cal O_{\calYdag /V}}\cHom_{V}\big(f^{-1}\cHom_{V}(\omega^{-1}_{\calXdag /V}, \cal O_{\calXdag /V}),
\cal O_{\calYdag /V}\big)\to
\\
\qquad\qquad\to \cHom_{V}\big(f^{-1}\cHom_{V}(\omega^{-1}_{\calXdag /V}, \cal O_{\calXdag /V}),
\omega_{\calYdag /V}\big)\simeq
\\\simeq\cHom_{V}(f^{-1}\omega_{\calXdag /V}, \omega_{\calYdag /V}).
\endlignes\postskip0pt$$ 
D'où un morphisme:
$$
\deuxlignes{0cm}{0cm}
\omega_{\calYdag /V}\otimes_{\cal O_{\calYdag /V}}\calDdag_{\calYdag \rightarrow \calXdag /V}
\otimes_{f^{-1}\cal O_{\calXdag /V}}f^{-1}\omega^{-1}_{\calXdag /V}\rightarrow\\
\to\cHom_{V}(f^{-1}\omega_{\calXdag /V}, \omega_{\calYdag /V}),
\endlignes$$ dont il
s'agit de voir que l'image est dans le module de transfert $\calDdag_{\calXdag \leftarrow \calYdag /V}$. Mais c'est une condition qui
porte
sur les réductions modulo $\goth m^s$, et qui est par construction satisfaite. De même, on construit le morphisme:
$$f^{-1}\omega_{\calXdag /V}\otimes_{f^{-1}\cal O_{\calXdag /V}}
\calDdag_{\calXdag \leftarrow \calYdag /V}\otimes_{\cal O_{\calYdag /V}}\omega^{-1}_{\calYdag /V}\rightarrow \calDdag_{\cal
Y\dag\rightarrow \calXdag /V}.$$ À partir de là, on voit aussitôt que ces morphismes sont des isomorphismes. 
\enddemo

\begin{coro} Le module $\omega_{\calYdag /V}\otimes_{\cal O_{\calYdag /V}}\calDdag_{\calYdag \rightarrow \calXdag /R}
\otimes_{f^{-1}\cal O_{\calXdag /V}}f^{-1}\omega^{-1}_{\calXdag /V}$ est naturellement
un $(f^{-1}\calDdag_{\calXdag /V}, \calDdag_{\calYdag /V})$-bimodule.
\end{coro}

\demo
Cela résulte  par transport de structure de la structure de $(f^{-1}\calDdag_{\calXdag /V}, \calDdag_{\calYdag /V})$-bimodule du module de transfert $\calDdag_{\calXdag\leftarrow \calYdag /V}$.
\enddemo
\begin{theo}\label{ind-tra''}Soit $f\dag: (Y,\cal O_{\cal
Y\dag/V})\rightarrow  (X,\cal O_{\calXdag /V})$ un morphisme de schémas $\dagger$-adiques lisses sur $V$. Le module  de transfert
$\calDdag_{\calXdag \leftarrow \calYdag /V}$ \em{ne dépend} que de la réduction modulo $\goth m$, $f: Y\rightarrow X$, de
$f\dag$.
\end{theo}
\demo
En vertu du lemme précédent, l'indépendance du module de transfert $\calDdag_{\calXdag\leftarrow \calYdag /V}$ du relèvement résulte alors de l'indépendance du module de transfert $\calDdag_{\cal
Y\dag\rightarrow \calXdag /V}$.
\enddemo
\begin{defi}Soit $f : Y\rightarrow X$ un morphisme de schémas lisses sur $k$, et soient
$\calYdag $ et $\calXdag $ des relèvements
 plats. On définit le module de transfert $\calDdag_{\calXdag \leftarrow \calYdag /V}$ comme le sous-bimodule  du 
  $(f^{-1}\calDdag_{\calXdag /V},
\calDdag_{\calYdag /V})$-bimodule
$\cHom_{V}(f^{-1}\omega_{\calXdag /V}, \omega_{\calYdag /V})$  des germes du module de transfert défini à l'aide d'un relèvement local du morphisme $f$. 
\end{defi}

\subsection{Le foncteur image inverse dans la catégorie  $ \Dmoins (\protect\Modd(\calDdag_{\Xdaginf /V}, \Sp))$}
On dispose de la catégorie abélienne  $ \Modd(\calDdag_{\Xdaginf /V}, \Sp)$ des modules à droite spéciaux, qui a suffisamment d'objets
d'injectifs et d'objets  plats et de sa catégorie dérivée $ \rmD(\Modd(\calDdag_{\Xdaginf /V}, \Sp))$. En considérant  les
modules de transfert
$\calDdag_{\calXdag \leftarrow \calYdag /V}$, on définit le foncteur image inverse $f^*_{\diff}$ par :
$$\cal N\dag_{\calXdag }\mapsto f^{-1}\cal N\dag_{\calXdag }\Lotimes_{f^{-1}\calDdag_{\calXdag /V}}\calDdag_{\calXdag \leftarrow
\calYdag /V}$$ puis le foncteur: 
$$   \Dmoins (\Modd(\calDdag_{\Xdaginf /V},\Sp))\rightarrow \Dmoins (\Modd(\calDdag_{
\Ydaginf /V},\Sp))\leqno f_{\diff}^* :=\bfL f_{\diff}^{*,0}:$$ pour un morphisme
$f : Y\to X$ de schémas lisses sur $k$ par factorisation comme dans le cas des modules à gauche spéciaux. Les propriétés de ce foncteur  sont toutes
parallèles aux propriétés du foncteur image inverse pour les modules à gauche spéciaux. On dispose ainsi et à titre d'exemple  d'un analogue au théorème \ref{inv-tri}.

\begin{theo}\label{inv-tri'}Soit un morphisme $f : Y\to X$ de schémas lisses  sur $k$. Alors, l'image inverse
$f^*_{\diff}\omega_{{\Xdaginf /V}}$ du fibré des formes différentielles séparées de degré maximum sur $ \Xdaginf $ est canoniquement isomorphe à
$\omega_{\Ydaginf /V}$ muni de sa structure canonique de $\calDdag_{\Ydaginf /V}$-module à droite spécial.
\end{theo}
La démonstration est
identique à celle de
\ref{inv-tri}.

\begin{Rema}Sous la condition $e<p^h(p-1)$, on peut reprendre les constructions précédentes pour construire le foncteur $f^{*}_{\diff,h}$, qui a les mêmes
propriétés que le foncteur $f^*_{\diff}$.
\end{Rema}

\section{Le foncteur image directe dans la catégorie  
$ \Dplus ((\calDdag_{\Xdaginf /K}, \Sp)\protect\Mod)$}
Pour définir le foncteur image directe il nous faut étudier d'abord quelques propriétés de finitude des opérateurs différentiels. 

\subsection{La platitude du module de transfert pour une immersion fermée}

\begin{lemm}\label{cri-jac}(Critère jacobien)
Soit  $i: Y\to X$ une immersion fermée de schémas affines lisses sur $k$. Soient $U$ un ouvert affine de $X$, $W$ sa trace sur $Y$, et $A\dag$ et $ B\dag$ des relèvements plats de $W $ et $U$. Si $u: A\dag\rightarrow B\dag$ est un
relèvement 
de l'inclusion de $W$ dans $U$,  
alors, \em{localement} sur $U$, il existe des éléments 
$x_1,\dots,x_r,x_{r+1},\dots,x_n$  de $A\dag$ tels que
$\set dx_1,\dots,dx_n/$ est une base locale  de $\Omega_{\cal
U\dag/V}$ et $\set du(x_1),\dots,du(x_r)/$ est une base locale de $\Omega_{\cal
W\dag/V}$. De plus, $\set x_{r+1},\dots, x_n/$ engendre l'idéal  noyau de $u$.
\end{lemm}

\demo Puisque  $W\to U$ est une immersion fermée, le morphisme $u: A\dag\rightarrow B\dag$ est {\bf surjectif} en vertu de la première partie du  théorème \ref{rel-pla}, d'où une suite exacte:
$$0\rightarrow \cal I_{\calWdag /V}\rightarrow\cal O_{\calUdag /V}\rightarrow\cal O_{\calWdag /V}\rightarrow0.$$ 
Soient 
$z_1,\dots,z_s$ des générateurs de l'idéal $\cal I_{\calWdag /V}$ et 
$\set dx_1,\dots,dx_n/$ une base locale de $\Omega_{\calUdag /V}$. Si $$0\rightarrow
\cal I_{W/k}\rightarrow\cal O_{ U/k}\rightarrow\cal O_{W/k}\rightarrow0$$ est la réduction modulo $\goth m$ de la suite précédente, les classes modulo $\goth m$ de
$z_1,\dots,z_s$ engendrent $\cal I_{W/k}$ et les classes de
$dx_1,\dots,dx_n$ forment une base locale de $\Omega_{ U/k}$. En vertu du critère jacobien algébrique ([EGA IV$_4$], 17.12.2), quitte à réindexer   $z_1,\dots,z_s$ et
$dx_1,\dots,dx_n$,   les classes de $z_{r+1},\dots,z_n$ engendrent l'idéal $\cal I_{W/k}$,  les classes de $dx_1,\dots,dx_r$ et de $dz_{r+1},\dots,dz_n$ forment
une base locale de
$\Omega_{ U/k}$, et les images des classes de $dx_1,\dots,dx_r$ forment une base locale de $\Omega_{
W/k}$. Les éléments $x_1,\dots,x_r,z_{r+1},\dots,z_n$ ont la propriété du lemme
parce que $\cal O_{\calWdag /V}$ est plat sur $V$.
\enddemo

\begin{defi}\label{couple-adapte}Étant donné un relèvement $u: A\dag\to B\dag$ d'une immersion fermée $Y\to X$ de schémas affines et lisses sur k, on dit qu'un système 
d'éléments  $z =\set z_1,\dots,z_q/$, $y=\set y_1,\dots,y_{n-q}/$ de $A\dag$ est adapté à $u$,
si $dz:= \set dz_1,\dots,dz_q/$, $dy:=\set dy_1,\dots,dy_{n-q}/$ forment une base de 
$\Omega_{A\dag/V}$, 
$du(y):=\set du(y_1),\dots,du(y_{n-q})/$ forment une base de $\Omega_{B\dag/V}$ et 
$\set z_1,\dots,z_q/$ engendre le noyau  de $u$.
\end{defi}
En vertu du critère jacobien \ref{cri-jac},  il existe, \em{\bf localement} sur $X,$ un système d'éléments  $x=(z,y)$ adapté à $u$. 

On rappelle que $\calD_{\calUdag /V}$ désigne le faisceau des opérateurs différentiels d'ordre fini et $\calD^{<\infty,h}_{\calUdag /V}$
le sous-faisceau des opérateurs différentiels d'ordre fini et d'échelon $h\geq0$.

\begin{coro}Le $\calD_{\calWdag /V}$-module à gauche $\cal O_{\calWdag /V}\otimes_{i^{-1}\cal O_{\calUdag /V}}i^{-1}\calD_{\cal
U\dag/V}$ est localement libre, et pour tout
$h\geq0$ le 
$\calD^{<\infty,h}_{\calWdag /V}$-module à gauche $\cal O_{\calWdag /V}\otimes_{i^{-1}\cal O_{\calUdag /V}}i^{-1}\cal
D^{<\infty,h}_{\calUdag /V}$ est localement libre.
\end{coro}

\demo En effet, si $x_1,\dots, x_r,x_{r+1}\dots,x_n$  sont des éléments  comme dans le 
critère jacobien  précédent, une base locale du $\calD_{\calWdag /V}$-module $\cal O_{\calWdag /V}\otimes_{i^{-1}\cal O_{\calUdag /V}}i^{-1}\cal
D_{\calUdag /V}$ est formée des tenseurs d'opérateurs $1\otimes\Delta_{r+1}^{\alpha_{r+1}}\cdots\Delta_n^{\alpha_n}$. Le même argument vaut pour
$\cal O_{\calWdag /V}\Otimes_{i^{-1}\cal O_{\calUdag /V}}i^{-1}\calD^{<\infty,h}_{\calUdag /V}$.
\enddemo

\begin{prop}\label{ind-tra'''} Soit  $i: Y\to X$ une immersion fermée de schémas   affines et lisses sur $k$ d'algèbres $B$ et $A$. Soient $A\dag$ et $ B\dag$ des relèvements plats et  $u:
A\dag\rightarrow B\dag$ un relèvement  du morphisme quotient
$A\rightarrow B$. Alors, le module de transfert $\calDdag_{\calXdag\leftarrow \calYdag /V}$ est un $i^{-1}\calDdag_{\calXdag /V}$-module à gauche {\bf localement} engendré par $dy_1\dots
dy_m\otimes u\otimes (dx_1\dots dx_n)^{*}$, où $dy_1\dots dy_m$, $dx_1\dots dx_n$ sont des bases locales  des fibrés des formes
différentielles séparées de $X$ et de $Y$.
\end{prop}

\demo En effet,  en vertu du théorème \ref{rg1}, le module de transfert $\calDdag_{\cal
Y\dag\rightarrow \calXdag /V}$ est engendré à droite par $u$ et l'on a l'égalité  :
$$dy_1\dots dy_m\otimes u \trans P\otimes
(dx_1\dots dx_n)^{*} = P(dy_1\dots dy_m\otimes u\otimes
(dx_1\dots dx_n)^{*})$$ comme morphismes:
$$\let\to\longrightarrow
i^{-1}\omega_{\calXdag/V}\stackrel{(dx)^*}{\to} i^{-1}\cal O_{\calXdag /V}\stackrel{u\circ\trans P}{\vrule height4pt width0pt\smash{\hbox to1cm{\rightarrowfill}}}\cal O_{\calYdag /V}\stackrel{dy}{\to} \omega_{\cal
Y\dag/V},$$ où $\trans P$ est le transposé de $P$ pour la forme $dx_1\dots dx_n$.
\enddemo

\begin{theo}\label{pla-tra}Soit $i : Y\rightarrow X$ une immersion fermée de schémas lisses  sur $k$, et soient $\calYdag $ et $\calXdag $ des
relèvements  plats. Les modules de transfert $\calDdag_{\cal
Y\dag\rightarrow \calXdag /V}$ et $\calDdag_{\calXdag\leftarrow \calYdag /V}$ sont {\bf plats} sur $\calDdag_{\calYdag /V}$.
\end{theo}

\demo La question est locale sur $Y$ et, par un choix de trivialisations locales des fibrés $\omega_{\calYdag /V}$ et $\omega_{\calXdag /V}$, il suffit de montrer que
le $\calDdag_{\calYdag /V}$-module à gauche $\calDdag_{\cal
Y\dag\rightarrow \calXdag /V}$ est plat.
Soient des ouverts affines assez petits $U $ et $W$ au-dessus desquels il existe  des éléments  $x_1,\dots,x_n$ comme dans le lemme du critère jacobien, et
soit $N$ un idéal  à droite de
$D\dag_{B\dag/V}$
de type fini.  Il faut montrer que le morphisme:
$$N\otimes_{D\dag_{B\dag/V}}\Gamma(W,\calDdag_{\cal
W\dag\rightarrow \calUdag /V})\rightarrow\Gamma(W,\calDdag_{\cal
W\dag\rightarrow \calUdag /V})$$ est injectif. L'anneau $D\dag_{B\dag/V}$ est filtré par $D^{\dagger,h}_{B\dag/V}$ et $\Gamma(W,\cal
D^{\dagger}_{\calWdag \rightarrow \calUdag /V})$ est filtré par les $D^{\dagger,h}_{A\dag/V}$-modules à droite  $\Gamma(W,\calD^{\dagger,h}_{\cal
W\dag\rightarrow \calUdag /V})$. 
Il 
suffit de montrer, pour $h$ assez grand, que:
$$N^h\otimes_{D^{\dagger,h}_{B\dag/V}}\Gamma(W,\calD^{\dagger,h}_{\cal
W\dag\rightarrow \calUdag /V})\rightarrow\Gamma(W,\calD^{\dagger,h}_{\cal
W\dag\rightarrow \calUdag /V})$$ est injectif, où $N^h$ est l'idéal engendré par les générateurs de $N$. Il suffit donc de montrer que $\Gamma(W,\calD^{\dagger,h}_{\cal
W\dag\rightarrow \calUdag /V})$ est plat sur $D^{\dagger,h}_{B\dag/V}$. Comme les anneaux 
$D^{\dagger,h}_{B\dag/V}, D^{\dagger,h}_{A\dag/V}$ sont  noethériens en vertu du théorème \ref{noe},  que l'idéal $\goth m$ est contenu dans la radical de l'anneau
$D^{\dagger,h}_{A\dag/V}$ ([Me$_3$], Thm. 4.1.1) et que $\Gamma(W,\calD^{\dagger,h}_{\calWdag \rightarrow \calUdag /V})$ est un $D^{\dagger,h}_{A\dag/V}$-module à droite de type fini en vertu du théorème \ref{rg1h},  on peut appliquer le  critère de platitude local \ref{cri-pla}. Pour cela, il faut s'assurer que les réductions modulo $\goth m^s$ du
$D^{\dagger,h}_{B\dag/V}$-module
$\Gamma(W,\calD^{\dagger,h}_{\calWdag \rightarrow \calUdag /V})$ sont plates. 
Mais ceci est une conséquence du corollaire précédent.
\enddemo

\subsection{Le théorème de changement de base pour une immersion fermée}
Le but de ce paragraphe est de montrer le théorème de changement de base pour une immersion fermée, qui nous servira à définir le foncteur image directe
dans le cas d'une immersion fermée.

\begin{prop}\label{iso-tra}Soient $\calUdag $ et $\cal U'{}\dag $ deux relèvements  plats d'un schéma $U$ affine et  lisse sur $k$.
Il existe alors un isomorphisme \em{canonique}
$$\calDdag_{\cal
U\dag\leftarrow \cal U'{}\dag /V}\otimes_{\calDdag_{\cal U'{}\dag /V}}\calDdag_{\cal
U'{}\dag \rightarrow \calUdag /V}\rightarrow \calDdag_{\calUdag /V}\leqno(*):$$ de $(\calDdag_{\calUdag /V},\calDdag_{\cal
U\dag/V})$-bimodules, et   un isomorphisme \em{canonique}
$$\calDdag_{\cal
U'{}\dag \rightarrow \calUdag /V}\otimes_{\calDdag_{\calUdag /V}}\calDdag_{\cal
U\dag\leftarrow \cal U'{}\dag /V}\rightarrow \calDdag_{\cal U'{}\dag /V}\leqno(**):$$ de $(\calDdag_{\cal U'{}\dag /V},\calDdag_{\cal
U'{}\dag /V})$-bimodules.
\end{prop}

\demo
Soit $u: A\dag\rightarrow A'{}\dag $ un relèvement de l'identité. Notons $u_0, u_n$  (l'indice $n$ rappelle que $u_n$ opère sur les $n$-formes différentielles), ses images dans $\calDdag_{\cal
U'{}\dag \rightarrow \calUdag /V}$ et $\calDdag_{\cal
U\dag\leftarrow \cal U'{}\dag /V}$. Le module de transfert $\calDdag_{\cal
U'{}\dag \rightarrow \calUdag /V}$ est un $(\calDdag_{\cal U'{}\dag /V},\calDdag_{\cal
U\dag/V})$-bimodule libre engendré par $u_0$, et le module de transfert $\calDdag_{\cal
U\dag\leftarrow \cal U'{}\dag /V}$ est un $(\calDdag_{\calUdag /V},\calDdag_{\cal
U'{}\dag /V})$-bimodule libre engendré par $u_n$. 

\begin{lemm}\label{uni}L'élément $u_n\otimes u_0$ de $\calDdag_{\cal
U\dag\leftarrow \cal U'{}\dag /V}\otimes_{\calDdag_{\cal U'{}\dag /V}}\calDdag_{\cal
U'{}\dag \rightarrow \calUdag /V}$ est indépendant de $u$.
\end{lemm}

\demo En effet, si $v$ est un autre relèvement, on a  les égalités 
$u= (uv^{-1})v=gv$  où $g$ est un élément de $G_{A'{}\dag }$. Alors, $u_n\otimes u_0= u_n\otimes gv_0= u_ng\otimes v_0$.
Il faut voir que $u_ng=v_n$ comme morphisme de $\omega_{A\dag/V}\to \omega_{A'{}\dag /V}$. Mais en vertu du théorème  \ref{act-drt},
l'action à droite de l'opérateur différentiel $g$ sur $\omega_{A'{}\dag /V}$ est égal au morphisme de restriction
$g^{-1}$ et $g^{-1}\circ u_n= v_n$.
\endsubdemo

On définit alors
le morphisme par linéarité: 
$$\calDdag_{\cal
U\dag\leftarrow \cal U'{}\dag /V}\otimes_{\calDdag_{\cal U'{}\dag /V}}\calDdag_{\cal
U'{}\dag \rightarrow \calUdag /V}\rightarrow \calDdag_{\calUdag /V}\leqno(*):$$ par 
$Pu_n\otimes u_0Q\mapsto PQ$. Le lemme suivant montre que le morphisme $(*)$ bien défini et  est injectif.

\begin{lemm}Pour tous opérateurs $P $ et $Q$ de $\calDdag_{\calUdag /V}$, on a l'égalité $Pu_n\otimes u_0Q=PQu_n\otimes u_0.$
\end{lemm}

\demo Il existe un unique opérateur $Q'$ de $\calDdag_{\cal U'{}\dag /V}$ tel que $Q'u_0= u_0Q$, comme éléments du module de transfert
$\calDdag_{\cal
U'{}\dag \rightarrow \calUdag /V}$. Il s'agit alors de voir que $u_nQ'= Qu_n$ comme éléments du module de transfert
$\calDdag_{\cal
U\dag\leftarrow \cal U'{}\dag /V}$, ou bien,   si le diagramme de gauche est commutatif, que le diagramme de droite est aussi commutatif:
$$\def\quad{\hskip0.5ex}
\matrix{A\dag&\stackrel{u_0}{\longrightarrow}&A'{}\dag \cr\downarrow\rscript Q&&\downarrow\rscript{ Q'}\cr A\dag&\stackrel{u_0}{\longrightarrow}&A'{}\dag }
\qquad\qquad \matrix{\omega_{A\dag/V}&\stackrel{u_n}{\longrightarrow}&\omega_{A'{}\dag /V}\cr
\downarrow \rscript{Q}&&\downarrow\rscript{ Q'}\cr
\omega_{A\dag/V}&\stackrel{u_n}{\longrightarrow}&\omega_{A'{}\dag /V}.}$$ La question est locale. Soit $x=(x_1,\dots,x_n)$ un système de coordonnées locales
 au-dessus d'un ouvert affine de $\calUdag $ et $x'=(x'_1,\dots,x'_n)=u(x)=(u(x_1),\dots,u(x_n))$ le système de coordonnées locales
au-dessus de l'ouvert affine de $\cal U'{}\dag $ obtenu par le morphisme $u$. Il s'agit de voir que si 
$Q'u_0= u_0Q$, alors $u_0 {\trans Q}=\trans Q'u_0$. Mais en vertu du théorème \ref{sym-tot}, si  l'opérateur $Q$ admet le développement $\sum_{\alpha}a_\alpha \Delta_x^\alpha,$ 
l'opérateur
$Q'$ admet le développement  $\sum_{\alpha}u(a_\alpha) \Delta_{x'}^\alpha,$  $\trans Q'= \sum_{\alpha}(-1)^{|\alpha|}\Delta_{x'}^\alpha u(a_\alpha)$  et  
$\trans Q'u_0 = \sum_{\alpha}(-1)^{|\alpha|}\Delta_{x'}^\alpha u(a_\alpha)u= u \circ(\sum_{\alpha}(-1)^{|\alpha|}\Delta_{x}^\alpha a_\alpha)= u_0{\trans Q}$.
\endsubdemo

\medskip
Le morphisme $(*)$  est alors un morphisme de bimodules surjectif. De même, on définit le morphisme $(**)$ qui est un isomorphisme.
\enddemo

\begin{coro}Soient $\calUdag $ et $\cal U'{}\dag $ deux relèvements  plats d'un schéma affine lisse sur $k$. Alors, les foncteurs:
$$\cal M'{}\dag _{\cal U'}\mapsto \calDdag_{\cal
U\dag\leftarrow \cal U'{}\dag /V}\Otimes_{\calDdag_{\cal U'{}\dag /V}}\cal M'{}\dag _{\cal
U'}\text{\qquad et\qquad }\cal M\dag_{\cal U}\mapsto \calDdag_{\cal
U'{}\dag \rightarrow \calUdag /V}\Otimes_{\calDdag_{\calUdag /V}}\cal M\dag_{\cal
U}$$ sont des équivalences de catégories \em{canoniquement} inverses l'une de l'autre entre la catégorie 
$ \calDdag_{\cal U'{}\dag /V}\Mod$
et la catégorie
$ \calDdag_{\calUdag /V}\Mod$.
\end{coro}

\demo C'est une conséquence directe des isomorphismes $(*) $ et $ (**)$ précédents.
\enddemo

\begin{prop}Soient $\calUdag $ et $\cal U'{}\dag $ deux relèvements  plats d'un schéma   affine lisse sur $k$. Alors, il existe un isomorphisme
canonique:
$$\calDdag_{\cal
U\dag\to \cal U'{}\dag /V}\simeq \calDdag_{\cal
U\dag\leftarrow \cal U'{}\dag /V}$$ de $(\calDdag_{\calUdag /V},\calDdag_{\cal U'{}\dag /V})$-bimodules.
\end{prop}

\demo
Soient $u: \cal O_{\cal U'{}\dag /V}\to \cal O_{\calUdag /V}$ un relèvement de l'identité et $P_{\cal
U\dag\to \cal U'{}\dag /V}$ un opérateur différentiel de $\calDdag_{\cal
U\dag\to \cal U'{}\dag /V}$. Alors, $P_{\cal
U\dag\to \cal U'{}\dag /V}$ admet la factorisation: $$\cal O_{\cal U'{}\dag /V}\stackrel{P}{\to}\cal O_{\cal U'{}\dag /V}\stackrel{u}{\to} \cal O_{\calUdag /V},$$ où $P$
est un opérateur différentiel. On lui associe l'opérateur $P_{\cal
U\dag\gets \cal U'{}\dag /V}$ composé: $$\omega_{\calUdag /V}\stackrel{u^{-1}}{\to}\omega_{\cal U'{}\dag /V}\stackrel{P}{\to} \omega_{\cal U'{}\dag /V},$$ lequel ne dépend
pas de la factorisation. En effet, si $v$ est un autre relèvement, on aura $u= vg$ où $g$ est un automorphisme d'algèbre qui se réduit à l'identité et donc un opérateur
différentiel
de $\cal O_{\cal U'{}\dag /V} $. En vertu du théorème \ref{act-drt}, on aura $gP\circ v^{-1}= P\circ (vg)^{-1}= P\circ u^{-1}$. Le morphisme ainsi défini est
un morphisme de
$(\calDdag_{\calUdag /V},\calDdag_{\cal U'{}\dag /V})$-bimodules qui est un isomorphisme.
\enddemo


Cela va nous permettre de définir les modules à gauche spéciaux à l'aide de l'image directe.

\begin{defi} Un $\calDdag_{\Xdaginf /V}$-module à gauche spécial $\calMdaginf $ est la donnée pour tout ouvert $\calUdag $ d'un $\calDdag_{\calUdag /V}$-module à gauche $\cal M\dag_{\calUdag }$ tel que $\cal M\dag_{\calUdag |W}=\cal M\dag_{\calUdag }|W$ pour tout $W\subset U$,
et  la donnée pour tout couple
$(\calWdag ,
\calUdag )$ avec
$r:W\hookrightarrow U$
d'un isomorphisme de $\calDdag_{\calWdag /V}$-modules à gauche 
$$\calDdag_{\calWdag \gets \calUdag |W/V}\otimes_{r^{-1}\calDdag_{\calUdag /V}}r^{-1}\cal M\dag_{\calUdag }\simeq \cal M\dag_{\cal
W\dag}\leqno (\diamond):$$ se réduisant au morphisme canonique dans le cas $(\calUdag |W, \calUdag )$ et satisfaisant aux conditions de transitivité
pour trois ouverts.
\end{defi}

Soient $i: Y\rightarrow X$ une immersion fermée de schémas affines lisses sur $k$,  $j:U\hookrightarrow X$ un ouvert affine de $X$ et $W$ sa
trace sur $Y$. Soient 
$\calYdag $, $\calXdag $, $\calUdag $ et  $\calWdag $ des relèvements  plats de $Y$, $X$, $U $ et $W$.

\begin{theo}\label{cha-bas}Soit $\cal N\dag_{\calYdag }$ un $\calDdag_{\calYdag /V}$-module à gauche. Il  
  existe alors un morphisme canonique de changement de base:
$$\DeuxLignes  
\calDdag_{\cal
U\dag\rightarrow \calXdag /V}\otimes_{{j^{-1}\calDdag_{\calXdag /V}}}j^{-1}i_*(\calDdag_{\calXdag\leftarrow \calYdag /V}\otimes_{\calDdag_{\calYdag /V}}\cal N\dag_{\calYdag })
\\
\longrightarrow\calDdag_{\cal
U\dag\leftarrow \calWdag /V}\otimes_{{\calDdag_{\calWdag /V}}}i_*(\calDdag_{\cal
W\dag\rightarrow \calYdag /V}\otimes_{j^{-1}\calDdag_{\calYdag /V}}j^{-1}\cal N\dag_{\calYdag }), 
\endlignes
$$
qui est un isomorphisme de
$\calDdag_{\calUdag /V}$-modules à gauche.
\end{theo}

\demo 
Nous allons d'abord construire un morphisme
$$\calDdag_{\cal
U\dag\rightarrow \calXdag /V}\Otimes_{{j^{-1}\calDdag_{\calXdag /V}}}j^{-1}i_*\calDdag_{\calXdag\leftarrow \calYdag /V}\longrightarrow i_*\mBig(\calDdag_{\cal
U\dag\leftarrow \calWdag /V}\Otimes_{{\calDdag_{\calWdag /V}}}\calDdag_{\cal
W\dag\rightarrow \calYdag /V}\mBig).\leqno(*):$$ Soit $P_{\calUdag \to\calXdag }\otimes j^{-1}Q_{\calXdag \gets \calYdag }$ un tenseur
élémentaire. En vertu de la proposition précédente, l'opérateur $P_{\calUdag \to\calXdag }$ donne naissance à un opérateur $P_{\calUdag \gets \calXdag|U}$ qui, composé avec
$j^{-1}Q_{\calXdag \gets \calYdag }$, donne un opérateur $R_{\calUdag \gets \calYdag |W}$. On définit l'image du tenseur $P_{\calUdag \to\calXdag}\otimes j^{-1}Q_{\calXdag \gets \calYdag }$ comme $R_{\calUdag \gets \calYdag |W}u_n\otimes u_0$ pour n'importe quel relèvement $u$ de l'identité de $W$.
En vertu  du lemme \ref{uni},
ce tenseur ne dépend pas du relèvement $u$. On obtient le morphisme $(*)$ de $(\calDdag_{\calUdag /V},j^{-1}\calDdag_{\calYdag /V})$-bimodules, qui est
un isomorphisme.

Soit maintenant $\cal N\dag_{\calYdag }$ un $\calDdag_{\calYdag /V}$-module à gauche. En tensorisant l'isomorphisme $(*)$ avec $j^{-1}\cal N\dag_{\cal
Y\dag}$ et tenant compte 
l'isomorphisme de
changement de base topologique  de 
$j^{-1}\calDdag_{\calXdag /V}$-modules à gauche:
$$j^{-1}i_*(\calDdag_{\calXdag\gets \calYdag /V}\otimes_{\calDdag_{\calYdag /V}}\cal N\dag_{\calYdag })\rightarrow i_*j^{-1}(\calDdag_{\calXdag\leftarrow \calYdag /V}\otimes_{\calDdag_{\calYdag /V}}\cal N\dag_{\calYdag })\,,$$ 
on trouve l'isomorphisme de changement de base
du théorème \ref{cha-bas}.
\enddemo

\subsection{Le foncteur image directe dans le cas d'une immersion fermée}

Soit $i: Y\rightarrow X$ une immersion fermée de schémas lisses sur $k$. On se propose de définir le foncteur image directe:
$$i^{\diff}_{*,0} :   (\calDdag_{\Ydaginf /V},\Sp)\Mod \rightarrow (\calDdag_{ \Xdaginf /V},\Sp)\Mod .$$

Soient $\calNdaginf $ un $\calDdag_{\Ydaginf /V}$-module à gauche  spécial et $\calUdag $  un ouvert affine du site
$\Xdaginf $. Si
$\calWdag $ est un relèvement de l'ouvert affine $Y\cap U$, la restriction $\cal N\dag_{\calWdag }$ est un $\calDdag_{\calWdag /V}$-module à
gauche, et on peut définir
$$\cal M\dag_{\calUdag }:= i_*\big(\calDdag_{\cal
U\dag\gets \calWdag /V}\otimes_{\calDdag_{\calWdag /V}}\cal N\dag_{\calWdag }\big).$$ En vertu du théorème de changement de
base pour une immersion fermée \ref{cha-bas}, le module
$\cal M\dag_{\calUdag }$ ne dépend pas à isomorphisme canonique près du relèvement $\calWdag $. Donc, en prenant un recouvrement par des ouverts
affines d'un ouvert
$\calUdag $, on définit un $\calDdag_{\calUdag /V}$-module à gauche $\cal M\dag_{\calUdag }$ donné localement qui ne dépend pas du recouvrement
choisi. D'autre part, si
$U'$ est un ouvert de $U$, on trouve que $\cal M\dag_{\calUdag |U'}=\cal M\dag_{\calUdag }|U'$.

Soit un 
couple d'objets du site 
$ (\cal U\dag_1,\cal U\dag_2)$ 
avec $r: U_1\hookrightarrow U_2$. Nous allons construire 
un isomorphisme de  $\calDdag_{\cal U\dag_1/V}$-modules à gauche:
$$\calDdag_{\cal U\dag_1\gets \cal U\dag_2|U_1}\otimes_{r^{-1}\calDdag_{\cal U\dag_2/V}}r^{-1}\cal M\dag_{\cal U\dag_2}\simeq \cal M\dag_{\cal
U\dag_1}\leqno (\diamond):$$ satisfaisant aux conditions de transitivité pour trois ouverts. 

\bigskip
Supposons d'abord $U_1$  affine et soit $\cal W\dag_1$ un relèvement de $U_1\cap Y$. Alors, $\cal M\dag_{\cal
U\dag_1}$ est par construction  canoniquement isomorphe à: 
$$i_*\big(\calDdag_{\cal
U\dag_1\gets \cal W\dag_1/V}\otimes_{\calDdag_{\cal W\dag_1/V}}\cal N\dag_{\cal W\dag_1}\big),$$ et 
$\cal M\dag_{\cal
U\dag_2|U_1}$ est par construction  canoniquement isomorphe à: $$i_*\big(\calDdag_{\cal
U\dag_2|U_1\gets \cal W\dag_1/V}\otimes_{\calDdag_{\cal W\dag_1/V}}\cal N\dag_{\cal W\dag_1}\big),$$ ce qui fournit l'isomorphisme
$(\diamond)$ dans ce cas-là par application du foncteur $\calDdag_{\cal U\dag_1\gets \cal U\dag_2|U_1/V}\otimes_{r^{-1}\calDdag_{\cal U\dag_2/V}}\!\! ?$, en tenant compte de l'isomorphisme canonique:
$$i^{-1}\calDdag_{\cal U\dag_1\gets \cal U\dag_2|U_1/V}\otimes_{i^{-1}r^{-1}\calDdag_{\cal U\dag_2/V}}i^{-1}\calDdag_{\cal
U\dag_2|U_1\gets \cal W\dag_1/V}\simeq \calDdag_{\cal U\dag_1\gets \cal W\dag_1/V}$$ 
parce que $\calDdag_{\cal U\dag_1\gets \cal U\dag_2|U_1/V}$ est un $r^{-1}\calDdag_{\cal U\dag_2/V}$-module à droite localement libre de rang $1$.

On obtient des isomorphismes locaux qui se recollent canoniquement pour fournir un isomorphisme
$(\diamond)$ dans le cas général.
Les
conditions de transitivité pour
$\calMdaginf $ sont conséquences de celles de $\calNdaginf $. Le morphisme $(\diamond)$ coïncide avec le morphisme canonique dans le cas du
couple $(\calUdag |U', \calUdag )$.
\begin{coro} Sous les conditions précédentes, le foncteur $\calUdag \fonct \cal M\dag_{\calUdag }$ est un $\calDdag_{\Xdaginf /V}$-module à gauche spécial
$\calMdaginf  $ sur le site
$\Xdaginf $.
\end{coro}

\begin{defi}Pour   une immersion fermée $i: Y\to X$ de schémas lisses sur $k$,  on définit l'image directe: $$i_{*,0}^{\diff}\calNdaginf :=\calMdaginf .$$ On obtient, en vertu du
théorème
\ref{pla-tra}, un foncteur {\bf exact} 
$$i^{\diff}_{*,0} :   (\calDdag_{\Ydaginf /V},\Sp)\Mod \rightarrow (\calDdag_{ \Xdaginf /V},\Sp)\Mod,\leqno i^{\diff}_{*,0} :  $$
qui se dérive  trivialement en un foncteur  exact de catégories triangulées:
$$ \rmD^*((\calDdag_{\Ydaginf /V},\Sp)\Mod)\rightarrow \rmD^*((\calDdag_{ \Xdaginf /V},\Sp)\Mod).\leqno i^{\diff}_{*}   :=\bfR  i^{\diff}_{*,0} :  $$
\end{defi}

\begin{coro}Soient  une immersion {\bf fermée} $i: Y\to X$ de schémas lisses sur $k$  et une immersion {\bf ouverte} $j: U\to X$. Alors, il existe un isomorphisme
canonique de foncteurs exacts de la catégorie $(\calDdag_{\Ydaginf /V},\Sp)\Mod$ vers la catégorie $(\calDdag_{\Udaginf /V},\Sp)\Mod$:
$$j^{*,0}_{\diff}\circ i_{*,0}^{\diff}\simeq i_{*,0}^{\diff}\circ j^{*,0}_{\diff},$$ et il existe un isomorphisme
canonique de foncteurs exacts de la catégorie $ \rmD^*((\calDdag_{\Ydaginf /V},\Sp)\Mod)$ vers la catégorie $ \rmD^*((\calDdag_{\Udaginf /V},\Sp)\Mod)$:
$$j^*_{\diff}\circ i_*^{\diff}\simeq i_*^{\diff}\circ j^*_{\diff}.$$
\end{coro}

\demo C'est la traduction du théorème de changement de base \ref{cha-bas}.
\enddemo
De même, en considérant un couple adapté \ref{couple-adapte} d'un relèvement  d'une immersion fermée de schémas lisses sur  $k$ on montre le changement de base pour une projection:
\begin{prop}\label{cha-pro} Soit un carré cartésien:
$$\def\quad{\hskip0.5ex}
\matrix{
Z\times_kY&\stackrel{i}{\too}& Z\times_kX\cr
\downarrow\rscript p&&\downarrow\rscript p \cr
Y&\stackrel{i}{\too}&X}
$$ de morphismes de schémas lisses sur $k$, où $i$ est une immersion fermée et $p$ est une projection. Il existe alors un isomorphisme
canonique de foncteurs exacts entre les catégories $(\calDdag_{\Ydaginf /V},\Sp)\Mod$ et  $(\calDdag_{(Z\times_kX)\daginf /V},\Sp)\Mod$:
$$p^{*,0}_{\diff}\circ i_{*,0}^{\diff}\simeq i_{*,0}^{\diff}\circ p^{*,0}_{\diff},$$ et il existe un isomorphisme
canonique de foncteurs exacts de la catégorie $ \rmD^*((\calDdag_{\Ydaginf /V},\Sp)\Mod)$ vers la catégorie $ \rmD^*((\calDdag_{(Z\times_kX)\daginf /V},\Sp)\Mod)$:
$$p^*_{\diff}\circ i_*^{\diff}\simeq i_*^{\diff}\circ p^*_{\diff}.$$
\end{prop}

\begin{prop}Soient $ Y\stackrel{i_{1}}{\to} X\stackrel{i_{2}}{\to} Z$ la composée  de deux immersions fermées de schémas lisses sur $k$. Alors, il existe un isomorphisme
canonique de foncteurs: $i_{2*}^{\diff}\circ i_{1*}^{\diff}\simeq (i_{2}\circ i_{1})^{\diff}_*$.
\end{prop}

\demo Les  foncteurs ${i_1}_{*,0}, {i_2}_{*,0}, ({i_2\circ i_1})_{*,0}$ sont exacts.  Pour $\calYdag $, $\calXdag $, $\calZdag $ des relèvements  plats,
alors on a un morphisme canonique  : 
$$i^{-1}_1\calDdag_{\calZdag \gets \calXdag /V}\otimes_{i^{-1}_1\calDdag_{\calXdag /V}}\calDdag_{\calXdag\gets
\calYdag /V}
\to\calDdag_{\calZdag \gets \calYdag /V}$$
de $(i_2\circ i_1)^{-1}\calDdag_{\calZdag /V}$-modules de type fini. En filtrant
ce morphisme par la filtration par les échelons, on montre que  c'est un isomorphisme. On obtient des isomorphisme locaux qui se recollent canoniquement.
\enddemo

\begin{Rema}On peut remplacer le faisceau $\calDdag_{X\daginf /V}$ dans la définition du foncteur image direct
$i_*^{\diff}$ par le faisceau $\calDdag_{X\daginf /K}$.
\end{Rema}

\subsection{Le foncteur image directe dans le cas d'une projection}
Pour  un couple $(Y,X)$ de schémas lisses sur $R_1$, on note $q $ et $ p$ les projections de $Y\times_{R_1}X$ sur le premier et le second facteur.

\begin{defi}Soit   $\calXdag $ un ouvert du site $ \Xdaginf $. On définit le site  $\Ydaginf \times\calXdag $ relatif à l'ouvert $\calXdag $
comme le sous-site de
$(Y\times_{R_1}X)\daginf$ dont les objets sont les ouverts de la forme $\calUdag \times_{R}\calXdag $ pour un ouvert $\calUdag $ du site 
$\Ydaginf $, dont les
morphismes sont de la forme $r\dag\times Id$ pour $r\dag$ un morphisme du site $\Ydaginf $ et muni de la topologie induite.

\end{defi} 
Sur le site $\Ydaginf \times\calXdag $ on a une description locale des faisceaux de modules analogue à celle des faisceaux de modules sur le site
$\Ydaginf $. En particulier, si $\cal L_{X}$ est un faisceau d'anneaux de Zariski sur $X$, on a un faisceau d'anneaux sur le site $\Ydaginf \times\calXdag$ dont la valeur  sur un ouvert $\calUdag \times_{R}\calXdag $ est  l'image inverse de $\cal L_{X}$ par la projection $U\times_{R_1}X\to X$ muni des restrictions
canoniques.
Notons $ (\Ydaginf \times\calXdag , p^{-1}\cal L_X)\Mod$ la catégorie des $p^{-1}\cal L_X$-modules à gauche sur le site
$\Ydaginf \times\calXdag $. Sur le site 
$\Ydaginf \times\calXdag $ on a aussi le faisceau $q^{-1}\cal G_{\Ydaginf }$, dont la valeur sur un ouvert 
$\calUdag \times_{R}\calXdag$ est le faisceau  $q^{-1}\cal G_{\calUdag }$  muni des restrictions canoniques.

\begin{prop} Le foncteur naturel:
$$ (Y\times_{R_1}X, p^{-1}\cal L_X)\Mod\to  (\Ydaginf \times\calXdag , p^{-1}\cal L_X)\Mod\leqno P_{Y\times_{R_1}X}:$$ induit une équivalence de catégories
entre la catégorie des faisceaux de Zariski sur le produit $Y\times_{R_1}X$, notée
 $ (Y\times_{R_1}X, p^{-1}\cal L_X)\Mod$, et la sous-catégorie des modules sur le site relatif $\Ydaginf \times\calXdag $
 dont l'action du faisceau de groupes $q^{-1}\cal G_{\Ydaginf }$ est {\bf triviale}.  Ce foncteur admet comme inverse canonique le foncteur  qui à un module
$\cal F_{\Ydaginf \times\calXdag }$ sur le site relatif associe le module dont la valeur sur $U\times_{R_1} X$,  si $U$ est affine, est :
$$\limproj_{\calUdag }\cal F_{\calUdag \times_{R}\calXdag },$$  la limite étant prise sur tous les ouverts qui relèvent
$U$.
\end{prop}

\demo
Remarquons d'abord que sous l'hypothèse de la proposition, pour deux ouverts affines $\cal U\dag_1, \cal U\dag_2$ qui relèvent le même ouvert, tous les
morphismes de restriction induisent le \em{même} isomorphisme :
$$\cal F_{\cal U\dag_1\times_{R}\calXdag }\simeq\cal F_{\cal U\dag_2\times_{R}\calXdag}.$$ La limite projective est bien définie. On obtient un faisceau donné localement sur un recouvrement de $Y\times_kX$ qui donne naissance
canoniquement à un faisceau de $p^{-1}\cal L_X$-modules sur $Y\times_{R_1}X$.
\enddemo


\bigskip
Sur le site relatif $\Ydaginf \times\calXdag $ on a le faisceau $q^{-1}\calDdag_{\Ydaginf /R}$, dont la valeur sur
l'ouvert $\calUdag \times_{R}\calXdag $ est le faisceau  $q^{-1}\calDdag_{\calUdag /R}$. Pour tout $q^{-1}\calDdag_{\Ydaginf /R}$-module à gauche $\calMdaginf $, on peut considérer le
faisceau sur le site relatif  
$\Ydaginf \times\calXdag $:
$$\calUdag \times_{R}\calXdag \fonct \cHom_{q^{-1}\calDdag_{\calUdag /R}}\big(q^{-1}\cal O_{\calUdag /R},\cal M\dag_{\cal
U\dag\times_{R}\calXdag }\big).$$

\begin{coro} Si $\cal M\dag_{(Y\times_{R_1}X)\daginf}$ est un $\calDdag_{(Y\times_{R_1}X)\daginf }$-module à gauche spécial, le foncteur: 
$$\calUdag \times_{R}\calXdag \fonct  \cHom_{q^{-1}\calDdag_{\calUdag /R}}\big(q^{-1}\cal O_{\calUdag /R},\cal M\dag_{\cal
U\dag\times_R\calXdag }\big)$$ définit un faisceau de \em{\bf Zariski} de $p^{-1}\calDdag_{\calXdag /R}$-modules à gauche sur le produit $ Y\times_{R_1}X$.
\end{coro}

\demo
En effet,  l'action d'un élément $g$ du groupe $q^{-1}\cal G_{\calUdag }$ sur un morphisme $\varphi$ est donnée par construction par $g\varphi
g^{-1}$, mais  l'action de $g$ sur $q^{-1}\cal O\dag_{\calUdag /R}$ et sur $\cal M\dag_{\cal
U\dag\times_{R}\calXdag }$ se fait à travers $q^{-1}\cal G_{\calUdag }\subseteq q^{-1}\calDdag_{\calUdag /R}$ et donc $g\varphi g^{-1}=\varphi
gg^{-1}=\varphi$. Le corollaire est conséquence de la proposition précédente.
\enddemo

\begin{notation}\label{nota}
Si $\calMdaginf $ est un $\calDdag_{(Y\times_{R_1}X)\daginf /R}$-module à gauche spécial  et $\calXdag $ est un ouvert 
du site 
$ \Xdaginf $,  on note $R_{\Ydaginf \times\calXdag }(\calMdaginf )$ sa restriction au site $\Ydaginf \times\calXdag $ et:
$$\cHom_{q^{-1}\calDdag_{\Ydaginf /R}}\mBig(q^{-1}\cal O_{\Ydaginf /R}, R_{\Ydaginf \times\calXdag }(\calMdaginf )\mBig)$$\glossary{$R_{\Ydaginf \times\calXdag }(\calMdaginf )$}le faisceau de {\bf Zariski} de
$p^{-1}\calDdag_{\calXdag /R}$-modules défini dans le corollaire précédent. Son image directe $p_{*}\cHom_{q^{-1}\calDdag_{\Ydaginf /R}}\big(q^{-1}\cal O_{\Ydaginf /R}, R_{\Ydaginf \times\calXdag }(\calMdaginf )\big)$ est donc un
$\calDdag_{\calXdag /R}$-module à gauche bien défini.\glossary{\vadjust{\kern-6pt}$\mathrigid0mu 
p_{*}\cHom\SubO{q^{-1}\calDdag_{\Ydaginf /R}}(q^{-1}\cal O_{\Ydaginf /R}, R_{\Ydaginf \times\calXdag }(\calMdaginf ))$ \hbox to4cm{}}
\end{notation}

Soit   $(\cal W\dag_1,\cal W\dag_2)$  un couple d'objets du site $ \Xdaginf $,  avec $r: W_1\hookrightarrow W_2$. Nous allons construire   un
isomorphisme de
$\calDdag_{\cal W\dag_1/V}$-modules  à gauche: 
$$
\scriptdp=6pt
\DeuxLignes  
(\diamond):\quad
\calDdag_{\cal W\dag_1\to \cal W\dag_2/R}\Otimes_{r^{-1}\calDdag_{\cal W\dag_2/R}}r^{-1} p_{*}\cHom\SubX {q^{-1}\calDdag_{\Ydaginf /R}}\mBig(q^{-1}\cal O_{\Ydaginf /R}, R_{\Ydaginf \times\cal W\dag_2}(\calMdaginf )\mBig)\simeq
\\
\simeq  p_{*}\cHom\SubX {q^{-1}\calDdag_{\Ydaginf /R}}
\mBig(q^{-1}\cal O_{\Ydaginf /R}, R_{\Ydaginf \times\cal W\dag_1}(\calMdaginf )\mBig)
\endlignes
$$ 
satisfaisant à la
condition de transitivité pour trois ouverts et les conditions portant sur les couples $(\cal W\dag_2|W_1, \cal W\dag_2)$. Le morphisme de projection
topologique suivant est un isomorphisme, parce que
$\calDdag_{\cal W\dag_1\to \cal W\dag_2/R}$ est un $r^{-1}\calDdag_{\cal W\dag_2/R}$-module à droite localement libre de rang $1$. On a donc l'isomorphisme:
$$\preskip0.5ex\postskip0.5ex\scriptdp=6pt
\DeuxLignes  
\calDdag_{\cal W\dag_1\to \cal W\dag_2/R}\Otimes_{r^{-1}\calDdag_{\cal W\dag_2/R}}r^{-1}p_{*}\cHom\SubX {q^{-1}\calDdag_{\Ydaginf /R}}\mBig(q^{-1}\cal O_{\Ydaginf /R}, R_{\Ydaginf \times\cal W\dag_2}(\calMdaginf )\mBig)\simeq
\\
\simeq p_*\mBig(p^{-1}\calDdag_{\cal W\dag_1\to
\cal W\dag_2/R}\Otimes_{p^{-1}r^{-1}\calDdag_{\cal W\dag_2/R}}r^{-1}\cHom\SubX {q^{-1}\calDdag_{\Ydaginf /R}}\big(q^{-1}\cal O_{\Ydaginf /R}, R_{\Ydaginf \times\cal W\dag_2}(\calMdaginf )\big)\mBig).
\endlignes$$ 
Il suffit de construire un
isomorphisme:
$$\preskip0.5ex\postskip0.5ex\scriptdp=6pt
\DeuxLignes  
p^{-1}\calDdag_{\cal W\dag_1\to\cal
W\dag_2/R}\Otimes_{p^{-1}r^{-1}\calDdag_{\cal W\dag_2/R}}r^{-1}\cHom\SubX {q^{-1}\calDdag_{\Ydaginf /R}}
\mBig(q^{-1}\cal O_{\Ydaginf /R},
R_{\Ydaginf \times\cal W\dag_2}(\calMdaginf )\mBig)\simeq
\\
\simeq \cHom\SubX {q^{-1}\calDdag_{\Ydaginf /R}}
\mBig(q^{-1}\cal O_{\Ydaginf /R}, R_{\Ydaginf \times\cal W\dag_1}(\calMdaginf )\mBig).
\endlignes$$
Les deux faisceaux précédents sur
$Y\times_kW_1$ sont donnés localement sur un recouvrement produit. Il suffit de construire   un isomorphisme de faisceaux donnés localement. Pour tout
ouvert $\calUdag $ du site $\Ydaginf ,$
le  morphisme produit tensoriel est un isomorphisme  parce que l'action de $q^{-1}\calDdag_{\calUdag /V}$ {\bf commute} avec celle de $p^{-1}\calDdag_{\cal W\dag_2/V}$:
$$\preskip1ex\postskip0ex\scriptdp=6pt
\DeuxLignes  
p^{-1}\calDdag_{\cal W\dag_1\to \cal
W\dag_2/R}\Otimes_{p^{-1}r^{-1}\calDdag_{\cal W\dag_2/R}}j^{-1}\cHom\SubX{q^{-1}\calDdag_{\calUdag /R}}\mBig(q^{-1}\cal O_{\calUdag /R}, \cal
M\dag_{\cal
U\dag\times_{R}\cal W\dag_2}\mBig)\simeq
\\
\simeq r^{-1}\cHom\SubX{q^{-1}\calDdag_{\calUdag /R}}\mBig(q^{-1}\cal O_{\calUdag /R}, p^{-1}\calDdag_{\cal W\dag_1\to \cal
W\dag_2/R}\Otimes_{p^{-1}r^{-1}\calDdag_{\cal W\dag_2/R}}\cal
M\dag_{\cal
U\dag\times_{R}\cal W\dag_2}\mBig).
\endlignes
$$ Il suffit de construire un isomorphisme:
$$p^{-1}\calDdag_{\cal W\dag_1\to\cal
W\dag_2/R}\otimes_{p^{-1}r^{-1}\calDdag_{\cal W\dag_2/R}}r^{-1}\cal
M\dag_{\cal
U\dag\times_{R}\cal W\dag_2}\simeq\cal
M\dag_{\cal
U\dag\times_{R}\cal W\dag_1}.$$ Mais  $p^{-1}\calDdag_{\cal W\dag_1\to \cal
W\dag_2/R}$ est un sous-faisceau de  $\calDdag_{\cal
U\dag\times_{R}\cal W\dag_1\to \cal
U\dag\times_{R}\cal W\dag_2/R}$ et le morphisme provient du morphisme du module spécial $\calMdaginf $. Pour montrer que c'est un
isomorphisme, on reprend le raisonnement de la démonstration du théorème \ref{iso-can}. Les compatibilités  proviennent des
compatibilités pour
$\calMdaginf $. Les conditions portant sur les couples induits sont évidentes.
Nous avons montré le théorème qui suit.
\begin{theo}Si $\calMdaginf $ est un $\calDdag_{(Y\times_{R_1}X)\daginf /R}$-module à gauche spécial, le
foncteur précédent:
$$\calXdag \fonct p_{*}\cHom_{q^{-1}\calDdag_{\Ydaginf /R}}\big(q^{-1}\cal O_{\Ydaginf /R}, R_{\Ydaginf \times\calXdag }(\calMdaginf )\big) $$ est un  $\calDdag_{ \Xdaginf /R}$-module à
gauche {\bf spécial}.
\end{theo}
\begin{defi}
Soit $\calMdaginf $ un  $\calDdag_{(Y\times_{R_1}X)\daginf /R}$-module à gauche spécial. On définit son image directe  $p_{*,0}^{\diff}\calMdaginf $ comme  le
$\calDdag_{\Xdaginf /R}$-module spécial placé en degré $-[\dim Y]$ dont la valeur  sur l'ouvert $\calUdag $ de $ \Xdaginf $ est: $$p_{*}\cHom_{q^{-1}\calDdag_{\Ydaginf /R}}\mBig(q^{-1}\cal O_{\Ydaginf /R}, R_{\Ydaginf \times\calUdag }(\calMdaginf )\mBig)[\dim Y].$$ 
\end{defi}

\begin{prop} L'application $\calMdaginf \mapsto  p_{*,0}^{\diff}\calMdaginf $ définit un foncteur covariant exact à gauche 
$$ (\calDdag_{(Y\times_{R_1}X)\daginf /R},\Sp)\Mod 
\rightarrow (\calDdag_{ \Xdaginf /R},\Sp)\Mod [\dim Y]\leqno p_{*,0}^{\diff}:$$ de 
$(\calDdag_{(Y\times_{R_1}X)\daginf /R},\Sp)\Mod$
vers la catégorie $(\calDdag_{ \Xdaginf /R},\Sp)\Mod [\dim Y]$ des modules placés en degré $-[\dim Y]$.
\end{prop}

\demo En effet, la construction du module  $p_{*,0}^{\diff}\calMdaginf $ se fait à l'aide de foncteurs exacts à gauche.
\enddemo

\begin{defi} Le foncteur
image directe
$$ (\calDdag_{(Y\times_kX)\daginf /R},\Sp)\Mod 
\rightarrow (\calDdag_{ \Xdaginf /R},\Sp)\Mod [\dim Y]\leqno p_{*,0}^{\diff}:$$
est un foncteur covariant exact à gauche, et se dérive donc en foncteur exact de catégories triangulées
$$\mathrigid2mu
\Dplus ((\calDdag_{(Y\times_kX)\daginf /R},\Sp)\Mod)\rightarrow \Dplus ((\calDdag_{
\Xdaginf /R},\Sp)\Mod).\leqno p_*^{\diff}:= \bfR p_{*,0}^{\diff} :$$
\end{defi}

\begin{Rema} Si $R	\to S$ est une extension d'anneaux, on peut remplacer les faisceaux $\cal O_{X\daginf /R}$ et $\calDdag_{X\daginf /R}$ dans la définition du foncteur image directe
$p_*^{\diff}$, par les faisceaux $\cal O_{X\daginf /S}:=\cal O_{X\daginf /R}\otimes_RS$ et
$\calDdag_{X\daginf /S}:=\calDdag_{X\daginf /R}\otimes_RS$. Le cas où $R=V$ et $S=K$ est particulièrement important.
\end{Rema}

\subsection{Le foncteur image directe dans le cas général}
Si $f: Y\to X$ est un morphisme de schémas lisses sur $k$ et si $X$ est {\bf séparé},
le morphisme graphe de $f$ est une immersion fermée $i_f: Y\to Y\times_kX.$
Pour simplifier la définition du foncteur image directe dans le cas général,  
nous supposons que $X$ est  {\bf séparé},  de sorte qu'un morphisme se factorise
par une immersion {\bf fermée} suivi d'une projection. 
Si $f$ est un morphisme de $Y\to X$, notons $f:= p\circ i_f$, où $i_f$ l'immersion graphe et où $p$ est la projection de $Y\times_kX$ sur $X$.
\begin{defi}\label{fon-dir}Soit $f: Y\to X$ un morphisme de schémas   lisses  sur $k$, avec  $X$ séparé. On définit le foncteur $f_*^{\diff}$ image directe sur $K$:
$$ \Dplus ((\calDdag_{\Ydaginf /K},\Sp)\Mod)
\longrightarrow \Dplus ((\calDdag_{ \Xdaginf /K},\Sp)\Mod)\leqno f_*^{\diff}:$$ 
comme le foncteur composé
$$\DeuxLignes  
\Dplus ((\calDdag_{\Ydaginf /K},\Sp)\Mod)\stackrel{i_{f*}^{\diff}}{\too} \Dplus ((\calDdag_{(Y\times_kX)\daginf /K},\Sp)\Mod)\stackrel{p_*^{\diff}}{\too}\\
\stackrel{p_*^{\diff}}{\too} \Dplus ((\calDdag_{ \Xdaginf /K},\Sp)\Mod),
\leqno{f_*^{\diff}:}\endlignes
$$ où $p\circ i_f$ est la
factorisation canonique de $f$.
\end{defi}

\begin{Rema} On ne peut pas définir le foncteur image directe par le module de transfert comme dans la théorie classique
parce qu'en général les schémas ne se relèvent pas. D'autre part, le foncteur image directe n'est pas en général  le foncteur dérivé d'un foncteur.
\end{Rema}

Si $U$ est un ouvert  de $Y$, notons $f_U$ la restriction de $f$ à $U$, et si $\calMdaginf $ est un complexe de 
$ \Dplus ((\calDdag_{\Ydaginf /K},\Sp)\Mod)$, notons $\cal H^{i}_{f}(\calMdaginf )$ le préfaisceau sur $Y$ à valeurs dans la catégorie
abélienne des 
$\calDdag_{ \Xdaginf /K}$-modules spéciaux qui à un ouvert $U$ associe le $i$-ème faisceau de cohomologie $ h^i(f^{\diff}_{U_*}\cal
M\daginf )$ du complexe
$f^{\diff}_{U_*}\calMdaginf $.

\begin{theo} Soit $\cal B$ un recouvrement de $Y$ par des ouverts. Il existe alors une suite spectrale dont le terme $\EEE_2$ est $\Rm
H^j(\cal B,\cal H^{i}_{f}(\calMdaginf ))$ et dont le terme $\EEE_\infty$ est le module bigradué associé à une filtration convenable du module  gradué
$h^{i+j}(f_*^{\diff}\calMdaginf )$.
\end{theo}

\demo
Soit $$i_{f*,0}^{\diff}\calMdaginf  = i_{f*}^{\diff}\calMdaginf \to \cal I^{\bullet}_{\inf} $$ une résolution $\calDdag_{(Y\times_kX)\daginf /V}$-injective spéciale de
$i_{f*}^{\diff}\calMdaginf$. Alors, le complexe  $p_{*,0}^{\diff}(\cal I^{\bullet}_{\inf})$ est isomorphe en catégorie dérivée  au complexe $f_*^{\diff}(\calMdaginf )$. 

\bigskip
Pour un ouvert 
$U$ de $Y$, on désignera par 
$p_U$ la projection
$U\times X\to X$.  Soient $\cal B$ un recouvrement de $Y$ par des ouverts et $\cal I_{\inf}$ un $\calDdag_{(Y\times_kX)\daginf /K}$-module  à gauche spécial.  Notons $\cal C^{\bullet}(\cal B, p_{*,0}^{\diff} , \cal I_{\inf})$ le complexe de Cech du recouvrement dont la valeur sur un ouvert $U$ est le $\calDdag_{X\daginf /K}$-module spécial $ p^{\diff}_{U,*,0}\, \cal I_{\inf}$, et notons $\cal C^{\bullet}(\cal B, p_{*,0}^{\diff} , \cal I^{\bullet}_{\inf})$ le bicomplexe de \v Cech du recouvrement~$\cal B$.

\medskip

Il suffit de montrer
que pour tout
$\calDdag_{(Y\times_kX)\daginf /K}$-module {\bf injectif spécial} $\cal I^j_{\inf}$ le complexe $\cal C^{\bullet}(\cal B, p^{\diff}_{*,0} ,\cal
I^j_{\inf})$ est une résolution de $p_{*,0}^{\diff}\cal I^j_{\inf}$. En effet, la suite spectrale du bicomplexe $\cal C^{\bullet}(\cal B, p_{*,0}^{\diff} , \cal I^{\bullet}_{\inf})$ montre alors le théorème.

La question est alors locale sur $X$.
Pour  un ouvert $\calWdag $ de $ \Xdaginf $,  la valeur du  module $p_{*,0}^{\diff}\cal
I^j_{\inf}[-\dim Y]$ sur l'ouvert $\calWdag $ est   :
$$p_{*}\cHom_{q^{-1}\calDdag_{\Ydaginf /K}}
\big(q^{-1}\cal O_{\Ydaginf /K}, R_{\Ydaginf \times\calWdag }(\cal I^j_{\inf})\big),$$ 
et la valeur du  module 
$p_{U*,0}^{\diff}\cal I^j_{\inf}[-\dim Y]$ sur l'ouvert $\calWdag $ est:
$$p_{U*}\cHom_{q^{-1}\calDdag_{\calUdag /K}}(q^{-1}\cal O_{\calUdag /K}, \cal I^j_{\calUdag \times_{V}\calWdag }).$$ 
Le complexe  de \v Cech: 
$$\cal C^{\bullet}\mBig(\cal B\times W,  \cHom_{q^{-1}\calDdag_{\Ydaginf /K}}
\big(q^{-1}\cal O_{\Ydaginf /K}, R_{\Ydaginf \times\calWdag }(\cal I^j_{\inf})\big)\mBig)$$ du recouvrement produit
$\cal B\times W$ à valeurs dans le faisceau de Zariski $\cHom_{q^{-1}\calDdag_{\Ydaginf /K}}
\big(q^{-1}\cal O_{\Ydaginf /K}, R_{\Ydaginf \times\calWdag }(\cal I^j_{\inf})\big)$, dont les termes locaux sont
les faisceaux $\cHom_{q^{-1}\calDdag_{\calUdag /K}}(q^{-1}\cal O_{\calUdag /K}, \cal I^j_{\calUdag \times_{V}\calWdag }),$ est une résolution de ce faisceau:
$$\DeuxLignes
0\to \cHom\Sub{13pt}{q^{-1}\calDdag_{\Ydaginf /K}}
\big(q^{-1}\cal O_{\Ydaginf /K}, R_{\Ydaginf \times\calWdag }(\cal I^j_{\inf})\big)\to \\\to
\cal C^{\bullet}\mBig(\cal B\times W, \cHom\Sub{13pt}{q^{-1}\calDdag_{\Ydaginf /K}}
\big(q^{-1}\cal O_{\Ydaginf /K}, R_{\Ydaginf \times\calWdag }(\cal I^j_{\inf})\big)\mBig).
\endlignes$$ 
Il suffit de voir que l'image directe reste une résolution:
$$\DeuxLignes
0\to p_*\cHom\Sub{13pt}{q^{-1}\calDdag_{\Ydaginf /K}}
\big(q^{-1}\cal O_{\Ydaginf /K}, R_{\Ydaginf \times\calWdag }(\cal I^j_{\inf})\big)\to \\\to
p_*\cal C^{\bullet}\mBig(\cal B\times W, 
\cHom\Sub{13pt}{q^{-1}\calDdag_{\Ydaginf /K}}
\big(q^{-1}\cal O_{\Ydaginf /K}, R_{\Ydaginf \times\calWdag }(\cal I^j_{\inf})\big)\mBig),
\endlignes
$$
qui, par construction, est  le complexe :
$$0\to p_{*,0}^{\diff}\cal I^j_{\inf}\to \cal C^{\bullet}(\cal B, p^{\diff}_{*,0} ,\cal
I^j_{\inf}).$$

Or, l'extension 
$q^{-1}\calDdag_{\calUdag /K}\to \calDdag_{\cal
U\dag\times_{V}\calWdag /K}$ est {\bf plate}, comme on peut le voir sur $V$ en considérant la filtration par les échelons et en appliquant le critère de platitude locale \ref{cri-pla}. Le $q^{-1}\calDdag_{\calUdag /K}$-module $\cal I^j_{\calUdag \times_{V}\calWdag }$ reste donc  {\bf injectif}. Il en
résulte que  le faisceau:
$$\cHom_{q^{-1}\calDdag_{\calUdag /K}}(q^{-1}\cal O_{\calUdag /K}, \cal I^j_{\calUdag \times_{V}\calWdag })$$ est {\bf flasque}. La résolution de \v Cech du
faisceau flasque:
$$\cHom_{q^{-1}\calDdag_{\Ydaginf /K}}\big(q^{-1}\cal O_{\Ydaginf /K}, R_{\Ydaginf \times\calWdag }(\cal I^j_{\inf})\big)$$ est une résolution par des faisceaux
flasques, et l'image  directe est une résolution.
\enddemo

Cette suite spectrale  ramène souvent les propriétés de l'image directe au cas où $Y$ est affine, ce qui est très utile dans la pratique.

\begin{defi}\label{tra-dro}Soit 
$f: Y\to X$ un morphisme de schémas lisses sur $k$, et soient $\calYdag $ et $\calXdag $ des relèvements  plats de
$Y$ et
$X$. On définit le module de transfert
$\calDdag_{\calXdag \gets \calYdag /V}$ comme le module de transfert induit par un relèvement local à la source et au but du morphisme, ce qui est légitime puisqu'en vertu du théorème \ref{ind-tra'} le module de transfert ne dépend pas du relèvement. On définit le module de transfert sur $K$ par: $$\calDdag_{\calXdag \gets
\calYdag /K} :=
\calDdag_{\calXdag \gets \calYdag /V}\otimes_VK.$$
\end{defi}

Nous allons étudier le rapport entre les foncteurs images directes dans le cas des relèvements des schémas. 

\begin{theo}\label{ima-com}Soit $f: Y\rightarrow X$ un morphisme de schémas 
lisses sur $k$, où $X$ est séparé. Soient $\calYdag $ et $\calXdag $ des
relèvements  plats, alors le diagramme suivant est {\bf commutatif}:
$$\def\quad{\hskip5pt}\matrix{f_*^{\diff}:=f_{*, \calYdag , \calXdag }^{\diff}:& \Db ((\calDdag_{\calYdag /K})\Mod)&\longrightarrow& \ \Db ((\calDdag_{ \calXdag /K})\Mod)\cr
&\downarrow P_{\calYdag }&&\downarrow P_{\calXdag }
\cr f_*^{\diff}:& \Db ((\calDdag_{\Ydaginf /K}, \Sp)\Mod)&\longrightarrow& \Db ((\calDdag_{\Xdaginf /K}, \Sp)\Mod),}$$\glossary{$f^*_{\diff}$}où le foncteur de la première ligne est le foncteur défini par le module de transfert:
$$\cal M\dag_{\calYdag }\fonct \bfR f_*(\calDdag_{\calXdag \gets \calYdag/K}\Lotimes_{\calDdag_{\calYdag /K}}\cal M\dag_{\calYdag })$$ et où les flèches
verticales sont les foncteurs de  prolongement canonique.
\end{theo}

\demo
Le foncteur de la première ligne est le foncteur: 
$$\cal M\dag_{\calYdag }\fonct \bfR f_*(\calDdag_{\calXdag \gets \calYdag/K}\Lotimes_{\calDdag_{\calYdag /K}}\cal M\dag_{\calYdag }),$$ ce qui nécessite que $\cal M\dag_{\calYdag }$ soit cohomologiquement borné. 
Dans le cas d'une immersion fermée, le théorème
est vrai par construction; même sur
$V$. Le cas d'une projection est plus délicat.  Le produit fibré
$\calYdag \times_{V}\calXdag $ existe et est un relèvement lisse de $Y\times_kX$ ([A-M$_2$]).

\begin{lemm}\label{tra-drt}Il existe un isomorphisme canonique:
$$q^{-1}\omega_{\calYdag /V}\otimes_{q^{-1}\calDdag_{\calYdag /V}}\calDdag_{\calYdag \times_{V}\calXdag /V}\to\calDdag_{\calXdag \gets \cal
Y\dag\times_{V}\calXdag /V}\,,$$
de $(p^{-1}\calDdag_{\calXdag /V},\calDdag_{\calYdag \times_{V}\calXdag/V})$-bimodules.
\end{lemm}

\demo
On définit le morphisme de $(p^{-1}\calDdag_{\calXdag /V},\calDdag_{\calYdag \times_{V}\calXdag /V})$-bimo\-dules:
$$q^{-1}\omega_{\calYdag /V}\otimes_{q^{-1}\calDdag_{\calYdag /V}}\calDdag_{\calYdag \times_{V}\calXdag /V}\to \cHom_{V}(p^{-1}\omega_{\calXdag /V},
\omega_{\calYdag \times_{V}\calXdag /V})$$ 
par $(dy\otimes P)\mapsto\big(dx\to (dx\otimes dy)P\big)$. 
Par réduction modulo $\goth m^s$,  on obtient un
opérateur différentiel dont le degré est localement  borné par une fonction linéaire en $s$. D'où le morphisme:
$$q^{-1}\omega_{\calYdag /V}\otimes_{q^{-1}\calDdag_{\calYdag /V}}\calDdag_{\calYdag \times_{V}\calXdag /V}\to\calDdag_{\calXdag \gets \cal
Y\dag\times_{V}\calXdag /V}.\leqno(*):$$ Le module $\calDdag_{\calXdag \gets \cal
Y\dag\times_{V}\calXdag /V}$ est un $\calDdag_{\cal
Y\dag\times_{V}\calXdag /V}$-module  à droite de type fini, et le conoyau du morphisme $(*)$ est donc  nul en vertu du lemme de Nakayama, parce que de type fini   
et de réduction modulo $\goth m$ nulle. Le morphisme $(*)$ est filtré par: 
$$q^{-1}\omega_{\calYdag /V}\otimes_{q^{-1}\calDdag_{\cal Y^{\dag,h}/V}}\calD^{\dagger,h}_{\calYdag \times_{V}\calXdag /V}\to\calD^{\dagger,h}_{\calXdag \gets
\calYdag \times_{V}\calXdag /V}.\leqno(*)_h:$$ Le module $\calD^{\dagger,h}_{\calXdag \gets \cal
Y\dag\times_{V}\calXdag /V}$ est un $\calD^{\dagger,h}_{\cal
Y\dag\times_{V}\calXdag /V}$-module  à droite de type fini, et le conoyau du morphisme $(*)_h$ est donc  nul en vertu du lemme de Nakayama, parce que de type
fini    et de réduction modulo $\goth m$ nulle. Le module $\calD^{\dagger,h}_{\calXdag \gets \cal
Y\dag\times_{V}\calXdag /V}$ est sans $\goth m$-torsion et le noyau du morphisme $(*)_h$ est nul modulo $\goth m^s$ pour tout $s\geq 1$.
Comme
$q^{-1}\omega_{\calYdag /V}\otimes_{q^{-1}\calD_{\cal Y^{\dag,h}/V}}\calD^{\dagger,h}_{\calYdag \times_{V}\calXdag /V}$ est séparé pour
topologie $\goth m$-adique comme module de type fini sur un anneau de Zariski (non commutatif), il en résulte que le noyau de $(*)_h$ est nul pour tout $h\geq0$, ce qui
entraîne que le noyau du morphisme  $(*)$ est nul.
\endsubdemo

En particulier, le lemme entraîne l'isomorphisme canonique:
$$q^{-1}\omega_{\calYdag /K}\otimes_{q^{-1}\calD_{\calYdag /K}}\calDdag_{\calYdag \times_{V}\calXdag /K}\simeq\calDdag_{\calXdag \gets \cal
Y\dag\times_{V}\calXdag /K}\,,$$ 
 de  $(p^{-1}\calDdag_{\calXdag /K},\calDdag_{\calYdag \times_{V}\calXdag /K})$-bimodules.
Mais en vertu du lemme \ref{Spen}, 
il existe un isomorphisme canonique  dans la catégorie dérivée des $q^{-1}\calD_{\calYdag /K}$-modules
à droite:
$$\bfR \cHom_{q^{-1}\calDdag_{\calYdag /K}}(q^{-1}\cal O_{\calYdag /K},q^{-1}\calDdag_{\calYdag /K})[\dim Y]
\simeq q^{-1}\omega_{\calYdag /K}.$$ 
Le fibré   $q^{-1}\cal O_{\calYdag /K}$ est un $q^{-1}\calDdag_{\calYdag /K}$-module à gauche parfait, en vertu du lemme \ref{Spen}.  Pour tout complexe $\calMdag $ de $ \Dplus \big((\calDdag_{\calYdag \times_{V}\calXdag/K})\Mod\big)$ on trouve des isomorphismes canoniques   dans la catégorie dérivée de la catégorie des
 $p^{-1}\calDdag_{\calXdag /K}$-modules à gauche:
$$
\DeuxLignes  
\bfR \cHom_{q^{-1}\calDdag_{\calYdag /K}}(q^{-1}\cal O_{\calYdag /K},\calMdag )[\dim Y]\simeq 
q^{-1}\omega_{\cal
Y\dag/K}
\smashtop{\Lotimes_{q^{-1}\calDdag_{\calYdag /K}}}\calMdag \simeq 
\\
\simeq\calDdag_{\calXdag \gets \cal
Y\dag\times_{V}\calXdag /K}\Lotimes_{\calDdag_{\calYdag \times_{V}\calXdag/K}}\calMdag .
\endlignes
$$  
En appliquant cet isomorphisme à une résolution injective $\calIdag $ de $\calMdag $, on trouve effectivement que
$p_*^{\diff}P_{\calYdag \times_{V}\calXdag }(\calMdag )$ est canoniquement isomorphe à $P_{\calXdag }\big(\bfR p_*(\calDdag_{\calXdag \gets \cal
Y\dag\times_{V}\calXdag /K}\Lotimes_{\calDdag_{\calYdag \times_{V}\calXdag/K}}\calMdag )\big).$

\bigskip
Dans le cas général, le théorème résulte du fait que l'image directe pour les modules spéciaux est, par définition, la composée du cas d'une  immersion fermée suivie du cas de
la projection, et de l'isomorphisme: 
$$i^{-1}_f\calDdag_{\calXdag \gets \cal
Y\dag\times_{V}\calXdag /K}\Lotimes_{i^{-1}_f\calDdag_{\cal
Y\dag\times_{V}\calXdag /K}}\calDdag_{ \cal
Y\dag\times_{V}\calXdag \gets \calYdag /K}\simeq\calDdag_{\calXdag \gets \cal
Y\dag/K}$$ conséquence du théorème \ref{acy-pro}.
\enddemo

\bigskip
Le foncteur image directe dans la catégorie des modules spéciaux sur $K$ sur le site infinitésimal $p$-adique prolonge,  comme il se doit, le foncteur image directe lorsque
les relèvements existent. Il s'applique aux morphismes de variétés qui se relèvent, courbes lisses, variétés abéliennes, grassmanniennes, etc. sans que les morphismes ne se relèvent.
\begin{Rema} C'est bien pour faire coïncider les deux foncteurs d'image directe que l'on a mis le décalage $[\dim Y]$ dans la définition
du foncteur $p_{*,0}^{\diff}$. 
\end{Rema}

\medskip
Il nous faut maintenant étudier  l'indépendance du foncteur  image directe $f_*^{\diff}$ de la factorisation du morphisme $f$. 

\begin{theo}\label{com-i-p}Soit une composition de morphismes de schémas lisses sur $k$: $Z\stackrel{i}{\to} Y\times_k X\stackrel{p}{\to} X$, où $i$ est une
immersion  fermée
et $p$ une projection. Si nous supposons que $p\circ i$ est une {\bf immersion  fermée}, alors le foncteur $p^{\diff}_{*,0}\circ i^{\diff}_{*,0}$ est {\bf exact}, et  il existe un isomorphisme canonique de foncteurs de la catégorie
$ (\calDdag_{Z\daginf /K},\Sp)\Mod$:
$$p_{*,0}^{\diff}\circ i_{*,0}^{\diff}\longrightarrow (p\circ i)^{\diff}_{*,0}$$ et  un isomorphisme canonique de foncteurs de 
$ \Db ((\calDdag_{Z\daginf /K},\Sp)\Mod))$:
$$p_*^{\diff}\circ i_*^{\diff}\longrightarrow (p\circ i)^{\diff}_*.$$
\end{theo}

\demo
Soit $z$ un point  de $Z$ et soit   $T$ un voisinage affine de $q(i(z))$, 
où $q$ est la projection de $Y\times_kX$ sur $Y$. 
Soit  $U$ un ouvert affine  voisinage de $p\circ i(z)$ dans $X$. Puisque  $p\circ i$ est  par hypothèse une {\bf immersion fermée},  quitte à remplacer
$U$ par un ouvert  affine plus petit, on  peut supposer que sa trace $W$ sur $Z$ est contenue dans l'image inverse de $T$
par $q\circ i$.  

Si $\calMdaginf $  est un 
$ \calDdag_{Z\daginf /K}$-module  spécial,  les restrictions à $U$ des complexes $p_*^{\diff}\circ i_*^{\diff}(\calMdaginf )$ et $(p\circ i)^{\diff}_*(\calMdaginf )$ ne dépendent par construction  que de la restriction
de $\calMdaginf $ à $W$ et du morphisme
$W\stackrel{i}{\to} T\times_k U\stackrel{p}{\to} U$, où $i$ est 
l'immersion  fermée induite 
et $p$ la projection induite.
Mais $W$, $ T $ et $ U$ sont affines et admettent des relèvements plats  $\calWdag , \calTdag , \calUdag $.
En vertu du théorème \ref{ima-com}, on a les isomorphismes canoniques:
\begingroup\arraycolsep2pt\mathrigid2mu\scriptwd1.3em
\begin{eqnarray*}
(p\circ i)^{\diff}_*(\cal M\dag_{\cal W})&\simeq& (p\circ i)_*(\calDdag_{\calUdag \gets \cal
W\dag/K}\otimes_{\calDdag_{\calWdag /K}}\cal M\dag_{\cal W}),
\\
i^{\diff}_*(\cal M\dag_{\cal W})&\simeq& i_*(\calDdag_{ \cal
T\dag\times_{V}\calUdag \gets \calWdag /K}\otimes_{\calDdag_{\calWdag /K}}\cal M\dag_{\cal W}),
\\
p^{\diff}_*\circ i^{\diff}_*(\cal M\dag_{\cal W})&\simeq& p_*\mBig(\calDdag_{\calUdag \gets \cal
T\dag\dtimes_{V}\calUdag /K}\LOtimes_{\calDdag_{\cal
T\dag\times_{V}\calUdag /K}}i_*\big(\calDdag_{ \cal
T\dag\dtimes_{V}\calUdag \gets \calWdag /K}\Otimes_{\calDdag_{\calWdag /K}}\cal M\dag_{\cal W}\big)\mBig).
\end{eqnarray*}
\endgroup

En vertu  du théorème \ref{acy-pro}, on a un isomorphisme canonique:
$$\DeuxLignes  
i^{-1}\calDdag_{\calUdag \gets \cal
T\dag\times_{V}\calUdag /K}\LOtimes_{i^{-1}\calDdag_{\cal
T\dag\times_{V}\calUdag /K}}\calDdag_{ \cal
T\dag\times_{V}\calUdag \gets \calWdag /K}\simeq
\\
\vrule height15pt width0pt
\simeq i^{-1}\calDdag_{\calUdag \gets \cal
T\dag\times_{V}\calUdag /K}\Otimes_{i^{-1}\calDdag_{\cal
T\dag\times_{V}\calUdag /K}}\calDdag_{ \cal
T\dag\times_{V}\calUdag \gets \calWdag /K}\simeq\calDdag_{\calUdag \gets \cal
W\dag/K}.
\endlignes
$$  
Cela montre que les restrictions à $\calUdag $  des complexes: $$p_*^{\diff}\circ i_*^{\diff}(\calMdaginf )  \quad {\rm et}\quad
(p\circ i)_*^{\diff}(\calMdaginf )$$  sont isomorphes. Mais le foncteur $(p\circ i)_{*,0}^{\diff}$ est exact, puisque par hypothèse
l'immersion $p\circ i$ est fermée, et donc le complexe
$p_{*}^{\diff}\circ i_{*}^{\diff}(\calMdaginf )$ est concentré cohomologiquement en degré zéro.
Le foncteur $p_{*,0}^{\diff}\circ i_{*,0}^{\diff}$ est alors  exact. 

Les $\calDdag_{X\daginf /K}$-modules à gauche spéciaux:
$$p_{*,0}^{\diff}\circ i_{*,0}^{\diff}(\calMdaginf )\quad \text{ et }\quad (p\circ i)_{*,0}^{\diff}(\calMdaginf )$$  sont localement isomorphes.
Il reste à voir que ces isomorphismes locaux sont compatibles. 

\medskip
Mais si $\cal T'$ est un autre relèvement de $T$,
les modules: $$i^{-1}\calDdag_{\calUdag \gets \cal
T\dag\times_{V}\calUdag /K}\otimes_{i^{-1}\calDdag_{\cal
T\dag\times_{V}\calUdag /K}}\calDdag_{ \cal
T\dag\times_{V}\calUdag \gets \calWdag /K}
\preskip2pt\postskip-5pt$$ et $$\preskip0pt
i^{-1}\calDdag_{\calUdag \gets \cal
T'{}\dag \times_{V}\calUdag /K}\otimes_{i^{-1}\calDdag_{\cal
T'{}\dag \times_{V}\calUdag /K}}\calDdag_{ \cal
T'{}\dag \times_{V}\calUdag \gets \calWdag /K}$$ sont canoniquement isomorphes. 
Cela montre que l'isomorphisme local entre les modules $p_{*,0}^{\diff}\circ i_{*,0}^{\diff}(\calMdaginf )$ et $(p\circ i)_{*,0}^{\diff}(\calMdaginf )$ ne dépend pas du relèvement de $T$, puis par construction, ne dépend pas des relèvements  de 
$U$ et de $W$. 

\medskip
Les $\calDdag_{X\daginf /K}$-modules $p_{*,0}^{\diff}\circ i_{*,0}^{\diff}(\calMdaginf )$ et $(p\circ i)_{*,0}^{\diff}(\calMdaginf )$
sont canoniquement isomorphes, les  foncteurs $p_{*,0}^{\diff}\circ i_{*,0}^{\diff}$ et $(p\circ i)_{*,0}^{\diff}$ sont canoniquement isomorphes et les foncteurs $p_{*}^{\diff}\circ i_{*}^{\diff}$ et $(p\circ i)_{*}^{\diff}$ sont aussi canoniquement isomorphes comme foncteurs dérivés de foncteurs exacts canoniquement isomorphes.
\enddemo

\begin{prop} Soit $ Y\times_k( X\times _kZ)\stackrel{p_{1}}{\too} X\times_k Z\stackrel{p_{2}}{\too} Z$ une composition de deux projections de schémas lisses sur $k$. Il existe alors un isomorphisme canonique de foncteurs: $p_{2*}^{\diff}\circ p_{1*}^{\diff}\simeq (p_{2}\circ p_{1})^{\diff}_*$ de $ \Dplus\!\big((\calDdag_{ (Y\times_kX\times_kZ)\daginf /K},\Sp)\Mod\big)$ vers $ \Dplus\!\big((\calDdag_{ Z\daginf /K},\Sp)\Mod\big)$.
\end{prop}

\demo
Voyons d'abord que 
la première projection  transforme module à gauche spécial injectif en  module à gauche spécial {\bf acyclique} pour la seconde projection. La question étant locale sur $Z$, on peut supposer qu'il existe un relèvement plat  $\cal Z\dag$ de $Z$. Supposons alors, dans un premier temps, qu'il existe aussi des relèvements  plats   $\calYdag$ et  $\calXdag$ de $Y$ et $X$ respectivement. Le morphisme canonique:
$$\DeuxLignes  
p_1^{-1}\calDdag_{\calZdag \get\calXdag \times_V \calZdag /K}\LOtimes_{p_1^{-1}\calDdag_{\calXdag \times_V \calZdag /K}}\calDdag_{\calXdag \times_V \calZdag \get\calYdag \times_V\calXdag \times_V \calZdag /K}\to \\
\vrule height15pt width0pt\to
p_1^{-1}\calDdag_{\calZdag \get\calXdag \times_V \calZdag /K}\Otimes_{p_1^{-1}\calDdag_{\calXdag \times_V \calZdag /K}}\calDdag_{\calXdag \times_V \calZdag \get\calYdag \times_V\calXdag \times_V \calZdag /K}
\endlignes$$
est un isomorphisme, parce que 
le $p_1^{-1}\calDdag_{\calXdag \times_V \calZdag /K}$-module à gauche:
$$\calDdag_{\calXdag \times_V \calZdag \get\calYdag \times_V\calXdag \times_V \calZdag /K}$$ est plat.
Le morphisme canonique:
$$ \let\times\dtimes
p_1^{-1}\calDdag_{\calZdag \get\calXdag \times_V \calZdag /K}\Otimes_{p_1^{-1}\calDdag_{\calXdag \times_V \calZdag /K}}\calDdag_{\calXdag \times_V \calZdag \get\calYdag \times_V\calXdag \times_V \calZdag /K}\to\calDdag_{ \calZdag \get\calYdag \times_V\calXdag \times_V \calZdag /K}$$ est un isomorphisme
de $\big((p_2\circ p_1)^{-1}\calDdag_{ \calZdag /K}, \calDdag_{ \calYdag \times_V\calXdag \times_V \calZdag /K}\big)$-bimodules, parce que $p_1$ et $p_2$ sont des projections. Cela entraîne que pour tout
complexe $\calMdag $ de $ \Db (\calDdag_{ \calYdag \times_V\calXdag \times_V \calZdag /K}\Mod)$, le morphisme canonique:
$$\let\times\dtimes\deuxlignes{0pt}{4pt}
\myeqno{0.7cm}(*): \mBig(p_1^{-1}\calDdag_{\calZdag \get\calXdag \times_V \calZdag /K}\LOtimes_{p_1^{-1}\calDdag_{\calXdag \times_V \calZdag /K}}\calDdag_{\calXdag \times_V \calZdag \get\calYdag \times_V\calXdag \times_V \calZdag /K}\mBig)
\LOtimes_{\calDdag_{\calYdag \times_V\calXdag \times_V \calZdag /K}}\calMdag \to
\\
\to\calDdag_{ \calZdag \get\calYdag \times_V\calXdag \times_V \calZdag /K}\LOtimes_{\calDdag_{\calYdag \times_V\calXdag \times_V \calZdag /K}}\calMdag 
\endlignes$$ est un isomorphisme dans la catégorie dérivée de 
$(p_2\circ p_1)^{-1}\calDdag_{ \calZdag /K}$-modules à gauche. L'isomorphisme $(*)$ donne l'isomorphisme dans la catégorie dérivée de la catégorie des $(p_2\circ p_1)^{-1}\calDdag_{ \calZdag /K}$-modules à gauche  :
$$
\DeuxLignes  
\myeqno{1cm}(**):\bfR \cHom\goodSub{5pt}{8pt}{-7mm}{(q\circ p_1)^{-1}\calDdag_{ \calXdag /K}}\mBig((q\circ p_1)^{-1}\cal O_{ \calXdag /K}, \bfR \cHom\SubX{q^{-1}\calDdag_{ \calYdag /K}}\big(q^{-1}\cal O_{ \calYdag /K}, \calMdag \big)\mBig)\to
\\
\to\bfR \cHom_{q^{-1}\calDdag_{ \calYdag \times_V\calXdag /K}}\big(q^{-1}\cal O_{\calYdag \times_V \calXdag /K},\calMdag  \big),\endlignes$$   parce que le faisceau structural  $\cal O_{\calTdag /K}$ d'un schéma $\dagger$-adique lisse  $\calTdag $ est un 
  $\calDdag_{\calTdag /K}$-module à gauche  parfait en vertu du lemme \ref{Spen}, où l'on désigne
par $q$ la projection sur le premier facteur.
En appliquant le foncteur $\bfR p_{1_*}$, on trouve un isomorphisme dans la catégorie dérivée de la catégorie  des  $p_2^{-1}\calDdag_{ \calZdag /K}$-modules à gauche: 
$$\DeuxLignes  
\myeqno{1.1cm}(\*\*\*):\bfR \cHom\goodSub{5pt}{8pt}{-5mm}{q^{-1}\calDdag_{ \calXdag /K}}\mBig(q^{-1}\cal O_{ \calXdag /K}, \bfR p_{1_*}\bfR \cHom\SubX{q^{-1}\calDdag_{ \calYdag /K}}\big(q^{-1}\cal O_{ \calYdag /K}, \calMdag\big)\mBig)\to
\\
\to\bfR p_{1_*}\bfR \cHom_{q^{-1}\calDdag_{ \calYdag \times_V\calXdag /K}}\big(q^{-1}\cal O_{\calYdag \times_V \calXdag /K},\calMdag  \big),\endlignes$$    parce que le faisceau structural $ \cal O_{\calTdag /K}$  d'un schéma $\dagger$-adique lisse  $\calTdag $ est un $\calDdag_{\calTdag /K}$-module à gauche  parfait.

\medskip
Si $\calIdag $ est un $\calDdag_{\calYdag \times_V\calXdag \times_V \calZdag /K}$-module à gauche  supposé  injectif,   il reste injectif  comme $q^ {-1}\calDdag_{ \calYdag \times_V\calXdag /K}$-modules à gauche
et comme $q^ {-1}\calDdag_{ \calYdag /K}$-module  à gauche. L'isomorphisme $(\*\*\*)$
devient de fait un isomorphisme  de $p_2^{-1}\calDdag_{ \calZdag /K}$-modules à gauche:
$$\DeuxLignes  
\myeqno{1.4cm}(\*\*\*\*):\bfR \cHom\goodSub{5pt}{8pt}{-7mm}{q^{-1}\calDdag_{ \calXdag /K}}\mBig(q^{-1}\cal O_{ \calXdag /K}, p_{1_*}\cHom\SubX{q^{-1}\calDdag_{ \calYdag /K}}\big(q^{-1}\cal O_{ \calYdag /K}, \calIdag \big)\mBig)\simeq
\\
\simeq p_{1_*}\cHom_{q^{-1}\calDdag_{ \calYdag \times_V\calXdag /K}}\big(q^{-1}\cal O_{\calYdag \times_V \calXdag /K},\calIdag  \big),\endlignes$$ parce que   le faisceau $\cHom_{q^{-1}\calDdag_{ \calYdag \times_V\calXdag /K}}(q^{-1}\cal O_{\calYdag \times_V \calXdag /K},\calIdag  )$ est flasque
et donc le faisceau $p_{1_*}\cHom_{q^{-1}\calDdag_{ \calYdag \times_V\calXdag /K}}(q^{-1}\cal O_{\calYdag \times_V \calXdag /K},\calIdag  )$ est aussi flasque. 

En appliquant le foncteur $\bfR p_{2_*}$ à la formule $(\*\*\*\*)$, on voit que  $p^{\diff}_{1_{*,0}}\calIdag $ est acyclique
pour le foncteur $p^{\diff}_{2_{*,0}}$ et, en même temps, que $p^{\diff}_{2_{*,0}}\circ p^{\diff}_{1_{*,0}}(\calIdag  )$ 
et $(p_{2_{*,0}}\circ p_{1_{*,0}})^{\diff}(\calIdag  )$ sont canoniquement isomorphes.

\medskip
L'isomorphisme $(\*\*\*\*)$  se recolle naturellement en un isomorphe  de $p_2^{-1}\calDdag_{ \calZdag /K}$-modules à gauche pour tout $Y$, qu'il admette ou non de relèvement plat global $\calYdag$.

\medskip
Dans le cas général, pour tous $Y$ et $X$,  si $\calIdaginf  $
est un 
$\calDdag_{(Y\times_kX\times_k Z)\daginf /K}$-module à gauche spécial et injectif, la méthode précédente montre,  que l'isomorphisme
$(\*\*\*\*)$ se recolle en un isomorphisme de complexes de $p_2^{-1}\calDdag_{ \calZdag /K}$-modules à gauche: 
$$\let\times\dtimes\mathrigid1mu\scriptspace0.5pt
\DeuxLignes  
\myeqno{1.7cm}(\*\*\*\*)_{\inf}:
\bfR \cHom\SubX{q^{-1}\calDdag_{ \Xdaginf /K}}\mBig(q^{-1}\cal O_{ \Xdaginf /K},  R_{X\daginf \times\calZdag}\big(p^{\diff}_{1_*}(\calIdaginf)\big)\mBig)\\
\vrule height12pt width0pt
\simeq p_{1_*}\cHom_{q^{-1}\calDdag_{(Y\times_kX)\daginf/K}}
\mBig(q^{-1}\cal O_{(Y\times_kX)\daginf/K}, R_{(Y\times_kX)\daginf \times\calZdag }(\calIdaginf  ) 
\mBig)
\endlignes$$
avec les notations 	\ref{nota}. Cet isomorphisme implique que $p^{\diff}_{1_{*,0}}(\calIdaginf  )$ est acyclique
pour le foncteur $p^{\diff}_{2_{*,0}}$.

\bigskip
En vertu de ce qui précède, $p^{\diff}_{2_*}(p^{\diff}_{1_*}(\calIdaginf  ))$ et $p^{\diff}_{2_*}\circ p^{\diff}_{1_*}(\calIdaginf  )$  sont deux $\calDdag_{Z\daginf /K}$-modules à gauche spéciaux en degré $-[\dim Y+\dim X]$, et  il suffit donc de montrer qu'ils sont canoniquement isomorphes. Par construction,  la restriction de $p^{\diff}_{2_{*,0}}(p^{\diff}_{1_{*,0}}(\calIdaginf  ))[-\dim Y-\dim X]$ à l'ouvert $\calZdag $ est le faisceau:
$$p_{2_*}\cHom_{q^{-1}\calDdag_{  \Xdaginf /K}}(q^{-1}\cal O_{\Xdaginf /K}, R_{\Xdaginf \times\calZdag }(p_{1*}^{\diff}\calIdaginf  ) )[-\dim Y],$$ alors que  la restriction à l'ouvert $\calZdag $ du faisceau 
$$(p_{2}\circ p_{1})_{*,0}^{\diff}(\calIdaginf  )[-\dim Y-\dim X]$$ est le faisceau:
$$(p_{2}\circ p_{1})_{*}\cHom_{q^{-1}\calDdag_{ (Y\times_k X)\daginf /K}}\big(q^{-1}\cal O_{ (Y\times_k X)\daginf /K}, R_{(Y\times_k X)\daginf \times\calZdag }(\calIdaginf  )\big).$$
Ces deux 
$\calDdag_{\calZdag /K}$-modules à gauche sont  des images directes de $p_2^{-1}\calDdag_{\calZdag /K}$-modules à gauche donnés localement  qui sont canoniquement isomorphes.\enddemo

\begin{coro} 
Soit  un morphisme $f:Y\to X$ de schémas lisses sur $k$ où $X$ est séparé, et soit  $ Y\stackrel{i}{\to} P\times_kX\stackrel{p}{\to} X$ une factorisation  
de $f$  par une immersion fermée suivie d'une projection, où $P$ est un schéma lisse sur $k$, et telle que le morphisme induit $Y\to Y\times_kP\times_kX$ est une immersion fermée. Alors, le foncteur 
 $f^{\diff}_*$ est canoniquement isomorphe au foncteur $p^{\diff}_*\circ i^{\diff}_*$.
 \end{coro}

\demo  On a  le diagramme commutatif :
$$\matrix{ &\decale{25pt}{Y\times_kX}\cr
\decale{-25pt}\nearrow&\searrow\cr 
Y\to Y\times_kP\times_kX\ \rlap{\hbox to1.3cm{\rightarrowfill}}\kern4mm &&\decale{25pt}{X,}\cr
\decale{-25pt}\searrow&\nearrow\cr 
&\decale{25pt}{P\times_kX}}$$ 
où le composé de deux
morphismes non tous deux horizontaux est soit un composé de deux projections soit un composé d'une immersion  fermée et d'une projection qui est une immersion  fermée. 
On est réduit aux résultats
précédents.

En particulier, si $X$ est séparé, pour une immersion fermée $Y\to X$,  le foncteur  image directe
est  canoniquement isomorphe au foncteur construit à l'aide de la factorisation canonique  et pour une projection $Y\times_kX\to X$  le foncteur  image directe
est  aussi canoniquement isomorphe au foncteur construit à l'aide de la factorisation canonique.
\enddemo

\begin{exemples}\begin{liste}
\item 1) Si $X$ est le point $\Spec(k)$, le complexe $f^{\diff}_*\cal O_{\Ydaginf /K}$\glossary{$f^{\diff}_*\cal
O_{\Ydaginf /K}$} est un complexe de $ \Dplus (K)
$ qui est à cohomologie de dimension finie sur $K$ lorsque $Y$ est une variété algébrique lisse en vertu du théorème de finitude \ref{fin}. Ce théorème a présenté de sérieuses difficultés pendant près d'un quart de siècle dans le cas d'une variété qui se relève [Me$_2$], qui est aussi à la base des progrès significatifs de la théorie $p$-adique.
\item 2) Plus généralement,  pour {\bf tout morphisme} $f$  de schémas  de type fini et lisses sur $k$, le complexe spécial $f^{\diff}_*\cal O_{\Ydaginf /K}$ sur la base $X$ est l'analogue $p$-adique du complexe de Gauss-Manin en caractéristique nulle, dont
on s'attend naturellement à ses  propriétés de finitude. Cela semble de nouveau présenter de sérieuses difficultés, sauf qu'aujourd'hui on est mieux équipé. Un cas intéressant  est celui d'un morphisme propre et lisse sur une courbe, où l'on
s'attend à ce que les images directes soient des fibrés $p$-adiques dont les rangs donnent les nombres de Betti $p$-adiques des fibres. 
\end{liste}
\end{exemples}

\begin{Remas}\begin{liste}\item 1) Sous la condition $e<p^h(p-1)$, on peut reprendre les constructions précédentes pour construire les foncteurs
$f^{\diff, h}_{*}$ sur $K$.
\item 2) Le module $q^{-1}\cal O_{\calYdag /V}$ n'est pas un $q^{-1}\calDdag_{\calYdag /V}$-module à gauche de présentation
finie, parce que l'idéal d'augmentation $(\Delta^\alpha, \alpha\in \Bbb N^n)$ n'est pas de type fini,  et on ne peut pas appliquer le raisonnement du théorème \ref{ima-com} sur $V$. 

\item 3) Mais dans le cas $h=0$, on peut si $e<p-1$  reprendre les constructions précédentes pour définir  les foncteurs
$f^{\diff,0}_{*}$ sur $V$, et 
on peut appliquer le raisonnement du théorème \ref{ima-com}, sur $V$,  pour montrer que 
$$\resetdisplay
R_{\calXdag }(f_{*}^{\diff,0}\cal M^{\dagger,0}_{\inf})\simeq \bfR f_*\big(\calD^{\dagger,0}_{\calXdag \gets \calYdag /V}\Lotimes_{\calD^{\dagger, 0}_{\calYdag /V}}R_{\calYdag }(\cal M^{\dagger,0}_{\inf})\big)$$ est un isomorphisme pour  des relèvements  de
$Y$ et de $X$.
\end{liste}
\end{Remas}

\subsection{Le foncteur image directe pour les modules à droite spéciaux}
En remplaçant le module de transfert $\calDdag_{\calXdag \leftarrow\calYdag /V}$ par $\calDdag_{\calYdag \rightarrow\calXdag /V}$, on définit
le foncteur $i_{*,0}^{\diff}$ pour la catégorie $ \Modd(\calDdag_{\Ydaginf /V}, \Sp)$  et, en remplaçant le $\calDdag_{Z\daginf /V}$-module à gauche spécial
$\cal O_{Z\daginf /V}$ par  le $\calDdag_{Z\daginf /V}$-module à droite spécial
$\omega_{Z\daginf /V}$, on définit le foncteur $p_{*,0}^{\diff}$ pour la catégorie $ \Modd(\calDdag_{(Y\times_kX)\daginf/V }, \Sp)$. On dispose donc, sur $K$ et quand le but est séparé, du
foncteur image direct $f_*^{\diff}$ pour les catégories de modules à droite spéciaux, qui a des propriétés parallèles à celles du   foncteur image directe des
modules à gauche.

\section{Le foncteur de cohomologie locale et ses compatibilités avec les foncteurs images directe et inverse}

\subsection{Le foncteur de cohomologie locale dans le cas des faisceaux}
Soit $X := (X, \cal O_{X/R_1})$ un schéma  lisse sur $R_1$,  et soit $i:Z\hookrightarrow X$ un fermé de l'espace topologique $X$, de complémentaire $j: X-Z\hookrightarrow X$.

\begin{prop}\label{suite-coh-loc}Soit $\cal F_{\inf}$ un faisceau de $R_{\Xdaginf }$-modules sur le site $ \Xdaginf $. Les données qui à un ouvert $\calUdag $
associent le faisceau
$\Gamma_Z(\cal F_{\calUdag })$ et le faisceau $j_*j^{-1}\cal F_{\calUdag }$ définissent alors  des faisceaux de $R_{\Xdaginf }$-modules sur le site
$\Xdaginf $, notés respectivement $\Gamma_Z(\cal F_{\inf})$ et $j_*^{\inf}j^{-1}_{\inf}\cal F_{\inf}$.
\end{prop}

\demo
Il s'agit de voir que si $\calWdag  $ et $ \calUdag $ sont deux ouverts de $\Xdaginf$ avec $r: W\subset U$, le morphisme de restriction $(\sharp)$ induit
des  isomorphismes:
$$\hss\let\bigotimes\otimes
\def\quad{\hskip0.2ex}\scriptwd1em\scriptspace0pt
\matrix{R[\cal G_{\calWdag \sto \calUdag }]\Otimes_{r^{-1}R[\cal G_{\calUdag }]}r^{-1}\Gamma_Z(\cal F_{\calUdag })&\to &R[\cal G_{\calWdag \sto \cal
U\dag}]\Otimes_{r^{-1}R[\cal G_{\calUdag }]}r^{-1}\cal F_{\calUdag }&\to &R[\cal G_{\calWdag \sto \calUdag }]\Otimes_{r^{-1}R[\cal G_{\cal
U\dag}]}r^{-1}j_*j^{-1}\cal F_{\calUdag }\cr
\downarrow\indiced\simeq&&\downarrow\indiced\simeq&&\downarrow\indiced\simeq\cr\noalign{\kern3pt}
\Gamma_Z(\cal F_{\calWdag })&\varto{0pt}{17ex}& \cal F_{\calWdag }&\varto{0pt}{17ex} &j_*j^{-1}\cal F_{\calWdag }}
\hss$$
qui sont  transitifs ({\it cf. }\ref{defi'}). Le morphisme canonique:
$$R[\cal G_{\calWdag \to \calUdag }]\Otimes_{r^{-1}R[\cal G_{\cal
U\dag}]}r^{-1}j_*j^{-1}\cal F_{\calUdag }\to j_*\mBig(j^{-1}R[\cal G_{\calWdag \to \calUdag }]\Otimes_{r^{-1}j^{-1}R[\cal G_{\cal
U\dag}]}r^{-1}j^{-1}\cal F_{\calUdag }\mBig)$$ est un isomorphisme, parce que $R[\cal G_{\calWdag \to \calUdag }]$ est un $r^{-1}R[\cal G_{\cal
U\dag}]$-module localement libre de rang $1$. Mais le morphisme:
$$j^{-1}R[\cal G_{\calWdag \to \calUdag }]\Otimes_{r^{-1}j^{-1}R[\cal G_{\cal
U\dag}]}r^{-1}j^{-1}\cal F_{\calUdag }\to j^{-1}\cal F_{\calWdag },\leqno (\sharp):$$ qui est le morphisme de restriction du couple $(\calWdag |j^{-1}W,
\calUdag |j^{-1}U)$, est un isomorphisme. Les familles $\Gamma_Z(\cal F_{\calUdag })$ et $j_*j^{-1}\cal F_{\calUdag }$ ont toutes les propriétés des faisceaux
de $R_{\Xdaginf }$-modules.
\enddemo

\begin{Rema}Dans le  raisonnement précédent on peut remplacer le faisceau $R_{\Xdaginf }$ par un faisceau d'anneaux $\cal A_{\Xdaginf }$ et définir
les foncteurs de cohomologie locale dans la catégorie $\cal A_{\Xdaginf }\Mod$.
\end{Rema}

\begin{defi}Si $Z$ est un fermé de $X$, on définit ainsi  deux foncteurs exacts
à gauche qu'on peut dériver  dans la catégorie $ \Dplus (R_{\Xdaginf })$ et obtenir le triangle de cohomologie
locale habituel:\glossary{$\bfR \Gamma_Z(\cal F_{\inf}), \bfR j_*^{\inf}j^{-1}_{\inf}\cal F_{\inf}$}
$$\bfR \Gamma_Z(\cal F_{\inf})\to\cal F_{\inf}\to\bfR j_*^{\inf}j^{-1}_{\inf}\cal F_{\inf}\to\cdot$$ 
\end{defi}

\subsection{Le foncteur de cohomologie locale dans le cas des modules spéciaux}
Si $\calMdaginf $ est un $ \calDdag_{\Xdaginf /R}$-module à gauche spécial, le même raisonnement montre que $\Gamma_Z(\cal
M\daginf )$ et
$j^{\inf}_*j^{-1}_{\inf}\calMdaginf $ sont des $ \calDdag_{\Xdaginf /R}$-modules à gauche spéciaux et définissent des foncteurs covariants exacts à gauche
qu'on peut dériver dans la catégorie  $ \Dplus ((\calDdag_{\Xdaginf /R}, \Sp)\Mod)$. 

\begin{defi}Si $Z$ est un fermé de $X$, on définit ainsi  deux foncteurs covariants exacts
à gauche qu'on peut dériver  dans la catégorie $ \Dplus ((\calDdag_{\Xdaginf /R}, \Sp)\Mod)$, et obtenir le triangle de cohomologie
locale habituel:
$$\bfR \Gamma_Z(\calMdaginf )\to\calMdaginf \to\bfR j_*^{\inf}j^{-1}_{\inf}\calMdaginf \to,$$ qui pour tout ouvert se restreint par construction
au complexe de cohomologie locale usuelle.
\end{defi}

\begin{prop}\label{com-chl}Soient $Z$ un fermé de $X$,  $\calMdaginf $ un complexe borné à gauche et $\calNdaginf $ un complexe borné à droite de modules à gauche spéciaux. Alors,
il existe un isomorphisme canonique de complexes de Zariski de la catégorie $ \rmD(R_{X})$:
\glossary{$\bfR \Gamma_Z(\calMdaginf ), \bfR j_*^{\inf}j^{-1}_{\inf}\calMdaginf $}
$$\bfR \cHom\SubX{\calDdag_{\Xdaginf /R},\Sp}\mBig(\calNdaginf , \bfR \Gamma_{Z}(\calMdaginf )\mBig)
\simeq 
\bfR \Gamma_{Z}\mBig(\bfR \cHom\SubX{\calDdag_{\Xdaginf /R},\Sp}
(\calNdaginf , \calMdaginf )\mBig).$$
\end{prop}

\demo Soit $\calMdaginf \to \calIdaginf  $ la
résolution  injective par des modules spéciaux injectifs de $\calMdaginf $ construite localement par la résolution injective de Godement. Le complexe 
$\Gamma_{Z}(\calIdaginf  )$ est une résolution injective par des modules spéciaux de $\bfR \Gamma_{Z}(\calMdaginf )$, et le complexe
$\cHom^{\bullet}_{\calDdag_{\Xdaginf /R}}\!\big(\calNdaginf ,\Gamma_{Z}(\calIdaginf  )\big)$ est alors une résolution du complexe  $\bfR \cHom_{\calDdag_{\Xdaginf /R},\Sp}\!\big(\calNdaginf , \bfR \Gamma_{Z}(\calMdaginf )\big)$. D'autre part, $\cHom^{\bullet}_{\calDdag_{\Xdaginf /R}}(\calNdaginf ,\calIdaginf  )$ est un complexe de faisceaux de Zariski flasques, de sorte que le complexe $\Gamma_Z\big(\cHom^{\bullet}_{\calDdag_{\Xdaginf /R}}(\calNdaginf ,\calIdaginf  )\big)$ est une résolution du complexe $\bfR \Gamma_{Z}\big(\bfR \cHom_{\calDdag_{\Xdaginf /R},\Sp}(\calNdaginf , \calMdaginf )\big)$. L'isomorphisme de la proposition provient
de l'isomorphisme canonique de faisceaux de Zariski:
$$\cHom_{\calDdag_{\Xdaginf /R}}\!\big(\calNdaginf ,\Gamma_{Z}(\calMdaginf )\big)\simeq \Gamma_{Z}\big(\cHom_{\calDdag_{\Xdaginf /R}}(\calNdaginf , \calMdaginf )\big)$$ pour deux modules à gauche spéciaux.
\enddemo

\begin{Rema}L'extension $V[\cal G_{\calXdag }]\to\cal O_{\calXdag /V}[\cal G_{\calXdag }]$ est plate, et les foncteurs de cohomologie locale
calculés dans  $ \Dplus ((V_{\Xdaginf })\Mod)$ et $ \Dplus ((\cal O_{\Xdaginf /V})\Mod)$ coïncident. Mais l'extension $\cal
O_{\calXdag /V}[\cal G_{\calXdag }]\to\calDdag_{\calXdag /V}$ est loin d'être plate, même sur $K$. Aussi, les foncteurs de cohomologie  locale  sont
distincts lorsque calculés dans les catégories: $$ \Dplus (\cal O_{\Xdaginf /K}\Mod )\hbox{ et } \Dplus ((\calDdag_{\Xdaginf /K}, \Sp)\Mod),$$ et leurs rapports
posent des problèmes non triviaux.
\end{Rema}

\begin{prop}\label{loc-ouv}Si $X$ est un schéma   séparé et  lisse sur $k$,  les foncteurs $\bfR j_*^{\inf}j^{-1}_{\inf}$ et $\bfR j_*^{\diff}j^*_{\diff}$
coïncident canoniquement dans la catégorie $ \Dplus ((\calDdag_{\Xdaginf /K}, \Sp)\Mod)$, où $j: X-Z\hookrightarrow X$ est l'inclusion canonique.
\end{prop}

\demo
Soit $\calMdaginf $ un complexe de $ \Dplus ((\calDdag_{\Xdaginf /K}, \Sp)\Mod)$. Si $\calUdag $ est un ouvert du site  $\Xdaginf$, la valeur du
complexe $\bfR j_*^{\inf}j^{-1}_{\inf}\calMdaginf $ sur l'ouvert $\calUdag $ est   par définition  le complexe $\bfR j_*j^{-1}\cal M\dag_{\calUdag }$.  Notons $j:U':=
(X-Z)\cap U\hookrightarrow U$ l'inclusion canonique. L'espace annelé  $\cal U'{}\dag := (U', j^{-1}\cal O_{\calUdag /V})$ est un relèvement de
$U'$, donc les modules de transfert
$\calDdag_{\cal U'{}\dag \to\cal
U\dag/V}$ et
$\calDdag_{\calUdag \gets \cal
U'{}\dag /V}$ sont définis et valent $\calDdag_{\cal U'{}\dag /V}$ par construction.
La valeur du complexe  $j^*_{\diff}\calMdaginf $ sur l'ouvert $\cal U'$  est $j^{-1}\cal M\dag_{\calUdag }$  et, en vertu du théorème \ref{ima-com}, la valeur du  complexe $j_*^{\diff}j^*_{\diff}\calMdaginf $ sur l'ouvert  $\calUdag $ est le complexe $\bfR j_*j^{-1}\cal M\dag_{\calWdag }$. En prenant une résolution injective de $\calMdaginf $, on voit
que les foncteurs $\bfR j_*^{\inf}j^{-1}_{\inf}$ et $\bfR j_*^{\diff}j^*_{\diff}$ sont canoniquement isomorphes.
\enddemo
\begin{Rema} La possibilité de définir le triangle distingué de cohomologie locale dans le site infinitésimal $\dagger$-adique dans le cas d'un relèvement,
a été l'une des principales raisons du succès de notre point de vue [M-N$_1$] dès le début de la théorie.\end{Rema}

\subsection{Le théorème de pureté pour un couple lisse}
Soit $i:Y\hookrightarrow X$ une immersion fermée de schémas lisses  sur $k$ et soit $ \calXdag $  un relèvement  plat de $ X$. On a alors le
théorème de pureté suivant.

\begin{theo}\label{pur-coh}Sous les conditions précédentes, les faisceaux 
$\cal H^i_{Y}(\cal O_{\calXdag /V})$ et $\cal H^i_{Y}(\calDdag_{\calXdag /V})$ sont nuls si
$i\neq\codim_XY$.
\end{theo}

\demo La question est locale. Commençons par le cas du faisceau structural $\cal O_{\calXdag /V}$.

\begin{lemm}\label{12.3-3}
Soit  $i:Y\hookrightarrow X$ une immersion fermée de schémas   lisses  sur $R_1$ et soit $\calXdag $ un relèvement 
plat de $X$. Alors, $\cal H^i_{Y}(\cal O_{\calXdag /R})=0$ si $i>\codim_XY$.
\end{lemm}

\demo 
La question étant  locale, on peut supposer que le schéma  $X$ est affine  et que l'idéal de $Y$ est engendré  par $x_1,\dots,x_q $.
Le complémentaire
$X-Y$ admet le recouvrement affine $X-V(x_i)$, qui est acyclique pour le faisceau $\cal O_{\calXdag /R}$ en vertu du théorème d'acyclicité \ref{acy}. Cela
montre que $\cal H^i_{Y}(\cal O_{\calXdag /R})=0$ si $i>\codim_XY$.
\endsubdemo

\begin{lemm}\label{122.3-3}
Soit  $i:Y\hookrightarrow X$ une immersion fermée de schémas   lisses  sur $k$ et soit $\calXdag $ un relèvement 
plat de $X$. Alors, $\cal H^i_{Y}(\cal O_{\calXdag /V})=0$ si $i<\codim_XY$.
\end{lemm}
\demo
La question est locale, et l'on peut supposer que $X$ est affine d'algèbre $\dagger$-adique $A\dag$. Soient $B\dag$ un relèvement lisse de $Y$,  $u: A\dag\to
B\dag$ un relèvement de $i$ et $\cal I_{\calYdag /V}$ le noyau de $u$.

Soit $x=(y,z)$ un couple adapté à $u$, au sens de \ref{couple-adapte}.
Notons $I_z:=I_{z_1,\dots,z_q}$
l'idéal engendré par $z_1,\dots,z_q$. On peut considérer le $D_{A\dag/V}$-module à gauche: 
$$\Rm\alg H^{\bullet}_Y(X, \cal O_{\calXdag /V}):= \lim_{\to k}\ext^{\bullet}_{A\dag}(A\dag/I^k_{z},A\dag) $$ donné par  les modules d'hypercohomologie du
complexe: $$\bfR \alg\Gamma_Y(\cal O_{\calXdag /V}):= \bfR \lim_{\to k}\cHom_{\cal O_{\calXdag /V}}(\cal O_{\calXdag /V}/\cal I^k_{\cal
Y\dag/V},\cal O_{\calXdag /V}).$$

\begin{prop}\label{alg-ana}\label{12.3-5}Il existe un isomorphisme de $D\dag_{A\dag/V}$-modules à gauche:
$$ D\dag_{A\dag/V}\otimes_{D_{A\dag/V}}\alg H^{\bullet}_Y(X, \cal O_{\calXdag /V})\simeq H^{\bullet}_Y(X, \cal O_{\calXdag /V}).$$
\end{prop}

Comme $\set z_1,\dots,z_q/$  forment une suite régulière dans $A\dag$, cela entraîne que les $V$-modules $\ext^j_{\cal O_{\calXdag /V}}( \cal O_{\calXdag/V}/\cal I^k_{\calYdag /V}, \cal O_{\calXdag /V})$ sont nuls pour
$j<q$, et la proposition implique donc  le lemme  \ref{122.3-3} précédent.

\medskip
Reste à montrer l'isomorphisme de la proposition. Si $\cal I\dag_1$ et $\cal I\dag_2$ sont deux idéaux de type fini de $\cal O_{\calXdag /V}$ et   $Y_1$,
$Y_2$  les supports de $\cal O_{\calXdag /V}/\cal I\dag_1$ et de $\cal O_{\calXdag /V}/\cal I\dag_2$, la suite de Mayer-Vietoris
topologique fournit un triangle distingué:
$$\preskip0.5ex\DeuxLignes
 \bfR \alg \Gamma\goodSub{8pt}{4pt}{-2.5mm}{Y_1\cap  Y_2}(X, \cal O_{\calXdag /V})
\to
\bfR \alg \Gamma_{ Y_1}(X, \cal O_{\calXdag /V})\oplus\bfR \alg 
\Gamma_{
Y_2}(X, \cal O_{\calXdag /V})\to\\
\to
\bfR \alg \Gamma\goodSub{8pt}{4pt}{-2.2mm}{Y_1\cup   Y_2}(X, \cal O_{\calXdag /V})\to
\endlignes\postskip-1ex$$ 
et  un triangle distingué: 
$$\preskip1em\DeuxLignes
 \bfR  \Gamma_{ Y_1\cap  Y_2}(X, \cal O_{\calXdag /V})
\to
\bfR  \Gamma_{Y_1}(X, \cal O_{\calXdag /V})\oplus\bfR  \Gamma_{
Y_2}(X, \cal O_{\calXdag /V})
\to\\\to
\bfR  \Gamma_{ Y_1\cup   Y_2}(X, \cal O_{\calXdag /V})\to\,\cdot
\endlignes$$ 
Prenant une résolution $\calDdag_{\calXdag /V}$-injective
de $\cal O_{\calXdag /V}$ qui reste donc $\cal O_{\calXdag /V}$-injective. Comme 
l'extension $D_{A\dag/V}\to D\dag_{A\dag/V}$ est plate
([Me$_3$], Thm. 4.2.7), on trouve un  morphisme de triangles:

\hbox{\kern2.5pt\vbox{\def\downarrow{\vardownarrow14pt}
\def\SubX{\goodSub{2pt}{3pt}{-1.5mm}}
\def\\{\cr\noalign{\kern-10pt}}\scriptwd=1.5em
\def\quad{\hskip4pt}
\mathrigid4mu
$$\hbox to\hsize{$\matrix{
 D\dag_{A\dag/V}\Otimes_{ D_{A\dag/V}}
 \bfR \alg \Gamma_{Y_1\cap  Y_2}(X, \cal O_{\calXdag /V})
\\
 \downarrow\cr
 \bfR \Gamma_{Y_1\cap  Y_2}(X, \cal O_{\calXdag /V})
}$\hss}$$
\vskip9pt$$\def\Gammas_#1{\Gamma\SubX{#1}}
\matrix{ 
\llap{$\to$\hvpath{-4 4 57 19 -3}\ }D\dag_{A\dag/V}\Otimes_{D_{A\dag/V}}
\bfR \alg \Gammas_{ Y_1}(X, \cal O_{\calXdag /V})
&\displaystyle\bigoplus& 
D\dag_{A\dag/V}\Otimes_{ D_{A\dag/V}}
\bfR \alg \Gammas_{ Y_2}(X, \cal O_{\calXdag /V})
\\
 \downarrow&&\downarrow\cr
\llap{$\to$\hvpath{-16 18 48 6 -3}\ }
\bfR  \Gamma_{ Y_1}(X, \cal O_{\calXdag /V})&\displaystyle\bigoplus&\bfR  \Gamma_{
 Y_2}(X, \cal O_{\calXdag /V})
 }$$
\vskip7pt
$$\hbox to\hsize{\hss$\matrix{ 
\llap{$\to$\hvpath{-5 4 58 19 -3}\ }
  D\dag_{A\dag/V}\Otimes_{D_{A\dag/V}}\bfR \alg \Gamma_{Y_1\cup   Y_2}(X, \cal O_{\calXdag /V})\cr
\noalign{\kern-7pt}\downarrow\cr\noalign{\kern3pt}
\llap{$\to$\hvpath{-18 19 50 5 -3}\ }
 \bfR  \Gamma_{Y_1\cup   Y_2}(X, \cal O_{\calXdag /V}).
 }$}$$}}

En utilisant ce morphisme de triangles et en raisonnant par récurrence sur le nombre d'équations, on est ramené à montrer la proposition dans le cas de
l'idéal engendré par
$z_1\cdots z_r$ pour $1\leq r\leq q$, ce qui est conséquence du lemme suivant.

\begin{lemm}\label{mero}\label{12.3-6}Sous les conditions précédentes, pour tout $r$ tel que $1\leq r\leq q$, le $D_{A\dag/V}$-module à gauche $A\dag[1/(z_1\cdots z_r)]$, resp.   le
$D\dag_{A\dag/V}$-module à gauche $(A\dag[1/(z_1\cdots z_r)])\dag$, est engendré par l'élément
$1/(z_1\cdots z_r)$ dont l'annulateur est l'idéal engendré
par les opérateurs différentiels: $$\Delta_1^{\alpha_1}z_1,\dots, \Delta_r^{\alpha_r}z_r,\Delta_{r+1}^{\alpha_{r+1}},\dots,
\Delta_n^{\alpha_n}$$  tels que  $\alpha_1>0,\dots,\alpha_r>0,\alpha_{r+1}>0,\dots,\alpha_n>0$.
\end{lemm}

\demo 
Commençons par $A\dag[1/(z_1\cdots z_r)]$, lequel  est engendré en tant que 
$D_{A\dag/V}$-module  à gauche, de façon évidente, par  $1/(z_1\cdots z_r)$. Si $P$ est un opérateur 
d'ordre fini qui annule $1/(z_1\cdots z_r)$, on peut le supposer   de la forme
$P= \sum_{\alpha\in \Bbb N^r} a_\alpha\Delta_1^{\alpha_1}\cdots \Delta_r^{\alpha_r}$. 
Considérons l'ordre lexicographique sur les $r$-uplets $(\alpha_1,\dots,\alpha_r)$ et notons $\exp(P)$ le plus grand indice des indices des coefficients
de $P$. La suite $z_1,\dots, z_r$ étant régulière dans $A\dag$, il s'ensuit 
que le  coefficient $a_{\exp(P)}=a_\alpha$  appartient à l'idéal $(z_i, \alpha_i\neq0)$. Modulo
l'idéal $(\Delta_1^{\alpha_1}z_1,\dots, \Delta_r^{\alpha_r}z_r,\alpha_1>0,\dots,\alpha_r>0)$, le polynôme  $P$ est alors égal à un opérateur dont l'exposant est strictement plus petit que $\exp(P)$.
Par récurrence, on obtient  que $P$ appartient à l'idéal $(\Delta_1^{\alpha_1}z_1,\dots, \Delta_r^{\alpha_r}z_r,\alpha_1>0,\dots,\alpha_r>0)$.

\medskip
Passons au cas d'ordre infini. Le $D\dag_{A\dag/V}$-module $(A\dag[1/(z_1\cdots z_r)])\dag$ est
aussi engendré, de façon évidente,  par $1/(z_1\cdots z_r)$. 
Si
$P$ un opérateur d'ordre infini tel que
$P(1/(z_1\cdots z_r))=0$, on peut supposer que $P= \sum_{\alpha\in \Bbb N^r} a_\alpha\Delta_1^{\alpha_1}\cdots \Delta_r^{\alpha_r}$. 
Écrivons $P$ sous la forme $\sum_{\alpha}\Delta_1^{\alpha_1}\cdots \Delta_r^{\alpha_r}a'_\alpha$. Considérons la réduction $P_s$ modulo $\goth m^s$ de $P$ qui annule
$1/(z_1\cdots z_r)$. Le raisonnement du cas où l'ordre est fini montre que chaque coefficient $a'_\alpha$ appartient modulo $\goth m^s$ à l'idéal
$(z_i, \alpha_i\neq0)$, et ces coefficients $a'_\alpha$ appartiennent donc à l'idéal  $(z_i, \alpha_i\neq0)$. En considérant une présentation de $A\dag$, 
on peut montrer, en vertu du théorème de continuité de la division par une base 
de division [M-N$_3$], qu'on peut écrire $a'_\alpha= \sum_{i,\alpha_i\neq0}b_{\alpha,i}z_i$, de sorte 
que $P_i:= \sum_{\alpha,\alpha_i\neq0}\Delta_r^{\alpha_r}b_{\alpha,i}$ est encore un opérateur différentiel et que $P=\sum_{i}P_iz_i$. 
En raisonnant par récurrence sur la longueur des monômes en $\Delta$ qui interviennent dans $P$, on trouve que $P$ appartient à l'idéal $(\Delta_1^{\alpha_1}z_1,\dots, \Delta_r^{\alpha_r}z_r,\alpha_1>0,\dots,\alpha_r>0)$.
Cela termine  la démonstration  du lemme \ref{12.3-6}, mais également  celles de la proposition  \ref{12.3-5} et du lemme \ref{122.3-3}.
\hfill$\scriptstyle\square$

\bigskip
Les mêmes  raisonnements   montrent que sous les conditions précédentes on a  la nullité:
$$ H^j_{Y\times_k A^n_k}(T^*X, \cal O_{(\cal T^*\cal X)\dag/V})= H^j_{Y\times_k A^n_k}(X\times_k A^n_k, \cal O_{(\cal T^*\cal X)\dag/V})=0,$$
où  $(\cal T\sp{*}\cal X)\dag:=(T\sp{*}X,\cal O_{(\cal T^*\cal X)\dag/V})$ est le fibré cotangent $\dagger$-adique de $X$,
pour $j\neq \codim_XY$. Mais en vertu du théorème du symbole total \ref{sym-tot}, on a
l'isomorphisme:
$$ H^j_{Y\times_k A^n_k}(X\times_k A^n_k, \cal O_{(\cal T^*\cal X)\dag/V})\simeq H^j_{Y}(X, \calDdag_{\calXdag /V}),$$
et $H^j_{Y}(X, \calDdag_{\calXdag /V})=0$ si $j\neq \codim_XY$.
On en déduit que: 
$$\cal H^j_{Y}(\calDdag_{\calXdag /V})=0$$ si $j\neq \codim_XY$, puisque ce faisceau est associé au préfaisceau précédent. Cela termine la preuve du théorème
\ref{pur-coh}.
\enddemo

\begin{coro}Le foncteur $\calUdag \fonct \cal H^{\codim_UY}_{Y}(\calDdag_{\calUdag /V})$ définit un $(\calDdag_{ \Xdaginf /V}, \calDdag_{ \Xdaginf /V})$-bimodule  $\cal H^{\codim_XY}_{Y}(\calDdag_{ \Xdaginf /V})$ sur le site infinitésimal $ \Xdaginf $.
\end{coro}

\demo En effet, si $r\dag : \calWdag \to \calUdag $ est un morphisme, on définit
$$ r^{-1}\cal H^{\codim_UY}_{Y}(\calDdag_{\calUdag /V})\to \cal H^{\codim_WY}_{Y}(\calDdag_{\calWdag /V})\leqno (\sharp)_{r\dag}:$$ 
comme le morphisme
induit par le morphisme de complexes 
$$ r^{-1}\bfR \Gamma_{Y}(\calDdag_{\calUdag /V})\to\bfR \Gamma_{Y}(\calDdag_{\calWdag /V})\,,\leqno (\sharp)_{r\dag}:$$
dont on voit en considérant les résolutions locales de \v Cech que c'est un morphisme $(\sharp)_{r\dag}$-linéaire à droite et à gauche. Ces morphismes sont transitifs de façon naturelle.
\enddemo

\begin{coro} Le foncteur: $$\calMdaginf \fonct \cal H^{\codim_XY}_{Y}(\calDdag_{ \Xdaginf /V})\otimes_{\calDdag_{ \Xdaginf /V}}\calMdaginf $$ est covariant exact à droite de la catégorie
$ (\calDdag_{ \Xdaginf /V}, \Sp)\Mod$ dans elle-même.
\end{coro}

\demo Il suffit de montrer que si $\calMdaginf $ est un ${\calDdag_{ \Xdaginf /V}}$-module à gauche spécial, le ${\calDdag_{ \Xdaginf /V}}$-module à gauche: $$\cal H^{\codim_XY}_{Y}(\calDdag_{ \Xdaginf /V})\otimes_{\calDdag_{ \Xdaginf /V}}\calMdaginf $$ est spécial. Mais par construction, un élément $g$ du groupe $\cal G_{\Xdaginf }$ agit par $gPg^{-1}\otimes gm=
gP\otimes m$.
\enddemo

\subsection{La comparaison entre la cohomologie locale et l'image directe dans le cas d'une immersion fermée. Le morphisme de foncteurs $ \Adj _*^i$}

L'analogue de la proposition \ref{loc-ouv}  pour une immersion fermée {\bf n'a pas lieu} en général et son étude est beaucoup plus délicate.
Nous allons utiliser le théorème de pureté pour construire un morphisme de foncteurs  $ \Adj _*^i$ entre les foncteurs $i_*^{\diff}i^*_{\diff}$ et $\bfR \Gamma_Y[\codim_XY]$ pour une
immersion fermée
$i: Y\to X$ de schémas lisses sur le corps
$k$. 
\medskip

Supposons que $X$ est affine et soit $u: A\dag\to B\dag$ un relèvement de $i$ et $x=(y, z)$ un système adapté à $u$.  
\begin{lemm}Pour $X$ affine 
{\bf assez petit}, 
le  $D_{\!A\dag/V}$-module à gauche $\Rm\alg H^{\codim_XY}_Y(X, \cal O_{\calXdag /V})$ est engendré par la classe $[1/z]=[1/(z_1\cdots z_q)]$
dont l'annulateur est l'idéal 
$$\big(z_1,\dots, z_q, \Delta_{q+1}^{\alpha_{q+1}},\dots,
\Delta_n^{\alpha_n}, \alpha_{q+1}>0,\dots, \alpha_{n}>0\big),$$
et le  $D\dag_{A\dag/V}$-module à gauche 
$ H^{\codim_XY}_Y(X, {\cal O_{\calXdag /V}}
)$ est
engendré par la classe $[1/z]=[1/(z_1\cdots z_q)]$ dont l'annulateur est l'idéal 
$$\big(z_1,\dots, z_q, \Delta_{q+1}^{\alpha_{q+1}},\dots,
\Delta_n^{\alpha_n}, \alpha_{q+1}>0,\dots, \alpha_{n}>0)\,.$$
\end{lemm}

\demo
En vertu de la proposition \ref{alg-ana}, il suffit de traiter le cas du module $\Rm\alg H^{\codim_XY}_Y(X, \cal O_{\calXdag /V})$. Il est évident que la
classe $[1/z]$ est un générateur. Soit un opérateur d'ordre fini $P(x, \Delta) := \sum_\alpha a_\alpha\Delta^\alpha$ qui annule $[1/z]$. On peut supposer
que
$\alpha= (\alpha_1,\dots,\alpha_q,0, \dots,0)$. Il nous faut montrer que $P(x, \Delta)$ appartient à l'idéal à gauche engendré par $z_1,\dots, z_q$.
L'équation $P(x, \Delta)([1/z])=0$ se traduit par une égalité dans $A\dag[1/z]$:
$$P(x, \Delta)(1/(z_1\cdots z_q))= \sum_{1\leq i\leq q}\frac{b_i}{(z_1\cdots\hat{z_i}\cdots
z_q)^{\alpha^i}}\,,\qquad \alpha^i_j\geq 1\,,$$ 
entre fonctions méromorphes, où les coefficients $b_i$ sont des fonctions de $A\dag$. Mais:
{\arraycolsep2pt\begin{eqnarray*}
\frac{b_i}{(z_1\cdots\hat{z_i}\cdots z_q)^{\alpha^i}}&=& b_i\Delta^{\alpha_1^i-1,\dots, \alpha^i_{q}-1}(1/(z_1\cdots\hat{z_i}\cdots z_q))\\
&=& 
b_i\Delta^{\alpha_1^i-1,\dots, \alpha^i_{q}-1}z_i(1/(z_1\cdots z_q)).
\end{eqnarray*}}
Il suffit de démontrer que l'opérateur 
$P(x, \Delta)-\sum_i b_i\Delta^{\alpha_1^i-1,\dots, \alpha^i_{q}-1}z_i$ appartient à l'idéal $(z_1,\dots, z_q)$. Mais on a déjà vu dans la démonstration
du lemme \ref{mero} que si un tel opérateur annule
$1/(z_1\cdots z_q)$ comme fonction méromorphe, il appartient à l'idéal 
$(\Delta^{\alpha_1}_1z_1,\dots, \Delta_q^{\alpha_r}z_q)$.
\enddemo

\begin{lemm}\label{sep} Pour $X$ affine {\bf assez petit} les $V$-modules: 
$$ H^{\codim_XY}_{Y}(X, \cal O_{\calXdag /V})\quad {\rm et }  \quad
H^{\codim_XY}_{Y}(X, \calDdag_{\calXdag /V})$$ sont {\bf séparés}  pour la topologie $\goth m$-adique.
\end{lemm}
\demo Commençons par le $V$-module  $ H^{\codim_XY}_{Y}\!(X, \cal O_{\calXdag })$. En vertu du lemme précédent,
il faut montrer que si un opérateur $P(x, \Delta) = \sum_{\alpha=0, \infty} a_\alpha\Delta^\alpha$ appartient à l'idéal: 
$$(z_1,\dots, z_q, \Delta_{q+1}^{\alpha_{q+1}},\dots,
\Delta_n^{\alpha_n}, \alpha_{q+1}>0,\dots, \alpha_{n}>0)$$ modulo $\goth m^s$ pour tout $s\geq 1$, alors il appartient à l'idéal: $$(z_1,\dots, z_q,
\Delta_{q+1}^{\alpha_{q+1}},\dots,
\Delta_n^{\alpha_n}, \alpha_{q+1}>0,\dots, \alpha_{n}>0).$$ On peut supposer que $P(x, \Delta) = \sum_{\alpha=0, \infty} a_\alpha\Delta^\alpha$, avec $\alpha_{q+1}=\cdots=\alpha_n=0$. 
On a alors, pour tout $s\geq1$, les égalités 
$$P(x, \Delta) = \sum_{\alpha=0, \infty} a_\alpha\Delta^{\alpha_1,\dots,\alpha_q}= 
\sum_{\gamma=0, \infty} b_{\gamma;s}\Delta^\gamma +\mBig( \sum_{i=1,q} Q_iz_i
+\sum_{\beta}Q_\beta\Delta^{\beta}\mBig),$$ 
avec des $n$-uplets $\beta$ en nombre fini tels que $ \beta_{q+1}>0,\dots, \beta_{n}>0$,
où les coefficients $b_{\alpha;s}$ sont dans $\goth m^sA\dag$. 

En vertu du théorème du symbole
total \ref{sym-tot}, on a les égalités: 
$$\DeuxLignes
\sum_{\alpha=0, \infty} a_\alpha\Delta^{\alpha_1,\dots,\alpha_q}=
\sum_{\alpha_1,\dots,\alpha_q}b_{\alpha_1,\dots,\alpha_q,0,\dots,0;s}\Delta^{\alpha_1,\dots,\alpha_q}+\\
+
\sum_{i}\sum_{\alpha_1,\dots,\alpha_q} q_{i,\alpha_1,\dots,\alpha_q,0,\dots,0}\Delta^{\alpha_1,\dots,\alpha_q}z_i,
\endlignes$$
où les coefficients
$q_{i,\alpha}$ sont ceux de $Q_i$. 
Si on écrit cette égalité sous la forme:
$$\DeuxLignes
\sum_{\alpha=0, \infty} \Delta^{\alpha_1,\dots,\alpha_q}a'_\alpha=
\sum_{\alpha_1,\dots,\alpha_q}\Delta^{\alpha_1,\dots,\alpha_q}b'_{\alpha_1,\dots,\alpha_q,0,\dots,0;s}+\\+
\sum_{i}\sum_{\alpha_1,\dots,\alpha_q} \Delta^{\alpha_1,\dots,\alpha_q}q'_{i,\alpha_1,\dots,\alpha_q,0,\dots,0}z_i,
\endlignes$$
on trouve que les coefficients
$a'_\alpha$ appartiennent à l'idéal $(z_1,\dots ,z_q)$ modulo $\goth m^s$ pour tout $s\geq1$, et qu'ils appartiennent donc à l'idéal $(z_1,\dots, z_q)$.

\bigskip
En vertu du théorème du symbole total \ref{sym-tot}, pour $X$ assez petit,  le $V$-module $ H^{\codim_XY}_{Y}(X, \calDdag_{\calXdag /V})$  est isomorphe au $V$-module $\Rm
H^{\codim_XY}_{Y\times_k A^n_k}(X\times_k A^n_k,
\cal O_{ (\cal T^*\cal X)\dag})$ qui est séparé pour la topologie $\goth m$-adique en vertu de ce qui précède. D'où le lemme \ref{sep}.
\enddemo

Supposons que $X$ est affine et soit $u: A\dag\to B\dag$ un relèvement de $i$ et $x=(y, z)$ un système adapté à $u$.  Nous avons vu en 
\ref{ind-tra} que le module de transfert
$\calDdag_{\calYdag \to\calXdag /V}$ est un $\calDdag_{\calXdag /V}$-module à droite engendré par $u$.
Et, nous avons vu en 
\ref{ind-tra'} que le module de transfert $\calDdag_{\calXdag \gets\calYdag /V}$ est un $\calDdag_{\calXdag /V}$-module à gauche engendré par
$dyu(dx)^{*}$, où $(dx)^*$ est la base duale de $dx= dx_1\cdots dx_n$. Un élément
$Pdyu(dx)^{*}$, considéré  comme morphisme de $i^{-1}\omega_{\calXdag /V}\to \omega_{\calYdag /V}$, opère comme $dyu\trans P(dx)^{*}$: $$i^{-1}\omega_{\calXdag/V}\stackrel{(dx)^*}{\to} i^{-1}\cal O_{\calXdag /V}\stackrel{u\circ\trans P}{\to}\cal O_{\calYdag /V}\stackrel{dy}{\to} \omega_{\cal
Y\dag/V}.$$ Un tenseur  élémentaire de $\calDdag_{\calXdag \gets
\calYdag /V}\otimes_{\calDdag_{\calYdag /V}}\calDdag_{\calYdag \to \calXdag /V}$ s'écrit comme $Pdyu(dx)^{*}\otimes uQ$. Au générateur
$dyu(dx)^{*}\otimes u$ du bimodule  on   associe la classe
$[1/z]=[1/(z_1\cdots z_q)]$ dans
$\cal H^q_{Y}(\calDdag_{\calXdag /V})$. On définit l'application:  
$$i_*\mBig(\calDdag_{\calXdag \gets \cal
Y\dag/V}\otimes_{\calDdag_{\calYdag /V}}\calDdag_{\calYdag \to \calXdag/V}\mBig)\to \cal H^q_{Y}(\calDdag_{\calXdag /V})\leqno(*)_{u,z,y}:$$ par linéarité par:
$$Pdyu(dx)^{*}\otimes uQ\mapsto P[1/z_1\cdots z_q]Q.$$ 

\begin{lemm}\label{bimod}L'application $(*)_{u,z,y}$ est un morphisme bien défini
de $(\calDdag_{\calXdag /V},\calDdag_{\calXdag /V})$-bimodules.
\end{lemm}\demo Il faut  montrer que si une somme finie $\sum Pdyu(dx)^{*}\otimes uQ$ est nulle dans le produit tensoriel,  alors la somme  $\sum P[1/z]Q$ est nulle dans la cohomologie
locale.  La question est locale. En vertu du lemme \ref{sep},  il suffit donc d'établir  que la réduction $\sum P_s[1/z_s]Q_s$ modulo $\goth m^s$ est
nulle pour
$s\geq1$. On est ramené ainsi à démontrer  le lemme \ref{bimod} modulo $\goth m^s$ pour tout $s\geq1$. 

\bigskip
Soient $P_s= \sum_{\alpha, \beta}\Delta_{y_s}^\alpha\Delta_{z_s}^\beta a_{\alpha, \beta} $ et $Q_s= \sum_{\alpha, \beta}b_{\alpha, \beta}
\Delta_{y_s}^\alpha\Delta_{z_s}^\beta$  des opérateurs différentiels tels que  $P_sdy_su_s(dx_s)^{*}\otimes u_sQ_s$ est nul. Si $a_{\alpha, \beta}$ ou
$b_{\alpha, \beta}$ appartiennent à l'idéal $(z_s)$, alors $$\Delta_{y_s}^\alpha\Delta_{z_s}^\beta a_{\alpha, \beta}dy_su_s(dx_s)^{*}\otimes u_sQ_s
\hbox{\quad et\quad }
P_sdy_su_s(dx_s)^{*}\otimes u_sb_{\alpha, \beta}\Delta_{y_s}^\alpha\Delta_{z_s}^\beta$$ sont nuls, et
$$\Delta_{y_s}^\alpha\Delta_{z_s}^\beta a_{\alpha, \beta}[1/z_s]Q_s
\hbox{\quad et\quad }
P_s[1/z_s]b_{\alpha, \beta}\Delta_{y_s}^\alpha\Delta_{z_s}^\beta$$ sont nuls aussi. On peut donc supposer que les images par $u_s$ des coefficients
$a_{\alpha,
\beta}$ et $b_{\alpha, \beta}$ sont non nuls. Le module de transfert $\calD_{ Y_s\to 
X_s/V_s}$ est un $\calD_{ Y_s/V_s}$-module à 
{gauche}
libre, les opérateurs $u_s\Delta^\beta_{z_s}$ formant une base,  et de même le module de transfert $\cal
D_{ X_s\gets  Y_s/V_s}$ est un $\calD_{ Y_s/V_s}$-module à 
{droite}
libre, les opérateurs $dy_su_s\Delta^\beta_{z_s}(dx_s)^*$ formant une base. On a alors  les
égalités:
$$
u_sQ_s = \sum_{\beta}Q_{s,\beta}u_s\Delta_{z_s}^\beta:= \sum_{\beta}\mBig(\sum_\alpha u_s(b_{\alpha, \beta})\Delta_{y_s}^\alpha\mBig)u_s\Delta_{z_s}^\beta
\postskip0pt$$ 
et
$$\DeuxLignes
P_sdy_su_s(dx_s)^{*}=\\= \sum_\beta dy_su_s\Delta^\beta_{z_s}(dx_s)^*P_{s,\beta}:= \sum_{\beta}dy_su_s\Delta^\beta_{z_s}(dx_s)^*\mBig(\sum_\alpha\Delta_{y_s}^\alpha
u_s(a_{\alpha, \beta})\mBig).
\endlignes$$
La nullité du tenseur $P_sdy_su_s(dx_s)^{*}\otimes u_sQ_s$ est équivalente à la nullité des produits $P_{s,\beta}Q_{s,\gamma}$.
Si l'on pose $\tilde P_{s,\beta}:= \sum_\alpha\Delta_{y_s}^\alpha
a_{\alpha, \beta}$ et $\tilde Q_{s,\gamma}:= \sum_\alpha b_{\alpha, \gamma}\Delta_{y_s}^\alpha$, on a les égalités:
$$P_s[1/z_s]Q_s= \sum_{\beta, \gamma}\Delta_{z_s}^\beta\tilde P_{s,\beta}[1/z_s]\tilde Q_{s,\gamma}\Delta_{z_s}^\gamma= \sum_{\beta,
\gamma}\Delta_{z_s}^\beta\tilde P_{s,\beta}\tilde Q_{s,\gamma}[1/z_s]\Delta_{z_s}^\gamma=0.$$
Un raisonnement similaire vaut pour une somme finie 
$$\sum P_sdy_su_s(dx_s)^{*}\otimes u_sQ_s.\eqno\enddemo$$

\begin{lemm}\label{12.4-4}Pour $u$ fixé le morphisme $(*)_{u}:= (*)_{u,z,y}$ ne dépend pas du couple adapté $(z, y)$.
\end{lemm}

\demo
Si $(z', y')$ est un autre couple, on a dans $\calDdag_{\calXdag \gets \cal
Y\dag/V}$ l'égalité :
$$dy'\otimes u\otimes (dx')^*= \det(\partial_{z'}z)dy\otimes u\otimes (dx)^*,$$ où $\partial_{z'}z$ est la matrice $\partial_{z'_i}z_j$. On doit montrer
l'égalité dans $\cal H^q_{Y}(\calDdag_{\calXdag /V})$:
$$[1/z']=\det(\partial_{z'}z)[1/z].$$ Il suffit de montrer  cette égalité dans $\cal H^q_{Y}(\cal
O_{\calXdag /V}).$ 

\begin{lemm}\label{12.4-5}Sous les conditions précédentes, on a l'égalité dans $\cal H^q_{Y}(\cal
O_{\calXdag /V})$: $$[1/z']=\det(\partial_{z'}z)[1/z].$$
\end{lemm}

\demo La question est de nature locale. Nous allons raisonner par récurrence sur $q\geq 1$. Supposons  que $q=1$, alors $z'= az $ et $z=a'z'$ ;  cela montre que
$\det(\partial_{z'}z)(1/z) = 1/z'+a\partial_{z'}a'$ et donc
$[1/z']=\det(\partial_{z'}z)[1/z]$.

\bigskip
Supposons que $q\geq2$ et que $z_1=z'_1$. Considérons la suite exacte:
$$0\to \cal O_{\calXdag /V}\stackrel{z_1}{\to} \cal O_{\calXdag /V}\to\cal O_{\calXdag /V}/z_1\to0.$$ En vertu du théorème de pureté \ref{pur-coh}, on
obtient une suite exacte:
$$0\to \cal H^{q-1}_{Y}(\cal O_{\calXdag /V}/z_1)\to\cal H^{q}_{Y}(\cal O_{\calXdag /V})\stackrel{z_1}{\to}\cal H^{q}_{Y}(\cal O_{\calXdag /V})\to0.$$
Notons $\bar z_2,\dots,\bar z_q$ les classes de $z_2,\dots,z_q$  dans le quotient par $z_1$, et de même notons  $\bar z'_2,\dots,\bar z'_q$ les classes de
$z'_2,\dots,z'_q$ dans le quotient par $z_1$. Comme
$z_1[1/z']=z_1\det(\partial_{z'}z)[1/z]=0$, il suffit de montrer en vertu de l'hypothèse de récurrence que  la classe
$[1/z]$ provient de la classe
$[1/\bar z_2\cdots \bar z_q]$ par le morphisme de connexion et que la classe $\det(\partial_{z'}z)[1/z]$ provient de la classe
$\det(\partial_{\bar z'}\bar z)[1/\bar z_2\cdots \bar z_q]$ par le morphisme de connexion. En considérant le complexe de \v Cech du recouvrement affine 
$\set U_{i},
i=1,\dots,q/$, resp. $\set U'_{i'}, i'=1,\dots,q/$, où $U_i$, resp.  $U'_{i'}$, est l'ouvert principal défini par la réduction modulo
$\goth m$ de $z_i$, resp. de $z'_{i'}$, de
$X-Y$ à valeurs dans le faisceau structural, on voit effectivement que tel est bien le cas.

\bigskip
Soient $z=(z_1,\dots,z_q)$ et $z'=(z'_1,\dots,z'_q)$ deux systèmes comme précédemment. Alors, nous allons montrer que localement il existe un entier $j$, avec $1\leq j\leq q,$ tel que $\set z'_j,z_2,\dots,z_q/$ engendre l'idéal  $\cal I_{\calYdag /V}$ noyau  de $u$. Considérons  la suite exacte de
$\cal O_{\calYdag /V}
$-modules localement libres de type fini:
$$0\to \cal I_{\calYdag /V}/\cal I_{\calYdag /V}^2\to \cal O_{\calYdag /V}\otimes_{i^{-1}\cal O_{\calXdag /V}}i^{-1}\Omega_{\calXdag /V}\to\Omega_{\cal
Y\dag/V}\to0$$ et sa réduction modulo $\goth m$:
$$0\to \cal I_{Y/k}/\cal I_{Y/k}^2\to \cal O_{ Y/k}\otimes_{i^{-1}\cal O_{ X/k}}i^{-1}\Omega_{ X/k}\to\Omega_{
Y/k}\to0.$$ Considérons les réductions $\set\bar z_2,\dots,\bar z_q/$ modulo $\goth m$ et leurs différentielles $\set d\bar z_2,\dots, d\bar z_q/$ au voisinage d'un
point. Alors, il existe un entier $j$ tel que les différentielles $d\bar z'_j, d\bar z_2,\dots, d\bar z_q$ forment une base du fibré conormal de $Y$
dans $X$ et, en vertu du lemme de Nakayama, les classes de $\bar z'_j,\bar z_2,\dots,\bar z_q$ forment une base locale du $\cal O_{Y/k}$-module $\cal
I_{Y/k}/\cal I_{Y/k}^2$. Cela montre que localement $\bar z_1$ est dans l'idéal 
$(\bar z'_j,\bar z_2,\dots,\bar z_q)$, et   que $ z_1$ est donc dans l'idéal
$ (z'_j, z_2,\dots,z_q)$.

\bigskip
Quitte à réindexer, on peut supposer que $j=1$. Notons  par $z''=(z''_1,z''_2,\dots,z''_q)$ l'idéal  $(z'_1,z_2,\dots,z_q)$. Alors, en vertu de ce qui
précède on a les égalités:
$$[1/z_1\cdots z_q]= \det\partial_z(z'')[1/z''_1\cdots z''_q]$$ et: 
$$[1/z''_1\cdots z''_q]= \det\partial_{z''}(z')[1/z'_1\cdots z'_q].$$ Mais comme $\det\partial_z(z'')\det\partial_{z''}(z')\equiv\det\partial_z(z')$ modulo
$\cal I_{\calYdag /V}$, on a bien l'égalité:
$$[1/z_1\cdots z_q]= \det\partial_z(z')[1/z'_1\cdots z'_q]\,,$$
ce qui termine la preuve des lemmes \ref{12.4-5} et \ref{12.4-4}.
\hskip0pt plus 1filll$\scriptstyle\square$\enddemo

\begin{lemm}Le morphisme $(*):= (*)_{u}$ ne dépend pas du relèvement $u$ de $i$.
\end{lemm}

\demo
Soit $u'$ un autre relèvement de $i$. En  vertu de la lissité,  on a la factorisation $u= u'g$, où $g$ est un élément du groupe $G_{A\dag}$. Soit $(z,y)$ un
couple adapté à $u$. Par construction, le couple $(z',y')$, avec 
$z':= g(z) $ et $y':= g(y)$,
est adapté à $u'$. Si $Pdyu(dx)^{*}\otimes uQ$ est un tenseur, on a l'égalité:
$$\DeuxLignes
Pdyu(dx)^{*}\otimes uQ= P\trans g\det(\partial_{z}z')dy'u'(dx')^{*}\otimes u'gQ=\\
=P\trans g(1)g^{-1}\det(\partial_{z}z')dy'u'(dx')^{*}\otimes u'gQ.\endlignes$$ 
Donc, son
image par le morphisme $(*)_{u'}$ construit à l'aide de $u'$ est égale à:
$$P\trans g(1)g^{-1}\det(\partial_{z}z')[1/z']gQ.$$ En vertu du lemme précédent (\ref{12.4-5}), on a l'égalité:
$$[1/z] = \det(\partial_{z}z')[1/z']$$ et l'égalité:$$g^{-1}\det(\partial_{z}z')[1/z']g= g^{-1}(\det(\partial_{z}z')[1/z]).$$ Mais en vertu du théorème \ref{act-drt},
on a l'égalité
$\trans g(1)= \det(\partial_{x}g(x))$. Or,  par construction:
$$\det(\partial_{x}g(x))g^{-1}\det(\partial_{z}z')\equiv1\,,$$ modulo le noyau de $u$. On a alors  l'égalité:
$$P[1/z]Q=P\trans g(1)g^{-1}\det(\partial_{z}z'[1/z'])gQ,$$ qui montre que le morphisme $(*)$ ne dépend pas du relèvement de $i$.
\enddemo

\goodbreak
Les lemmes précédents  légitiment la définition qui suit.

\begin{defi}\label{adj-loc} Soient $i:Y\hookrightarrow X$ une immersion fermée de schémas lisses  sur $k$
et $\calXdag $ et $\calYdag $ des relèvement
 plats. On définit le morphisme 
$ \Adj _*^i(\calDdag_{\calXdag /V})$  de $(\calDdag_{\calXdag /V},\calDdag_{\calXdag /V})$-bimodules:
$$\DeuxLignes
i_*\mBig(\calDdag_{\calXdag \gets \cal
Y\dag/V}\halfsmash{\Otimes_{\calDdag_{\calYdag /V}}}
\calDdag_{\calYdag \to \calXdag/V}\mBig)\to\cal H^{\codim_XY}_{Y}(\calDdag_{\calXdag /V})
\simeq\\
\simeq\bfR \Gamma_Y(\calDdag_{\calXdag /V})[\codim_XY],
\leqno{ \Adj _*^i(\calDdag_{\calXdag /V}):}
\endlignes
$$ comme  le morphisme
$(*)$  précédent, défini localement.
\end{defi}

\begin{theo} Pour $s\geq1$, la réduction modulo $\goth m^s$ du morphisme $ \Adj _*^i(\calDdag_{\calXdag /V})$
$$\scriptwd1.2em\DeuxLignes
i_*\mBig(\calD_{ X_s\gets 
Y_s/V_s}\Otimes _{\calD_{ Y_s/V_s}}\calD_{ Y_s\to 
X_s/V_s}\mBig)\to\cal H^{\codim_XY}_{Y}(\cal
D_{ X_s/V_s})\simeq\\
\simeq
\bfR \Gamma_Y(\calD_{ X_s/V_s})[\codim_XY]
\leqno{ \Adj _{*,s}^i(\calD_{X_s/V_s}):}
\endlignes$$
est un isomorphisme.
\end{theo}

\demo
La question est locale. Soit un triplet $(u,z,y)$ comme  précédemment, où $x=(z,y)$  est adapté à $u$, et soit $(u_s,z_s,y_s)$ sa réduction.
Le morphisme $ \Adj _{*,s}^i$ est défini par: 
$$P_sdy_su_s(dx_s)^{*}\otimes u_sQ_s\mapsto P_s[1/z_s]Q_s.$$ On a un isomorphisme de $\calD_{
X_s/V_s}$-modules à droite :
$$\cal H^q_{Y}(\cal
O_{ X_s/V_s})\otimes_{\cal
O_{ X_s/V_s}}\cal
D_{ X_s/V_s}\simeq\cal H^q_{Y}(\cal
D_{ X_s/V_s}).$$ Tout tenseur s'écrit de manière unique comme somme finie:
$$\sum \delta_\alpha\otimes\Delta_{x_s}^\alpha,$$ où $\delta_\alpha$ est un élément de $\cal H^q_{Y}(\cal
O_{ X_s/V_s})$. L'action à gauche est définie par:
$$\Delta^\beta_{x_s}(\delta_\alpha\otimes\Delta_{x_s}^\alpha)= \sum_{0\leq k\leq
\beta}\Delta_{x_s}^{\beta-k}\delta_\alpha\otimes\Delta_{x_s}^k\Delta_{x_s}^\alpha.$$ D'autre part, tout élément $\delta_\alpha$ est de la forme
$P_\alpha([1/z_s])$ pour un opérateur différentiel $P_\alpha$.  Mais:
$$\Delta_{x_s}^\alpha([1/z_s]\otimes 1)= \sum_{0\leq k\leq \alpha}\Delta_{x_s}^{\alpha-k}[1/z_s]\otimes\Delta_{x_s}^k.$$ Cela montre par récurrence sur la
longueur de
$\alpha$ que $ \Adj _{*,s}^i(\calD_{X_s/V_s})$ est un morphisme surjectif.

Comme on a l'égalité par construction:
 $$dy_su_s(dx_s)^{*}\otimes u_sb_\alpha\Delta_{y_s}^\alpha=dy_su_sb_\alpha\Delta_{y_s}^\alpha(dx_s)^{*}\otimes u_s,$$
tout tenseur  s'écrit comme une somme finie de tenseurs $\sum P_\alpha dy_su_s(dx_s)^{*}\otimes u_s\Delta^\alpha_{z_s}$ où $P_\alpha$ est un opérateur
différentiel.
L'image du tenseur:  
$$P_\alpha dy_su_s(dx_s)^{*}\otimes u_s\Delta^\alpha_{z_s}$$ par le morphisme $ \Adj _{*,s}^i$ 
est égale à: $$P_\alpha([1/z_s]\otimes\Delta^\alpha_{z_s}).$$ Tout opérateur différentiel $P_\alpha$ s'écrit de manière unique comme une somme
$\sum_{\beta}P_{\alpha,\beta}\Delta_{z_s}^\beta$ où $P_{\alpha,\beta}$ est un opérateur différentiel en $\Delta_{y_s}^\alpha$.  Soit:
$$\sum_{\alpha,\beta}P_{\alpha,\beta}\Delta_{z_s}^\beta dy_su_s(dx_s)^{*}\otimes u_s\Delta^\alpha_{z_s}$$ un tenseur dont l'image par $ \Adj ^i_{*,s}$ est
nulle. Nécessairement, $P_{0,0}([1/z_s]\otimes1)$ est nul, l'élément  $P_{0,0}$ appartient à l'idéal $(z_1,\dots, z_q)$, et cela entraîne que le tenseur: 
$$P_{0,0} dy_su_s(dx_s)^{*}\otimes u_s1= dy_su_s(dx_s)^{*}\otimes u_sP_{0,0}$$ est nul. Par récurrence sur l'ordre en $\Delta_{z_s}$ on voit que le tenseur:
$$\sum_{\alpha,\beta}P_{\alpha,\beta}\Delta_{z_s}^\beta dy_su_s(dx_s)^{*}\otimes u_s\Delta^\alpha_{z_s}$$ est nul. Le morphisme $ \Adj _{*,s}^i(\cal
D_{X_s/V_s})$  est injectif.
\enddemo

\begin{coro}Si $\cal M_s$ est un complexe de $ \Dplus (\calD_{ X_s/V_s})$ à cohomologie $\cal O_{\cal X_s/V_s}$-quasi-cohérente, alors le morphisme:
$$
\DeuxLignes  
i_*\mBig(\calD_{ X_s\gets 
Y_s/V_s}\Otimes_{\calD_{ Y_s/V_s}}\calD_{ Y_s\to \cal
X_s/V_s}\mBig)\Lotimes_{\calD_{ X_s/V_s}}\cal M_s\to\\\to\bfR \Gamma_Y(\calD_{ X_s/V_s})\Lotimes_{\cal
D_{ X_s/V_s}}\cal M_s[\codim_XY]
\to\bfR \Gamma_Y(\cal M_s)[\codim_XY]
\endlignes
$$
est un isomorphisme.
\end{coro}

\demo C'est une conséquence du théorème précédent et du lemme du way-out foncteur.
\enddemo

 Nous n'utiliserons pas dans la suite de cet article le théorème précédent et son corollaire.

\bigskip
À présent, nous allons construire, pour une
immersion fermée
$i: Y\to X$ de schémas lisses sur le corps
$k$, le morphisme de foncteurs  $ \Adj _*^i$ 
entre les foncteurs $i_*^{\diff}i^*_{\diff}$ et $\bfR \Gamma_Y[\codim_XY]$,
annoncé au début de cette section.  On rappelle que, sous les conditions  du théorème de pureté
\ref{pur-coh}, on a:
$$\bfR \Gamma_Y(\calDdag_{\calXdag /V})[\codim_XY]\simeq\cal H^q_{Y}(\calDdag_{\calXdag /V})\,,\hbox{ avec }q=\codim _XY\,.$$


\begin{lemm}\label{adjp}Si $\cal P\daginf $ est un $\calDdag_{\Xdaginf /V}$-module à gauche spécial et {\bf plat}, il existe un morphisme:
$$i_{*,0}^{\diff}i_{\diff}^{*,0}\cal P\daginf \to \cal H^{\codim_XY}_{Y}(\calDdag_{ \Xdaginf /V})\otimes_{\calDdag_{ \Xdaginf /V}}\cal P\daginf $$ de $\calDdag_{\Xdaginf /V}$-modules à gauche spéciaux.
\end{lemm}

\demo
Soit $\calUdag $ un ouvert affine du site $ \Xdaginf $ et soit $\calWdag $ un relèvement  plat de la trace $W$ de $U$ dans $Y$. 
Le morphisme $ \Adj _*^i(\calDdag_{\calUdag /V})$ pour le couple $(\calUdag , \calWdag )$ fournit un morphisme de $\calDdag_{\cal
U\dag/V}$-modules à gauche par produit tensoriel:
$$\scriptwd=1.5em
i_*\mBig(\calDdag_{\calUdag \gets \cal
W\dag/V}\Otimes_{\calDdag_{\calWdag /V}}\calDdag_{\calWdag \to\cal
U\dag/V}\mBig)\Otimes_{\calDdag_{\calUdag /V}}\cal P\dag_{\calUdag }\to\cal H^{\codim_XY}_Y(\calDdag_{\cal
U\dag/V})\Otimes_{\calDdag_{\calUdag /V}}\cal P\dag_{\calUdag }.$$
Par construction, le complexe de
gauche est la valeur de
$i_{*,0}^{\diff}i_{\diff}^{*,0}\cal P\daginf $ sur l'ouvert $\calUdag $.   Il faut voir que
ce morphisme  commute aux morphismes du site. 
Soient  $r\dag:\cal U'{}\dag \to \calUdag $ un morphisme du site, $u$ un relèvement de l'inclusion $ W\to U$ et 
$\cal W'{}\dag $ le produit fibré  $\cal
U'{}\dag \times_{\calUdag }\calWdag $, de sorte qu'on a un diagramme commutatif:
$$\def\quad{\hskip0.5ex}\matrix{\cal W'{}\dag &\stackrel{r\dag}{\too}&\calWdag \cr \noalign{\kern2pt}
\lscript{u'}\downarrow&&\lscript{u}\downarrow\cr\noalign{\kern-2pt}
\cal U'{}\dag &\stackrel{r\dag}{\too}&\calUdag. }$$ Si $(z,y)$
est un couple adapté à $u$,  le couple $(z',y')$, avec  
$z':= r^*z $ et $ y':= r^*y$,
 est par construction adapté au relèvement $u'$. Il en résulte  un diagramme commutatif:
$$\def\quad{\hskip0.5ex}
\matrix{r^{-1}i_*\mBig(\calDdag_{\calUdag \gets \cal
W\dag/V}\otimes_{\calDdag_{\calWdag /V}}\calDdag_{\calWdag  \to\cal
U\dag/V}\mBig)&\too&r^{-1}\cal H^{\codim_XY}_Y(\calDdag_{\calUdag /V})\cr\downarrow&&\downarrow\cr i_*\mBig(\calDdag_{\cal U'{}\dag \gets \cal
W'{}\dag /V}\otimes_{\calDdag_{\cal W'{}\dag /V}}\calDdag_{\cal W'{}\dag \to \cal
U'{}\dag /V}\mBig)&\too&\cal H^{\codim_XY}_Y(\calDdag_{\cal U'{}\dag /V}),}$$ où les lignes sont les morphismes $ \Adj _*^i$ et les colonnes sont les
morphismes:
$$Pdyu(dx)^{*}\otimes uQ\mapsto r^*Pr^{*-1}dy'u'(dx')^{*}\otimes u'r^*Qr^{*-1}$$ et   $$P[1/z]Q\mapsto r^*Pr^{*-1}[1/z']r^*Qr^{*-1}.$$ Par produit
tensoriel, on obtient que le morphisme commute aux restrictions. On obtient comme cela un morphisme défini localement.
\enddemo

\begin{lemm}Si $X$ est {\bf séparé} et si $\cal P\daginf $ est un complexe $\calDdag_{\Xdaginf /V}$-modules à gauche spéciaux et {\bf plats}, il existe un
morphisme de la catégorie
$ \rmD((\calDdag_{ \Xdaginf /V}, \Sp)\Mod)$:
$$\cal H^{\codim_XY}_{Y}(\calDdag_{ \Xdaginf /V})\otimes_{\calDdag_{ \Xdaginf /V}}\cal P\daginf [-\codim_XY]\to \cal P\daginf .$$ 
\end{lemm}

\demo L'extension $V_{\Xdaginf }\to \calDdag_{ \Xdaginf /V}$ est plate, et un injectif sur $\calDdag_{ \Xdaginf /V}$ reste donc  injectif sur $V_{\Xdaginf }$. 
Prenant une résolution injective de $\calDdag_{ \Xdaginf /V}$ par des $\calDdag_{ \Xdaginf /V}$-modules à droite, on réalise 
le morphisme:
$$\cal H^{\codim_XY}_{Y}(\calDdag_{ \Xdaginf /V})\otimes_{\calDdag_{ \Xdaginf /V}}\cal P\daginf [-\codim_XY]\to \cal P\daginf $$ dans la catégorie $ \rmD((V_{\Xdaginf }\Mod ))$. Rien ne dit alors que c'est un
morphisme de la catégorie $ \rmD((\calDdag_{ \Xdaginf /V}, \Sp)\Mod)$, car a priori  on ne sait pas construire une résolution de $\calDdag_{ \Xdaginf /V}$ par des bimodules qui soit injective à droite. C'est pour construire une résolution par des bimodules que l'hypothèse de séparation
intervient. Soit $B_\alpha$ un recouvrement de $X$ par des ouverts affines au-dessus desquels $Y$ est défini par des équations $z_1,\dots,z_q$. Notons
$j_{\alpha,i}: B_{\alpha,i}:= B_\alpha-V(z_i)\to X$ l'inclusion canonique. Considérons le complexe de \v Cech $C^{\bullet}(\cal B_{\alpha,i}, \calDdag_{ \Xdaginf /V})$:
$$0\to\calDdag_{ \Xdaginf /V}\to \prod_{\alpha,i}j_{\alpha,i_*}j_{\alpha,i}^{-1}\calDdag_{ \Xdaginf /V}\to\cdots $$ En vertu du théorème de pureté \ref{pur-coh} et du corollaire d'acyclicité \ref{acy-dif}, si $X$ est
\em{séparé}, le complexe
$C^{\bullet}( \cal B_{\alpha,i}, \calDdag_{ \Xdaginf /V})$ est une résolution 
par des $(\calDdag_{ \Xdaginf /V}, \calDdag_{ \Xdaginf /V})$-bimodules
de 
$\cal H^{\codim_XY}_{Y}(\calDdag_{ \Xdaginf /V})[-\codim_XY]$.  Le complexe tronqué $\sigma_{\leq \codim_XY}C^{\bullet}(\cal B_{\alpha,i}, \calDdag_{ \Xdaginf /V})$ est une résolution du bimodule  $\cal H^{\codim_XY}_{Y}(\calDdag_{ \Xdaginf /V})[-\codim_XY]$, et le morphisme de complexes qui provient du morphisme canonique $\sigma_{\leq \codim_XY}C^{\bullet}(\cal B_{\alpha,i}, \calDdag_{ \Xdaginf /V})\to \calDdag_{ \Xdaginf /V}$, à savoir:
$$\sigma_{\leq \codim_XY}C^{\bullet}(\cal B_{\alpha,i}, \calDdag_{ \Xdaginf /V})\otimes_{\calDdag_{ \Xdaginf /V}}\cal P\daginf \to\cal P\daginf \,,$$ réalise le morphisme du lemme.
\enddemo

Soit  maintenant $\calMdaginf $  un complexe de $ \Db ((\calDdag_{\Xdaginf /V}, \Sp)\Mod)$. Il admet une résolution $\cal P\daginf $ par des
modules spéciaux plats,
d'où résulte par le lemme précédent  un morphisme canonique:
$$\cal H^{\codim_XY}_{Y}(\calDdag_{ \Xdaginf /V})\Lotimes_{\calDdag_{ \Xdaginf /V}, \Sp}\calMdaginf \to \calMdaginf [\codim_XY].$$  Comme $\cal H^{\codim_XY}_{Y}(\calDdag_{ \Xdaginf /V})$ est à support dans $Y,$ ce morphisme se factorise nécessairement par $\bfR \Gamma_Y(\cal
M\daginf )\to \calMdaginf $. Il en résulte  le morphisme:
$$\cal H^{\codim_XY}_{Y}(\calDdag_{ \Xdaginf /V})\Lotimes_{\calDdag_{ \Xdaginf /V}, \Sp}\calMdaginf \to\bfR \Gamma_Y(\calMdaginf )[\codim_XY],$$ qui composé avec  morphisme du lemme \ref{adjp} donne le
morphisme:\glossary{$ \Adj _*^i(\calMdaginf )$}
$$i_*^{\diff}i_{\diff}^*\calMdaginf \to \bfR \Gamma_Y(\calMdaginf )[\codim_XY].\leqno  \Adj _*^i(\calMdaginf ):$$

\medskip\begin{defi}\label{adj}Soit $i:Y\hookrightarrow X$ une immersion fermée de schémas lisses et {\bf séparés} sur $k$. On définit le morphisme $ \Adj _*^i(\calMdaginf )$
pour tout complexe $\calMdaginf $
de la catégorie 
$ \Db ((\calDdag_{\Xdaginf /V}, \Sp)\Mod)$ comme le morphisme précédent  de la catégorie  $ \Dplus ((\calDdag_{\Xdaginf /V}, \Sp)\Mod)$:
$$i_{*}^{\diff}i_{\diff}^{*}\calMdaginf \to \bfR \Gamma_Y(\calMdaginf )[\codim_XY].\leqno  \Adj _*^i(\calMdaginf ):$$
\end{defi}
\begin{Rema}Le morphisme $ \Adj _*^i(\calDdag_{\calXdag /V})$ \em{n'est pas}  un isomorphisme, même sur $K$. Mais
nous allons montrer qu'il induit un isomorphisme dans le cas géométrique $\calMdaginf = \cal O_{\Xdaginf /V}$.
\end{Rema}

\begin{theo} \label{coh-loc}Si $X$ est un schéma lisse et séparé sur $k$, le morphisme $\Adj _*^i(\cal O_{\Xdaginf /V}):$
$$i_*^{\diff}i_{\diff}^*\cal O_{\Xdaginf /V}\to \bfR \Gamma_Y(\cal O_{\Xdaginf /V})[\codim_XY]\simeq \cal H^{\codim_XY}_Y(\cal O_{\Xdaginf /V})$$ est un
isomorphisme de $\calDdag_{\Xdaginf /V}$-modules à gauche  spéciaux.
\end{theo}

\demo  En vertu du théorème \ref{inv-tri}, on a un isomorphisme canonique de $\calDdag_{\Ydaginf /V}$-modules à gauche spéciaux:
$$\cal O_{\Ydaginf /V}\simeq i_{\diff}^*\cal O_{\Xdaginf /V}.$$ 
Il suffit donc de montrer
que le morphisme: $$i_*^{\diff}\cal O_{\Ydaginf /V}\to \cal H^{\codim_XY}_Y(\cal O_{\Xdaginf /V})$$ est un isomorphisme de
$\calDdag_{\Xdaginf /V}$-modules à gauche spéciaux. La question est locale,  et l'on peut supposer que $X$ est affine.
Soient $\calXdag $ un relèvement  plat de $X$ et $\calYdag $ un relèvement 
 plat de $Y$. Par construction, il suffit de montrer que le morphisme:
$$\calDdag_{\calXdag \gets \cal
Y\dag/V}\otimes_{\calDdag_{\calYdag /V}}\cal O_{\calYdag /V}\to \cal H^{\codim_XY}_Y(\cal O_{\calXdag /V})$$ est un isomorphisme de $\calDdag_{\calXdag /V}$-modules à gauche. Mais le $\calDdag_{\calXdag /V}$-module  $\cal H^{\codim_XY}_Y(\cal O_{\calXdag /V})$ est localement de type fini, dont les générateurs locaux sont donnés par la classe $[1/z]$
pour  un système  $z= (z_1,\dots,z_q)$ de relèvements des équations de $Y$. D'autre part, le $\calDdag_{\calXdag /V}$-module  à gauche $\calDdag_{\calXdag \gets \cal
Y\dag/V}\otimes_{\calDdag_{\calYdag /V}}\cal O_{\calYdag /V}$ est engendré localement par $dyu(dx)^{*}\otimes1$ pour un relèvement $u$ de l'immersion
$i$ et un système $x = (y,z)$ adapté à $u$. On voit que le morphisme $ \Adj _*^i(\cal O_{\calXdag /V})$ envoie 
le générateur $dyu(dx)^{*}\otimes1$ au générateur $[1/z]$ dont l'annulateur annule aussi $dyu(dx)^{*}\otimes1$. 
Le morphisme $ \Adj _*^i(\cal O_{\calXdag /V})$ est nécessairement un isomorphisme.
\enddemo

\bigskip

De la même façon,  on construit le morphisme $ \Adj _*^i(\omega_{\Xdaginf /V})$ et on obtient le théorème qui suit. 

\begin{theo} \label{coh-loc'}Le morphisme $\Adj _*^i(\omega_{\Xdaginf /V}):$
$$i_*^{\diff}i_{\diff}^*\omega_{\Xdaginf /V}\to \bfR \Gamma_Y(\omega_{\Xdaginf /V})[\codim_XY]\simeq\cal H^{\codim_XY}_Y(\omega_{\Xdaginf /V})$$ est un
isomorphisme de $\calDdag_{\Xdaginf /V}$-modules  à droite spéciaux.
\end{theo}
\begin{Rema} Les résultats précédents valent sur l'extension 
$V\to K$ et nous aurons surtout à les utiliser dans ce cas-là dans l'article présent. \end{Rema}

\subsection{La comparaison entre la cohomologie locale et l'image directe dans 
le cas d'une immersion fermée. Le morphisme de foncteurs $ \Adj _i^*$}

\begin{lemm}\label{Tors=0}Soient $i:Y\hookrightarrow X$ une immersion fermée de schémas lisses  sur $k$, et $\calYdag $ et $\calXdag $  des relèvements
 plats de
$Y$ et $X$. Alors, les modules $\cTor_j^{\calDdag_{\calXdag /V}}(\calDdag_{\calYdag \to \calXdag/V},\calDdag_{\calXdag \gets\cal
Y\dag/V})$ sont nuls pour $j\neq \codim_XY$.
\end{lemm}

\demo
La question étant  locale, on peut supposer que $X$ est affine. Soit $u$ un relèvement de l'inclusion $i$ et soit $x=(z_1,\dots,z_q,y)$ un système adapté à $u$.
Considérons le complexe  de Koszul $K(\calDdag_{\calXdag /V},z_1,\dots,z_q)$, où l'action des $z$ se fait à gauche, augmenté  vers $\calDdag_{\cal
Y\dag\to \calXdag /V}$ par $P\mapsto uP$.  
C'est est une résolution par des  $ \calDdag_{\calXdag /V}$-modules 
à droite de $\calDdag_{\calYdag \to\calXdag /V},$ 
en vertu du théorème \ref{rg1} et du lemme \ref{ima-invd}. 
Il suffit maintenant de démontrer que le complexe $K(\calDdag_{\calXdag /V},z_1,\dots,z_q)\otimes_{\calDdag_{\calXdag /V}}\calDdag_{\calXdag\gets\calYdag /V}$ est concentré en degré $-\codim_XY$. Mais c'est un complexe de $\cal O_{\calXdag /V}$-modules à gauche, et  l'application
symbole total:
$$\hbox{\Large$\sigma$}\mBig(\sum_\alpha(-1)^{|\alpha|}\Delta_x^\alpha  b_\alpha\mBig):= \sum_\alpha(-1)^{|\alpha|}\xi^\alpha  b_\alpha$$ est un isomorphisme $\cal O_{\calXdag/V}$-linéaire, en vertu du théorème du symbole total \ref{sym-tot}. Il suffit donc de démontrer que le complexe: $$K\big((\cal O_{\calXdag/V}[\xi_1,\dots,\xi_n])\dag,z_1,\dots,z_q\big)\Otimes_{(\cal O_{\calXdag /V}[\xi_1,\dots,\xi_n])\dag}\big(\cal O_{\calYdag /V}[\xi_1,\dots,\xi_n]\big)\dag$$ est
concentré en degré
$-\codim_XY$, mais comme les objets de ce complexe sont sans $\goth m$-torsion et de type
fini sur $(\cal O_{\calXdag /V}[\xi_1,\dots,\xi_n])\dag$,
suffit de montrer que la réduction modulo $\goth m$ de ce complexe est concentrée en degré
$-\codim_XY$. Or, cela résulte du fait que l'immersion de $Y\times_k A_k^n$ dans le fibré cotangent $T^*X$ est une immersion régulière de codimension $\codim_XY$.
\enddemo

\medskip
Nous allons construire un isomorphisme canonique de $(\calDdag_{\calYdag /V},\calDdag_{\calYdag /V})$-bimodules:
$$\calDdag_{\calYdag /V}\simeq\cTor_{\codim_XY}^{\calDdag_{\calXdag /V}}(\calDdag_{\calYdag \to \calXdag/V},\calDdag_{\calXdag \gets\cal
Y\dag/V}).$$

\begin{lemm}\label{Tor=D}Soit $u$ un relèvement de l'immersion $i$ et soit $x=(z,y)$ un système de fonctions adapté à $u$. Alors, le morphisme $P\mapsto P\circ dyu(dx)^{*}$ est un
isomorphisme de $\calDdag_{\calYdag /V}$-modules à droite de $\calDdag_{\calYdag /V}$ sur le noyau $\Ker(z,\calDdag_{\calXdag \gets\cal
Y\dag/V})$ du morphisme: 
$$\matrix{
\calDdag_{\calXdag \gets\calYdag /V}&\too&(\calDdag_{\calXdag \gets\cal
Y\dag/V})^q\cr\noalign{\kern2pt}
Q&\longmapsto& (z_1Q,\dots,z_qQ).
}$$
\end{lemm}

\demo
Comme les actions de $\calDdag_{\calYdag /V}$ et de $z_1,\dots,z_q$ commutent, l'image de $\calDdag_{\calYdag /V}$ par
le morphisme $P\mapsto P\circ dyu(dx)^{*}$ est contenue dans le noyau $\Ker(z,\calDdag_{\calXdag \gets\cal
Y\dag/V})$. Le module de transfert $\calDdag_{\calXdag \gets\calYdag /V}$ est un $\calDdag_{\calXdag /V}$-module à gauche engendré par $dyu(dx)^{*}$. Un élément $Pdyu(dx)^{*}$ opère
par $$f(x)dxPdyu(dx)^{*}= u(\!\trans P(f(x))dy,$$ et $z_iPdyu(dx)^{*}$ opère par  $f(x)dxz_iPdyu(dx)^{*}= u(\trans P(z_if(x))dy$. Si on développe $\trans P=\sum_\alpha
a_{\alpha,\alpha'}\Delta_y^\alpha\Delta_z^{\alpha'}$,  on trouve que $z_iPdyu(dx)^{*}$ est nul pour $i=1,\dots,q$  si seulement si les coefficients
$a_{\alpha,\alpha'}$ pour
$\alpha'$ non nul appartiennent au noyau  de $u$. Cela entraîne  que le morphisme $\calDdag_{\calYdag /V}\to \Ker(z,\calDdag_{\calXdag \gets\cal
Y\dag/V})$ est surjectif.
 D'autre part, si l'image d'un opérateur $\sum_{\alpha}a_\alpha\Delta_y^\alpha$ est nulle,   l'action sur les monômes
$y^\alpha$ de cet opérateur est nulle, et l'opérateur est donc  nul. Le morphisme du lemme est injectif.
\enddemo

En vertu du lemme précédent, pour tout triplet $(u,z,y)$ on dispose d'un isomorphisme de $\calDdag_{\calYdag /V}$-modules à droite:
$$\calDdag_{\calYdag /V}\simeq\cTor_{\codim_XY}^{\calDdag_{\calXdag /V}}(\calDdag_{\calYdag \to \calXdag/V},\calDdag_{\calXdag \gets\cal
Y\dag/V}).\leqno(*):$$ Soit $v$ un autre relèvement de $i$ qui donne la présentation de $\calDdag_{\calXdag /V}$-modules à droite
$$\calDdag_{\calXdag /V}\stackrel{v}{\to}\calDdag_{\calYdag \to \calXdag/V}.$$ Mais $u= vg$ pour un élément $g$ du groupe $\cal G_{\calXdag }$ dont l'action fournit un diagramme commutatif:
$$\def\quad{\hskip0pt}\hskip-0.5cm
\matrix{&&K(\calDdag_{\calXdag /V},z_1,\dots,z_q)
\smash{\Otimes_{\calDdag_{\calXdag /V}}}\calDdag_{\calXdag \gets\cal
Y\dag/V}&&\cr&\nearrow&&\nwarrow\cr
\hdecale{0.5cm}{\calDdag_{\calYdag /V}}
&&g\ \raise2pt\hbox{\smash{\vardownarrow35pt}}\kern1em
&&
\hskip-2cm\cTor_{\codim_XY}^{\calDdag_{\calXdag /V}}(\calDdag_{\calYdag \to \calXdag/V},\calDdag_{\calXdag \gets\cal
Y\dag/V})\,.&\cr\noalign{\kern2pt}
&\searrow&&\swarrow\cr
&&K(\calDdag_{\calXdag /V},z'_1,\dots,z'_q)\Otimes_{\calDdag_{\calXdag /V}}\calDdag_{\calXdag \gets\cal
Y\dag/V}&&\cr
}$$ 
Cette description montre que le morphisme $(*)$ précédent ne dépend pas du triplet choisi $(u,z,y)$ et fournit donc  un isomorphisme canonique global
de $\calDdag_{\calYdag /V}$-modules à droite. En inversant les rôles de $\calDdag_{\calXdag \gets\cal
Y\dag/V}$ et de $\calDdag_{\calYdag \to\calXdag/V}$, on montre que $(*)$ est un isomorphisme canonique global
de $\calDdag_{\calYdag /V}$-modules à gauche et donc le morphisme $(*)$ est un isomorphisme de bimodules.

\begin{defi}Soit $i:Y\hookrightarrow X$ une immersion fermée de schémas  lisses sur $k$, et $\calYdag $ et $\calXdag $  des relèvements plats de
$Y$ et $X$. On définit l'isomorphisme $ \Adj ^*_i(\calDdag_{\calYdag }/V)$  de $(\calDdag_{\calYdag /V},\calDdag_{\cal
Y\dag/V})$-bimodules:
$$\TroisLignes  
\leqno{ \Adj ^*_i(\calDdag_{\calYdag /V}):}
\calDdag_{\calYdag /V}[\codim_XY]\to
\\
\to\cTor_{\codim_XY}^{\calDdag_{\calXdag /V}}(\calDdag_{\calYdag \to \calXdag/V},\calDdag_{\calXdag \gets\cal
Y\dag/V})[\codim_XY]\simeq
\\\simeq
\calDdag_{\calYdag \to \calXdag/V}\Lotimes_{\calDdag_{\calXdag /V}}\calDdag_{\calXdag \gets\cal
Y\dag/V}
\endlignes$$ comme l'isomorphisme $(*)$ précédent.
\end{defi}

\begin{coro}\label{adj-imm}Soit
$i:Y\hookrightarrow X$ une immersion fermée de schémas lisses et séparés sur $k$. Pour tout complexe $\calMdaginf $ de la catégorie  $ \Dmoins (\calDdag_{\Ydaginf /V}, \Sp)$ il existe un isomorphisme canonique dans la catégorie  $ \Dmoins (\calDdag_{\Ydaginf /V}, \Sp)$:\glossary{$ \Adj ^*_i(\calMdaginf )$}
$$ \calMdaginf [\codim_XY]\simeq i^*_{\diff}i_*^{\diff}\calMdaginf .\leqno  \Adj ^*_i(\calMdaginf ):$$
\end{coro}

\demo
Soit $W$ un ouvert affine de $Y$ qui est la trace d'un ouvert affine $U$ de $X$. Soient $\calWdag $ et $\calUdag $ des relèvements de $W$ et de $U$. 
Si $\calMdaginf $ est un $\calDdag_{\Ydaginf /V}$-module spécial, la valeur du complexe $i^*_{\diff}i_*^{\diff}\calMdaginf $ sur l'ouvert    $\calWdag $ est par construction  le complexe:
$$\mBig(\calDdag_{\calWdag \to \cal
U\dag/V}\Lotimes_{\calDdag_{\calUdag /V}}\calDdag_{\calUdag \gets\cal
W\dag/V}\mBig)\Lotimes_{\calDdag_{\calWdag /V}}\cal M_{\calWdag }$$  dont la cohomologie
est concentrée, en vertu des lemmes \ref{Tors=0} et \ref{Tor=D}, en degré $-\codim_XY$, et dont la valeur est 
le faisceau $$\cTor_{\codim_XY}^{\calDdag_{\calUdag /V}}(\calDdag_{\calWdag \to \cal
U\dag/V},\calDdag_{\calUdag \gets\cal
W\dag/V})\otimes_{\calDdag_{\calWdag /V}}\cal M\dag_{\calWdag }.$$ Il en résulte  un isomorphisme de $\calDdag_{\calYdag /V}$-modules à gauche induit
par  l'isomorphisme 
$ \Adj ^*_i(\calDdag_{\calYdag /V})$:
$$\cal M\dag_{\calWdag }\to \cTor_{\codim_XY}^{\calDdag_{\calXdag /V}}(\calDdag_{\calYdag \to \calXdag/V},\calDdag_{\calXdag \gets\cal
Y\dag/V})\otimes_{\calDdag_{\calWdag /V}}\cal M\dag_{\calWdag },$$ qui commute aux restrictions du site $\Ydaginf $. On a comme cela un 
isomorphisme de 
$\calDdag_{\Ydaginf /V}$-modules spéciaux:
$$\calMdaginf [\codim_XY]\simeq i^*_{\diff}i_*^{\diff}\calMdaginf, \leqno  \Adj ^*_i(\calMdaginf ):$$ qui se dérive trivialement.
\enddemo

On en déduit  le corollaire suivant qui donne un cas où la cohomologie locale commute avec  l'image inverse.

\begin{coro}Soit $i:Y\hookrightarrow X$ une immersion fermée de schémas lisses  sur $k$. Il existe un isomorphisme canonique 
$\calDdag_{\Ydaginf /V}$-modules à gauche spéciaux:
$$ \cal O_{\Ydaginf /V}\simeq i_{\diff}^*\bfR \Gamma_Y(\cal O_{\Xdaginf /V}).$$
\end{coro}

\demo C'est une conséquence de l'isomorphisme $ \Adj _*^i(\cal O_{\Xdaginf /V})$ et de l'isomorphisme du théorème \ref{coh-loc}.
\enddemo

Autrement dit, les modules spéciaux à gauche: 
$$\cal O_{\Ydaginf /V}\text{\quad et\quad} \cal H^{\codim_XY}_Y(\cal O_{\Xdaginf /V})$$ se correspondent par les foncteurs
image directe et inverse dans le cas d'une immersion fermée. 

De même, les modules spéciaux à droite $\omega_{\Ydaginf /V}$ \glossary{$\omega_{\Xdaginf /V}$}et $\cal H^{\codim_XY}_Y(\omega
_{\Xdaginf /V})$ se correspondent par les foncteurs image directe et inverse dans le cas d'une immersion fermée.
\begin{Rema} Dans le cas du site infinitésimal $p$-adique formel la correspondance précédente n'a pas lieu, aussi
le module $i^{\diff}_*\cal O_{Y^\wedge /V}$ se substitue avantageuseemnt au faisceau de cohomologie locale de $Y$ à valeur dans le faisceau structural $\cal O_{X^\wedge /V}$.\end{Rema}
\section{La fonctorialité de la cohomologie de de Rham $p$-adique et le foncteur de dualité}
Dans ce paragraphe, nous allons démontrer que la cohomologie de de Rham $p$-adique sur le corps des fractions $K$ varie
fonctoriellement de façon contravariante pour les morphismes $f$. La fonctorialité sur $V$ nécessite des changements dans les définitions. On note encore $Z$ le produit $Y\times_kX$.

\subsection{La compatibilité des foncteurs image  inverse  topologique et différentielle}

\begin{theo}\label{dr-dir} Soit $f: Y\to X$ un morphisme de schémas  lisses et {\bf séparés} sur $k$, et soit $\calMdaginf $ un complexe
de \relax{$ \Db ((\calDdag_{\Xdaginf /K}, \Sp)\Mod)$}. Il existe un morphisme $\dR(f/K,\calMdaginf )$ fonctoriel en $\calMdaginf $:
$$f^{-1}\bfR \cHom\goodSub{10pt}{6pt}{-6mm}{\calDdag_{\Xdaginf /K},\Sp}(
\cal O_{\Xdaginf /K},\calMdaginf )\to\bfR \cHom\goodSub{10pt}{6pt}{-6mm}{\calDdag_{\Ydaginf /K},\Sp}(
\cal O_{\Ydaginf /K},f^*_{\diff}\calMdaginf )$$ de complexes de $ \rmD(K_Y)$.\glossary{$\dR(f/K,\calMdaginf )$}
\end{theo}

\demo
Comme $f^{-1}= i^{-1}\circ p^{-1}$ et $f^*_{\diff}= i_{\diff}^*\circ p^*_{\diff}$, il suffit de considérer le cas d'une projection puis celui d'une immersion fermée. Le foncteur
$p^{*,0}_{\diff}$ étant exact, le cas d'une projection est élémentaire. Il suffit de construire un morphisme entre faisceaux de Zariski:
$$p^{-1}\cHom_{\calDdag_{\Xdaginf /K}}(
\calNdaginf ,\calMdaginf )\to\cHom_{\calDdag_{Z\daginf /K}}(
p^{*,0}_{\diff}\calNdaginf ,p^{*,0}_{\diff}\calMdaginf ),$$ si $\calNdaginf $ et $ \calMdaginf $ sont deux $\calDdag_{\Xdaginf /K}$-modules à
gauche  spéciaux. Si
$Y$ et
$X$ sont affines et $\calYdag $ et $\calXdag $ sont des relèvements $\dagger$-adiques  plats, le morphisme:
$$\DeuxLignes  
p^{-1}\cHom_{\calDdag_{\calXdag /K}}\big(
\cal N\dag_{\calXdag },\cal M\dag_{\calXdag }\big)\to
\\
\to\cHom\goodSub{10pt}{6pt}{-6mm}{\calDdag_{\calYdag \times_V\calXdag /K}}\mBig(
\calD_{\calYdag \times_V\calXdag \to\calXdag /K}\Otimes_{p^{-1}\calDdag_{\calXdag /K}}p^{-1}\cal N\dag,\calD_{\calYdag \times_V\calXdag\to\calXdag /K}\Otimes_{p^{-1}\calDdag_{\calXdag /K}}p^{-1}\calMdag \mBig)
\endlignes$$ 
provient de la fonctorialité du produit tensoriel. Ces morphismes
locaux se recollent de façon naturelle pour fournir un morphisme global dans le cas d'une projection. Prenant une résolution injective $\calMdaginf \to
\calIdaginf  $ de $\calMdaginf $ par des $\calDdag_{\calXdag /K}$-modules spéciaux puis une résolution injective de $p^*_{\diff}\cal
I\daginf $ par des
$\calDdag_{Z\daginf /K}$-modules spéciaux, on trouve le morphisme du théorème dans le cas d'une projection. En fait, dans le cas d'une projection
le résultat est vrai sur $V$, même sans l'hypothèse de séparation. 

\medskip 
Le morphisme du théorème précédent est défini
chaque fois que le foncteur $f^{*,0}_{\diff}$ est exact.

\medskip
Le cas d'une immersion ouverte est élémentaire mais le cas d'une immersion
fermée  où le foncteur image inverse n'est plus exact est nettement plus profond, délicat et nécessite l'extension $V\to K$ et l'hypothèse de séparation. Pour cela, nous allons d'abord définir le foncteur de dualité pour la catégorie des modules spéciaux et montrer la compatibilité du foncteur dualité avec le morphisme image directe pour une immersion fermée.

\subsection{Le foncteur dualité et le théorème de dualité pour une immersion fermée}
Soit $i: Y\to X$ une immersion fermée
de schémas lisses sur
$k$. 

\begin{lemm}\label{dr-coh}   Si $\cal P\daginf $ est un $\calDdag_{\Xdaginf /V}$-module à gauche spécial et {\bf plat},  le faisceau 
$\cTor_{j}^{\calDdag_{\Ydaginf /V}, \Sp}(\omega_{Y\dag_{\inf/V}}, i_{\diff}^{*,0}\cal P\daginf )$ est alors   nul pour $j\neq 0$, et  il existe un isomorphisme
canonique de faisceaux de Zariski:
$$\omega_{Y\dag_{\inf/V}}\otimes_{\calDdag_{\Ydaginf /V}}i_{\diff}^{*,0}\cal P\daginf \simeq \cal H^{\codim_XY}_{Y}(\omega_{
\Xdaginf /V})\otimes_{\calDdag_{\Xdaginf /V}}\cal P\daginf .$$  Pour tout complexe $\calMdaginf $ de la catégorie $ \Dmoins ((\calDdag_{\Xdaginf /V}, \Sp)\Mod)$ il existe un isomorphisme canonique de la catégorie $ \Dmoins (V_Y)$:
$$\omega_{Y\dag_{\inf/V}}\LOtimes_{\calDdag_{\Ydaginf /V}, \Sp}i_{\diff}^*\calMdaginf \simeq \bfR \Gamma_{Y}(\omega_{
\Xdaginf /V})\LOtimes_{\calDdag_{\Xdaginf /V}, \Sp}\calMdaginf [\codim_XY].$$
\end{lemm}

\demo La première assertion étant de nature locale, on peut supposer que $X$ est affine et que $\calXdag $, $\calYdag $ sont des relèvements plats de $X$
et $Y$. Le complexe
$\omega_{Y\dag_{\inf/V}}\Lotimes_{\calDdag_{\Ydaginf /V},\Sp}i_{\diff}^*\cal P\daginf $ 
se représente par le complexe 
$$\omega_{\cal
Y\dag/V}\Lotimes_{\calDdag_{\calYdag /V}}\mBig(\calDdag_{\calYdag \to \calXdag /V}\Lotimes_{i^{-1}\calDdag_{ \calXdag /V}}i^{-1}\cal P\dag_{\calXdag}\mBig).$$ 
En vertu du théorème \ref{coh-loc'}, il existe un isomorphisme canonique:
$$i_*(\omega_{\cal
Y\dag/V}\Lotimes_{\calDdag_{\calYdag /V}}\calDdag_{\calYdag \to \calXdag /V})\simeq \bfR \Gamma_{Y}(\omega_{\calXdag/V})[\codim_XY].$$ 
Mais en vertu du théorème de pureté \ref{pur-coh}, on a l'isomorphisme: 
$$\bfR \Gamma_{Y}(\omega_{\calXdag/V})[\codim_XY]\simeq \cal H^{\codim_XY}_Y(\omega_{\calXdag/V}).$$ Comme $\cal P\dag_{\calXdag }$ est plat sur $\calDdag_{\calXdag /V}$  par hypothèse,
les faisceaux de torsion $\cTor_{j}^{\calDdag_{\calYdag /V}, \Sp}(\omega_{\calYdag /V}, i_{\diff}^{*,0}\cal P\dag_{\calXdag })$ sont nuls pour $j\neq 0$ et, pour $j=0$, on a l'isomorphisme:
$$\omega_{\calYdag /V}\otimes_{\calDdag_{\calYdag /V}}i_{\diff}^{*,0}\cal P\dag_{\calXdag }\simeq \cal H^{\codim_XY}_{Y}(\omega_{\calXdag/V})\otimes_{\calDdag_{\calXdag /V}}\cal P\dag_{\calXdag }\,,$$ qui se globalise naturellement dans le cas non affine.

\bigskip
Prenant une résolution  de $\calMdaginf $ par des $\calDdag_{\Xdaginf /V}$-modules à gauche spéciaux plats, on trouve l'isomorphisme canonique de complexes de Zariski:
$$\omega_{Y\daginf /V}\LOtimes_{\calDdag_{\Ydaginf /V}, \Sp}i_{\diff}^*\calMdaginf \simeq \bfR \Gamma_{Y}(\omega_{
\Xdaginf /V})\LOtimes_{\calDdag_{\Xdaginf /V}, \Sp}\calMdaginf [\codim_XY].
\eqno\endsubdemo$$

\begin{coro}\label{com-dr}
 Si $\calNdaginf $ est  $\calDdag_{\Ydaginf /K}$-module à gauche spécial, on a un isomorphisme:
$$\omega_{Y\daginf /V}\Lotimes_{\calDdag_{\Ydaginf /V}, \Sp}\calNdaginf \simeq \bfR \Gamma_{Y}(\omega_{
\Xdaginf /V})\Lotimes_{\calDdag_{\Xdaginf /V}, \Sp}i^{\diff}_*\calNdaginf .$$\end{coro}
\demo C'est une conséquence de la proposition précédente et de l'isomorphisme du corollaire \ref{adj-imm}.


\begin{lemm}Soit $X$ un  schéma lisse sur $k$ et soit $\calMdaginf $ un $\calDdag_{\Xdaginf /V}$-module à gauche spécial. Alors, le faisceau
$\cHom_{\calDdag_{\Xdaginf /V}}(\calMdaginf , \calDdag_{\Xdaginf /V})$ est
un $\calDdag_{\Xdaginf /V}$-module à droite spécial. Autrement dit, le foncteur de dualité:
$$\calMdaginf \fonct \cHom_{\calDdag_{\Xdaginf /V}}(\calMdaginf , \calDdag_{\Xdaginf /V})$$ est  contravariant et exact à
gauche de la catégorie $ (\calDdag_{\Xdaginf /V}, \Sp)\Mod$ des modules à gauche spéciaux dans la catégorie $ \Modd(\calDdag_{\Xdaginf /V}, \Sp)$ des modules à droite spéciaux.
\end{lemm}

\demo Par construction, l'action géométrique $(\sharp)$ de $\cal G_{\Xdaginf }$ sur $\cHom_{\calDdag_{\Xdaginf /V}}(\cal
M\daginf ,
\calDdag_{\Xdaginf /V})$ se fait à travers le morphisme $\Inv: \cal G_{\Xdaginf }\to \calDdag_{\Xdaginf /V}$ qui associe
à $g$ l'opérateur différentiel $g^{-1}$.
\endsubdemo

\begin{prop}Si $X$ est {\bf séparé} sur $k$, le foncteur de dualité:  
$$\calMdaginf \to 
\cHom_{\calDdag_{\Xdaginf /V}}(\calMdaginf , \calDdag_{\Xdaginf /V})$$ se dérive en un foncteur 
de la catégorie dérivée $\Rm
\rm\Dmoins ((\calDdag_{\Xdaginf /V}, \Sp)\Mod)$ vers la catégorie dérivée 
$ \Dplus (\Modd(\calDdag_{\Xdaginf /V}, \Sp))$, foncteur que l'on note:
$$\calMdaginf \fonct \bfR \cHom^{\Sp}_{\calDdag_{\Xdaginf /V}}(\calMdaginf , \calDdag_{\Xdaginf /V}).$$
\end{prop}

\demo Il faut montrer qu'un $\calDdag_{\Xdaginf /V}$-module à gauche spécial admet une résolution à gauche
par des modules à gauche spéciaux qui sont acycliques pour le foncteur dualité. 

Soient $\calUdag $ un ouvert du site  $\Xdaginf $ et $j: W\to U$ 
un ouvert affine
de $U$. Le $\calDdag_{\calUdag /V}$-module à droite $\cHom_{\calDdag_{\calUdag /V}}(j_!\calDdag_{\calUdag |W/V}, \calDdag_{\calUdag /V})$
est isomorphe à $j_*\calDdag_{\calUdag |W/V}$. Pour que le $\calDdag_{\calUdag /V}$-module à gauche $j_!\calDdag_{\calUdag |W/V}$ soit acyclique
pour le foncteur de dualité, il suffit que les images directes $R^lj_*\calDdag_{\calUdag |W/V}$ soient nulles pour $l\geq1$. Pour cela et en vertu du
corollaire \ref{acy-dif}, il suffit que $W$ soit affine  assez petit et que la trace sur $W$ d'un ouvert affine {\bf reste affine} : et c'est cela qui demande
l'hypothèse de
\em{séparation}. Si donc $U$ est {\bf séparé}, tout $\calDdag_{\calUdag /V}$-module à gauche est quotient d'un module acyclique pour le
foncteur de dualité, et la méthode de la démonstration du théorème \ref{exi-pla} montre que tout $\calDdag_{\Xdaginf /V}$-module à gauche
spécial admet une résolution à gauche par des modules à gauche spéciaux qui sont acycliques pour le foncteur dualité.
\endsubdemo

\begin{notation} \label{def-dua}On note $ \DDD^\vee  _{\calDdag_{\Xdaginf /V}}$ le foncteur de dualité décalé:\glossary{$ \DDD^\vee  _{\calDdag_{\Xdaginf /V}}(\calMdaginf )$}
$$ \DDD^\vee  _{\calDdag_{\Xdaginf /V}}(\calMdaginf ):= \bfR \cHom^{\Sp}_{\calDdag_{\Xdaginf /V}}(\calMdaginf , \calDdag_{\Xdaginf /V})[\dim X].$$
\end{notation}
Le foncteur de dualité est aussi défini sur $K$.
\begin{prop}\label{mor-tra}Si $X$ est séparé, il existe un morphisme canonique de la catégorie  dérivée des modules à droite $ \Dplus (\Modd(\calDdag_{\Xdaginf /K}, \Sp))$
$$\omega_{
\Xdaginf /K}\to \bfR \Gamma_{Y}(\omega_{
\Xdaginf /K})[2\codim_XY].\leqno E^\vee_i:$$
\end{prop}

\demo
Soit: 
$$  \bfR \Gamma_{Y}(\cal O_{
\Xdaginf /K})\to \cal O_{
\Xdaginf /K},\leqno E_i:$$ 
le morphisme canonique de la catégorie $ \Db ((\calDdag_{
\Xdaginf /K},\Sp)\Mod)$. En dualisant sur $K$, on trouve un morphisme:
$$\bfR \cHom^{\Sp}\SubV{\calDdag_{\Xdaginf /K}}(\cal O_{
\Xdaginf /K},\calDdag_{\Xdaginf /K})\to\bfR \cHom^{\Sp}\SubV{\calDdag_{\Xdaginf /K}}(\bfR \Gamma_{Y}(\cal O_{
\Xdaginf /K}),\calDdag_{\Xdaginf /K}) .$$ Le complexe de gauche est canoniquement isomorphe à $\omega_{
\Xdaginf /K}[-\dim X]$. Le théorème  suivant permet de montrer  que le complexe de droite est canoniquement isomorphe au complexe:  $$\bfR \Gamma_{Y}(\omega_{
\Xdaginf /K})[2\codim_XY-\dim X],$$ ce qui fournit le morphisme $E^\vee_i$ de la proposition.
\endsubdemo

\begin{theo}\label{dua-imm}Soit $X$ un schéma  lisse et séparé sur $k$ et soient $i:Y\hookrightarrow X$ un sous-schéma  fermé et lisse et $\calMdaginf $ un complexe
de la catégorie $ \Dmoins ((\calDdag_{\Ydaginf /V},\Sp)\Mod)$. Alors, il existe un morphisme fonctoriel canonique de dualité de la catégorie $ \Dplus (\Modd(\calDdag_{\Xdaginf /V},\Sp))$:
$$ i^{\diff}_*\DDD^\vee  _{\calDdag_{\Ydaginf /V}}(\calMdaginf )\to \DDD^\vee  _{\calDdag_{\Xdaginf /V}}(i^{\diff}_*\calMdaginf ).
\leqno  \DDD^i_*(\calMdaginf ):$$
De plus, 
le morphisme de dualité est un isomorphisme  pour tout complexe parfait $\calMdaginf $.
\glossary{$ \DDD^i_*(\calMdaginf )$}
\end{theo}

\demo
Soient $\calUdag $ un ouvert affine de $ \Xdaginf $, $W=U\cap Y$ et $\calWdag $ un relèvement de $W$. 

La restriction à $\calUdag $ du
complexe 
\smashbot{$i^{\diff}_*\bfR \cHom^{\Sp}_{\calDdag_{\Ydaginf /V}}(\calMdaginf , \calDdag_{\Ydaginf /V})$} est le complexe:
$$\bfR \cHom_{\calDdag_{\calWdag /V}}(\cal M\dag_{\calWdag }, \calDdag_{\calWdag /V})\otimes_{\calDdag_{\cal
W\dag/V}}\calDdag_{\calWdag \to\calUdag /V}\,.$$

La restriction à $\calUdag $ du
complexe $\bfR \cHom^{\Sp}_{\calDdag_{\Xdaginf /V}}(i^{\diff}_*\calMdaginf , \calDdag_{\Xdaginf /V})$ est par définition le complexe de
modules à droite: 
$$\bfR \cHom_{\calDdag_{\calUdag /V}}(i_*\calDdag_{\calUdag \leftarrow \calWdag /V}\otimes_{\calDdag_{\cal
W\dag/V}}\cal M_{\calWdag }, \calDdag_{\cal
U\dag/V})\,.$$
Soit $\cal P_{\calWdag }^{\bullet}\to \cal M\dag_{\calWdag }$ une résolution à gauche de $\cal M\dag_{\calWdag }$ par des $\calDdag_{\calWdag /V}$-modules à gauche qui  sont somme directe de modules élémentaires. Alors, on a un isomorphisme dans la catégorie dérivée:
$$\bfR \cHom_{\calDdag_{\calWdag /V}}(\cal M\dag_{\calWdag }, \calDdag_{\calWdag /V})\simeq \cHom_{\calDdag_{\calWdag /V}}(\cal P_{\calWdag }^{\bullet}, \calDdag_{\calWdag /V}).$$

Le morphisme du produit tensoriel  fournit un morphisme canonique:
$$
\cHom\SubX{\calDdag_{\calWdag /V}}(\cal P_{\calWdag }^{\bullet}, \calDdag_{\calWdag /V})\Otimes_{\calDdag_{\cal
W\dag/V}}\calDdag_{\calWdag \to\calUdag /V}\to
\cHom\SubX{\calDdag_{\calWdag /V}}(\cal P_{\calWdag }^{\bullet}, \calDdag_{\calWdag \to\calUdag /V}).$$
Par fonctorialité, on a un morphisme canonique:
$$\DeuxLignes
i_*\cHom\SubX{\calDdag_{\calWdag /V}}(\cal P_{\calWdag }^{\bullet}, \calDdag_{\calWdag \to\calUdag /V})\to\\\to
\cHom\SubX{\calDdag_{\calUdag /V}}(i_*\calDdag_{\calUdag \leftarrow \cal
W\dag/V}\Otimes_{\calDdag_{\calWdag /V}}\cal P_{\calWdag }^{\bullet},i_*\calDdag_{\calUdag \leftarrow \cal
W\dag/V}\Otimes_{\calDdag_{\calWdag /V}} \calDdag_{\cal W\to \calUdag /V}).
\endlignes
$$

Le morphisme précédent  suivi du morphisme
$ \Adj _*^i(\calDdag_{\calUdag /V})$ fournit un morphisme:
$$\TroisLignes
i_*\cHom\SubX{\calDdag_{\calWdag /V}}(\cal P_{\calWdag }^{\bullet}, \calDdag_{\calWdag /V})\Otimes_{\calDdag_{\cal
W\dag/V}}\calDdag_{\calWdag \to\calUdag /V}\to
\\
\to \cHom\SubX{\calDdag_{\calUdag /V}}\mBig(i_*\calDdag_{\calUdag \leftarrow \cal
W\dag/V}\Otimes_{\calDdag_{\calWdag /V}}\cal P_{\calWdag }^{\bullet},\cal H_Y^{[\codim_XY]}(\calDdag_{\calUdag /V})\mBig)\to
\\
\to \bfR \cHom\SubX{\calDdag_{\calUdag /V}}\mBig(i_*\calDdag_{\calUdag \leftarrow \cal
W\dag/V}\Otimes_{\calDdag_{\calWdag /V}}\cal P_{\calWdag }^{\bullet},\cal H_Y^{[\codim_XY]}(\calDdag_{\calUdag /V})\mBig)
\endlignes
$$

 Le morphisme canonique:  $$\cal H_Y^{[\codim_XY]}(\calDdag_{\calUdag /V})\to \calDdag_{\cal
U\dag/V}[\codim_XY]$$  fournit un morphisme canonique:
$$\DeuxLignes
i_*\cHom\sub{4pt}{\calDdag_{\calWdag /V}}(\cal P_{\calWdag }^{\bullet}, \calDdag_{\calWdag /V})\Otimes_{\calDdag_{\cal
W\dag/V}}\calDdag_{\calWdag \to\calUdag /V}\to
\\\to
\bfR \cHom\sub{4pt}{\calDdag_{\calUdag /V}}\mBig(i_*\calDdag_{\calUdag \leftarrow \calUdag /V}\Otimes_{\calDdag_{\calWdag /V}}\cal P_{\calWdag }^{\bullet},\calDdag_{\cal
U\dag/V}\mBig)[\codim_XY]
\endlignes$$
D'où un morphisme canonique:
$$\DeuxLignes
i_*\bfR \cHom\sub{4pt}{\calDdag_{\calWdag /V}}(\cal M_{\calWdag }, \calDdag_{\calWdag /V})\Otimes_{\calDdag_{\cal
W\dag/V}}\calDdag_{\calWdag \to\calUdag /V}
\to\\\to \bfR \cHom\sub{4pt}{\calDdag_{\calUdag /V}}(i_*\calDdag_{\calUdag \leftarrow \calUdag /V}\Otimes_{\calDdag_{\calWdag /V}}\cal M\dag_{\calWdag },\calDdag_{\cal
U\dag/V})[\codim_XY],
\endlignes$$
 qui est le morphisme de dualité $ \DDD^i_*(\cal M\dag_{\calWdag })$.

En considérant la résolution canonique  $\cal P^\bullet_{\inf}$ de $\cal
M\daginf $ par des $\calDdag_{\Ydaginf /V}$-modules à gauche spéciaux 
acycliques pour le foncteur de dualité, on voit que le morphisme local de
dualité  est la restriction d'un morphisme global de dualité
$ \DDD^i_*(\calMdaginf )$ parce que dans le cas d'un module élémentaire, de la forme $P_{\calUdag  !}(\calDdag_{\calUdag /V})$ pour un ouvert affine $\calUdag $, les deux membres du morphisme de dualité sont concentrés en 
degré $0$.

De plus, le morphisme $ \DDD^i_*(\calDdag_{\calWdag /V})$:
$$i_*\calDdag_{\calWdag \to\calUdag /V}\to \bfR \cHom_{\calDdag_{\calUdag /V}}(i_*\calDdag_{\calUdag \leftarrow \calWdag /V},\calDdag_{\cal
U\dag/V})[\codim_XY]$$ est un isomorphisme en vertu du lemme \ref{ima-invd}. Cela montre la seconde assertion du lemme de dualité. 
\endsubdemo

L'isomorphisme du dualité précédent a aussi  lieu sur l'extension $V\to K$.

\begin{coro}Soient un schéma  $X$  lisse et séparé sur $k$ et $Y$ un sous-schéma lisse sur $k$. On a alors un isomorphisme canonique de la catégorie $\Rm
\Dmoins (\Modd(\calDdag_{\Xdaginf /K},\Sp))$:
$$\DeuxLignes
\bfR \Gamma_{Y}(\omega_{
\Xdaginf /K})[2\codim_XY-\dim X]\simeq
\\\simeq
\bfR \cHom^{\Sp}\SubV{\calDdag_{\Xdaginf /K}}(\bfR \Gamma_{Y}(\cal O_{
\Xdaginf /K}),\calDdag_{\Xdaginf /K}).
\endlignes
$$
\end{coro}

\demo
En effet, on applique l'isomorphisme de dualité au fibré trivial  $\cal O_{\Xdaginf /K}$, en tenant compte des isomorphismes des théorèmes \ref{coh-loc} et
\ref{coh-loc'}.
\endsubdemo


L'isomorphisme
du corollaire induit le morphisme de la proposition \ref{mor-tra}.

En vertu de la proposition  \ref{mor-tra}, on obtient par produit tensoriel un morphisme de la catégorie $ \rmD(K_X)$:
$$\omega_{X\dag_{\inf/K}}\ \LOtimes_{\calDdag_{\Xdaginf /K},\Sp}\ \calMdaginf \to \bfR \Gamma_{Y}(\omega_{
\Xdaginf /K})\ \LOtimes_{\calDdag_{\Xdaginf /K}, \Sp}\ \calMdaginf [2\codim_XY].$$ Mais le complexe
de gauche est canoniquement isomorphe au complexe de de Rham décalé vers la droite: 
$$\bfR \cHom\SubX{\calDdag_{\Xdaginf /K,\Sp}}(\cal O_{
\Xdaginf /K},\calMdaginf )[\dim X],$$ alors qu'en vertu du lemme \ref{dr-coh} le complexe de droite est canoniquement isomorphe 
au complexe de de Rham décalé vers la gauche: $$\bfR \cHom_{\calDdag_{\Ydaginf /K,\Sp}}(\cal O_{\Ydaginf /K},i_{\diff}^*\calMdaginf )[\dim X].$$ 
 Cela fournit le morphisme du théorème \ref{dr-dir} dans le cas d'une immersion fermée.
\enddemo

\begin{prop}\label{comp} Soient $i: Y\to X$ une immersion fermée de schémas affines lisses sur $k$ et $i\dag: \calYdag \to \calXdag $ un relèvement de
$i$. Alors, le morphisme: $$\omega_{
\Xdaginf /K}\to \bfR \Gamma_{Y}(\omega_{
\Xdaginf /K})[2\codim_XY]\leqno E^\vee_i:$$ est représenté par le morphisme de complexes défini par $i\dag$:
$$i^{-1}\mBig(\Omega^{\bullet}_{\calXdag /K}\Otimes_{\cal O_{\calXdag /K}}\calDdag_{\calXdag /K}\mBig)[\dim X]\to 
\mBig(\Omega^{\bullet}_{\calYdag /K}\Otimes_{\cal
O_{\calYdag /K}}\calDdag_{\calYdag \to\calXdag /K}\mBig)[\dim X].$$
\end{prop}

\demo
Voyons d'abord que $u:= i^*$ définit un morphisme de complexes:
$$\calDdag_{\calXdag \gets\calYdag /K}\Otimes_{\calDdag_{\calYdag /K}}\Sp^{\bullet}(\cal O_{\calYdag /K})\to \cal H^{\codim_XY}_{Y}(\calDdag_{\calXdag/K})\Otimes_{\calDdag_{\calXdag /K}}\Sp^{\bullet}(\cal O_{\calXdag /K})\leqno u:$$ 
par $$P\otimes1\otimes \eta_1\wedge\cdots\wedge\eta_k\mapsto { \Adj }^i_*(P\otimes u)\otimes\trans\, u(\eta_1\wedge\cdots\wedge\eta_k),$$ où $\trans\, u$ est la
différentielle de
$u$ et ${ \Adj }^i_*(P\otimes u):= { \Adj }^i_*(\calDdag_{\calXdag /K})(P\otimes u)$. Ce morphisme représente le morphisme
$$i^{\diff}_*i^*_{\diff}\cal O_{\calXdag /K}\to \cal H^{\codim_XY}_{Y}(\calDdag_{\calXdag/K})\LOtimes_{\calDdag_{\calXdag /K}}\Sp^{\bullet}(\cal O_{\calXdag /K})$$

\vskip0pt plus 1ex
\noindent Considérons alors le diagramme:

\penalty10000
\vskip0pt plus 1ex
\nobreak
\vbox{$$\def\quad{\hskip0.5ex}
\def\cHom_#1{\mathop{\cal H\!\it om}\nolimits\sub{4pt}{#1}}
\matrix{
i^{-1}(\Omega^{\bullet}_{\calXdag /K}\otimes_{\cal O_{\calXdag /K}}\calDdag_{\calXdag /K})
&\simeq&
i^{-1}\cHom_{\calDdag_{\calXdag /K}}\big(\Sp^{\bullet}(\cal O_{\calXdag /K}),\calDdag_{\calXdag /K}\big)
\cr
\downarrow\indiced{ u:=i^*}&&\downarrow\cr
\Omega^{\bullet}_{\calYdag /K}\Otimes_{\cal O_{\cal
Y\dag/K}}\calDdag_{\calYdag \to\calXdag /K}
&\simeq&
\cHom_{\calDdag_{\calYdag /K}}\big(\Sp^{\bullet}(\cal O_{\calYdag /K}),\calDdag_{\calYdag \to \calXdag /K}\big)
&\hskip1.2cm\cr}$$
\vskip20pt
$$
\matrix{
\hskip8mm&
\hmsmash{\llap{$\to$\hvpath{-5 7 103 23 -6}\ }\cHom\SubO{\calDdag_{\calXdag /K}}\big(\cal H^{\codim_XY}_Y(\calDdag_{\calXdag /K})\Otimes_{\calDdag_{\calXdag /K}}\Sp^{\bullet}(\cal O_{\calXdag /K})\,,\, \cal
H^{\codim_XY}_{\,Y}(\calDdag_{\calXdag /K})\big)}\cr
\noalign{\kern-12pt}&
\downarrow\indiced u\cr \noalign{\kern-2pt}
&
\hmsmash{\llap{$\to$\hvpath{-10 20 104 8 -4}\ }\cHom_{\calDdag_{\calXdag /K}}\big(\calDdag_{\calXdag \gets\calYdag /K}\Otimes_{\calDdag_{\calYdag /K}}\Sp^{\bullet}(\cal O_{\calYdag /K})\,,\,\cal
H^{\codim_XY}_{\,Y}(\calDdag_{\calXdag /K})\big),}
}\postskip0ex$$ 
}
\vskip0pt plus 1ex
\noindent qui est par construction \em{commutatif}.  
D'autre part, le diagramme suivant est aussi commutatif:
\vskip0pt plus 1ex

\penalty10000
\vbox{\def\quad{\hskip0.5ex}$$\scriptwd1.5em
\def\calH{\cal H}
\def\cHom_#1{\mathop{\cal H\!\it om}\nolimits\Sub{0pt}{#1}}
\vcenter{\vskip5pt$$\kern-1.3cm\matrix{
\cHom_{\calDdag_{\calXdag /K}}\big(\calH^{\codim_XY}_{\,Y}(\calDdag_{\calXdag /K})\Otimes_{\calDdag_{\calXdag /K}}\Sp^{\bullet}(\cal O_{\calXdag /K})\,,\, \calH^{\codim_XY}_{\,Y}(\calDdag_{\calXdag /K})\big)
\cr
\noalign{\kern-10pt}\downarrow \cr
\cHom_{\calDdag_{\calXdag /K}}\big(\calDdag_{\calXdag \gets\cal
Y\dag/K}\Otimes_{\calDdag_{\calYdag /K}}\Sp^{\bullet}(\cal O_{\calYdag /K})\,,\, 
\calH^{\codim_XY}_{\,Y}(\calDdag_{\calXdag /K})\big)
\cr}$$
\vskip20pt
$$\hskip12pt\matrix{
\llap{$\to$\hvpath{-4 8 105 24 -4}\ }
\bfR \cHom_{\calDdag_{\calXdag /K}}\big(\calH^{\codim_XY}_{\,Y}(\calDdag_{\calXdag /K})
\Otimes_{\calDdag_{\calXdag /K}}\Sp^{\bullet}(\cal O_{\calXdag /K})\,,\, 
\calH^{\codim_XY}_{\,Y}(\calDdag_{\calXdag /K})\big)
\cr
\noalign{\kern-10pt}\hskip10pt \downarrow \cr
\llap{$\to$\hvpath{-12 23 104 9 -9}\ }
\bfR \cHom_{\calDdag_{\calXdag /K}}\big(\calDdag_{\calXdag \gets\calYdag /K}\Otimes_{\calDdag_{\calYdag /K}}\Sp^{\bullet}(\cal O_{\calYdag /K})\,,\, 
\calH^{\codim_XY}_{\,Y}(\calDdag_{\calXdag /K})\big)
}$$
\vskip20pt$$\kern30pt\matrix{
\llap{$\to$\hvpath{-4 7 108 25 -2}\ }
\bfR \cHom_{\calDdag_{\calXdag /K}}\big(\calH^{\codim_XY}_Y(\calDdag_{\calXdag /K})\Otimes_{\calDdag_{\calXdag /K}}\Sp^{\bullet}(\cal O_{\calXdag /K})\,,\, 
\calDdag_{\calXdag /K}\big)[\codim_XY]\cr
\noalign{\kern-10pt}\hskip10pt \downarrow \cr
\llap{$\to$\hvpath{-12 22 108 11 -7}\ }
\bfR \cHom_{\calDdag_{\calXdag /K}}\big(\calDdag_{\calXdag \gets\calYdag /K}\Otimes_{\calDdag_{\calYdag /K}}\Sp^{\bullet}(\cal
O_{\calYdag /K})\,,\, 
\calDdag_{\calXdag /K}\big)[\codim_XY]\cr
}$$}\belowdisplayskip0pt$$} 
\vskip0pt plus 2pt
\noindent ce qui montre que le morphisme: 
$$i^{-1}\omega_{
\Xdaginf /K}\to \bfR \Gamma_{Y}(\omega_{
\Xdaginf /K})[2\codim_XY]\postskip1ex\preskip1ex$$ 
est représenté par le morphisme de complexes défini par $i\dag$:
$$\preskip1ex\ecartlignes0pt
\TroisLignes
\vrule depth15pt width0pt
i^{-1}\omega_{\calXdag /K}\simeq i^{-1}\big(\Omega^{\bullet}_{\calXdag /K}\Otimes_{\cal O_{\calXdag /K}}\calDdag_{\calXdag /K}\big)[\dim X]
\stackrel{u}{\too}\\
\hskip5em\stackrel{u}{\too}
\big(\Omega^{\bullet}_{\calYdag /K}\Otimes_{\cal O_{\calYdag /K}}\calDdag_{\calYdag \to\calXdag /K}[\dim X]\big)
\simeq \\%
\simeq \cal H_Y^{\codim_XY}(\omega_{\calXdag /K})[\codim_XY].
\endlignes
\eqno\enddemo$$

\begin{coro}\label{fro-deR}Soient $f: Y\to X$ un morphisme de schémas  lisses et {\bf séparés} sur $k$, $f\dag:  \calYdag \to \calXdag $ un relèvement de $f$ et $\cal
M\daginf $ un $\calDdag_{\Xdaginf /K}$-module à gauche spécial  et acyclique pour $f^*_{\diff}$. Alors, le morphisme 
$$
f^{-1}\bfR \cHom\goodSub{8pt}{4pt}{-6mm}{\calDdag_{\Xdaginf /K},\Sp}\big(
\cal O_{\Xdaginf /K},\calMdaginf \big)\to\bfR \cHom\goodSub{8pt}{4pt}{-6mm}{\calDdag_{\Ydaginf /K},\Sp}(
\cal O_{\Ydaginf /K},f^*_{\diff}\calMdaginf )
$$
est représenté par le morphisme de complexes:
$$f^{-1}\big(\cal M\dag_{\calXdag }\otimes_{\cal O_{\calXdag /K}}\Omega^{\bullet}_{\calXdag /K}\big)\to f^{*,0}_{\diff}(\cal M\dag_{\calXdag})\otimes_{\cal O_{\calYdag /K}}\Omega^{\bullet}_{\calYdag /K}.$$
\end{coro}

\demo La méthode du complexe de Spencer du théorème \ref{ind-rel}  montre le  cas d'une projection. Le cas    d'une immersion ouverte est élémentaire. Le cas d'une immersion fermée est conséquence, précisément, de la compatibilité \ref{comp}. Cela montre 
que le morphisme \ref{dr-dir} recolle les morphismes locaux que l'on peut construire à partir de relèvements. En particulier,
le corollaire s'applique à chaque fois que le foncteur $f^{*,0}_{\diff}$ est exact.
\enddemo

\begin{Rema}Les morphismes locaux précédents sont des morphismes de la catégorie dérivée, qui n'est pas un champ, et donc ne se recollent pas a priori, ce qui explique la nécessité et l'utilité du foncteur 
de dualité
pour définir le morphisme \ref{dr-dir} dans le cas d'une immersion fermée. Dans la pratique, le corollaire permet d'expliciter localement le morphisme précédent et c'est ce qui est utile.
\end{Rema}
La même démonstration montre la compatibilité sur $K$ du foncteur de dualité avec le foncteur image inverse pour une projection.
\begin{theo}\label{dua-pro} Soient $p: Y\times_kX\to X$ une projection de schémas lisses et séparés sur $k$, 
 et $\calMdaginf $ un complexe
de la catégorie $ \Dmoins ((\calDdag_{\Xdaginf /K},\Sp)\Mod)$. Il existe alors un morphisme fonctoriel canonique de dualité de la catégorie $ \Dplus (\Modd(\calDdag_{(Y\times_kX)\daginf /K},\Sp))$:
$$ p_{\diff}^*\DDD^\vee  _{\calDdag_{\Xdaginf /K}}(\calMdaginf )\to \DDD^\vee  _{\calDdag_{(Y\times_kX)\daginf /K}}(p_{\diff}^*\calMdaginf ).
\leqno  \DDD_p^{*}(\calMdaginf ):$$
De plus, 
le morphisme de dualité est un isomorphisme  pour tout complexe parfait $\calMdaginf $.
\end{theo}
\subsection{La fonctorialité  de la cohomologie de de Rham $p$-adique}
Nous allons déduire la fonctorialité de la cohomologie des résultats précédents.

\begin{coro}\label{com-inv}Soit $f:Y\to X$ un morphisme de schémas lisses et séparés  sur $k$ et soit $\calMdaginf $ un complexe de $ \Db ((\calDdag_{\Xdaginf /K}, \Sp)\Mod)$. Alors, il existe un morphisme canonique  $\DR(f/K,\calMdaginf ):$
$$\bfR \hom\SubX{\calDdag_{\Xdaginf /K},\Sp}\big(
\cal O_{\Xdaginf /K},\calMdaginf \big)\to\bfR \hom\SubX{\calDdag_{\Ydaginf /K},\Sp}\big(
\cal O_{\Ydaginf /K},f^*_{\diff}\calMdaginf \big),$$ qui est en fonctoriel en $\calMdaginf $.\glossary{$\DR(f/K,\calMdaginf )$}
\end{coro}
C'est une conséquence du théorème \ref{dr-dir} et du fait que la cohomologie globale est l'hypercohomologie pour la topologie de Zariski du complexe de de Rham
$p$-adique local. 
\bigskip

Spécialisant au cas géométrique $\calMdaginf = \cal O_{\Xdaginf /K}$ et tenant compte du théorème \ref{inv-tri}, on trouve la fonctorialité de la cohomologie de de
Rham
$p$-adique pour un morphisme.

\begin{theo}\label{fon-coh}Soit $f:Y\to X$  un morphisme de schémas lisses et séparés  sur $k$. Il existe alors    un  morphisme canonique $\DR(f/K):$
$$\bfR \hom\SubX{\calDdag_{\Xdaginf /K},\Sp}(
\cal O_{\Xdaginf /K},\cal O_{\Xdaginf /K})
\to
\bfR \hom\SubX{\calDdag_{\Ydaginf /K},\Sp}(
\cal O_{\Ydaginf /K},\cal O_{\Ydaginf /K}),$$ qui est transitif  pour
deux morphismes.\glossary{$\DR(f/K)$}
\end{theo}
\demo
Si $Y\stackrel{f}{\to}X\stackrel{g}{\to} Z$ sont deux morphismes 
entre schémas séparés et lisses sur $k$, le théorème \ref{inv-tri} fournit des isomorphismes canoniques:
$$\mathrigid2mu f^{*}_{\diff}(g^{*}_{\diff}(\cal O_{Z\daginf /K}))\simeq f^{*}_{\diff}(\cal O_{\Xdaginf /K})\simeq\cal O_{\Ydaginf /K}\simeq(g\circ f)^{*}_{\diff}(\cal O_{Z\daginf /K}).$$

 Reste à voir que le foncteur $\DR(-/K)$ est compatible avec la composition des morphismes. Considérons le diagramme commutatif où les morphismes $i_f, i_g$ sont les morphismes graphes de $f$ et $g$ et les morphismes $p$ sont les projections naturelles :
$$\catcode`\!=12
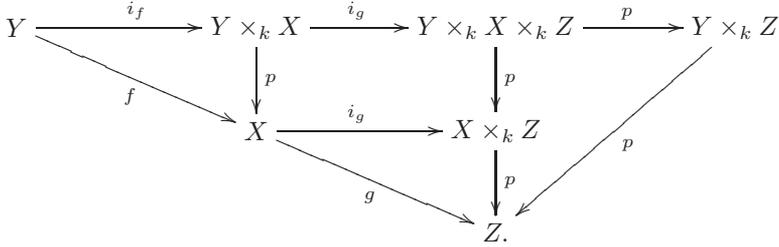
\xymatrix@!C{Y\ar[r]^{i_f}\ar[dr]_f &
Y\times_k X\ar[d]^p\ar[r]^{i_g\hskip0.5cm}&
Y\times_k X\times_k Z\ar[d]^p \ar[r]^{\hskip0.3cm p}
&Y\times_k Z\ar[ddl]^p\\
&X\ar[r]^{i_g\hskip0.5cm}\ar[dr]_g&X\times_k Z\ar[d]^p\\
&&Z.}$$
 Par construction, le foncteur  $\DR(-/K)$ est compatible avec la composition de deux immersions fermées ou de deux projections. Reste à voir la compatibilité des morphismes du carré cartésien, c'est-à-dire l'égalité: $$DR(i_g/K)\circ \DR(p/K)= DR(p/K)\circ \DR(i_g/K).\leqno (*):$$ On a  un  premier   diagramme commutatif:
$$\def\quad{\hskip0.2ex}\mathrigid1mu
\scriptspace0pt\def\Sub{\goodSub{0pt}{7pt}{-9mm}}
\matrix{\DR(X\times_k Z/K)&\to&\bfR \hom\Sub{\calDdag_{(X\times_k Z)\daginf /K},\Sp}\big(
\bfR \Gamma_{i_g(X)}(\cal O_{(X\times_k Z)\daginf /K})\,,\cal O_{(X\times_k Z)\daginf /K}\big)&
\cr\noalign{\kern-12pt}
\downarrow&&\downarrow
\cr
\DR(Y\timesk X\timesk Z/K)&\to&\bfR \hom\Sub{\calDdag_{(Y\timesk X\timesk Z)\daginf /K},\Sp}\big(\bfR \Gamma_{i_g(Y\timesk X)}(\cal O_{(Y\timesk X\timesk Z)\daginf /K})\,,\cal O_{(Y\timesk X\timesk Z)\daginf /K}\big)\,,}$$ 
puis  un diagramme commutatif dont les morphismes horizontaux sont des isomorphismes:

\vtop{\kern-10pt\def\quad{\hskip0.2ex}\mathrigid1mu
\scriptspace0pt\def\Sub{\goodSub{0pt}{7pt}{-7mm}}
$$\hbox to\hsize{$
\matrix{\bfR \hom\Sub{\calDdag_{(X\times_k Z)\daginf /K},\Sp}\big(
\bfR \Gamma_{i_g(X)}(\cal O_{(X\timesk Z)\daginf /K})\,,\,\cal O_{(X\timesk Z)\daginf /K}\big)\cr\noalign{\kern-12pt}
\downarrow\cr
\bfR \hom\Sub{\calDdag_{(Y\times_k X\times_k Z)\daginf /K},\Sp}\big(\bfR \Gamma_{i_g(Y\timesk X)}(\cal O_{(Y\timesk X\timesk Z)\daginf /K})\,,\,\cal O_{(Y\timesk X\timesk Z)\daginf /K}\big)}
$\hss}$$
\vskip3mm$$
\hbox to\hsize{\hss$\matrix{
\llap{\def\addendpath{$\!\leftarrow$}$\leftarrow$\hvpath{-4 6 67 25 -16}}
\bfR \hom\Sub{\calDdag_{(X\times_k Z)\daginf /K},\Sp}\big(\cal H^{\dim Z}_{i_g(X)}(\cal O_{(X\timesk Z)\daginf /K})\,,\,\cal H^{\dim Z}_{i_g(X)}(\cal O_{(Y\timesk X)\daginf /K})\big)
\cr\noalign{\kern-12pt}
\downarrow
\cr
\llap{\def\addendpath{$\!\leftarrow$}$\leftarrow$\hvpath{-4 21 73 10 -8}\ }
 \bfR \hom\Sub{\calDdag_{(Y\times_k X\times_k Z)\daginf /K},\Sp}\big(\cal H^{\dim Z}_{i_g(Y\timesk X)}(\cal O_{(Y\timesk X\timesk Z)\daginf /K})\,,\,\cal H^{\dim Z}_{i_g(Y\timesk X)}(\cal O_{(Y\timesk X\timesk Z)\daginf /K})\big)}
$}$$}
\noindent et enfin un diagramme: 
$$\hss\def\quad{\hskip0.ex}\mathrigid0mu
\scriptspace0pt\def\Sub{\goodSub{0pt}{7pt}{-7mm}}
\matrix{\bfR \hom\Sub{\calDdag_{(X\times_k Z)\daginf /K},\Sp}\big(
\cal H^{\dim Z}_{i_g(X)}(\cal O_{(X\times_k Z)\daginf /K})\mkern1mu,\mkern1.7mu \cal H^{\dim Z}_{i_g(X)}(\cal O_{(X\times_k Z)\daginf /K})\big)&\hdecale{-3mm}\gets \DR(X/K)
\cr\noalign{\kern-7pt}\downarrow&\downarrow
\cr\noalign{\kern5pt}
\def\timesk{\goodtimesk{3mu}{5mu}}\bfR \hom\Sub{\calDdag_{(Y\timesk X\timesk Z)\daginf /K},\Sp}\big(\mkern-1mu\cal H^{\dim Z}_{i_g(Y\timesk X)}(\cal O_{(Y\timesk X\timesk Z)\daginf /K})\mkern1mu,\mkern1.7mu\cal H^{\dim Z}_{i_g(Y\timesk X)}(\cal O_{(Y\timesk X\timesk Z)\daginf /K})\mkern-2mu\big)&\gets \DR(Y\timesk \mkern1 mu X/K),}\hss$$
dont les morphismes horizontaux sont des isomorphismes et dont la commutativité entraîne la compatibilité $(*)$.

En effet, soit: $$\cal O_{\Xdaginf/K}\to \cal I^{\bullet}_{X\daginf}$$ une résolution par des $\calDdag_{X\daginf /K}$-modules à gauche spéciaux et 
$$ p^{*,0}_{\inf} \cal I^{\bullet}_{X\daginf}\to \cal J^{\bullet}_{(Y\timesk X)\daginf}$$ une résolution par des $\calDdag_{(Y\times_k X)\daginf /K}$-modules à gauche spéciaux, alors le dernier diagramme  se représente par le diagramme {\bf commutatif} de complexes d'espaces vectoriels sur $K$ dont les morphismes horizontaux sont des isomorphismes en vertu du corollaire \ref{com-dr}:
$$\hss\def\quad{\hskip0.1ex}\mathrigid0mu
\scriptspace0pt\def\Sub{\goodSub{0pt}{7pt}{-7mm}}\matrix{
\hom\SubX{\calDdag_{(X\times_k Z)\daginf /K},\Sp}\big(
\cal H^{\dim Z}_{i_g(X)}(\cal O_{(X\times_k Z)\daginf /K})\mkern1mu,\mkern1.7mu i^{\diff}_{g,*}\cal I^{\bullet}_{X\daginf}\big)&\hdecale{-0.3cm}\gets \hom\SubX{\calDdag_{X\daginf/K},\Sp}(
\cal O_{X\daginf/K}\mkern1mu,\mkern1.7mu \cal I^{\bullet}_{X\daginf})
\cr\noalign{\kern-5pt}\downarrow&\downarrow
\cr
\def\timesk{\goodtimesk{2mu}{4mu}}\hom\Sub{\calDdag_{(Y\times_k X\times_k Z)\daginf /K},\Sp}\big(\mkern-2mu\cal H^{\dim Z}_{i_g(Y\timesk X)}(\cal O_{(Y\timesk X\timesk Z)\daginf/K})\mkern1mu,\mkern1.7mui^{\diff}_{g,*}\cal J^{\bullet}_{(Y\timesk X)\daginf}\big)
&\gets \hom\Sub{\calDdag_{(Y\times_k X)\daginf /K},\Sp}\big(\cal O_{(Y\timesk X)\daginf /K}\mkern1mu,\mkern1.7mu\cal J^{\bullet}_{(Y\timesk X)\daginf}\big),}
\hss$$ 
où le morphisme
vertical de gauche provient de l'isomorphisme de changement de base pour une projection \ref{cha-pro}. 

La commutativité de ce diagramme implique l'égalité: $$DR(i_g/K)\circ \DR(p/K)= DR(p/K)\circ \DR(i_g/K)$$ et donc l'égalité:
$$\DR(f/K)\circ\DR(g/K)= \DR(g\circ f/K).\eqno\enddemo$$

\medskip

Le théorème suivant résume ce qui précède.

\begin{theo}\label{fon-deR}La correspondance $\DR(-/K):$
$$\matrix{
\Sms (k)&\fonct&\Dplus (K)\hfill\cr
X&\fonct& \DR(X/K):=\bfR \hom\SubX{\calDdag_{\Xdaginf /K},\Sp}(
\cal O_{\Xdaginf /K},\cal O_{\Xdaginf /K})
}$$ est un {\bf foncteur} contravariant de la catégorie $ \Sms (k)$  des schémas lisses et
séparés sur $k$  vers la catégorie dérivée $ \Dplus (K)$, qui étend le foncteur de Monsky-Washnitzer \ref{fon-aff} du cas affine.\glossary{$\DR(-/K)$}
\end{theo}

\begin{Rema}Sous la condition $e<p^h(p-1)$, tous les arguments précédents se transposent pour le couple $\cal G_{\Xdaginf }\hookrightarrow \cal
D^{\dagger, h}_{\Xdaginf /V}$.  On n'obtient rien de nouveau pour la cohomologie  sur  $K$  quand $h\neq0$. 
Mais dans le cas
$h=0$, on obtient déjà la fonctorialité de la cohomologie sur $V$:
\end{Rema}

\begin{theo}\label{fon-ent}Si $e<p-1$, alors la correspondance  $\DR_0(-/V):$ 
$$\matrix{\Sms (k)&\fonct& \Dplus (V)\hfill\cr
X&\fonct& \bfR \hom\SubX{\calD^{\dagger, 0}_{\Xdaginf /V},\Sp}(
\cal O_{\Xdaginf /V},\cal O_{\Xdaginf /V})}$$ est un {\bf foncteur} contravariant de la catégorie $ \Sms (k)$ des schémas lisses et
séparés sur $k$ dans la catégorie dérivée $ \Dplus (V)$.
\end{theo}

\demo En effet, dans ce cas on peut considérer le complexe de Spencer $\Sp^{\bullet,0}(\cal O_{\Xdaginf /V})$ \glossary{$\Sp^{\bullet,0}(\cal O_{\Xdaginf /V})$}d'échelon zéro qui se construit en remplaçant 
le faisceau $\calD^{\dagger}_{\Xdaginf /K}$ par le faisceau $\calD^{\dagger, 0}_{\Xdaginf /V},$  qui est une résolution du faisceau structural,
qui est alors un module de présentation finie. 
A partir de là tous les arguments se transposent.
\enddemo

\begin{Rema}
D'autre part, si $e<p-1$, la cohomologie de Rham $p$-adique d'échelon zéro du site infinitésimal formel doit être canoniquement
isomorphe à la cohomologie cristalline sur $V$, et cet isomorphisme doit être fonctoriel.
\end{Rema}
\subsection{Le théorème de factorisation $p$-adique de la fonction Zêta d'une variété 
algébrique lisse sur un corps fini}
Nous allons appliquer la fonctorialité précédente, en fait dans une situation plus simple,  au morphisme de Frobenius d'une variété algébrique lisse sur un corps fini, morphisme plat mais non lisse, et en déduire l'expression cohomologique $p$-adique de sa fonction Zêta.
\begin{defi}\label{fro-coh} Supposons  que $k$ est un corps fini $\Bbb F_q$ et notons $f: =fr$ le morphisme de Frobenius relatif à $\Bbb F_q$ d'une variété algébrique lisse $X$  sur $\Bbb F_q$.
On définit l'endomorphisme ${\bf \cal F}$ de complexes de Zariski par :
$$\def\sub{\goodSub{15pt}{6pt}{-6mm}}\ecartlignes5pt
\DeuxLignes
\bfR \cHom\sub{\calDdag_{\Xdaginf /K},\Sp}(\cal O_{\Xdaginf /K},\cal O_{\Xdaginf /K})
\to\\\to
\bfR \cHom\sub{\calDdag_{\Xdaginf /K},\Sp}(\cal O_{\Xdaginf /K},fr^*_{\diff}\cal O_{\Xdaginf /K})
\simeq
\bfR \cHom\sub{\calDdag_{\Xdaginf /K},\Sp}(\cal O_{\Xdaginf /K},\cal O_{\Xdaginf /K}),
\endlignes
$$
où le premier morphisme est celui du théorème  \ref{dr-dir} et le second provient
de l'isomorphisme canonique $fr^*_{\diff}\cal O_{\Xdaginf /K}\simeq \cal O_{\Xdaginf /K}$ du théorème \ref{inv-tri}.
On définit l'endomorphisme de Frobenius de la cohomologie de de Rham $p$-adique,   et l'on note \glossary{$\bfF^{\bullet}$}$$ H^{i}_{\DR}(X/K, \cal O_{\Xdaginf /K})\to  H^{i}_{\DR}(X/K, \cal O_{\Xdaginf /K})\,,\leqno {\bfF}^{i}:$$ 
le morphisme induit en cohomologie par le morphisme $\bfF$:
$$\def\sub{\goodSub{10pt}{6pt}{-5mm}}
\ecartlignes5pt\DeuxLignes
\bfR \hom\sub{\calDdag_{\Xdaginf /K},\Sp}(\cal O_{\Xdaginf /K},\cal O_{\Xdaginf /K})
\to
\bfR \hom\sub{\calDdag_{\Xdaginf /K},\Sp}(\cal O_{\Xdaginf /K},fr^*_{\diff}\cal O_{\Xdaginf /K})
\simeq\\\simeq
\bfR \hom\sub{\calDdag_{\Xdaginf /K},\Sp}(\cal O_{\Xdaginf /K},\cal O_{\Xdaginf /K}).
\endlignes
$$\end{defi}

\begin{lemm}\label{sui-lon}
Si $X_1\cup X_2$ est un recouvrement ouvert de $X,$ le triangle de Mayer-Vietoris de \ref{may-vie}   induit une suite exacte longue de cohomologie de $K$-espaces vectoriels munis de l'endomorphisme  $\bfF^\bullet$:
$$\preskip1ex\ecartlignes6pt
\TroisLignes 
 \cdots \to H^{\bullet}_{\DR}(X/K, \cal O_{\Xdaginf /K})
\to\\\to  H^{\bullet}_{\DR}(X_1/K, \cal O_{(X_1)\daginf /K})\oplus  H^{\bullet}_{\DR}(X_2/K, \cal O_{(X_2)\daginf /K})
\to\\
\to  H^{\bullet}_{\DR}(X_{12}/K, \cal O_{(X_{12})\daginf /K})\to H^{\bullet+1}_{\DR}(X/K, \cal O_{\Xdaginf /K})\to \cdots
\endlignes
$$
\end{lemm}

\demo
Soit $\bfR \cHom_{\calDdag_{\Xdaginf /K},\Sp}(\cal O_{\Xdaginf /K},\cal O_{\Xdaginf /K})\to \cal I^{\bullet}_X$ une résolution injective par un  complexe de Zariski de faisceaux d'espaces vectoriels sur $K$. Alors, l'endomorphisme $\bf \cal F$ se prolonge en un endomorphisme $\tilde {\cal F}$ du complexe $\cal I^{\bullet}_X$. On obtient alors un endomorphisme d'une  {\bf suite exacte courte} de complexes d'espaces vectoriels sur $K$, parce que les termes du complexe $\cal I^{\bullet}_X$ sont {\bf flasques}:
$$\matrix{0&\to&\Gamma(X,\cal I^{\bullet}_X)&\to&\Gamma(X_1,\cal I^{\bullet}_X)\oplus\Gamma(X_2,\cal I^{\bullet}_X)&\to\Gamma(X_{12},\cal I^{\bullet}_X)&\to 0\phantom{\,.}\cr
&&\downarrow&&\downarrow&\downarrow
\cr 0&\to&\Gamma(X,\cal I^{\bullet}_X)&\to&\Gamma(X_1,\cal I^{\bullet}_X)\oplus\Gamma(X_2,\cal I^{\bullet}_X)&\to\Gamma(X_{12},\cal I^{\bullet}_X)&\to 0\,.}$$ La suite longue de cohomologie fournit le lemme. 

\begin{coro} Les endomorphismes de Frobenius  $\bfF^{\bullet}$ sur la cohomologie de de Rham $p$-adique sont {\bf bijectifs}.\end{coro}
\demo La suite longue de cohomologie du lemme précédent ramène
la bijectivité de $\bfF^{\bullet}$ pour $X$ aux cas de $X_1$, $X_2$  et  de $X_1\cap X_2$. En raisonnant sur le nombre  d'ouverts affines d'un recouvrement de $X$,  on se ramène au cas affine,  traité par Monsky dans [M$_2$], sachant que {\bf localement} les morphismes de Frobenius sont induits par un morphisme de complexes de de Rham. 
\enddemo

\begin{defi} Si la variété  $X$ est purement de dimension $\dim X$, on définit le   nombre de Lefschetz $L(X)$  comme la somme alternée des traces de l'automorphisme gradué  de la cohomologie $q^{\dim X}({\bfF}^\bullet)^{-1}$. \end{defi}
Le nombre $L(X)$  est bien défini comme un élément de $K$ en vertu du théorème de finitude \ref{fin} [Me$_2$].
\begin{coro}\spaceskip0.85ex plus2pt 
On a l'égalité $L(X)=N_1(X)$, où  $N_1(X)$  est le nombre de points $\Bbb F_q$-rationnels de $X$.\end{coro}
\demo
En considérant un recouvrement de $X$ par des ouverts affines et en raisonnant par récurrence sur le nombre des ouverts
affines du recouvrement, la suite longue de cohomologie
du triangle  de Mayer-Vietoris
réduit la formule des traces $N_1(X)= L(X) $ au cas affine, traité dans [M$_2$],  puisque les deux nombres sont additifs pour deux ouverts. 
\enddemo
Soit $V\to V'$ une extension finie  d'anneaux de valuation discrète complets de corps résiduels parfaits. Notons  
$(k', K', X')$  le triplet obtenu à l'aide de ce changement de base
à partir du triplet $(k, K, X)$.
\begin{lemm} Soit $V\to V'$ une telle extension. Alors, le morphisme canonique:
$$\scriptwd1.5em\def\SubX{\goodSub{13pt}{4pt}{-6mm}}
\bfR \hom\SubX{\calDdag_{\Xdaginf /K},\Sp}(
\cal O_{\Xdaginf /K},\cal O_{\Xdaginf /K})\Otimes_KK'\to\bfR \hom\SubX{\calDdag_{{X'}\daginf /K'},\Sp}(
\cal O_{{X'}\daginf /K'},\cal O_{{X'}\daginf /K'})$$ est un isomorphisme. \end{lemm}
\demo Dans le cas affine, la cohomologie de de Rham $p$-adique commute au changement de base fini. La suite de Mayer-Vietoris   ramène le cas général au cas affine.
\enddemo

On en déduit la première démonstration, à notre connaissance, de la factorisation $p$-adique de la fonction Zêta d'une variété algébrique lisse sur un corps fini  pour une variété algébrique lisse sur un corps fini, sans aucune autre restriction, ce qui constitue le véritable test de nos méthodes:

\begin{theo}\label{fac-zet}Supposons  que $k$ est un corps fini $\Bbb F_q$ et 
que $f=fr$ est le morphisme de Frobenius  relatif à $\Bbb F_q$ d'une variété algébrique 
$X$ lisse purement de dimension $\dim X$. 
On a alors la factorisation de la fonction Zêta de $X$:
$$Z(X,t)= {\displaystyle\prod_{i\,\rm impair} P_{p,i}(X)\mBig/\displaystyle
\prod_{i\,\rm pair} P_{p,i}(X)}$$ où $P_{p,i}(X):= \det(1-q^{\dim X}({\bfF}^i)^{-1}t) \in K[t]$ avec $0\leq i\leq 2\dim X,$ est le polynôme caractéristique 
de $q^{\dim X}({\bfF}^i)^{-1}$ agissant sur  la cohomologie $ H^{i}_{\DR}(X/K):= H^{i}_{\DR}(X/K, \cal O_{\Xdaginf /K})$.
\end{theo}
\demo
Pour $m\geq 1$, soit $V\to V_m$ l'extension non ramifiée relevant l'extension $\Bbb F_q\to \Bbb F_{q^m}$. Par le lemme précédent,
on obtient un isomorphisme:
$$ H^{i}_{\DR}(X/K, \cal O_{\Xdaginf /K})\otimes_KK_m\simeq H^{i}_{\DR}(X_m/K_m, \cal O_{{(X_m)\daginf} /K_m})$$
de $K_m$-espaces vectoriels compatible avec  l'action de $({\bfF}^i)^m$. D'où l'égalité:
$$\sum_i(-1)^i \Tr\big ((q^{\dim X}({\bfF}^i)^{-1})^{m},  H^{i}_{\DR}(X/K, \cal O_{\Xdaginf /K})\big)= L(X_m)= N_m$$  où $N_m$ est le nombre de points
de $X$ à valeurs dans $\Bbb F_{q^m}$, qui implique la factorisation de la fonction Zêta, par le raisonnement habituel, en comparant les dérivées logarithmiques des deux membres du théorème.
\enddemo
Comme nous l'avons dit dans l'introduction,
cette factorisation était connue auparavant dans le cas affine ([M$_2$]),
complété par le théorème de finitude [Me$_2$], et dans le cas propre et lisse 
par la théorie cristalline [B$_2$].
\begin{Rema}
Le morphisme de Frobenius {\bf n'est pas défini} comme morphisme de complexes sur le complexe de de Rham $\Omega^{\bullet}_{\Xdaginf/K}$ du site infinitésimal, parce que l'action de Frobenius {\bf ne commute pas} à  l'action du groupe
$\cal G_{\Xdaginf}$. L'action de Frobenius n'est définie sur ce complexe qu'en {\bf catégorie dérivée}, parce que ses faisceaux de cohomologies
$\cal Ext^{\bullet}_{\calDdag_{\Xdaginf /K}}(
\cal O_{\Xdaginf /K},\cal O_{\Xdaginf /K})$ sont isomorphes aux faisceaux de Zariski $\cal Ext^{\bullet}_{\calDdag_{\Xdaginf /K},\Sp}(
\cal O_{\Xdaginf /K},\cal O_{\Xdaginf /K})$. On ne voit pas a priori comment définir globalement l'action de Frobenius sur la cohomologie  à partir du complexe $\Omega^{\bullet}_{\Xdaginf/K}$
sans passer par la catégorie des complexes spéciaux. C'est le point.\end{Rema}
\section{La suite exacte de Gysin en cohomologie de de Rham $p$-adique et la classe de cohomologie d'un cycle}
Dans ce paragraphe   nous allons établir la suite exacte de Gysin en cohomologie de de Rham $p$-adique sur $K$ pour un couple $Y\subset X$ lisse sur $k$,
généralisant la suite  exacte de Gysin  pour un couple de variétés affines lisses ([M$_1$], [Me$_2$]), et nous allons construire un morphisme gradué entre le groupe des cycles de $X$ et l'algèbre graduée de la cohomologie de de Rham $p$-adique de $X$. Nous suivons la méthode, introduite  dans le cas affine dans  
[Me$_2$], qui consiste d'abord à montrer un résultat  de compatibilité des images directes différentielle et topologique pour le foncteur de de Rham local.

\subsection{La compatibilité des foncteurs images directes topologique et différentielle}
Nous allons d'abord démontrer le théorème de compatibilité entre le foncteur image directe des modules spéciaux et le foncteur image direct  du
complexe de de  Rham local.

\begin{theo}\label{ima-dr}Soit $f : Y\to X$ un morphisme de schémas  lisses   sur $k$,
où $X$ est séparé sur $k$, et soit $\calMdaginf $ un complexe de la catégorie
$ \Dplus ((\calDdag_{\Ydaginf /K}, \Sp)\Mod)$. Alors, il existe un isomorphisme canonique de la catégorie $ \Dplus (K_X)$:
$$\DeuxLignes
\bfR \cHom_{\calDdag_{\Xdaginf /K}, \Sp}(
\cal O_{\Xdaginf /K},f_*^{\diff}\calMdaginf )[\dim X]
\simeq\\\simeq
\bfR f_*\bfR \cHom_{\calDdag_{\Ydaginf /K}, \Sp}(
\cal O_{\Ydaginf /K},\calMdaginf )[\dim Y].\endlignes$$ 
\end{theo}


Par construction, le foncteur image directe des modules spéciaux $f_*^{\diff}$ est le composé $p_*^{\diff}\circ i_*^{\diff}$ du foncteur image directe dans le cas d'une
immersion fermée et dans le cas d'une projection. Le foncteur $ \bfR f_*$ est le composé $ \bfR p_*\circ  \bfR i_*$. Pour monter le
théorème précédent il suffit donc de le montrer dans le cas d'une immersion 
fermée puis dans le cas d'une projection. 


\begin{lemm}Soit  $i: Y\to X$ une immersion fermée de schémas lisses sur $k$,  et soient $\cal P\daginf $ un $\calDdag_{\Xdaginf /V}$-module à
droite spécial plat et $\calMdaginf $ un $\calDdag_{\Ydaginf /V}$-module à gauche spécial. Il existe un isomorphisme canonique de
projection de faisceaux de Zariski de 
$V$-modules sur $X$:
$$\cal P\daginf \otimes_{\calDdag_{\Xdaginf /V}}i_{*,0}^{\diff}\calMdaginf \to i_{*}\big(i^{*,0}_{\diff}\cal P\daginf \otimes_{\calDdag_{\Ydaginf /V}}\calMdaginf \big).$$
\end{lemm}

\demo
Soient $U$ un ouvert affine de $X$ et $W$ sa trace sur $Y$. Soit $\calUdag $ un relèvement de $U$ et $\calWdag $ un relèvement de $W$.
Par définition, la restriction à $U$ du faisceau de Zariski $\cal P\daginf \otimes_{\calDdag_{\Xdaginf /V}}i_{*,0}^{\diff}\calMdaginf $ est le
faisceau
$\cal P\dag_{\calUdag }\otimes_{\calDdag_{\calUdag /V}}i_*(\calDdag_{\calUdag \leftarrow \calVdag }\otimes_{\calDdag_{\cal
W\dag/V}}\cal M\dag_{\calWdag })$. De même, par construction, la restriction à $U$ du faisceau de Zariski 
$i_*(i^{*,0}_{\diff}\cal P\daginf \otimes_{\calDdag_{\Ydaginf /V}}\cal
M\daginf )$ est le faisceau
$i_*(i^{-1}\cal P\dag_{\calUdag }\otimes_{i^{-1}\calDdag_{\calUdag /V}}\calDdag_{\calUdag \leftarrow \calWdag }\otimes_{\calDdag_{\cal
W\dag/V}}\cal M\dag_{\calWdag })$. Le morphisme du lemme est le morphisme naturel de projection:
$$\cal P\dag_{\calUdag }\Otimes_{\calDdag_{\calUdag /V}}
i_*\big(\calDdag_{\calUdag \leftarrow \calWdag }\Otimes_{\calDdag_{\cal
W\dag/V}}\cal M\dag_{\calWdag }\big)
\to i_*\big(i^{-1}\cal P\dag_{\calUdag }\Otimes_{i^{-1}\calDdag_{\calUdag /V}}\calDdag_{\calUdag \leftarrow \cal
W\dag}\Otimes_{\calDdag_{\calWdag /V}}\cal M\dag_{\calWdag }\big),$$ lequel  est un isomorphisme puisque $i$ est une immersion fermée. Cet isomorphisme de projection est invariant par l'action du groupe $\cal
G_{\calUdag }$, et   se recolle en un isomorphisme de faisceaux de $V$-modules sur~$X$. 
\enddemo

Comme le foncteur $i_{*,0}^{\diff}$ est
exact, l'isomorphisme précédent se dérive de façon évidente en un isomorphisme de projection:
$$\calNdaginf \Lotimes_{\calDdag_{\Xdaginf /V}, \Sp}i_*^{\diff}\calMdaginf \to i_*(i^*_{\diff}\calNdaginf 
\Lotimes_{\calDdag_{\Ydaginf /V},\Sp}\calMdaginf )\,,$$ pour un complexe $\calNdaginf $ de la catégorie $ \Db (\Modd(\calDdag_{\Xdaginf /V}, \Sp))$ et un complexe $\calMdaginf $ de la catégorie $ \Dplus ((\calDdag_{\Ydaginf /V}, \Sp)\Mod)$. Appliquons l'isomorphisme
précédent au cas $\calNdaginf =
\omega_{\Xdaginf /V}$ et tenons compte de l'isomorphisme du théorème \ref{inv-tri'}. On trouve l'isomorphisme de projection:
$$ \omega_{\Xdaginf /V}\Lotimes_{\calDdag_{\Xdaginf /V}, \Sp}i_*^{\diff}\calMdaginf \to i_*
(\omega_{\Ydaginf /V}\Lotimes_{\calDdag_{\Ydaginf /V}, \Sp}\calMdaginf )$$ 
et, en particulier, l'isomorphisme sur $K$:
$$\omega_{\Xdaginf /K}\Lotimes_{\calDdag_{\Xdaginf /K}, \Sp}i_*^{\diff}\calMdaginf \to i_*
(\omega_{\Ydaginf /K}\Lotimes_{\calDdag_{\Ydaginf /K},\Sp}\calMdaginf )\,.$$ 
Mais 
le complexe: $$\omega_{\Xdaginf /K}\Lotimes_{\calDdag_{\Xdaginf /K}, \Sp}i_*^{\diff}\calMdaginf $$ est canoniquement isomorphe
au complexe: $$\bfR \cHom_{\calDdag_{\Xdaginf /K}, \Sp}(
\cal O_{\Xdaginf /K},i^{\diff}_*\calMdaginf )[\dim X],$$ et la formule de projection précédente fournit l'isomorphisme du théorème 
\ref{ima-dr} dans le cas d'une immersion fermée.

\bigskip
Considérons le cas d'une projection $p: Y\times_k X\to X$. Posons pour simplifier les notations  $Z:= Y\times_k X$. Supposons d'abord
que $Y$ et $X$ sont affines et soient $\calYdag $ et $\calXdag $ des relèvements plats sur $V$ de $Y$ et de $X$. Alors, 
$p_*^{\diff}\calMdaginf $ est canoniquement isomorphe au complexe: 
$$ \bfR p_*(\calDdag_{\calXdag \gets \cal
Y\dag\times_{V}\calXdag /K}\Lotimes_{\calDdag_{\calYdag \times_{V}\calXdag/K}}\cal M\dag_{\calYdag \times_{V}\calXdag}),$$
et le complexe de de Rham de $\calMdaginf $ est canoniquement isomorphe au complexe:
$$\omega_{\calYdag \times_{V}\calXdag/K}\Lotimes_{\calDdag_{\calYdag \times_{V}\calXdag/K}}\cal M\dag_{\calYdag \times_{V}\calXdag}[-\dim Z].$$ 
Le morphisme de projection topologique:
$$\DeuxLignes  
\omega_{\calXdag /K}\LOtimes_{\calDdag_{\calXdag /K}}\bfR p_*\big(\calDdag_{\calXdag \gets \cal
Y\dag\times_{V}\calXdag /K}\LOtimes_{\calDdag_{\calYdag \times_{V}\calXdag/K}}\cal M\dag_{\calYdag \times_{V}\calXdag}\big)\to
\\
\to \bfR p_*\big(p^{-1}\omega_{\calXdag /K}\LOtimes_{p^{-1}\calDdag_{\calXdag /K}}\calDdag_{\calXdag \gets \cal
Y\dag\times_{V}\calXdag /K}\LOtimes_{\calDdag_{\calYdag \times_{V}\calXdag/K}}\cal M\dag_{\calYdag \times_{V}\calXdag}\big).\endlignes
$$ 
est un \em{isomorphisme}, parce que le $\calDdag_{\calXdag /K}$-module à droite  $\omega_{\calXdag /K}$ admet
une résolution finie par des $\calDdag_{\calXdag /K}$-modules localement libres de type fini. Mais en vertu du théorème \ref{inv-tri'}, l'image
inverse: 
$$p^*_{\diff}\omega_{\calXdag /K}:= p^{-1}\omega_{\calXdag /K}\Lotimes_{p^{-1}\calDdag_{\calXdag /K}}\calDdag_{\calXdag \gets \cal
Y\dag\times_{V}\calXdag /K}$$ est canoniquement isomorphe à $\omega_{\cal
Y\dag\times_{V}\calXdag /K}$, ce qui fournit la compatibilité  du théorème \ref{ima-dr} dans ce cas-là.


\bigskip
Supposons que $X$ est affine et soit $\calXdag $ un relèvement $V$-plat  de $X$. Le complexe $p_*^{\diff}\calMdaginf $ se représente par
construction par le complexe: 
$$p_{*}\cHom_{q^{-1}\calDdag_{\Ydaginf /K}}\big(q^{-1}\cal O_{\Ydaginf /K},R_{\Ydaginf \times\calXdag}( \cal I^{\bullet}_{Z\daginf}  )\big)[\dim Y],$$ où $ \cal I^{\bullet}_{Z\daginf} $ est une résolution de $\calMdaginf $ par des $\calDdag_{Z_{\inf}/K}$-modules à gauche
spéciaux injectifs. La formule de projection
fournit encore un isomorphisme:
$$\DeuxLignes  
\omega_{\calXdag /K}\LOtimes_{\calDdag_{\calXdag /K}}\bfR p_*\cHom\SubX{q^{-1}\calDdag_{\Ydaginf /K}}
\big(q^{-1}\cal O_{\Ydaginf /K}, R_{\Ydaginf \times\calXdag}( \cal I^{\bullet}_{Z\daginf}  )\big)[\dim Y]\simeq
\\
\simeq\bfR p_*\mBig(p^{-1}\omega_{\calXdag /K}\LOtimes_{p^{-1}\calDdag_{\calXdag /K}}\cHom\SubX{q^{-1}\calDdag_{\Ydaginf /K}}\big(q^{-1}\cal O_{\Ydaginf /K}, R_{\Ydaginf \times\calXdag}( \cal I^{\bullet}_{Z\daginf}  )\big)\mBig)[\dim Y].
\endlignes$$ Mais si $Y$ est un ouvert affine  
 et si $\calIdaginf  $ est un $\calDdag_{Z_{\inf}/K}$-module spécial injectif,   le faisceau:
$$\cTor_i^{p^{-1}\calDdag_{\calXdag /K}}\mBig(p^{-1}\omega_{\calXdag /K}, \cHom_{q^{-1}\calDdag_{\Ydaginf /K}}\big(q^{-1}\cal O_{\Ydaginf /K}, R_{\Ydaginf \times\calXdag}(\calIdaginf  )\big)\mBig)$$ est en vertu du cas précédent nul si $i\neq \dim Y$, et l'on a l'isomorphisme local:
$$\DeuxLignes
\cTor_{\dim Y}^{p^{-1}\calDdag_{\calXdag /K}}\mBig(p^{-1}\omega_{\calXdag /K}, \cHom\Sub{5pt}{q^{-1}\calDdag_{\Ydaginf /K}}\big(q^{-1}\cal O_{\Ydaginf /K}, R_{\Ydaginf \times\calXdag}(\calIdaginf  )\big)\mBig)
\simeq\\\simeq \cHom\Sub{5pt}{\calDdag_{Z\daginf /K}}\big(\cal O_{Z\daginf /K}, R_{\Ydaginf \times\calXdag}(\calIdaginf  )\big),
\endlignes
$$ qui se recolle naturellement en un isomorphisme
global au-dessus de $Y\times X$,  fournissant ainsi  l'isomorphisme du théorème dans le cas où $X$ est affine.

\bigskip
Le cas où $X$ est affine montre que si $\calIdaginf  $ est un $\calDdag_{Z_{\inf}/K}$-module à gauche spécial et injectif,   le faisceau de Zariski
$\cTor_i^{\goodsmash{1}{0.4}{\calDdag_{\Xdaginf /K}, \Sp}}(\omega_{\Xdaginf /K}, p_{*,0}^{\diff}\calIdaginf  )$ est nul si $i\neq \dim Y$, et l'on a l'isomorphisme local:
$$\cTor_{\dim Y}^{\calDdag_{\Xdaginf /K}, \Sp}(\omega_{\Xdaginf /K}, p_{*,0}^{\diff}\calIdaginf  )\simeq p_*\cHom_{\calDdag_{Z\daginf /K}, \Sp}(\cal O_{Z\daginf /K}, \calIdaginf  ),$$ qui se recolle naturellement en isomorphisme
global au-dessus de $ X$\glossary{$\cTor_{\dim Y}^{\calDdag_{\Xdaginf /K}, \Sp}(\omega_{\Xdaginf /K}, p_*^{\diff}\calIdaginf  )$}.  En prenant une résolution de $\calMdaginf $ par des $\calDdag_{Z_{\inf}/K}$-modules spéciaux injectifs, on obtient un isomorphisme de complexes  qui représente l'isomorphisme du théorème \ref{ima-dr} dans le cas
d'une projection.

\subsection{La suite exacte  de Gysin en cohomologie de de Rham $p$-adique}
Nous allons   établir d'abord  l'analogue $p$-adique de l'isomorphisme de Thom, puis    en déduire formellement la suite de Gysin.

\begin{theo}\label{sui-gys} Soit une immersion fermée $i: Y\to X$  de schémas lisses et séparés sur $k$, de complémentaire $j: U\to X$. Alors, il existe des
isomorphismes canoniques:
$$ H_{\DR}^{\bullet -2\codim_XY}(Y/K)\simeq H_{\DR,Y}^{\bullet }(X/K) $$ dont on déduit  une suite exacte longue de Gysin:
$$ \to H_{\DR}^{\bullet -2\codim_XY}(Y/K)\to H_{\DR}^{\bullet }(X/K)\to H_{\DR}^{\bullet }(U/K)\to \cdot$$
\end{theo}

\demo  En vertu du théorème  \ref{coh-loc}, et c'est là le point crucial, on a un isomorphisme
canonique de la catégorie $ \Dplus ((\calDdag_{\Xdaginf /K}, \Sp)\Mod)$:
$$i^{\diff}_*\cal O_{\Ydaginf /K}\simeq \bfR \Gamma_Y(\cal O_{\Xdaginf /K})[\codim_XY].$$ En vertu du théorème \ref{ima-dr} appliqué à $\cal
M\daginf =
\cal O_{\Ydaginf /K}$ dans le cas d'une immersion fermée, on trouve  un isomorphisme canonique de complexes de Zariski: 
$$\DeuxLignes \bfR \cHom_{\calDdag_{\Xdaginf /K}, \Sp}(
\cal O_{\Xdaginf /K},\bfR \Gamma_Y(\cal O_{\Xdaginf /K}))
\simeq  \\\simeq
\bfR \cHom_{\calDdag_{\Ydaginf /K}, \Sp}(
\cal O_{\Ydaginf /K},\cal O_{\Ydaginf /K})[-2\codim_X Y].
\endlignes
$$ 
L'hypercohomologie du membre de droite calcule  par définition les espaces $ H_{\DR}^{\bullet
-2\codim_XY}(Y/K)$.

En vertu de la proposition
\ref{com-chl}, on a un isomorphisme canonique de complexes de Zariski:
$$\DeuxLignes
\bfR \cHom_{\calDdag_{\Xdaginf /K}, \Sp}\big(
\cal O_{\Xdaginf /K},\bfR \Gamma_Y(\cal O_{\Xdaginf /K})\big)
\simeq\\\simeq
\bfR \Gamma_Y\big(\bfR \cHom_{\calDdag_{\Xdaginf /K}, \Sp}(
\cal O_{\Xdaginf /K},\cal O_{\Xdaginf /K})\big).\endlignes$$ 
L'hypercohomologie du membre de droite calcule  par définition les espaces $ H_{\DR,Y}^{\bullet
}(X/K)$, ce qui fournit les isomorphismes de Thom.

Soient $j: U\to X$ le complémentaire de $Y$ et  le triangle de cohomologie locale:
$$\bfR \Gamma_Y(\cal O_{\Xdaginf /K})\to\cal O_{\Xdaginf /K}\to\bfR j^{\inf}_*j^{-1}_{\inf}\cal O_{\Xdaginf /K}\to\ ,$$
qui est un triangle distingué de la catégorie $ \Dplus ((\calDdag_{\Xdaginf /K}, \Sp)\Mod)$.
En vertu de la proposition \ref{com-chl}, on a un isomorphisme canonique:
$$\DeuxLignes
\bfR \hom_{\calDdag_{\Xdaginf /K}, \Sp}(
\cal O_{\Xdaginf /K},\bfR j^{\inf}_*j^{-1}_{\inf}\cal O_{\Xdaginf /K})
\simeq\\\simeq\bfR \hom_{\calDdag_{\Udaginf /K}, \Sp}(
\cal O_{\Udaginf /K},\cal O_{\Udaginf /K}).\endlignes
$$ 
On en déduit le  triangle distingué de la catégorie $ \Dplus (K)$:
$$\TroisLignes  
\bfR \hom_{\calDdag_{\Ydaginf /K}, \Sp}(
\cal O_{\Ydaginf /K},\cal O_{\Ydaginf /K})[-2\codim_X Y]
\to\\\to
\bfR \hom_{\calDdag_{\Xdaginf /K}, \Sp}(
\cal O_{\Xdaginf /K},\cal O_{\Xdaginf /K})\to
\\
\to\bfR \hom_{\calDdag_{\Udaginf /K}, \Sp}(
\cal O_{\Udaginf /K},\cal O_{\Udaginf /K})\to\ ,
\endlignes$$
dont la suite longue de cohomologie fournit la suite exacte de Gysin pour le couple $Y\subset X$:
$$ \cdots \to H_{\DR}^{\bullet -2\codim_XY}(Y/K)\to H_{\DR}^{\bullet }(X/K)\to H_{\DR}^{\bullet }(U/K)\to \cdots\eqno\enddemo$$

\medskip
\begin{Rema} Le lecteur prendra garde à ce que la suite exacte de Gysin n'a pas lieu pour la cohomologie de de Rham $p$-adique du site infinitésimal formel.\end{Rema}
\subsection{La classe de cohomologie d'un cycle}
Soient $X$ un schéma séparé et lisse sur $k$ et $Z$ un sous-$k$-schéma intègre de $X$. Notons $\sing(T)$  le lieu singulier du schéma {\bf réduit} associé à un schéma  $T$.

\begin{theo}Si  le corps de base est parfait, sous les conditions précédentes le morphisme de restriction:
$$ H_{\DR}^{2\codim_XZ}(X/K)\to H_{\DR}^{2\codim_XZ }(X-\sing(Z)/K)$$ est bijectif.
\end{theo}

\demo
Comme $Z$ est localement de type fini  sur $k$, le lieu singulier  $\sing(Z)$ est fermé en vertu du critère de  Zariski,  et la suite $\sing(Z)\supset(\sing(\sing Z))\dots$ est strictement décroissante. Si le corps de base $k$ est parfait, le lieu  non singulier coïncide avec le lieu  lisse  sur $k$ et $Z-\sing(Z)$ est lisse sur $k$. 
En considérant la filtration $X-\sing(Z)\subset X-\sing(\sing(Z))\dots $, on se ramène à supposer que $\sing(Z)$ est lisse sur 
$k$. 

Pour un sous-schéma  fermé $Y$ de $X$, lisse sur $k$, on a
$ H_{\DR, Y}^{l}(X/K)=0$ pour $l< 2\codim_XY$. En effet, $ H_{\DR, Y}^{\bullet}(X/K)$ est par définition l'hypercohomologie du complexe de Zariski: 
$$\ecartlignes-3pt\DeuxLignes
\bfR \Gamma_Y\big( \bfR \cHom_{\calDdag_{\Xdaginf /K}, \Sp}(
\cal O_{\Xdaginf /K},\cal O_{\Xdaginf /K})\big)
\simeq\\\simeq  \bfR \cHom_{\calDdag_{\Xdaginf /K}, \Sp}\big(
\cal O_{\Xdaginf /K},\bfR \Gamma_Y(\cal O_{\Xdaginf /K})\big),
\endlignes$$
qui est isomorphe, en vertu des résultats précédents, au complexe: $$  \bfR \cHom_{\calDdag_{\Ydaginf /K}, \Sp}(
\cal O_{\Ydaginf /K},\cal O_{\Ydaginf /K})[-2\codim_X Y].$$ La suite exacte de cohomologie locale du couple $\sing(Z)\subset X$:
$$\ecartlignes5pt
\DeuxLignes
 H_{\DR,\sing(Z)}^{2\codim_XZ}(X/K)\to H_{\DR}^{2\codim_XZ}(X/K)\to H_{\DR}^{2\codim_XZ }(X{-}\sing(Z)/K)\to\\\to H_{\DR, \sing(Z)}^{2\codim_XZ+1}(X/K)\to\endlignes
\postskip-0,ex$$ 
montre le théorème, puisque 
$$ 2\codim_XZ<2\codim_XZ+1<2\codim_X\sing(Z)\,.\eqno
\enddemo$$
\smallskip

\begin{defi} Soit $X$ un schéma séparé et lisse sur $k$. \begin{liste} \item 1) 
Si $Y$ est un sous-schéma fermé de $X$ lisse sur $k$,
on définit la classe de cohomologie $c\ell(Y)$ de $Y$ comme l'image du morphisme identique 
de $\cal O_{\Ydaginf /K}$ par le morphisme canonique:
$$\ecartlignes-10pt\DeuxLignes
 \hom\goodSub{5pt}{4pt}{-7mm}{\calDdag_{\Ydaginf /K}, \Sp}(
\cal O_{\Ydaginf /K},\cal O_{\Ydaginf /K}) =: H_{\DR}^{0}(Y/K)\simeq H_{\DR, Y}^{2\codim_XY }(X/K) \to\\\to H_{\DR}^{2\codim_XY }(X/K).\endlignes$$
\item 2) Si  le corps de base est parfait, on définit la classe de cohomologie $c\ell(Z)$ d'un sous-schéma
$Z$ intègre comme l'image inverse par l'isomorphisme du théorème précédent de la classe de  cohomologie  de
$Z-\sing(Z)$ dans $ H_{\DR}^{2\codim_XZ}(X-\sing(Z)/K)$.
\end{liste}
\end{defi}
La classe $c\ell(Z)$ est donc un élément bien défini de $ H_{\DR}^{2\codim_XZ}(X/K)$, dont l'image
dans $ H_{\DR}^{2\codim_XZ}(X-Z/K)$ est nulle. On obtient donc par linéarité un morphisme gradué:
$$ C^\bullet(X)\to H_{\DR}^{2\bullet}(X/K)\leqno c\ell:$$ entre le groupe des cycles de $X$ et le $K$-espace gradué de la cohomologie de de Rham $p$-adique de $X$. Naturellement, on s'attend à ce que l'intersection des cycles soit compatible avec la composition des classes de cohomologie.

\section{Le lemme de Poincaré et la formule de Künneth pour l'espace affine}
Soient $X$ un schéma lisse sur $k$ et $p: X\times_k A^n_k\to X$ la projection de l'espace affine au-dessus de $X$ sur la base $X$. 
Nous allons déduire  du théorème de compatibilité \ref{ima-dr} le
lemme de Poincaré, qui est aussi un cas particulier de la formule de Künneth.

\begin{theo}Soit $\calMdaginf $ un $\calDdag_{\Xdaginf /K}$-module à gauche spécial   {\bf parfait}. Il existe alors un isomorphisme canonique:
$$\calMdaginf \simeq p^{\diff}_*p^*_{\diff}\calMdaginf [-n]$$ de complexes spéciaux.
\end{theo}

\demo
Par récurrence, on peut supposer que $n=1$, que   $t$ est une coordonnée sur $A^1_k$ et enfin que $\cal A^1_V$ est le relèvement associé d'algèbre $(V[t])\dag$.  
Cela entraîne que, pout tout ouvert $\calUdag $, les images directes supérieures topologiques $R^lp_*p^*_{\diff}\cal M\dag_{\calUdag }$ sont nulles pour
$l\geq1$. En effet, localement sur la base, $p^*_{\diff}\calMdaginf $ admet une résolution finie par des sommes finies du module $\calDdag_{\cal
U\dag\times_V\cal A^1_V\to\cal A^1_V/K}$ qui  est acyclique pour le foncteur $p_*$, en vertu du théorème du symbole total et du théorème d'acyclicité. Cela
entraîne que, pour tout ouvert $\calUdag $, le complexe $p^{\diff}_*p^*_{\diff}\cal M\dag_{\calUdag }[-1]$ se représente par:
$$0\to p_*p^*_{\diff}\cal M\dag_{\calUdag }\stackrel{\partial_t}{\to} p_*p^*_{\diff}\cal M\dag_{\calUdag }\otimes_{\cal O_{\cal
U\dag\times_V\cal A^1_V/K}}\Omega_{\cal
U\dag\times_V\cal A^1_V/\calUdag /K}\to0.$$
On en déduit un morphisme:
$$\cal M\dag_{\calUdag }\to p^{\diff}_*p^*_{\diff}\cal M\dag_{\calUdag }[-1],$$ et il s'agit d'établir que c'est un quasi-isomorphisme. Le complexe
précédent n'a de la cohomologie qu'en degré zéro, parce l'action de $\partial_t$ sur $(A\dag[t])\dag_K$ est surjective, et il est facile de voir
que le morphisme précédent induit un isomorphisme  entre $\cal M\dag_{\calUdag }$ et le noyau du morphisme du complexe précédent. Ces isomorphismes locaux
se recollent en un isomorphisme de modules spéciaux entre $\calMdaginf $ et la cohomologie de degré zéro du complexe concentré cohomologiquement 
en degré zéro $p^{\diff}_*p^*_{\diff}\calMdaginf [-1]$.
\enddemo

\begin{coro}Sous les conditions précédentes, l'isomorphisme précédent induit des isomorphismes canoniques:
$$ H^{\bullet}_{\DR}(X/K, \calMdaginf )\simeq H^{\bullet}_{\DR}(X\times_k A^n_k/K, p^*_{\diff}\calMdaginf ).$$
\end{coro}

\demo C'est une conséquence du théorème \ref{ima-dr}, dans le cas d'une projection, et du théorème précédent.
\enddemo

\begin{coro}Sous les conditions précédentes, l'isomorphisme précédent induit des isomorphismes canoniques:
$$ H^{\bullet}_{\DR}(X/K, \cal O_{\Xdaginf /K})\simeq H^{\bullet}_{\DR}(X\times_k A^n_k/K, \cal O_{(X\times_kA^n_k)\daginf /K}).$$
\end{coro}

\demo C'est une conséquence du théorème \ref{inv-tri} et du corollaire précédent, sachant que la résolution de Spencer du lemme \ref{Spen} montre que 
$\cal O_{\Xdaginf /K}$ est parfait. Autrement dit, le foncteur de de Rham $p$-adique est homotopiquement invariant, comme on dit aujourd'hui dans la littérature consacrée aux motifs.
\enddemo


\catcode`@=11
\def\IT#1,{{\it #1},}
\def\BF#1 {{\bf #1} }
\def\@biblabel#1{\hss[#1]}
\def\auteurs=#1\livre{\def\AU{\ignorespaces#1\unskip, }\livre}
\def\livre=#1\article{\def\LI{\setbox0=\hbox{\ignorespaces#1\unskip}%
\ifdim\wd0>0pt{\it\ignorespaces#1\unskip\/}\fi}\article}
\def\article=#1\divers{\def\AR{{\rm\ignorespaces#1\unskip}}\divers}
\def\divers=#1\par{\def\EX{\ignorespaces#1\unskip.}{\AU}{``\LI\AR''}. \EX}
\def\bib[#1]{\bibitem[#1]{#1}\def\AU{}\def\LI{}\def\AR{}\def\EX{}\ignorespaces}

\footnotesize

\closeoutglossary

\def\indexname{Liste des notations}
\def\glossaryentry#1#2{\par\vskip0.5ex\noindent
\hangafter1\hangindent2em #1\cleaders\hbox{.}\hfill #2}
\let\indexentry\glossaryentry

\begin{theindex}\parskip0pt
\baselineskip=1\baselineskip plus1pt minus 1pt
\input \jobname.glo
\end{theindex}

\end{document}